# Assouad Dimension and Fractal Geometry

Jonathan M. Fraser

For Rayna and Dylan

# Contents

# Acknowledgements

I am grateful to the *Leverhulme Trust* and the *EPSRC* for financial support during the writing of this book and for supporting much of the research which I document. I held a *Leverhulme Trust Research Fellowship* (RF-2016-500) from 2016-2018 for a project entitled *Fractal Geometry and Dimension Theory* which concerned many problems relating to the Assouad dimension. In particular, the papers [89, 92] arose directly from that Fellowship. I also held a *Leverhulme Trust Research Project Grant* (RPG-2019-034) from September 2019 for a project entitled *New Perspectives in the Dimension Theory of Fractals.* This project concerns many aspects of this book, most specifically the Assouad spectrum. I held an *EPSRC Standard Grant* (EP/R015104/1) from April 2018 for a project entitled *Fourier Analytic Techniques in Analysis and Geometry.* I was also financially supported by the EPSRC during my PhD, which included the writing of the paper [88], and during my time as a Research Fellow in Warwick, where I worked on the paper [96], for example.

I am grateful to Cambridge University Press for their assistance during the writing and publishing process. I thank Tom Harris for carefully answering many questions, especially early on, and Suresh Kumar for technical assistance.

I was fortunate to have many friends and colleagues read drafts of this book. I am especially grateful to Kenneth Falconer for reading the entire manuscript and making numerous detailed suggestions. I also thank Amlan Banaji, Stuart Burrell, Haipeng Chen, Antti Käenmäki, István Kolossváry, Lawrence Lee, Juha Lehrbäck, Chris Miller, Tuomas Orponen, Pablo Shmerkin, Liam Stuart, and Sascha Troscheit for making helpful suggestions. I am also very grateful to Tom Coleman, Ailsa Fraser, and Iain Fraser for thoroughly proof reading the book. I thank

Sascha Troscheit for expertly producing the percolation sets depicted in Figure 9.6.

Finally, I thank my family. Most of all my wife, Rayna, for her support, encouragement and patience during the writing process. The book is dedicated to her and our son, Dylan, born 30th November 2018. Maybe he will read this book someday. I am also extremely grateful to my mum, Ailsa, my dad, Iain, my sister, Cara, and her partner, Tom, for their support, especially relating to childcare! Those reading closely will observe that Cara and Tom did not help with the proof reading.

# Preface

I first encountered the Assouad dimension on Wednesday 20th April 2011 whilst attending an EPSRC workshop on *Dynamical Systems and Dimension* hosted by the University of Warwick. James Robinson gave a talk entitled *Assouad dimension, cube slicing, and the dynamics on finite-dimensional attractors*, which I remember included a discussion of the fact that the maximal volume of the intersection of the unit cube in $\mathbb{R}^d$ with an affine hyperplane is $\sqrt{2}$ (counter-intuitively independent of $d$). I was intrigued by the Assouad dimension, and surprised that I had not heard of it before, especially given my interest in the box and Hausdorff dimensions. Following the talk, I found two papers on the topic, one by Lars Olsen which gave a direct proof of the fact that self-similar sets satisfying the open set condition have equal Hausdorff and Assouad dimensions [219], and one by John Mackay which computed the Assouad dimension of certain self-affine sets [194]. Coincidentally both of these papers were also published in 2011. Another significant paper on the topic, which I found shortly after the papers of Mackay and Olsen, is an article by Jouni Luukkainen [192]. This article established many of the basic properties of the Assouad dimension, but its main focus was in proving a 'Szpilrajn Theorem' for Assouad dimension: the topological dimension of a separable metric space $X$ is the infimum of the Assouad dimension of metric spaces $Y$ such that $X$ and $Y$ are homeomorphic. Szpilrajn proved this with Assouad dimension replaced by Hausdorff dimension in 1937 [263].[1]

Olsen posed two (classically fractal) questions concerning the Assouad dimension which at the time lay only at the fringes of the fractal geometry dialogue, see [219, Questions 1.3 and 1.4]. The first question asked if the Assouad and Hausdorff dimensions could be distinct for self-similar sets (obviously requiring failure of the open set condition) and the second asked for the Assouad dimension of the famous self-affine sets introduced by Bedford and McMullen in the 1980s [26, 210]. The second of these questions was answered in Mackay's paper and I was able to answer the first one (in the affirmative) [88]. All I needed was an example, which I will describe later in this book in Theorem 7.2.1, but this example led to fruitful collaboration with James Robinson, Eric Olson and Alexander Henderson, in which we precisely described the Assouad dimension of all self-similar sets in $\mathbb{R}$, see Theorem 7.2.4. This work was completed in 2014 and published the following year [96]. I was hooked and ever since then I have spent a lot of time investigating the Assouad dimension, learning about its subtle, often counter-intuitive, properties and exploring its connections with dimension theory, fractal geometry and wider mathematics.

[1] Szpilrajn later changed his name to Marczewski to avoid persecution in Nazi Germany.

The Assouad dimension is rapidly becoming part of the mainstream dialogue in fractal geometry and the need for this book is highlighted by the fact that it is not mentioned in some of the most important and influential books in the field. For example, it is not mentioned in the books by Falconer [64, 69, 70, 72], Mattila [205, 207], Bishop and Peres [34], Edgar [60], Mandelbrot [197], or Barnsley [22, 23]. In many ways I see this book as a companion to, for example, Kenneth Falconer's *Fractal Geometry: Mathematical Foundations and Applications* [70], which presents a detailed and comprehensive treatment of dimension theory from a fractal geometry perspective. Many of the examples and problems I consider will have a similar flavour to those in [70].

The Assouad dimension is mentioned in — and is a central focus of — some books closely related to fractal geometry. In particular, James Robinson's *Dimensions, Embeddings, and Attractors* [240], and John Mackay and Jeremy Tyson's *Conformal dimension: Theory and application* [195]. Robinson's book focuses on the general embedding theory of dynamical systems, where the Assouad dimension is a key technical tool with many applications. We will discuss some examples in this direction in Section 13.3. Mackay and Tyson's book focuses on the problem of lowering dimension (usually Hausdorff or Assouad) via quasi-symmetric maps. This is a subtle and challenging problem and the Assouad dimension highlights many key features of this exploration. We will touch on this area in Chapter 12. Despite the important role played by dimension theory, these books [240, 195] do not have what I would call *classical fractal geometry* at their heart, but rather encounter fractal ideas on a journey motivated by a different set of problems. As such, this book serves a different purpose and will consider more *classically fractal* questions in the context of Assouad dimension, such as the dimension theory of iterated functions systems (self-similar, self-affine sets, other dynamical invariants) and geometric constructions (projections, products, slices, distance sets). As well as developing the theory in the context of fractal geometry, I will also explore numerous applications to problems in areas such as embedding theory, arithmetic geometry, Diophantine approximation, probability theory, and functional analysis. There will inevitably be some overlap with [240, 195] — for example, Mackay and Tyson's work on weak tangents and Assouad dimension has been central to recent progress in fractal geometry — but I will attempt to keep repetition to a minimum. Another purpose of this book is to consider various natural variations on the Assouad dimension, including its natural dual the lower dimension, the quasi-Assouad dimension, the Assouad spectrum (which I introduced with my first PhD student, Han Yu) and the corresponding dimensions of measures. I aim to present these ideas in a unified way and to inspire future explorations in these directions. As such this book also contains several new results as I attempt to present a unified and comprehensive theory.

The book is not meant to be a comprehensive survey and I apologise in advance for omitting a detailed discussion of much interesting work. Where appropriate I have tried to include a large set of references and to indicate where more details could have been included. Often results are not presented in their most general form. For example, I almost exclusively restrict my attention to sets and measures in Euclidean space despite most of the theory applying in more general metric spaces. My goal here is to present the key ideas in as simple a framework as possible. Again, where possible, I will indicate where more general results can be found. Elementary real analysis will be assumed and some familiarity with fractal geometry and dimension theory would be beneficial, but should not be required. Some of my favourite introductory books on analysis include: Howie [138], Rudin [244] and Stewart-Tall [257].

# PART ONE

# THEORY

# 1
# Fractal geometry and dimension theory

In this introductory chapter we briefly discuss the history and development of fractal geometry and dimension theory. We introduce and motivate some important concepts such as Hausdorff and box dimension. As part of this discussion we encounter covers and packings, which are central notions in dimension theory, and introduce the dimension theory of measures.

## 1.1 The emergence of fractal geometry

A *fractal* can be described as an object which exhibits interesting features on a large range of scales, see Figure 1.1. In pure mathematics, the Sierpiński triangle, the middle third Cantor set, the boundary of the Mandelbrot set, and the von Koch snowflake are archetypal examples and, in 'real life', examples include the surface of a lung, the horizon of a forest, and the distribution of stars in the galaxy. The fractal story began in the nineteenth century with the appearance of a multitude of strange examples exhibiting what we now understand as fractal phenomena. These included Weierstrass' example of a continuous nowhere differentiable function, Cantor's construction of an uncountable set with zero length, and Brown's observations on the path taken by a piece of pollen suspended in water (Brownian motion). During the first half of the twentieth century the mathematical foundations for fractal geometry were laid down by, for example, Besicovitch, Bouligand, Hausdorff, Julia, Marstrand and Sierpiński, and the theory was unified and popularised by the extensive writings of Mandelbrot in the 1970s, for example [197]. It was Mandelbrot who coined the term 'fractal', derived from the Latin *fractus* meaning 'broken'. Since then the subject has grown and developed as a self-contained discipline in pure mathematics touching on many other subjects, such as dynamical systems, geometric measure theory, analysis (real and complex), topology, number theory, probability theory and harmonic analysis. However, the importance of fractals is not restricted to abstract mathematics, with many naturally occurring physical phenomena exhibiting a fractal structure, such as graphs of random processes, percolation models and fluid turbulence. The mathematical challenge is to understand the mechanisms which generate and underpin fractal behaviour and to develop robust and rich theories concerning the geometric properties that such objects possess.

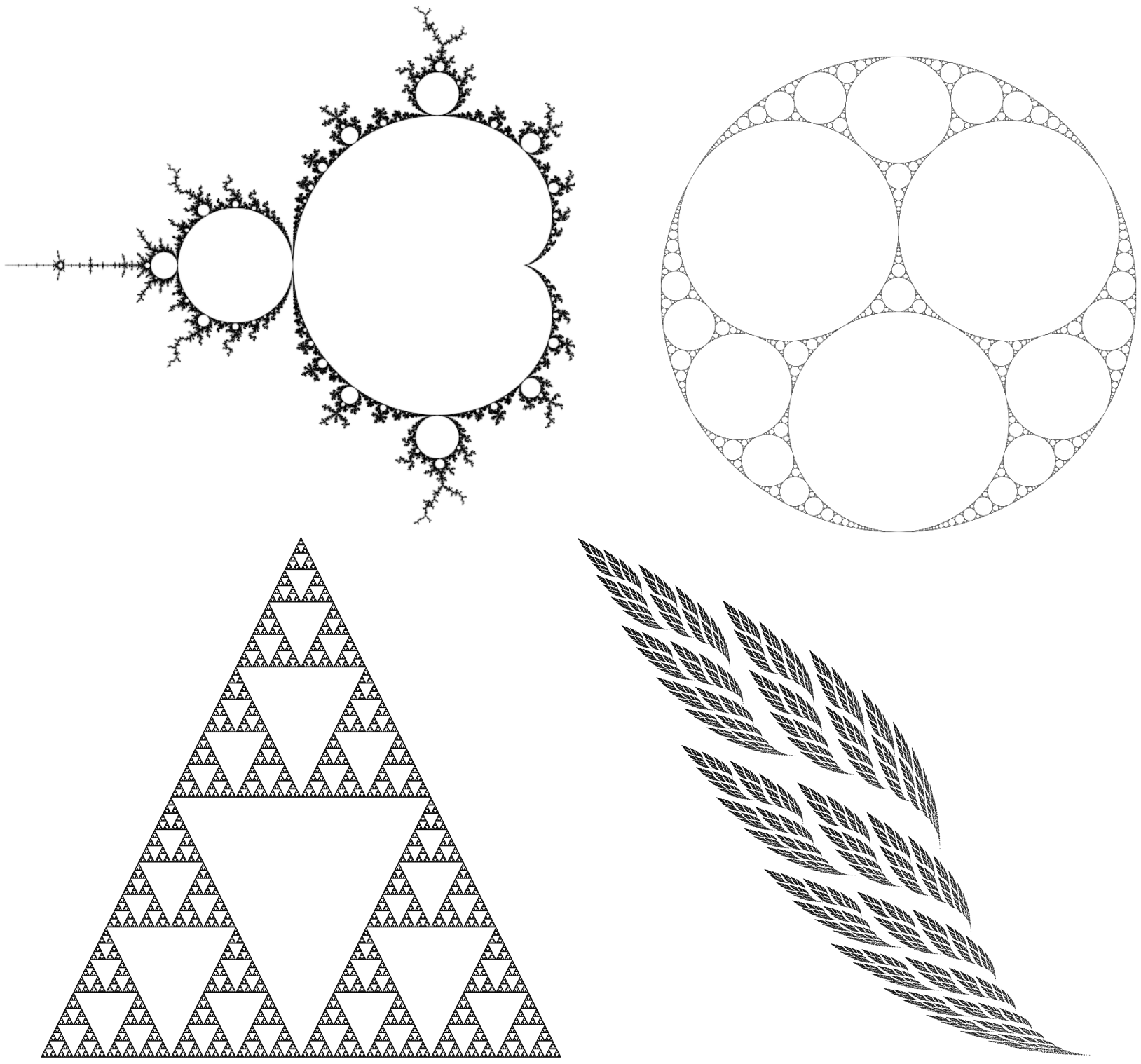

Figure 1.1 *Four fractals.* From top left moving clockwise: the boundary of the Mandelbrot set, the Apollonian circle packing, a (self-affine) leaf and the Sierpiński triangle.

## 1.2 Dimension theory

At the heart of fractal geometry lies *dimension theory*, the subject dedicated to understanding how to define, interpret, understand, and calculate dimensions of sets in Euclidean space or more general metric spaces. A *dimension* is a (usually non-negative real) number which gives geometric information concerning how the object in question fills up space on small scales. There are many distinct notions of dimension and one of the joys (and central components) of the subject is in understanding how these notions relate to each other, and how their behaviour compares in different settings or when applied to different families of examples.

A natural approach to dimension theory is to quantify how large a set is at a given scale by considering optimal *covers* by balls whose diameter is related to the scale. More precisely, given a scale $r > 0$, a finite or countable collection of sets $\{U_i\}_i$ is called an *$r$-cover* of a set $F$ if each of the sets $U_i$ has diameter less than or equal to $r$, and $F$ is contained in the union $\bigcup_i U_i$, see Figure 1.2. Throughout the book we write $|U| = \sup_{x,y \in U} |x-y|$ for the diameter of a non-empty set $U \subseteq \mathbb{R}^d$. Understanding how to find covers of a set at small scales underpins much of dimension theory and often the 'covering strategy' is specific to the setting, sometimes driven by dynamical invariance or *a priori* knowledge of another, related, set. This book is dedicated to a thorough analysis of the Assouad dimension and some of its natural variants. However, we will often attempt to put our discussion in a wider context for which we require other notions.

The Hausdorff dimension is arguably the most well-studied and important notion of fractal dimension. It was introduced by Hausdorff in 1918 [129], greatly developed by Besicovitch, and is considered extensively in many of the important books on fractal geometry, such as [34, 69, 70, 205]. It is defined in terms of Hausdorff measure, which can be viewed as a natural extension of Lebesgue measure to non-integer dimensions. Given $s \geqslant 0$ and $r > 0$, the *$r$-approximate $s$-dimensional Hausdorff measure* of a set $F \subseteq \mathbb{R}^d$ is defined by

$$\mathcal{H}^s_r(F) = \inf\left\{ \sum_i |U_i|^s : \{U_i\}_i \text{ is a countable } r\text{-cover of } F \right\}$$

and the *$s$-dimensional Hausdorff (outer) measure* of $F$ is $\mathcal{H}^s(F) = \lim_{r\to 0} \mathcal{H}^s_r(F)$. The limit exists because the sequence $\mathcal{H}^s_r(F)$ increases as $r$ decreases, but it may be infinite. The measures $\mathcal{H}^s$ can now be used to identify the *critical exponent* or *dimension* in which it is most appro-priate to consider $F$. First consider the square $[0,1]^2$, which has infinite length (length measures objects which are much smaller, such as line segments or smooth curves), and zero volume (volume measures objects which are much bigger, such as cubes or spheres). However, the *area* of the square is positive and finite (it hardly matters that the precise area is 1), demonstrating that the natural measure to use when considering squares is area, that is, 2-dimensional Lebesgue measure. It is no coincidence that we think of the square as a 2-dimensional object. Since we have continuum many Hausdorff measures to choose from, this leads naturally to the *Hausdorff dimension* of $F$ being defined as

$$\dim_{\mathrm{H}} F = \inf\Big\{s \geqslant 0 : \mathcal{H}^s(F) = 0\Big\} = \sup\Big\{s \geqslant 0 : \mathcal{H}^s(F) = \infty\Big\}.$$

It is a useful exercise to show that these two expressions for the Hausdorff dimension actually coincide. The value of the Hausdorff measure in the critical dimension, that is, $\mathcal{H}^{\dim_{\mathrm{H}} F}(F)$ is often rather hard to compute exactly and can be any value in $[0,\infty]$.

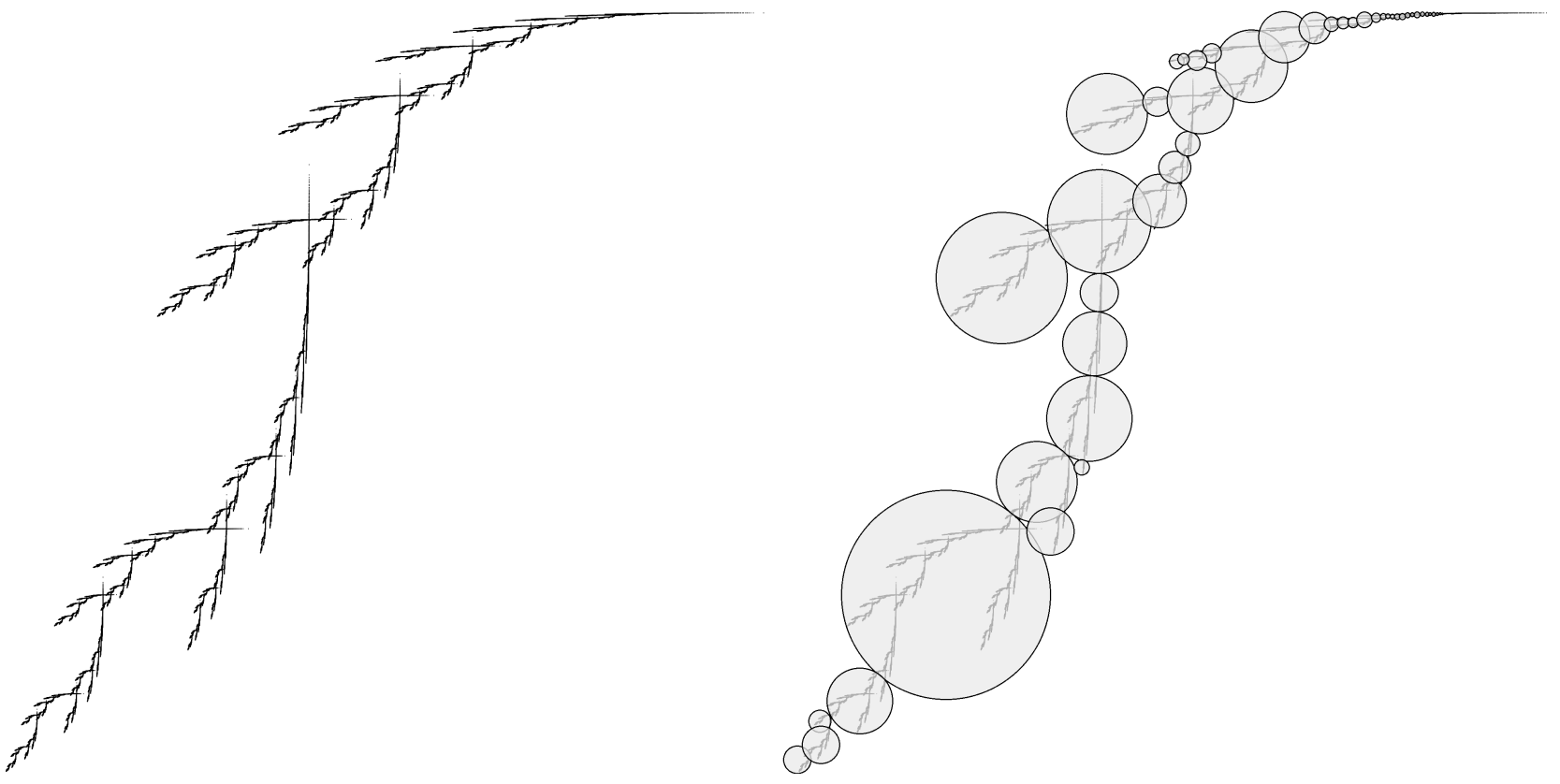

Figure 1.2 Left: a self-affine fractal. Right: a covering of the self-affine fractal using balls of arbitrarily varying radii. Understanding such covers leads to calculation of the Hausdorff dimension.

A less sophisticated, but nevertheless very useful, notion of dimension is box dimension. The *lower and upper box dimensions* of a non-empty bounded set $F \subseteq \mathbb{R}^d$ are defined by

$$\underline{\dim}_{\mathrm{B}} F = \liminf_{r\to 0} \frac{\log N_r(F)}{-\log r} \qquad \text{and} \qquad \overline{\dim}_{\mathrm{B}} F = \limsup_{r\to 0} \frac{\log N_r(F)}{-\log r},$$

respectively, where $N_r(F)$ is the smallest number of open sets required for a $r$-cover of $F$, see Figure 1.3. If $\underline{\dim}_{\mathrm{B}} F = \overline{\dim}_{\mathrm{B}} F$, then we call the common value the *box dimension* of $F$ and denote it by $\dim_{\mathrm{B}} F$. Note that, unlike the Hausdorff dimension, the box dimension is usually only defined for bounded sets since $N_r(F) = \infty$ for any unbounded set.

Notice that for a bounded set $F$, $r > 0$ and $s \geqslant 0$,

$$\mathcal{H}^s_r(F) \leqslant r^s N_r(F)$$

which immediately gives $\dim_{\mathrm{H}} F \leqslant \underline{\dim}_{\mathrm{B}} F \leqslant \overline{\dim}_{\mathrm{B}} F$.

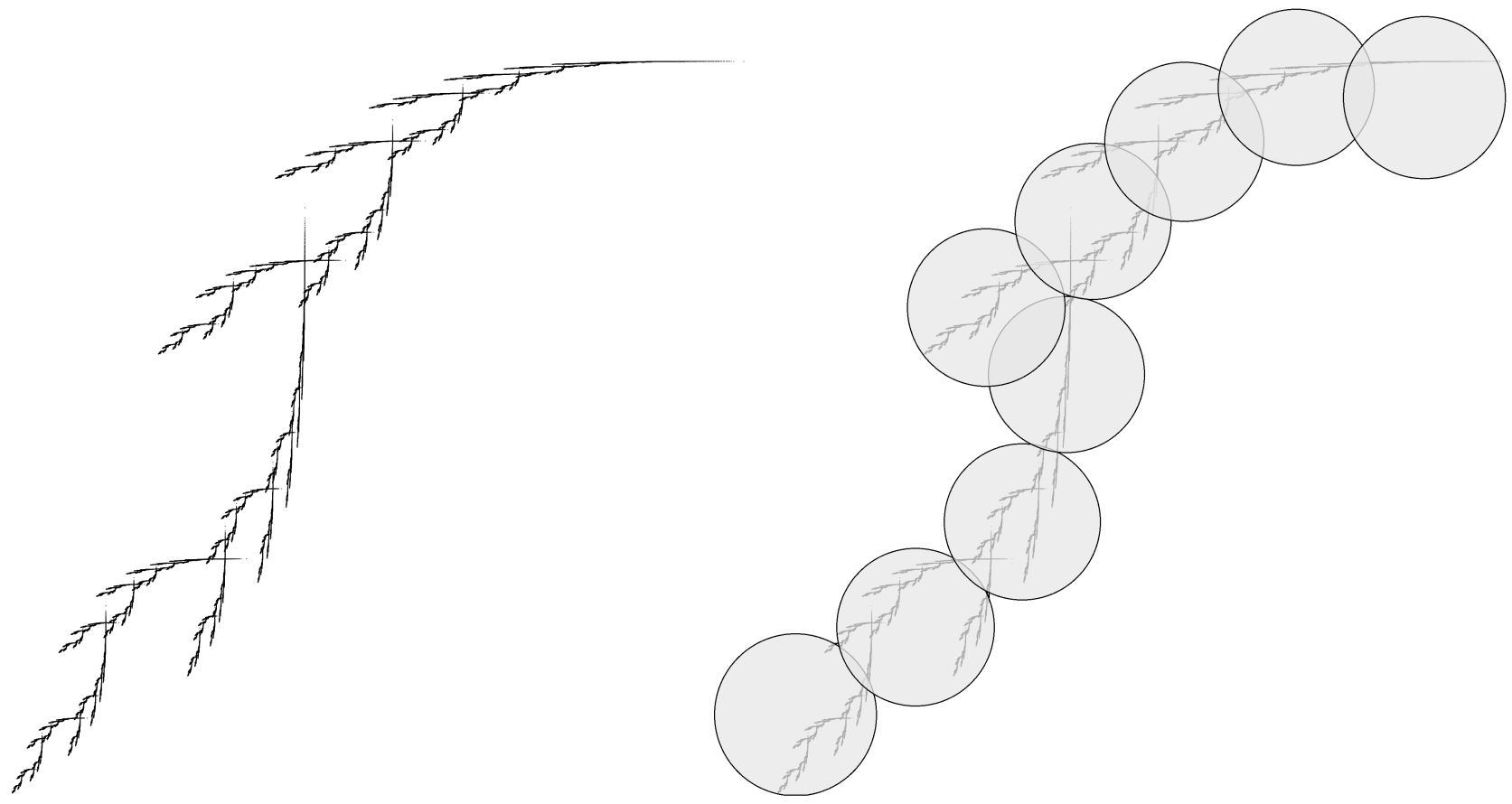

Figure 1.3 Left: the self-affine set from Figure 1.2. Right: a covering using balls of constant radii. Understanding such covers leads to calculation of the box dimensions.

One final notion, which we will mention less frequently, is the packing dimension. This can be defined by a suitable modification of the upper box dimension designed to make it countably stable, see Section 2.4. For $F \subseteq \mathbb{R}^d$, the *packing dimension* of $F$ is defined by

$$\dim_{\mathrm{P}} F = \inf \left\{ \sup_i \overline{\dim}_{\mathrm{B}} F_i \ : \ F = \bigcup_i F_i \right\}. \tag{1.1}$$

This definition works perfectly well for unbounded sets if we assume the $F_i$ are bounded. Moreover, one immediately gets $\dim_{\mathrm{H}} F \leqslant \dim_{\mathrm{P}} F \leqslant \overline{\dim}_{\mathrm{B}} F$. The packing dimension turns out to be a natural 'dual notion' to the Hausdorff dimension, where packings are used instead of covers. The usual formulation of packing dimension first defines packing measure, as a dual to the Hausdorff measure, and then packing dimension in the natural way. It was first introduced by Claude Tricot in 1982 [265].

We write $B(x,r) = \{y \in \mathbb{R}^d : |y - x| \leqslant r\}$ to denote the closed ball with centre $x \in \mathbb{R}^d$ and radius $r > 0$. A collection of balls $\{B(x_i, r)\}_i$ is called a (centred) *r-packing* of a set $F$ if the balls are pairwise disjoint and, for all $i$, $x_i \in F$. A related notion is that of $r$-separated sets. A set $X \subseteq F$ is called an *r-separated* subset of $F$ if each pair of distinct points from $X$ are separated by distance at least $r$. If $r_i = r$ for all $i$, then the set of centres $\{x_i\}_i$ of balls from an $r$-packing form a $2r$-separated subset of $F$. It is an elementary but instructive exercise to prove that if one replaces $N_r(F)$ in the definition of upper and lower box dimensions with any of:

(i) the maximum number of balls in an $r$-packing of $F$,

(ii) the maximum cardinality of an $r$-separated subset of $F$,

(iii) the number of $r$-cubes in an axes oriented mesh which intersect $F$,

(iv) the minimum number of cubes of sidelength $r$ required to cover $F$

then one obtains the same values for the box dimensions, see Figure 1.2. See [70, Section 2.1] for a more detailed discussion of this, along with some direct calculations.

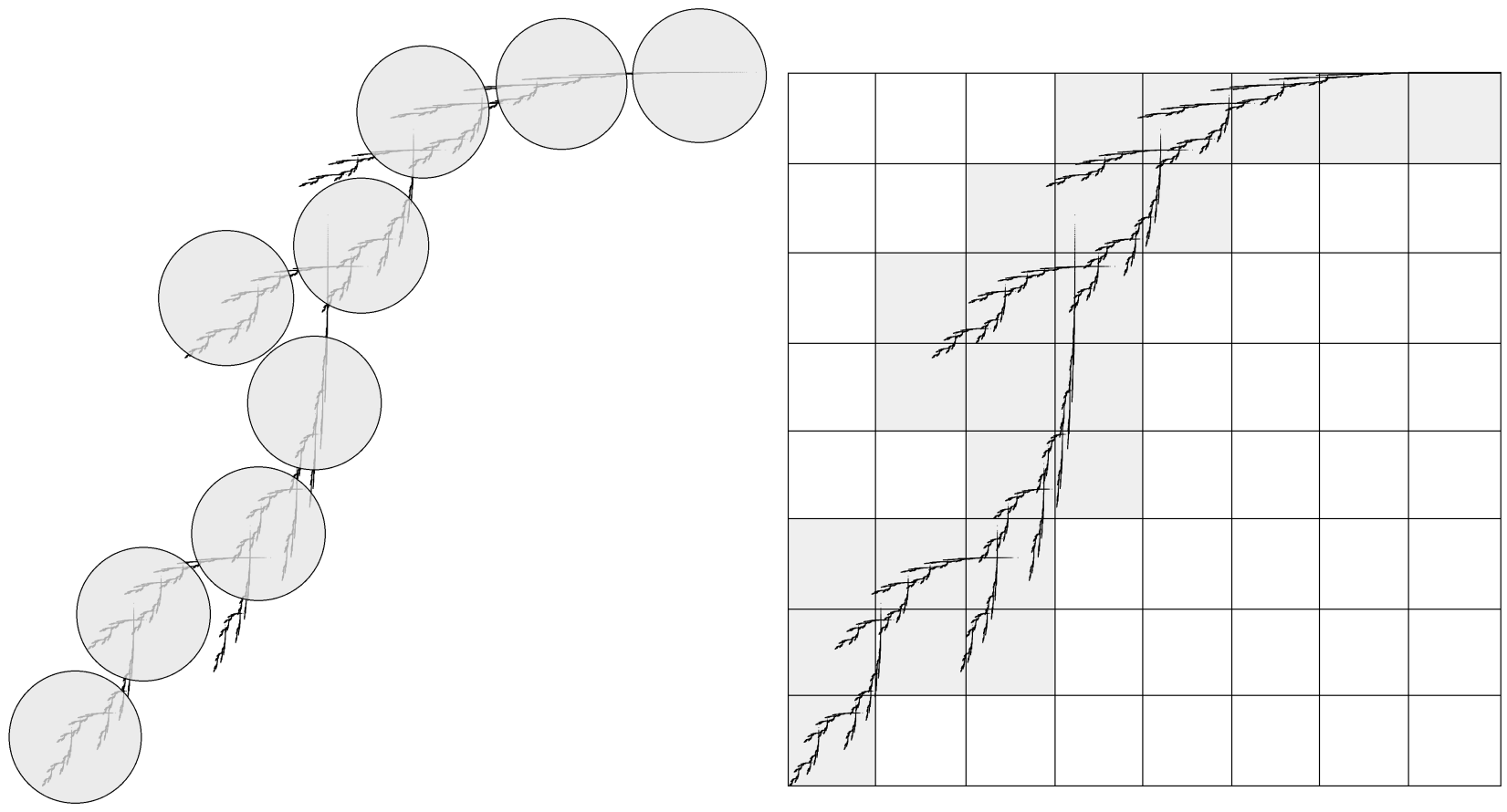

Figure 1.4 A packing of the self-affine set from Figure 1.2 using balls of constant radii centred in $F$ (left) and a mesh of squares imposed on $F$ with the squares intersecting $F$ shown in grey (right).

### 1.2.1 Dimension theory of measures

An important aspect of dimension theory is the interplay between the dimensions of sets and the dimensions of measures, see, for example, the mass distribution principle, Lemma 3.4.2. For this we need analogous notions of dimension for measures. The *(lower) Hausdorff dimension* of a Borel measure $\mu$ on $\mathbb{R}^d$ is

$$\dim_{\mathrm{H}} \mu = \inf\{\dim_{\mathrm{H}} E : \mu(E) > 0\}.$$

Thus the dimension of a measure is conveniently expressible in terms of dimensions of sets which 'see' the measure. We write

$$\mathrm{supp}(\mu) = \{x \in \mathbb{R}^d : \mu(B(x,r)) > 0 \text{ for all } r > 0\} \tag{1.2}$$

for the *support* of $\mu$, which is necessarily a closed set, and we say $\mu$ is *fully supported* on $F \subseteq \mathbb{R}^d$ if $\mathrm{supp}(\mu) = F$ and *supported* on $F \subseteq \mathbb{R}^d$ if $\mathrm{supp}(\mu) \subseteq F$. Straight from the definition one has, for $\mu$ supported on $F$,

$$\dim_{\mathrm{H}} \mu \leqslant \dim_{\mathrm{H}} F.$$

In fact, for Borel sets $F$,

$$\dim_{\mathrm{H}} F = \sup\{\dim_{\mathrm{H}} \mu : \ \mathrm{supp}(\mu) \subseteq F\}. \tag{1.3}$$

This follows by finding compact subsets $E \subseteq F$ with positive and finite $s$-dimensional Hausdorff measure for all $s < \dim_{\mathrm{H}} F$, see [70, Theorem 4.10]. Similarly, the *(lower) packing dimension* of a Borel measure $\mu$ on $\mathbb{R}^d$ is

$$\dim_{\mathrm{P}} \mu = \inf\{\dim_{\mathrm{P}} E : \mu(E) > 0\}.$$

Again one has, for $\mu$ supported on $F$,

$$\dim_{\mathrm{P}} \mu \leqslant \dim_{\mathrm{P}} F$$

and, for Borel sets $F$,

$$\dim_{\mathrm{P}} F = \sup\{\dim_{\mathrm{P}} \mu : \ \mathrm{supp}(\mu) \subseteq F\}.$$

This final result follows by a similar approach, this time due to Joyce and Preiss [149]. The box dimension of a measure is a less well-developed concept. We formulate a definition in Section 4.2 following [75], which is partially motivated by the Assouad spectrum, see Section 3.3.

# 2
# The Assouad dimension

In this chapter we define the Assouad dimension, which is the central notion of the book. We discuss its origins in Section 2.3 and establish many of its basic properties in Section 2.4 such as stability under Lipschitz mappings and monotonicity. These are compared with the basic properties of the Hausdorff and box dimensions.

## 2.1 The Assouad dimension and a simple example

If the Hausdorff dimension provides fine, but global, geometric information, then the Assouad dimension provides coarse, but local, geometric information. The key difference between the Assouad dimension and the dimensions discussed in Chapter 1 is that only a small part of the set is considered at any one time. This is what gives it its 'local quality' and what leads to many of its interesting features, see Figure 2.1

The *Assouad dimension* of a non-empty set $F \subseteq \mathbb{R}^d$ is defined by

$$\dim_{\mathrm{A}} F = \inf \Bigg\{ \alpha : \text{ there exists a constant } C > 0 \text{ such that, for all } 0 < r < R \text{ and } x \in F, \; N_r\big(B(x,R) \cap F\big) \leqslant C\left(\frac{R}{r}\right)^{\alpha} \Bigg\}.$$

Recall that $N_r(E)$ is the smallest number of open sets required for an $r$-cover of a bounded set $E$. Note that we can replace $N_r$ in the definition of the Assouad dimension with any of the standard covering or packing functions, see the discussion in Chapter 1, and still obtain the same value for the dimension. For example, $N_r(E)$ could denote the number of cubes in an $r$-mesh oriented at the origin which intersect $E$ or the maximum cardinality of an $r$-separated subset of $E$. We also obtain an equivalent definition if the ball $B(x, R)$ is taken to be open or closed (although we usually think of it as being closed) or if it is replaced by a cube of sidelength $R$ centred at $x$.

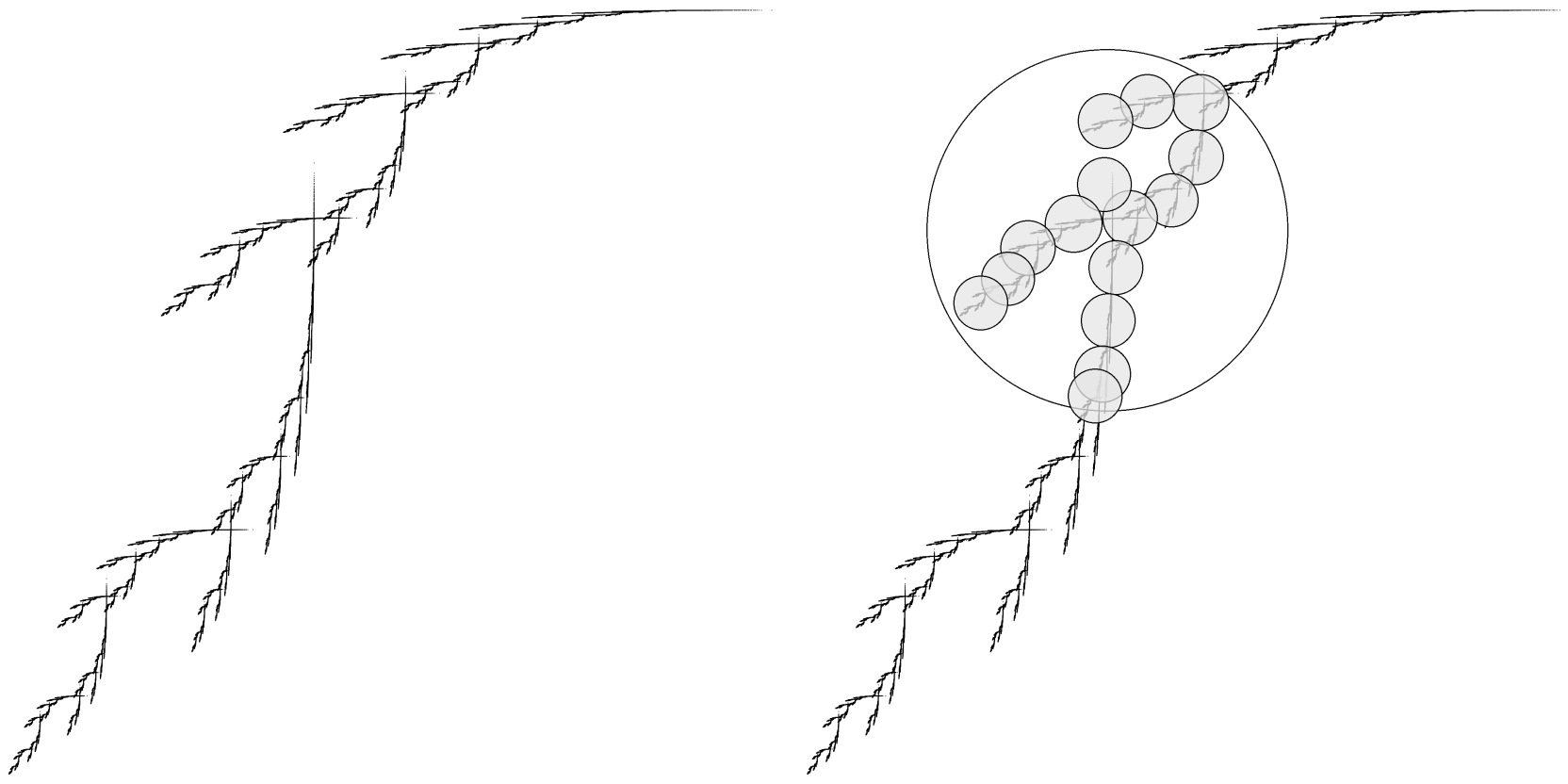

Figure 2.1 Left: the self-affine set from Figure 1.2. Right: A covering of a ball intersected with a 'thick part' of the self-affine set from Figure 1.2 using balls of smaller radii. Understanding how large such covers have to be leads to calculation of the Assouad dimension.

Before we move on let us consider a simple but fundamental example, see Figure 2.1. This example serves to demonstrate how inhomogeneity across a set can cause the box and Assouad (and Hausdorff) dimensions to be distinct.

**Theorem 2.1.1** For $F = \{0\} \cup \{1/n : n \in \mathbb{N}\}$,

$$\dim_{\mathrm{A}} F = 1,$$

$$\dim_{\mathrm{B}} F = 1/2,$$

and

$$\dim_{\mathrm{H}} F = 0.$$

*Proof* To prove that $\dim_{\mathrm{A}} F = 1$ it suffices to find a constant $c > 0$ and a sequence of points $x_n \in F$ and scales $0 < r_n < R_n$ such that $R_n/r_n \to \infty$ and for all $n$

$$N_{r_n}(B(x_n, R_n) \cap F) \geqslant cR_n/r_n.$$

To this end, let $x_n = 0$, $R_n = 1/n$ and $r_n = 1/n^2$. For $k \geqslant n$,

$$1/k - 1/(k+1) = \frac{1}{k(k+1)} < 1/n^2 = r_n$$

and so finding an $r_n$-cover of $B(0, R_n) \cap F$ is at least as hard as finding an $r_n$-cover of the interval $[0, R_n]$. More precisely,

$$N_{r_n}(B(0, R_n) \cap F) \geqslant 2^{-1}n = 2^{-1}R_n/r_n \tag{2.1}$$

as required.

To prove that $\dim_B F = 1/2$ we must control $N_r(F)$ from above and below for all small $r > 0$. Given $r \in (0, 1)$, let $n$ be the largest integer satisfying $1/(n(n+1)) > r$. It follows that any $r$-cover of $F$ must contain at least one set for each $m \leqslant n$ and, similar to (2.1), covering $[0, 1/n] \cap F$ is at least as hard as covering $[0, 1/n]$. Therefore

$$N_r(F) \geqslant N_r([0, 1/n] \cap F) \geqslant 2^{-1}n \geqslant 4^{-1}r^{-1/2}$$

and

$$N_r(F) \leqslant N_r([0, 1/n]) + n \leqslant 3r^{-1/2}$$

and we have thus bounded $N_r(F)$ above and below by an expression of the form $cr^s$ for all small $r$. This is sufficient to prove that $\dim_B F = s$ noting that the precise constants $c$ are not important.

Finally, for the Hausdorff dimension we note that $\dim_H X = 0$ for *any* countable set, which makes it rather different from the box and Assouad dimension. Write $X = \{x_1, x_2, \dots\} \subseteq \mathbb{R}^d$ and, given $r, \varepsilon > 0$, let $U_i = B(x_i, 2^{-i}r)$ for $i \geqslant 1$. Since $\{U_i\}_i$ is an $r$-cover of $X$,

$$\mathcal{H}_r^\varepsilon(X) \leqslant (2r)^\varepsilon \sum_{i \geqslant 1} 2^{-i\varepsilon} = (2r)^\varepsilon \frac{2^{-\varepsilon}}{1 - 2^{-\varepsilon}} \to 0$$

as $r \to 0$, which proves that $\mathcal{H}^\varepsilon(X) = 0$ for all $\varepsilon > 0$ and therefore $\dim_H X = 0$. □

Figure 2.2 The set $F = \{0\} \cup \{1/n : n \in \mathbb{N}\}$.

## 2.2 A word or two on the definition

The first thing one notices about the definition of the Assouad dimension is that it appears complicated. This is underlined by the fact that we required three lines to properly express the definition![1] However, despite its apparent complexity, the underpinning idea is very simple. In fact, here are five distinct ways of defining the Assouad dimension in one line:

$$\dim_{\mathrm{A}} F = \sup \{\dim_{\mathrm{H}} E \ : \ E \text{ is a weak tangent of } F\}$$

$$\dim_{\mathrm{A}} F = \inf \{\alpha \ : \ F \text{ is } \alpha\text{-homogeneous}\}$$

$$\dim_{\mathrm{A}} F = \inf \{\dim_{\mathrm{A}} \mu \ : \ \mu \text{ is a Borel measure supported on } F\}$$

$$\dim_{\mathrm{A}} F = \inf \left\{ \alpha : (\exists C > 0) \ (\forall 0 < r < R) \ \sup_{x \in F} N_r\big(B(x,R) \cap F\big) \leqslant C\left(\frac{R}{r}\right)^{\alpha} \right\}$$

$$\dim_{\mathrm{A}} F = \inf \left\{ \alpha : \text{there exists } C > 0 \text{ such that, for all } 0 < r < R \text{ and } x \in F,\ N_r(B(x,R) \cap F) \leqslant C\left(\frac{R}{r}\right)^{\alpha} \right\}.$$

Humour aside, each of these expressions demonstrates an important principle. The first is a beautiful fact, which was first explicitly stated by Käenmäki, Ojala and Rossi [154], but has origins in the work of Furstenberg [114, 115, 116] and Bishop and Peres [34]. It states that one may express the Assouad dimension entirely in terms of the Hausdorff dimension and entirely at the level of tangents. See Chapter 5 for more details on this important result, especially Theorem 5.1.3. The second expression is how one often explains the notion to an inexperienced party: *suppose we can cover any $R$-ball in a set $F$ with at most a constant multiple of $(R/r)^{\alpha}$ many $r$-balls. This is natural since if $F$ is $\mathbb{R}^d$ then one readily covers any $R$-ball by at most $5 \cdot 2^d (R/r)^d$ many $r$-balls. We call a set with this property $\alpha$-homogeneous and clearly this becomes harder to satisfy as $\alpha$ decreases. The Assouad dimension is simply the infimum of $\alpha$ such that the set is $\alpha$-homogeneous.* The third expression brings in the natural interplay between dimensions of sets and dimensions of measures. It is well known that the Hausdorff dimension of a set is the *supremum* of the Hausdorff dimension of Borel measures which it supports, recall (1.3). This duality between sets and measures is a powerful concept in dimension theory and the fact that one has a parallel dualism for Assouad dimension is appealing — and turns out to be useful as well. Of course, it depends on how one defines the Assouad dimension of a measure. This concept, initially called the upper regularity dimension, was introduced by Käenmäki, Lehrbäck and Vuorinen [153], and we will discuss it in detail in Chapter 4, especially Theorem 4.1.3. Hidden in this statement is also the seminal fact that a doubling metric space always carries a doubling measure. The fourth expression demonstrates the power of good and efficient notation and the fifth expression is illustrative of the fact that any mathematical statement, theorem, or proof, can be reduced to one line if one has a long enough line or a small enough font!

[1] "The problem with the Assouad dimension is that it is impossible to write the definition in one line." - Nick Sharples (in jest), Manchester Dynamics Seminar, October 2015.

We conclude this section with some technical remarks on the definition. In several instances in the literature, for example [219], the following subtly different definition of the Assouad dimension is given:

$$\inf \left\{ \alpha : \text{ there exist constants } C > 0, \rho > 0 \text{ such that, for all } 0 < r < R < \rho \text{ and } x \in F, \quad N_r\big(B(x,R) \cap F\big) \leqslant C\left(\frac{R}{r}\right)^{\alpha} \right\}. \tag{2.2}$$

It is easy to check that this definition and our definition coincide for all bounded $F \subseteq \mathbb{R}^d$ but they may differ drastically for unbounded sets.[2] For example, our definition assigns the integers dimension 1, whereas the alternative definition (2.2) assigns them dimension 0. Whether one wants the integers to have dimension 1 or 0 is perhaps a matter of taste but at the core of the decision is whether you want to look on large *and* small scales for the 'thickest part' of the space, or whether you only want to consider behaviour on small scales. We propose that the integers should have Assouad dimension 1 and therefore reject the alternative definition (2.2). This has several theoretical advantages, such as allowing the Assouad dimension to be invariant under Möbius transformations and inversions such as $z \mapsto 1/z$, see Section 9.3. It also means we do not have to specify whether tangents are obtained by 'zooming out' or 'zooming in' on the fractal. Moreover, the definition we use, see the beginning of Chapter 2, was also adopted by Robinson [240] and Mackay and Tyson [195].

In an attempt to simplify the definition, several people have asked us if the constant $C$ is really necessary. Specifically, does it not get dwarfed by $(R/r)^{\alpha-\dim_A F}$ for all $\alpha > \dim_A F$ since we should clearly only be interested in sequences of scales such that $R/r \to \infty$? However, if the constant $C$ is omitted from the definition, then one does not obtain a sensible notion of dimension. In particular, the finite set $\{0, 1\}$ would have dimension 1 since if we choose $x = 0$, $R = 1$ and $r = 1/2$, then $N_r\big(B(x,R) \cap F\big) = 2 = (R/r)$ and so we cannot choose any $\alpha < 1$ to satisfy the required condition.

[2] It was Chris Miller who first drew my attention to this inconsistency back in 2014.

Finally, we point out that the upper box dimension may be expressed in a similar way to the Assouad dimension, which partially serves to motivate both definitions and also highlights the potential for interplay between the two notions. Indeed, for bounded $F \subseteq \mathbb{R}^d$,

$$\overline{\dim}_{\mathrm{B}} F \;=\; \inf\Bigg\{ \alpha : \text{ there exists a constant } C > 0 \text{ such that, for all } 0 < r < 1 \text{ and } x \in F,\; N_r\big(B(x,1) \cap F\big) \;\leqslant\; C\left(\frac{1}{r}\right)^{\alpha} \Bigg\}.$$

The reason that one may simplify this expression (for example, the constant $C$ can be removed) is that there is only one scale involved. Therefore, the upper box dimension may be expressed as an upper limit as this single scale tends to 0.

## 2.3 Some history

Although interesting in its own right, and as a notion alongside Hausdorff and box dimension in fractal geometry, the Assouad dimension first came to prominence in other fields. Perhaps most notably, it plays an important role in quasi-conformal geometry, see [130, 192, 195], and embedding theory, see [222, 221, 240, 6, 7]. However, there has been an explosion of activity concerning the Assouad dimension in the fractal geometry literature in the last few years. We put this down to the publication of the books by Robinson [240] and Mackay and Tyson [195], as well as interesting research articles published around the same time which consider the Assouad dimension in classically fractal settings, such as [192, 194, 219, 153, 88, 115].

The Assouad dimension takes its name from French mathematician Patrice Assouad, who used the concept extensively during his doctoral studies in the 1970s as a tool in embedding problems, see [6, 7, 8]. This work culminated with the seminal *Assouad Embedding Theorem*, see Theorem 13.1.3, which served to motivate the Assouad dimension and popularise it throughout the mathematics community. However, the notion did not quite originate with Assouad and, in fact, appeared (independently) at least three times before Assouad's work in the 1970s.

The earliest reference we are aware of is Bouligand's 1928 paper [35], which is best known for the origin of the box dimension, or 'Minkowski-Bouligand dimension'. Discussion here is a little vague, but certainly there is an attempt to consider the volume of the $r$-neighbourhood of $B(x, R) \cap F$ for sets $F \subseteq \mathbb{R}^d$. Bouligand observes that this may lead to a different notion of dimension than if the volume of the $r$-neighbourhood of the *whole* of $F$ is considered (this is the box dimension) — not least since this dimension would depend on the point $x$. Bouligand writes[3] "One must thus give, instead of a dimensional number, a dimensional order, and also give to the notion of dimensional order a local character". He also referred to sets as *isodimensional* if the 'local dimensional order' is the same at all points. In modern day language we might think of this as sets with equal Assouad and lower dimension.

The next appearance we are aware of is in Larman's paper [179] where it is referred to as the *dimensional number* and denoted by dim-$n$. Among many other things, Larman's work contains the result that, for fixed $p > 0$, $\dim_{\mathrm{A}}\{1/n^p : n \in \mathbb{N}\} = 1$ and the fact that the lower dimension (Larman's *minimal dimensional number* $m$.dim-$n$) of a compact set is bounded above by its Hausdorff dimension, see Theorem 3.4.3.

The Assouad dimension also appears implicitly in Furstenberg's early work [114]. This was made more explicit in a celebrated article from 2008 [115], see also [116]. We briefly describe some of the results from [115] to highlight the connection with Assouad dimension. Let $F \subseteq [0, 1]^d$ be a non-empty compact set and say that a non-empty compact set $E \subseteq [0, 1]^d$ is a *miniset* of $F$ if, for some $c \geqslant 1$ and $t \in \mathbb{R}^d$, $E \subseteq cF + t$. Moreover, $E$ is called a *microset* of $F$ if it is a limit of a sequence of minisets of $F$ in the Hausdorff metric. A family $\mathcal{G}$ of non-empty compact subsets of $[0, 1]^d$ is called a *gallery* if it satisfies:

(i) $\mathcal{G}$ is closed with respect to the Hausdorff metric, see (5.1),
(ii) for each $E \in \mathcal{G}$, every miniset of $E$ is also in $\mathcal{G}$.

Given a non-empty compact set $F \subseteq [0, 1]^d$, the collection of all microsets of $F$ forms a gallery which Furstenberg denoted by $\mathcal{G}_F$. Then the $*$-dimension of $F$ is defined by

$$\dim^* F = \sup\{\dim_{\mathrm{H}} E : E \in \mathcal{G}_F\}.$$

[3] I am quoting from the English translation of [35] found in [61]. The translation is due to Ilan Vardi.

In [115, Theorem 5.1], Furstenberg proved the following.

**Theorem 2.3.1** Let $\mathcal{G}$ be a gallery. Let

$$\Delta(\mathcal{G}) = \limsup_{k\to\infty} \frac{1}{k} \log_2 \left( \sup_{X\in\mathcal{G}} \#\{Q \in \mathcal{D}_k : X \cap Q \neq \varnothing\} \right),$$

where $\mathcal{D}_k$ denotes the collection of half-open dyadic intervals of side length $2^{-k}$ and $\log_2$ is the base-2 logarithm. Then there exists a set $A \in \mathcal{G}$ such that

$$\dim_{\mathrm{H}} A = \Delta(\mathcal{G}).$$

It is easy to see that $\dim_{\mathrm{H}} X \leqslant \dim_{\mathrm{A}} X \leqslant \Delta(\mathcal{G})$ for any $X \in \mathcal{G}$, and therefore, for any gallery $\mathcal{G}$,

$$\sup\{\dim_{\mathrm{H}} X : X \in \mathcal{G}\} = \sup\{\dim_{\mathrm{A}} X : X \in \mathcal{G}\} = \Delta(\mathcal{G}),$$

and, moreover, both these suprema are obtained. In particular, for any non-empty compact $F$

$$\dim_{\mathrm{A}} F = \dim^* F = \Delta(\mathcal{G}_F)$$

and there exists a microset $E$ of $F$ such that $\dim_{\mathrm{H}} E = \dim_{\mathrm{A}} F$. Using our terminology, the gallery $\mathcal{G}_F$ of a compact set $F$ is the set of all (subsets of) weak tangents and the coincidence of the $*$-dimension and the Assouad dimension is manifest in Theorem 5.1.3. See Section 5 for more details on this topic.

## 2.4 Basic properties: *the greatest of all dimensions*

Following [70], in this section we describe some basic properties which we might hope for a 'dimension' to satisfy and then decide which of these properties are satisfied by the Assouad dimension. To put these results in context we will frequently compare the behaviour of the Assouad dimension with that of the box, packing and Hausdorff dimensions.

The following is a list of basic properties which a given notion of dimension dim may satisfy:

- *Monotonicity:* dim is said to be monotone if $E \subseteq F \Rightarrow \dim E \leqslant \dim F$ for all $E, F \subseteq \mathbb{R}^d$.
- *Finite stability:* dim is said to be finitely stable if $\dim(E \cup F) = \max\{\dim E, \dim F\}$ for all $E, F \subseteq \mathbb{R}^d$.
- *Countable stability:* dim is said to be countably stable if $\dim \bigcup_i F_i = \sup_i \dim F_i$ for all countable collections $\{F_i\}_i$ of subsets of $\mathbb{R}^d$.
- *Lipschitz stability:* dim is said to be Lipschitz stable if $\dim T(F) \leqslant \dim F$ for all $F \subseteq \mathbb{R}^d$ and all Lipschitz maps $T$ on $\mathbb{R}^d$. Recall that a map $T : F \to \mathbb{R}^d$ is *Lipschitz* if there exists a constant $K > 0$ such that, for all $x, y \in F$,
$$|T(x) - T(y)| \leqslant K|x - y|.$$
- *Bi-Lipschitz stability:* dim is said to be stable under bi-Lipschitz maps if $\dim T(F) = \dim F$ for all $F \subseteq \mathbb{R}^d$ and all bi-Lipschitz maps $T$ on $F$. Recall that a map $T : F \to \mathbb{R}^d$ is *bi-Lipschitz* if it is Lipschitz with a Lipschitz inverse (defined on $T(F)$).
- *Stability under closure:* dim is said to be stable under closure if $\dim F = \dim \overline{F}$ for all $F \subseteq \mathbb{R}^d$. Here and throughout $\overline{F}$ denotes the closure of a set $F$, which is the intersection of all closed sets containing $F$.
- *Open set property:* dim is said to satisfy the open set property if for any bounded open set $U \subseteq \mathbb{R}^d$, $\dim U = d$.

The proof of the following lemma is left as an exercise, although proofs can be found in [192]. The proof is straightforward, but the unfamiliar reader will find it useful to prove it carefully, paying particular attention to the precise definitions involved.

**Lemma 2.4.1** The Assouad dimension is monotone, finitely stable, stable under closure and satisfies the open set property.

The fact that the Assouad dimension is stable under closure and satisfies the open set property immediately shows that the Assouad dimension is *not* countably stable. For example,

$$\dim_A([0,1] \cap \mathbb{Q}) = \dim_A[0,1] = \dim_A(0,1) = 1$$

but the set $[0,1] \cap \mathbb{Q}$ is a countable union of singletons, each of which has Assouad dimension 0.

As an example, and partly to warm us up for future sections, we include a full proof that the Assouad dimension is stable under bi-Lipschitz maps. There will be examples later in the book which demonstrate that the Assouad dimension is *not* stable under Lipschitz maps, see Theorem 3.4.12, Corollary 7.2.2 and also [192, Example A.6 2].

**Lemma 2.4.2** If $F \subseteq \mathbb{R}^d$ is non-empty and $T : \mathbb{R}^d \to \mathbb{R}^d$ is bi-Lipschitz, then $\dim_A T(F) = \dim_A F$.

*Proof* Since $T$ is bi-Lipschitz, we can find a constant $K > 1$ such that for all $x, y \in \mathbb{R}^d$

$$K^{-1}|x-y| \leqslant |T(x)-T(y)| \leqslant K|x-y|. \tag{2.3}$$

Let $s > \dim_A F$, and fix an arbitrary point $x \in T(F)$, and arbitrary scales $0 < r < R$. In search of an $r$-cover of $B(x,R) \cap T(F)$, consider $T^{-1}(B(x,R) \cap T(F))$ which is necessarily contained in the set $B(T^{-1}(x), RK) \cap F$ by (2.3). By the definition of the Assouad dimension of $F$ there exists a constant $C > 0$, depending only on $s$ and $F$, such that we may cover $B(T^{-1}(x), RK) \cap F$ by fewer than

$$C\left(\frac{RK}{rK^{-1}}\right)^s = CK^{2s}\left(\frac{R}{r}\right)^s$$

balls of diameter $rK^{-1}$. Denote such an optimal cover by $\{B_i\}_i$ and note that $\{T(B_i)\}_i$ is an $r$-cover of $B(x,R) \cap T(F)$ and therefore

$$N_r(B(x,R) \cap T(F)) \leqslant CK^{2s}\left(\frac{R}{r}\right)^s$$

which proves that $\dim_A T(F) \leqslant s$ and since $s > \dim_A F$ was chosen freely we can upgrade this bound to $\dim_A T(F) \leqslant \dim_A F$. Finally, note that we can obtain the reverse inequality by using $T^{-1}$ in place of $T$ and interchanging the roles of $F$ and $T(F)$. □

The fact that Assouad dimension is not stable under Lipschitz maps is already manifest in the above argument. It was essential in the proof of Lemma 2.4.2 that $T$ was bi-Lipschitz, even to establish the bound in one direction. On the other hand, to prove that $\overline{\dim}_B T(F) \leqslant \overline{\dim}_B F$ we only require that $T$ is Lipschitz.

Another important fact to keep in mind is the relationships between the different notions of dimension. It turns out that the Assouad dimension is the greatest of all dimensions (at least of all the dimensions we consider).[4]

[4] "The Assouad dimension is the largest reasonable dimension one can define using coverings" - Antti Käenmäki, St Andrews Analysis Seminar, May 2018.

| Property | $\dim_{\mathrm{H}}$ | $\dim_{\mathrm{P}}$ | $\underline{\dim}_{\mathrm{B}}$ | $\overline{\dim}_{\mathrm{B}}$ | $\dim_{\mathrm{A}}$ |
|---|---|---|---|---|---|
| Monotone | ✓ | ✓ | ✓ | ✓ | ✓ |
| Finitely stable | ✓ | ✓ | × | ✓ | ✓ |
| Countably stable | ✓ | ✓ | × | × | × |
| Lipschitz stable | ✓ | ✓ | ✓ | ✓ | × |
| Bi-Lipschitz stable | ✓ | ✓ | ✓ | ✓ | ✓ |
| Stable under closure | × | × | ✓ | ✓ | ✓ |
| Open set property | ✓ | ✓ | ✓ | ✓ | ✓ |

Table 2.1 *This table summarises which properties are satisfied by which dimensions. Proofs of these facts for the Hausdorff, packing and box dimensions, can be found in [70, Chapters 2-3]. Many basic properties of the Assouad dimension are considered in [192]. See Table 3.1 which includes a similar analysis for different notions of dimension.*

**Lemma 2.4.3** For any set $F \subseteq \mathbb{R}^d$

$$\dim_{\mathrm{H}} F \leqslant \dim_{\mathrm{P}} F \leqslant \dim_{\mathrm{A}} F$$

and, if $F$ is bounded, then

$$\begin{array}{ccccc} & & \dim_{\mathrm{P}} F & & \\ & \leqslant & & \leqslant & \\ \dim_{\mathrm{H}} F & & & & \overline{\dim}_{\mathrm{B}} F \quad \leqslant \dim_{\mathrm{A}} F. \\ & \leqslant & & \leqslant & \\ & & \underline{\dim}_{\mathrm{B}} F & & \end{array}$$

*Proof* The relationships between the Hausdorff, packing and box dimensions are proved in [70, Chapters 2-3]. As such we omit the proofs of these facts, but encourage the reader to consider them as exercises. The fact that, for bounded $F \subseteq \mathbb{R}^d$, $\overline{\dim}_{\mathrm{B}} F \leqslant \dim_{\mathrm{A}} F$ follows immediately from the definitions by choosing $R = |F|$. More precisely, for $\varepsilon > 0$ there exists $C \geqslant 1$ such that for all $0 < r < |F|$ and $x \in F$,

$$N_r(F) = N_r(B(x, |F|) \cap F) \leqslant C \left( \frac{|F|}{r} \right)^{\dim_{\mathrm{A}} F + \varepsilon}$$

which proves that $\overline{\dim}_{\mathrm{B}} F \leqslant \dim_{\mathrm{A}} F + \varepsilon$ and letting $\varepsilon \to 0$ yields the desired inequality. This inequality passes to packing dimension via (1.1)

since

$$\begin{aligned}\dim_{\mathrm{P}} F &= \inf\{\sup_i \overline{\dim}_{\mathrm{B}} F_i \ : \ F = \bigcup_i F_i, \ F_i \text{ bounded}\} \\ &\leqslant \inf\{\sup_i \dim_{\mathrm{A}} F_i \ : \ F = \bigcup_i F_i, \ F_i \text{ bounded}\} \\ &\leqslant \dim_{\mathrm{A}} F\end{aligned}$$

by monotonicity. □

There is no general relationship between the lower box and packing dimensions.

# 3
# Some variations on the Assouad dimension

In this chapter we introduce several variants of the Assouad dimension, which will also play a key role in this book. Most importantly, the lower dimension, see Section 3.1, is the natural dual to the Assouad dimension, and the Assouad spectrum, see Section 3.3, is a function which interpolates between the quasi-Assouad dimension and the upper box dimension. The quasi-Assouad dimension is closely related to the Assouad dimension and is introduced in Section 3.2. When we say the Assouad spectrum 'interpolates' we mean that it is a continuous function of $\theta \in (0, 1)$ which tends to the upper box dimension as $\theta \to 0$ and tends to the quasi-Assouad dimension as $\theta \to 1$. Therefore, the Assouad spectrum provides more information than the upper box and quasi-Assouad dimension considered in isolation and serves to explore the gap between these dimensions, should there be one. In Section 3.4 we develop some of the basic properties of these notions and in subsequent chapters we will attempt to present a unified theory of the Assouad dimension together with all of its natural variants whenever possible.

## 3.1 The lower dimension

An important theme in dimension theory is that dimensions often come in pairs. For example, the Hausdorff and packing dimension are a natural pair, as are the upper and lower box dimensions. See Section 10.1 for a particularly transparent example of the role of dimension pairs in the context of product sets. The natural partner of the Assouad dimension is the *lower dimension*, which was introduced by Larman [179] where it was called the *minimal dimensional number*.[1] We adopt the name 'lower dimension' following an important paper written by Bylund and Gudayol [40], which will be revisited in Section 4. Just as the Assouad dimension seeks to quantify the 'thickest' part of the set in question, the lower dimension identifies the 'thinnest' part, that is, the part of the set which is easiest to cover, see Figure 3.1. This leads to the lower dimension having some strange properties, such as failing to be monotonic.[2] The *lower dimension* of $F$ is defined by

$$\dim_{\mathrm{L}} F \;=\; \sup\Bigg\{\,\alpha:\ \text{there exists a constant } C>0 \text{ such that, for all } 0<r<R\leqslant |F| \text{ and } x\in F,\ N_r\big(B(x,R)\cap F\big) \;\geqslant\; C\left(\frac{R}{r}\right)^{\alpha}\Bigg\}.$$

Comparing this with the definition of the Assouad dimension, it is clear that the lower dimension is a dual notion. However, it is worth pointing out the one asymmetry in the definition which is the inclusion of the diameter $|F|$ (which may be infinity) as an upper bound on allowable scales $r, R$. Sometimes this is included in the definition of the Assouad dimension, but it leads to an equivalent definition, since for $R > |F|$ it is just as easy to cover $B(x, R) \cap F$ as it is to cover $B(x, |F|) \cap F$. It is crucial, however, to include this bound in the definition of the lower dimension as otherwise all bounded sets would have dimension 0.

If a set has some exceptional points around which it is distributed very sparsely in comparison with the rest of the set, then the lower dimension will reflect this. For instance, sets with isolated points always have lower dimension equal to 0. However, for sets with some degree of homogeneity, such as attractors of iterated function systems (IFSs) or dynamically invariant sets, the lower dimension is very suitable and reveals easily interpreted information about the set. See Theorems 6.3.1 and 7.4.2, as well as the papers [148, 260]. Since the Assouad and lower dimensions are extremal, the difference between them provides a crude measure of inhomogeneity. For example, self-similar sets satisfying the open set condition and, more generally, Ahlfors-regular sets, have equal Assouad and lower dimensions, indicating that these sets are as homogeneous as possible, see Section 6.4. However, for more complicated sets, such as

[1] It was referred to as the *lower Assouad dimension* by Käenmäki, Lehrbäck and Vuorinen [153] and Tuomas Sahlsten used to call it the *uniformity dimension*.
[2] 'Thickest' is monotonic but 'thinnest' is not.

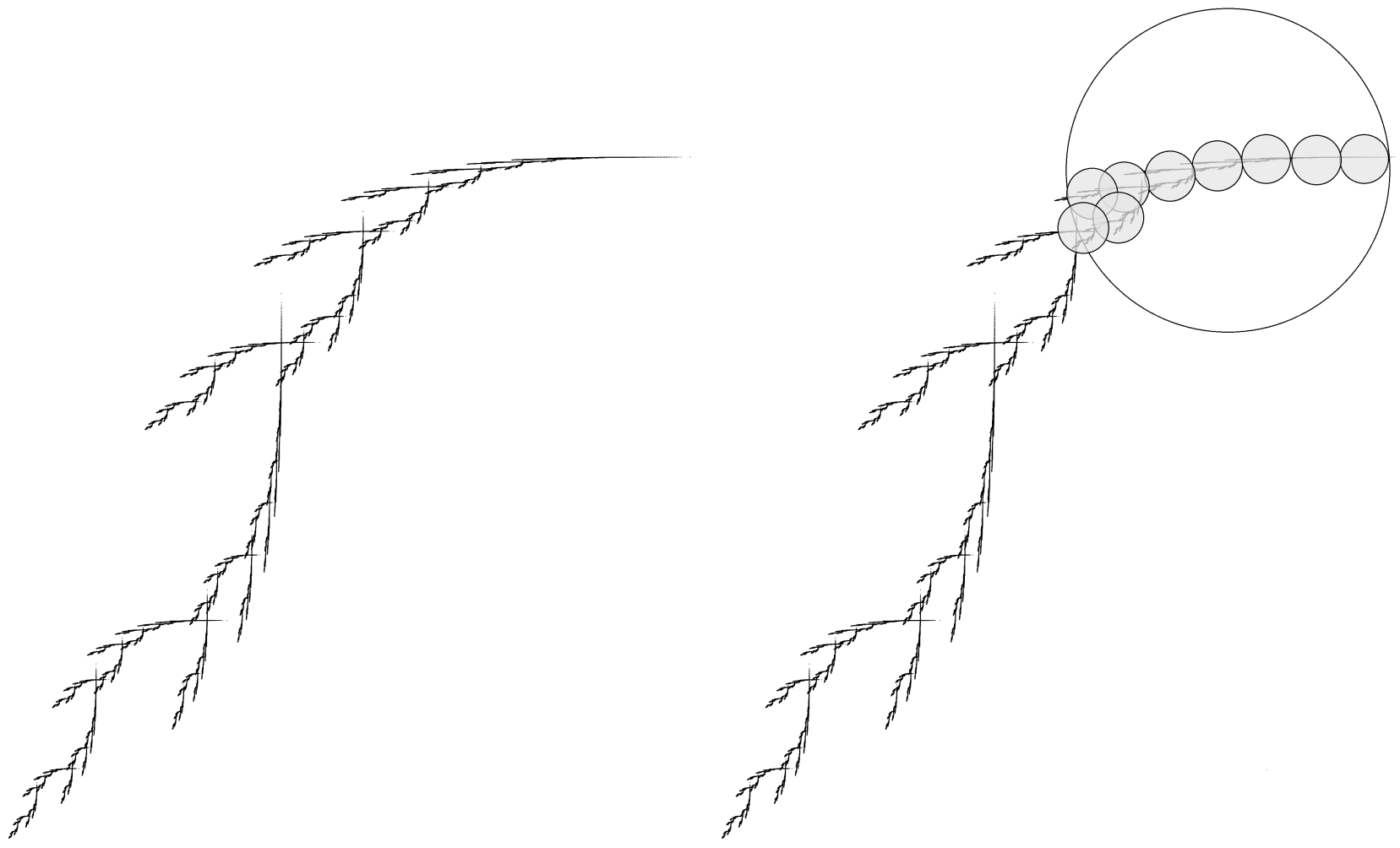

Figure 3.1 Left: the self-affine set from Figure 1.2. Right: a covering of a ball intersected with a 'thin part' of the self-affine set from Figure 1.2 using balls of smaller radii. Understanding how small such covers can be leads to calculation of the lower dimension.

*self-affine sets*, *self-similar sets with overlaps*, and *limit sets of Kleinian groups*, the quantities can be, and often are, different.

A natural trick which 'forces' the lower dimension to be monotonic is to define the *modified lower dimension* as

$$\dim_{\mathrm{ML}} F = \sup\{\dim_{\mathrm{L}} E : E \subseteq F\}.$$

This dimension has not yet received much attention in the literature. However, it turns out to have some interesting applications in, for example, Diophantine approximation, see Section 14.2.

## 3.2 The quasi-Assouad dimension

The *quasi-Assouad* dimension, introduced by Lü and Xi [191] much more recently than the Assouad and lower dimensions, is defined by

$\dim_{\mathrm{qA}} F = \lim_{\delta \to 0} h_F(\delta)$ where, for $\delta \in (0,1)$,

$$h_F(\delta) = \inf \Bigg\{ \alpha : \text{ there exists a constant } C > 0 \text{ such that, for all } 0 < r \leqslant R^{1+\delta} < 1 \text{ and } x \in F, \; N_r\big(B(x,R) \cap F\big) \leqslant C\left(\frac{R}{r}\right)^{\alpha} \Bigg\}.$$

The quasi-Assouad dimension leaves an 'exponential gap' between the scales $r$ and $R$, which in some settings can make a huge difference — see, for example, Section 9.4.

One of the motivations behind the definition of quasi-Assouad dimension is that it is preserved under quasi-Lipschitz maps. Given compact sets $X, Y \subseteq \mathbb{R}^d$, a bijection $T : X \to Y$ is called *quasi-Lipschitz* if

$$\frac{\log |T(x) - T(y)|}{\log |x - y|} \to 1$$

uniformly as $|x - y| \to 0$. The concept of quasi-Lipschitz equivalence was introduced in [276]. Clearly any bi-Lipschitz map is quasi-Lipschitz. Lü and Xi [191] proved that for any quasi-Lipschitz map $T$ and $F \subseteq \mathbb{R}^d$

$$\dim_{\mathrm{qA}} T(F) = \dim_{\mathrm{qA}} F.$$

This property is shared by the Hausdorff, box and packing dimensions, but *not* by the Assouad dimension. We will revisit this concept later, see Corollary 3.4.15. There is also an analogously defined quasi-lower dimension, which we will introduce formally below, see (3.3).

## 3.3 The Assouad spectrum

Motivated by the desire to obtain more nuanced information about the scaling structure of a fractal set, Fraser and Yu introduced the *Assouad spectrum* [111], see also the survey [94]. The idea is to understand which pairs of scales $r < R$ give rise to the extreme behaviour captured by the Assouad (and quasi-Assouad) dimension. To each $\theta \in (0,1)$, one associates the appropriate analogue of the Assouad dimension with the restriction that the two scales $r$ and $R$ used in the definition satisfy $R = r^{\theta}$. The resulting 'dimension spectrum' (as a function of $\theta$) thus gives finer geometric information regarding the scaling structure of the set and, in some precise sense, interpolates between the upper box dimension and the (quasi-)Assouad dimension. The Assouad spectrum is generally better behaved than the Assouad dimension and so this interpolation can lead to a clearer understanding of the Assouad dimension considered in isolation. More precisely, given $\theta \in (0,1)$ we define

$$\dim_{\mathrm{A}}^{\theta} F \;=\; \inf\Bigg\{ \alpha : \text{ there exists a constant } C>0 \text{ such that, for all } 0<r<1 \text{ and } x \in F, \; N_r\big(B(x,r^\theta)\cap F\big) \;\leqslant\; C\left(\frac{r^\theta}{r}\right)^{\alpha} \Bigg\},$$

and the dual notion

$$\dim_{\mathrm{L}}^{\theta} F \;=\; \sup\Bigg\{ \alpha : \text{ there exists a constant } C>0 \text{ such that, for all } 0<r<1 \text{ and } x \in F, \; N_r\big(B(x,r^\theta)\cap F\big) \;\geqslant\; C\left(\frac{r^\theta}{r}\right)^{\alpha} \Bigg\}.$$

We refer to $\dim_{\mathrm{A}}^{\theta} F$ as the *Assouad spectrum* and $\dim_{\mathrm{L}}^{\theta} F$ as the *lower spectrum* and wish to understand these concepts as functions of $\theta$. Many properties of these dimension spectra were considered in [111, 95] and further examples were considered in [112, 108]. Sometimes it is convenient to consider scales expressed as $R^{1/\theta} < R$ rather than $r < r^\theta$, that is, the smaller scale as a function of the bigger scale rather than the bigger scale as a function of the smaller scale. This was how the definition was originally expressed in [111] but the definition we give here is clearly equivalent. Since we require $r < r^\theta$ for $\theta \in (0,1)$, it is necessary to restrict our attention to small scale behaviour, that is, $r < 1$, which is another notable difference between the Assouad spectrum and quasi-Assouad dimension in comparison with the Assouad dimension. As discussed above, this makes no difference for bounded sets, but for unbounded sets it leads to rather different behaviour. For example, for all $\theta \in (0,1)$, $\dim_{\mathrm{A}}^{\theta} \mathbb{Z} = \dim_{\mathrm{qA}} \mathbb{Z} = 0 < 1 = \dim_{\mathrm{A}} \mathbb{Z}$.

We will study many fundamental properties of the Assouad and lower spectra in the following sections, but we close this introductory discussion by emphasising that these spectra lie in between the Assouad and lower dimensions and the box dimensions. That is, for bounded $F \subseteq \mathbb{R}^d$ and all $\theta \in (0,1)$,

$$\dim_{\mathrm{L}} F \leqslant \dim_{\mathrm{L}}^{\theta} F \leqslant \underline{\dim}_{\mathrm{B}} F \leqslant \overline{\dim}_{\mathrm{B}} F \leqslant \dim_{\mathrm{A}}^{\theta} F \leqslant \dim_{\mathrm{A}} F. \tag{3.1}$$

See Lemma 3.4.1 and Lemma 3.4.4 for more refined estimates.

### 3.3.1 First observations on the Assouad spectrum

One of the most important and useful properties of the Assouad and lower spectra is that they are continuous in $\theta \in (0,1)$. This will follow from the following theorem, which has other useful consequences as well. This proposition unifies [111, Propositions 3.4 and 3.7], see also [42, Proposition 2.1].

**Theorem 3.3.1** For any set $F \subseteq \mathbb{R}^d$ and $0 < \theta_1 < \theta_2 < 1$,

$$\left(\frac{1-\theta_2}{1-\theta_1}\right) \dim_{\mathrm{A}}^{\theta_2} F \leqslant \dim_{\mathrm{A}}^{\theta_1} F$$

$$\leqslant \left(\frac{1-\theta_2}{1-\theta_1}\right) \dim_{\mathrm{A}}^{\theta_2} F + \left(\frac{\theta_2-\theta_1}{1-\theta_1}\right) \dim_{\mathrm{A}}^{\theta_1/\theta_2} F,$$

and

$$\left(\frac{1-\theta_2}{1-\theta_1}\right) \dim_{\mathrm{L}}^{\theta_2} F \leqslant \dim_{\mathrm{L}}^{\theta_1} F$$

$$\leqslant \left(\frac{\theta_1-\theta_1\theta_2}{\theta_2-\theta_1\theta_2}\right) \dim_{\mathrm{L}}^{\theta_2} F + \left(\frac{\theta_2-\theta_1}{\theta_2-\theta_1\theta_2}\right) \dim_{\mathrm{A}}^{\theta_1/\theta_2} F.$$

*Proof* First consider the Assouad spectrum. Let $0 < \theta_1 < \theta_2 < 1$,

$$s < \dim_{\mathrm{A}}^{\theta_2} F < t$$

and

$$\dim_{\mathrm{A}}^{\theta_1/\theta_2} F < u.$$

For $r \in (0,1)$, we clearly have

$$\sup_{x \in F} N_r(B(x, r^{\theta_2}) \cap F) \leqslant \sup_{x \in F} N_r(B(x, r^{\theta_1}) \cap F).$$

Therefore, by applying the definition of the Assouad spectrum directly, we can find arbitrarily small $r > 0$ satisfying

$$\sup_{x \in F} N_r(B(x, r^{\theta_1}) \cap F) \geqslant \left(\frac{r^{\theta_2}}{r}\right)^s = \left(r^{\theta_1-1}\right)^{s\left(\frac{1-\theta_2}{1-\theta_1}\right)}$$

which proves

$$\dim_{\mathrm{A}}^{\theta_1} F \geqslant s\left(\frac{1-\theta_2}{1-\theta_1}\right)$$

and letting $s \to \dim_{\mathrm{A}}^{\theta_2} F$ yields the desired lower bound. Moving on to the upper bound for the Assouad spectrum, for $r \in (0,1)$,

$$\sup_{x \in F} N_r(B(x, r^{\theta_1}) \cap F) \leqslant \sup_{x \in F} N_{r^{\theta_2}}(B(x, r^{\theta_1}) \cap F) \sup_{x \in F} N_r(B(x, r^{\theta_2}) \cap F)$$

which can be seen by first covering $B(x, r^{\theta_1}) \cap F$ by $r^{\theta_2}$-balls and then covering each of these by $r$-balls. Therefore, for all $r \in (0,1)$,

$$\sup_{x \in F} N_r(B(x, r^{\theta_1}) \cap F) \leqslant \left(\frac{r^{\theta_1}}{r^{\theta_2}}\right)^u \left(\frac{r^{\theta_2}}{r}\right)^t = \left(r^{\theta_1 - 1}\right)^{u\left(\frac{\theta_2 - \theta_1}{1 - \theta_1}\right) + t\left(\frac{1 - \theta_2}{1 - \theta_1}\right)}$$

which proves

$$\dim_{\mathrm{A}}^{\theta_1} F \leqslant u\left(\frac{\theta_2 - \theta_1}{1 - \theta_1}\right) + t\left(\frac{1 - \theta_2}{1 - \theta_1}\right)$$

and letting $t \to \dim_{\mathrm{A}}^{\theta_2} F$ and $u \to \dim_{\mathrm{A}}^{\theta_1/\theta_2} F$ yields the desired upper bound.

For the lower spectrum, let

$$s < \dim_{\mathrm{L}}^{\theta_2} F$$

and

$$v < \dim_{\mathrm{L}}^{\theta_1/\theta_2} F.$$

For $r \in (0,1)$,

$$\inf_{x \in F} N_r(B(x, r^{\theta_2}) \cap F) \leqslant \inf_{x \in F} N_r(B(x, r^{\theta_1}) \cap F).$$

Therefore, by applying the definition of the lower spectrum, for all $r \in (0,1)$,

$$\inf_{x \in F} N_r(B(x, r^{\theta_1}) \cap F) \geqslant \left(\frac{r^{\theta_2}}{r}\right)^s = \left(r^{\theta_1 - 1}\right)^{s\left(\frac{1 - \theta_2}{1 - \theta_1}\right)}$$

which proves

$$\dim_{\mathrm{L}}^{\theta_1} F \geqslant s\left(\frac{1 - \theta_2}{1 - \theta_1}\right)$$

and letting $s \to \dim_{\mathrm{L}}^{\theta_2} F$ yields the desired lower bound.

The upper bound is a little more complicated. Let

$$\dim_{\mathrm{L}}^{\theta_2} F < t$$

and

$$\dim_{\mathrm{A}}^{\theta_1/\theta_2} F < u.$$

Then, for $r \in (0,1)$,

$$\inf_{x \in F} N_r(B(x, r^{\theta_1}) \cap F)$$

$$\leqslant \inf_{x \in F} N_{r^{\theta_1/\theta_2}}(B(x, r^{\theta_1}) \cap F) \sup_{x \in F} N_r(B(x, r^{\theta_1/\theta_2}) \cap F)$$

which comes by first finding an $r^{\theta_1}$-ball which minimises

$$\inf_{x \in F} N_{r^{\theta_1/\theta_2}}(B(x, r^{\theta_1}) \cap F),$$

covering it efficiently by $r^{\theta_1/\theta_2}$-balls, and then covering each of these by $r$-balls. By applying the definition of the Assouad and lower spectrum, we can find arbitrarily small $r > 0$ satisfying

$$\inf_{x \in F} N_r(B(x, r^{\theta_1}) \cap F) \leqslant \left( \frac{r^{\theta_1}}{r^{\theta_1/\theta_2}} \right)^t \left( \frac{r^{\theta_1/\theta_2}}{r} \right)^u$$

$$= \left( \frac{r^{\theta_1}}{r} \right)^{t\left( \frac{\theta_1 - \theta_1/\theta_2}{\theta_1 - 1} \right) + u\left( \frac{\theta_1/\theta_2 - 1}{\theta_1 - 1} \right)}$$

which proves

$$\dim_{\mathrm{L}}^{\theta_1} F \leqslant t \left( \frac{\theta_1 - \theta_1\theta_2}{\theta_2 - \theta_1\theta_2} \right) + u \left( \frac{\theta_2 - \theta_1}{\theta_2 - \theta_1\theta_2} \right)$$

and letting $t \to \dim_{\mathrm{L}}^{\theta_2} F$ and $u \to \dim_{\mathrm{A}}^{\theta_1/\theta_2} F$ yields the desired upper bound. □

Using similar ideas, one may establish various other bounds for the spectra. We omit further details since the estimates above serve our purpose. See [111, Propositions 3.4 and 3.7] and, in particular, [42, Proposition 2.1] which establishes a better lower bound for the lower spectrum than the one we give here.

The most appealing thing about Theorem 3.3.1 is its applications. The first and most important of these is continuity of the spectra.

**Corollary 3.3.2** For a fixed non-empty set $F \subseteq \mathbb{R}^d$, the functions $\theta \mapsto \dim_{\mathrm{A}}^{\theta} F$ and $\theta \mapsto \dim_{\mathrm{L}}^{\theta} F$ are continuous in $\theta \in (0,1)$.

One immediately sees that for any set $F$ and all $\theta \in (0,1)$, $\dim_{\mathrm{A}}^{\theta} F \leqslant \dim_{\mathrm{qA}} F$. Another useful consequence of Theorem 3.3.1 is that if the Assouad spectrum reaches the quasi-Assouad dimension, then it is constant from that point on. We will see that the Assouad spectrum is not necessarily non-decreasing across the whole range, so this information is *a priori* useful.

**Corollary 3.3.3** Fix a non-empty set $F \subseteq \mathbb{R}^d$. If $\dim_{\mathrm{A}}^{\theta} F = \dim_{\mathrm{qA}} F$ for some $\theta \in (0,1)$, then

$$\dim_{\mathrm{A}}^{\theta'} F = \dim_{\mathrm{qA}} F$$

for all $\theta' \in [\theta, 1)$.

*Proof* This follows from the upper bound on $\dim_{\mathrm{A}}^{\theta_1} F$ from Theorem 3.3.1, by choosing $\theta_1 = \theta$, $\theta_2 = \theta'$ and bounding $\dim_{\mathrm{A}}^{\theta_1/\theta_2} F$ above by $\dim_{\mathrm{qA}} F$. □

As we will see, a typical behaviour for the Assouad spectrum is to be convex and increasing on $(0, \rho)$ and constantly equal to the quasi-Assouad dimension on $[\rho, 1)$ for some $\rho \in (0,1)$. Another useful consequence of Theorem 3.3.1 is that, if $\rho$ is known, then we get a simple lower bound for the increasing part of the spectrum. One reason we mention this is that $\rho$ turns up again in other contexts, see, for example, Section 17.7.

**Corollary 3.3.4** For a fixed non-empty set $F \subseteq \mathbb{R}^d$, suppose

$$\rho = \inf\{\theta \in (0,1) : \dim_{\mathrm{A}}^{\theta} F = \dim_{\mathrm{qA}} F\}$$

exists and $\rho \in (0,1)$. Then for $\theta \in (0, \rho)$

$$\dim_{\mathrm{A}}^{\theta} F \geqslant \left(\frac{1-\rho}{1-\theta}\right) \dim_{\mathrm{qA}} F.$$

*Proof* This follows from Theorem 3.3.1 with $\theta_1 = \theta < \theta_2 = \rho$. □

The spectra exhibit stronger regularity properties than just continuity. In fact, [111, Corollary 3.7] shows that the spectra are Lipschitz on every closed interval strictly contained in $(0,1)$ and are therefore differentiable almost everywhere by Rademacher's Theorem. That said, the set of points for which the spectra are not differentiable may be dense by the following theorem proved in [95, Theorem 2.5].

**Theorem 3.3.5** Let $f : [0,1] \to [0,1]$ be continuous, concave, non-decreasing and satisfy $f(0) > 0$ and $f(\theta) \leqslant f(0)(\theta + 1)$ for all $\theta \in [0,1]$. Then there exists a compact set $F \subseteq [0,1]$ such that

$$\dim_{\mathrm{A}}^{\theta} F = f(\theta)$$

for all $\theta \in (0,1)$.

The strategy to prove Theorem 3.3.5 is as follows. First argue that, given any $0 < s < t \leqslant 1$, there exists a compact set $E \subseteq [0,1]$ with

$$\dim_{\mathrm{A}}^{\theta} E = \min\left\{\frac{s}{1-\theta}, t\right\}.$$

The specific sets used in [95] are Moran constructions, see also [112]. With these examples in place, enumerate the rationals in $(0,1)$ by $q_i$ and let $E_i$ be the set $E$ above with $t = f(q_i)$ and $s = f(q_i)(1-q_i)$. Finally let

$$F = \overline{\bigcup_i E_i^*}$$

where $E_i^*$ is a scaled and translated copy of $E_i$ such that the sets $E_i^*$ are separated and $F$ is bounded. If $\dim_{\mathrm{A}}^{\theta}$ was countably stable, then we could conclude

$$\dim_{\mathrm{A}}^{\theta} F = \sup_i \dim_{\mathrm{A}}^{\theta} E_i = f(\theta).$$

However, $\dim_{\mathrm{A}}^{\theta}$ is *not* countably stable in general and so more care is needed. Essentially what is needed is that the sets $E_i^*$ are scaled such that their diameter decays very quickly. This means that, for a given pair of scales $r, r^{\theta}$, only one of the sets $E_i$ is 'seen'. For the details, see [95].

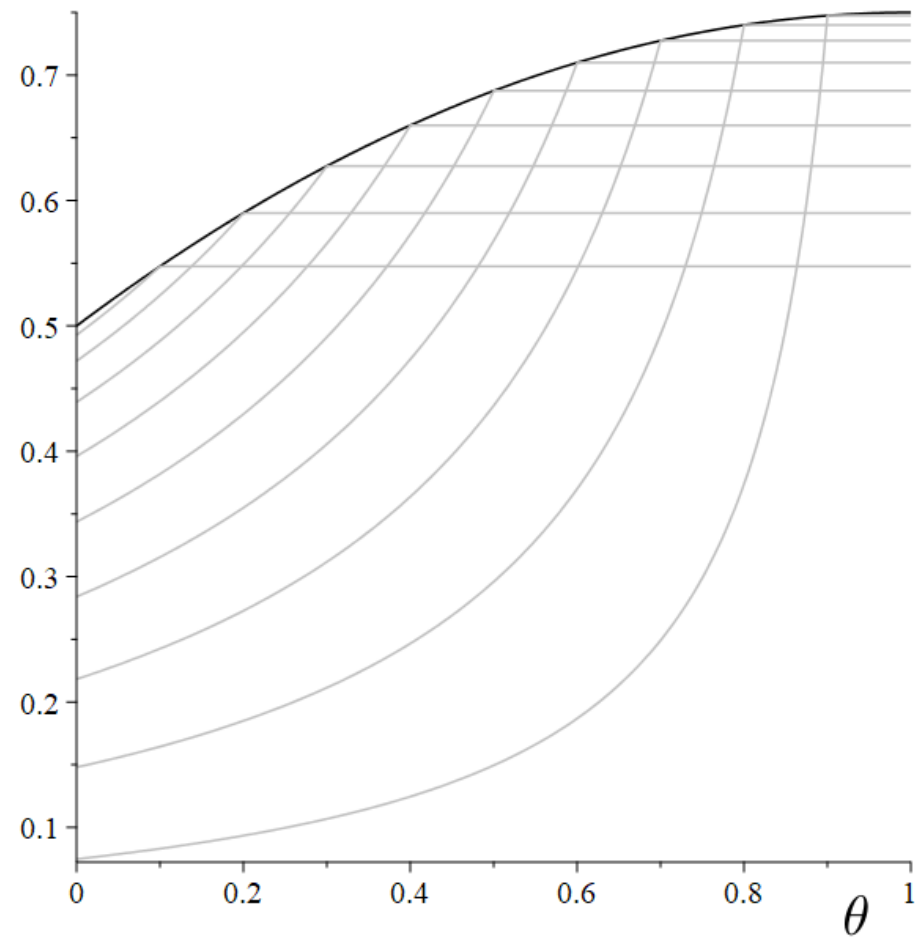

Figure 3.2 The sketch proof of Theorem 3.3.5 in action. The function $f(\theta) = (1+\theta)/2 - (\theta/2)^2$ (black) is approximated from below by 'ghost functions' $\dim_{\mathrm{A}}^{\theta} E_i$ (grey) for $q_i = i/10$ $(i = 1, \ldots, 9)$.

In many of the examples we consider, the spectra turn out to be piecewise convex analytic and constant in a neighbourhood of $\theta = 1$, but Theorem 3.3.5 says that these properties are not necessarily satisfied. Most examples also turn out to be monotonic, but non-monotonic examples are possible, see Section 3.4.4. An example was constructed in [111, Section 8] where the spectra switch between strictly increasing and strictly decreasing infinitely many times.

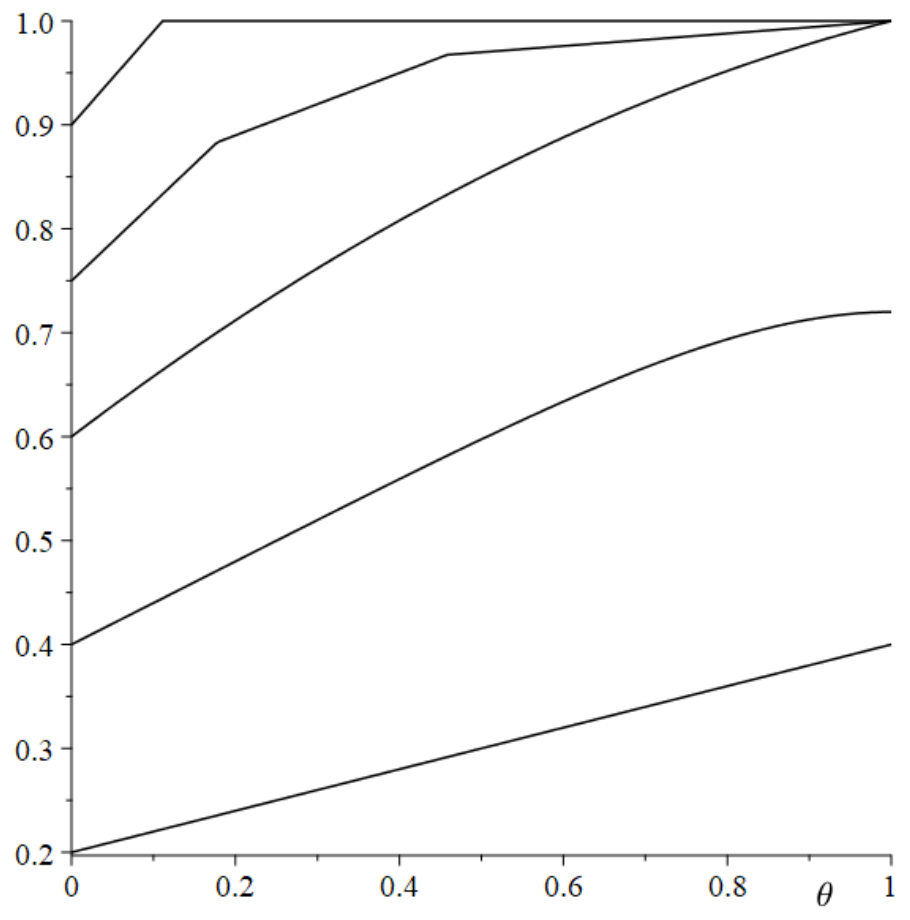

Figure 3.3 Five examples of Assouad spectra which may be constructed using Theorem 3.3.5. Starting from the top, we see the following features which are not seen elsewhere in the book: (i) maximum of affine and constant; (ii) piecewise affine with three pieces; (iii) strictly increasing, strictly concave and smooth; (iv) strictly increasing and differentiable but with derivative approaching 0 at $\theta = 1$; (v) affine but not constant.

### 3.3.2 A tale of two spectra

It is perhaps equally natural to consider different spectra where instead of fixing $r = R^{1/\theta}$, one simply requires $r \leqslant R^{1/\theta}$. This leads to what we will refer to as the *upper spectrum*[3], $\overline{\dim}_{\mathrm{A}}^{\theta} F$. More precisely, given

[3] During a research visit to Waterloo in 2018, Kathryn Hare, Kevin Hare, Sascha Troscheit and I briefly referred to $\dim_{\mathrm{A}}^{\theta} F$ as the 'Scottish spectrum' and to $\overline{\dim}_{\mathrm{A}}^{\theta} F$ as the 'Canadian spectrum'.

$\theta \in (0,1)$ we define

$$\overline{\dim}_{\mathrm{A}}^{\theta} F \;=\; \inf \Bigg\{ \alpha : \text{ there exists a constant } C > 0 \text{ such that, for all } 0 < r \leqslant R^{1/\theta} < R < 1 \text{ and } x \in F, \; N_r\big(B(x,R) \cap F\big) \;\leqslant\; C\left(\frac{R}{r}\right)^{\alpha} \Bigg\}$$

noting that $\overline{\dim}_{\mathrm{A}}^{\theta} F$ is equal to $h_F(\delta)$ in the definition of the quasi-Assouad dimension, for $\delta = 1/\theta - 1$.

It follows immediately from the definition that

$$\dim_{\mathrm{A}}^{\theta} F \leqslant \overline{\dim}_{\mathrm{A}}^{\theta} F \leqslant \dim_{\mathrm{qA}} F \leqslant \dim_{\mathrm{A}} F$$

and, moreover, the upper spectrum is non-decreasing in $\theta$, whereas the Assouad spectrum is not necessarily non-decreasing. In particular, the notions are distinct. However, it turns out that the upper spectrum is determined entirely by the Assouad spectrum. The following result was proved in [95] and we follow the proof given there.

**Theorem 3.3.6** Let $F \subseteq \mathbb{R}^d$. Then, for all $\theta \in (0,1)$,

$$\overline{\dim}_{\mathrm{A}}^{\theta} F = \sup_{0<\theta' \leqslant \theta} \dim_{\mathrm{A}}^{\theta'} F.$$

*Proof* We give the proof in the case where $F$ is bounded. The unbounded case requires a delicate argument for which we refer the reader to [95]. We only need to prove

$$\sup_{0<\theta' \leqslant \theta} \dim_{\mathrm{A}}^{\theta'} F \geqslant \overline{\dim}_{\mathrm{A}}^{\theta} F$$

since the other direction is trivial.

Let $\theta \in (0,1)$, $s = \overline{\dim}_{\mathrm{A}}^{\theta} F$ which we may assume is strictly positive, and $0 < \varepsilon < s$. By definition we can find sequences $x_i, r_i, R_i$ ($i \geqslant 1$) such that $x_i \in F$, $0 < r_i \leqslant R_i^{1/\theta} < R_i < 1$, $(r_i/R_i) \to 0$ and

$$N_{r_i}\left(B(x_i, R_i) \cap F\right) \geqslant \left(\frac{R_i}{r_i}\right)^{s-\varepsilon}. \tag{3.2}$$

We can assume that $R_i \to 0$ since otherwise (using (3.1), see Lemma 3.4.1) $\dim_{\mathrm{A}}^{\theta} F \geqslant \overline{\dim}_{\mathrm{B}} F \geqslant s - \varepsilon$, which is sufficient.

For each $i$, let $\theta_i$ be defined by $r_i = R_i^{1/\theta_i}$, noting that $0 < \theta_i \leqslant \theta$ for all $i$. Using compactness of $[0,\theta]$ to extract a convergent subsequence, we may assume that $\theta_i \to \theta' \in [0,\theta]$ and, by taking a further subsequence if necessary, we may assume that $|\theta_i - \theta'| < \delta$ for all $i$ where $\delta > 0$ can be chosen arbitrarily. We may also assume that the sequence $\theta_i$ is either non-increasing or strictly increasing.

First suppose that $\theta' = 0$. Since $F$ is assumed to be bounded

$$N_{R_i^{1/\theta_i}}(F) \geqslant N_{R_i^{1/\theta_i}}(B(x_i, R_i) \cap F) \geqslant \left(\frac{R_i}{R_i^{1/\theta_i}}\right)^{s-\varepsilon} \geqslant \left(\frac{1}{R_i^{1/\theta_i}}\right)^{(s-\varepsilon)(1-\delta)}$$

by (3.2). Note that the final inequality uses the fact that $\delta > \theta_i$, which follows since $\delta > |\theta_i - \theta'| = \theta_i$. Then (using (3.1), see Lemma 3.4.1) $\dim_{\mathrm{A}}^{\theta} F \geqslant \overline{\dim}_{\mathrm{B}} F \geqslant (s-\varepsilon)(1-\delta)$. Since $\delta > 0$ and $\varepsilon > 0$ can be chosen arbitrarily small, this yields the desired result.

From now on suppose $\theta' > 0$. If the sequence $\theta_i$ is non-increasing, then $\theta' \leqslant \theta_i$, and therefore $R_i^{1/\theta_i} \geqslant R_i^{1/\theta'}$, for all $i$. It follows that

$$\begin{aligned} N_{R_i^{1/\theta'}}(B(x_i, R_i) \cap F) &\geqslant N_{R_i^{1/\theta_i}}(B(x_i, R_i) \cap F) \\ &\geqslant \left(\frac{R_i}{R_i^{1/\theta_i}}\right)^{s-\varepsilon} \quad \text{by (3.2)} \\ &= \left(\frac{R_i}{R_i^{1/\theta'}}\right)^{\left(\frac{1-1/\theta_i}{1-1/\theta'}\right)(s-\varepsilon)} \\ &\geqslant \left(\frac{R_i}{R_i^{1/\theta'}}\right)^{\frac{\theta'(1-\theta'-\delta)}{(\theta'+\delta)(1-\theta')}(s-\varepsilon)} \end{aligned}$$

where the final inequality uses

$$\frac{1-1/\theta_i}{1-1/\theta'} = \frac{\theta'(1-\theta_i)}{\theta_i(1-\theta')} \geqslant \frac{\theta'(1-\theta'-\delta)}{(\theta'+\delta)(1-\theta')}$$

which holds since $\theta_i \leqslant \theta' + \delta$. This yields $\dim_{\mathrm{A}}^{\theta'} F \geqslant \frac{\theta'(1-\theta'-\delta)}{(\theta'+\delta)(1-\theta')}(s-\varepsilon)$ and, since $\delta > 0$ can be chosen arbitrarily small (after fixing $\theta'$), we obtain $\dim_{\mathrm{A}}^{\theta'} F \geqslant s - \varepsilon$.

On the other hand, if $\theta_i$ is strictly increasing, then $\theta' > \theta_i$ for all $i$. Taking another subsequence if necessary we can also assume that $\theta_i > \theta'/2$ for all $i$. Covering by $R_i^{1/\theta'}$-balls and then covering each $R_i^{1/\theta'}$- ball by $R_i^{1/\theta_i}$-balls we obtain

$$\begin{aligned} &N_{R_i^{1/\theta_i}}\left(B(x_i,R_i)\cap F\right) \\ &\qquad \leqslant N_{R_i^{1/\theta'}}\left(B(x_i,R_i)\cap F\right)\left(\sup_{z\in\mathbb{R}^d} N_{R_i^{1/\theta_i}}\left(B\left(z,R_i^{1/\theta'}\right)\right)\right) \\ &\qquad \leqslant N_{R_i^{1/\theta'}}\left(B(x_i,R_i)\cap F\right)c(d)\left(\frac{R_i^{1/\theta'}}{R_i^{1/\theta_i}}\right)^d \end{aligned}$$

where $c(d)\geqslant 1$ is a constant depending only on $d$. Therefore

$$\begin{aligned} N_{R_i^{1/\theta'}}\left(B(x_i,R_i)\cap F\right) &\geqslant c(d)^{-1}N_{R_i^{1/\theta_i}}\left(B(x_i,R_i)\cap F\right)R_i^{(1/\theta_i-1/\theta')d} \\ &\geqslant c(d)^{-1}\left(\frac{R_i}{R_i^{1/\theta_i}}\right)^{s-\varepsilon}R_i^{(1/\theta_i-1/\theta')d} \qquad \text{by (3.2)} \\ &= c(d)^{-1}\left(\frac{R_i}{R_i^{1/\theta'}}\right)^{\left(\frac{1-1/\theta_i}{1-1/\theta'}\right)(s-\varepsilon)+\left(\frac{1/\theta_i-1/\theta'}{1-1/\theta'}\right)d} \\ &\geqslant c(d)^{-1}\left(\frac{R_i}{R_i^{1/\theta'}}\right)^{s-\varepsilon-\frac{\delta d}{(1-\theta')\theta'/2}}, \end{aligned}$$

where the final inequality uses the coefficient bounds

$$\left(\frac{1-1/\theta_i}{1-1/\theta'}\right)\geqslant 1$$

which holds since $\theta'>\theta_i$ and

$$\left(\frac{1/\theta_i-1/\theta'}{1-1/\theta'}\right)=-\left(\frac{\theta'-\theta_i}{(1-\theta')\theta_i}\right)\geqslant-\frac{\delta}{(1-\theta')\theta'/2}$$

which holds since $0<\theta'-\theta_i\leqslant\delta$ and $\theta_i>\theta'/2$. It follows that $\dim_{\mathrm{A}}^{\theta'}F\geqslant s-\varepsilon-\frac{\delta d}{(1-\theta')\theta'/2}$ and since $\delta>0$ can be chosen arbitrarily small (after fixing $\theta'$) we obtain $\dim_{\mathrm{A}}^{\theta'}F\geqslant s-\varepsilon$ as before. Since $\varepsilon>0$ was arbitrary it follows that

$$\sup_{0<\theta'\leqslant\theta}\dim_{\mathrm{A}}^{\theta'}F\geqslant s$$

completing the proof. □

One of the benefits of Theorem 3.3.6 is that it allows us to focus future study on the Assouad spectrum rather than the upper spectrum which could have *a priori* contained new information in its own right. Not only does the Assouad spectrum contain strictly more information than the upper spectrum, but it is also easier to work with since the family of scales is 1-parameter, rather than 2-parameter.

A theoretically significant corollary to Theorem 3.3.6, is that we obtain the interpolation result which motivated the introduction of the Assouad spectrum in the first place, albeit with Assouad dimension replaced by quasi-Assouad dimension. In many (but certainly not all) cases of interest, the quasi-Assouad and Assouad dimensions coincide and so genuine interpolation between the upper box and Assouad dimension is achieved.

**Corollary 3.3.7** Let $F \subseteq \mathbb{R}^d$. Then $\dim_A^\theta F \to \dim_{qA} F$ as $\theta \to 1$.

Theorem 3.3.6 only directly implies that $\limsup_{\theta \to 1} \dim_A^\theta F = \dim_{qA} F$, but the fact that the limit of $\dim_A^\theta F$ as $\theta \to 1$ exists follows from [111, Remark 3.9], see [95, Section 3.2].

The analogous 'tale of two spectra' problem for the lower spectrum was considered in [41, 42]. The quasi-lower dimension, $\dim_{qL} F$, was introduced in [41] and in [42] it was proved that

$$\dim_{qL} F = \lim_{\theta \to 1} \dim_L^\theta F \tag{3.3}$$

provided $\dim_L F > 0$. The assumption that $\dim_L F > 0$ was removed in [127, Theorem A.2], which allows us to take (3.3) as our definition of quasi-lower dimension for any set $F \subseteq \mathbb{R}^d$.

### 3.3.3 Recovering the interpolation

The Assouad spectrum was introduced to understand the 'gap' in between the upper box and Assouad dimensions. This was partially achieved since, for bounded $F \subseteq \mathbb{R}^d$,

(i) $\dim_A^\theta F \to \overline{\dim}_B F$ as $\theta \to 0$,
(ii) $\dim_A^\theta F$ is continuous in $\theta \in (0, 1)$,
(iii) $\dim_A^\theta F \to \dim_{qA} F$ as $\theta \to 1$.

In many cases of interest $\dim_{qA} F = \dim_A F$ and so the full range $[\overline{\dim}_B F, \dim_A F]$ is 'witnessed' by the Assouad spectrum. However, if $\dim_{qA} F < \dim_A F$, then there is a gap and the desired interpolation is not achieved. An approach for 'recovering the interpolation' was outlined in [111], which we briefly explain here.

Let $\phi : [0,1] \to [0,1]$ be an increasing continuous function such that $\phi(R) \leqslant R$ for all $R \in [0,1]$. The *$\phi$-Assouad dimension*, introduced in [111], is defined by

$$\dim_{\mathrm{A}}^{\phi} F \;=\; \inf \Bigg\{ \alpha : \text{ there exists a constant } C > 0 \text{ such that, for all } 0 < r \leqslant \phi(R) \leqslant R < 1 \text{ and } x \in F, \; N_r\big(B(x,R) \cap F\big) \;\leqslant\; C\left(\frac{R}{r}\right)^{\alpha} \Bigg\}.$$

If $\phi(R) = R^{1/\theta}$, then $\dim_{\mathrm{A}}^{\phi} F$ recovers the upper Assouad spectrum (which can be expressed entirely in terms of the Assouad spectrum), and if $\phi(R) = R$, then $\dim_{\mathrm{A}}^{\phi} F = \dim_{\mathrm{A}} F$ (for bounded $F$). Therefore, one may recover the desired interpolation by identifying precise conditions on $\phi$ which guarantee $\dim_{\mathrm{A}}^{\phi} F = \dim_{\mathrm{A}} F$. Often $\dim_{\mathrm{A}}^{\theta} F = \dim_{\mathrm{A}} F$ for some $\theta \in (0,1)$, in which case the threshold for witnessing the Assouad dimension is provided by the function $\phi(R) = R^{1/\theta}$. The *$\phi$-Assouad dimension* has been considered in detail[4] by García, Hare, and Mendivil [120, 121] and Troscheit [268]. A particular case of interest is Mandelbrot percolation, see Section 9.4 and Theorem 9.4.4.

Many interesting results pertaining to the $\phi$-Assouad dimension, and several worthwhile variants, were given in [120]. We highlight some here but refer the reader to [120] for more detail.

**Theorem 3.3.8**

(i) If $\frac{\log \phi(R)}{\log R} \to 1$ as $R \to 0$, then, for all sets $F \subseteq \mathbb{R}^d$,

$$\dim_{\mathrm{A}}^{\phi} F \geqslant \dim_{\mathrm{qA}} F.$$

(ii) If $\inf_{R>0} \frac{\phi(R)}{R} > 0$, then, for all bounded sets $F \subseteq \mathbb{R}^d$,

$$\dim_{\mathrm{A}}^{\phi} F = \dim_{\mathrm{A}} F.$$

(iii) There exist sets $F \subseteq \mathbb{R}^d$ such that $\dim_{\mathrm{qA}} F < \dim_{\mathrm{A}} F$ and for all $s \in [\dim_{\mathrm{qA}} F, \dim_{\mathrm{A}} F]$, there exists $\phi$ (with the properties described above) such that $\dim_{\mathrm{A}}^{\phi} F = s$.

[4] Slightly different notation is used in these papers. Instead of writing $\phi(R)$ they write $R^{1+\Phi(R)}$. Therefore some care is required when transferring results from one set of notation to another.

## 3.4 Basic properties: *revisited*

Before we lose track of how many definitions of dimension we are considering, we briefly return to the topic of basic properties and collect together some observations on the variants of the Assouad dimension we have encountered.

### 3.4.1 Relationships between dimensions

Perhaps the most important facts to establish are how these dimensions relate to each other. This was already touched on in (3.1) above. The following should be compared with Lemma 2.4.3 and was observed in [111].

**Lemma 3.4.1** For a bounded set $F \subseteq \mathbb{R}^d$,

$$\overline{\dim}_{\mathrm{B}} F \;\leqslant\; \dim_{\mathrm{A}}^{\theta} F \;\leqslant\; \dim_{\mathrm{qA}} F \;\leqslant\; \dim_{\mathrm{A}} F$$

and

$$\begin{array}{ccccc} & & \dim_{\mathrm{ML}} F & & \\ & \rotatebox{30}{$\leqslant$} & & \rotatebox{-30}{$\leqslant$} & \\ \dim_{\mathrm{L}} F & & & & \underline{\dim}_{\mathrm{B}} F. \\ & \rotatebox{-30}{$\leqslant$} & & \rotatebox{30}{$\leqslant$} & \\ & & \dim_{\mathrm{L}}^{\theta} F & & \end{array}$$

*Proof* Most of these inequalities follow easily from the definitions and their proofs are left as exercises. The details can be found in [111]. We prove $\overline{\dim}_{\mathrm{B}} F \leqslant \dim_{\mathrm{A}}^{\theta} F$ here, which is not quite immediate. We may assume $\overline{\dim}_{\mathrm{B}} F > 0$ and let $0 < t < \overline{\dim}_{\mathrm{B}} F < s$. By covering $F$ with $r^{\theta}$-balls and then covering each of these $r^{\theta}$-balls with $r$-balls we obtain

$$N_r(F) \leqslant N_{r^{\theta}}(F) \left( \sup_{x \in F} N_r(B(x, r^{\theta}) \cap F) \right).$$

By the definition of upper box dimension, there exist arbitrarily small $r > 0$ such that

$$\sup_{x \in F} N_r(B(x, r^{\theta}) \cap F) \geqslant \frac{N_r(F)}{N_{r^{\theta}}(F)} \geqslant \frac{r^{-t}}{r^{-\theta s}} = \left( \frac{r^{\theta}}{r} \right)^{\frac{\theta s - t}{\theta - 1}}$$

which proves that $\dim_{\mathrm{A}}^{\theta} F \geqslant \frac{t - \theta s}{1 - \theta}$. Since $s$ and $t$ can be made arbitrarily close to $\overline{\dim}_{\mathrm{B}} F$, the lower bound follows. □

In general the lower dimension and lower spectrum are not comparable to the Hausdorff dimension. For example, for all $\theta \in (0,1)$,

$$\dim_{\mathrm{L}} \mathbb{Q} = \dim_{\mathrm{L}}^{\theta} \mathbb{Q} = 1 > 0 = \dim_{\mathrm{H}} \mathbb{Q}$$

and

$$\dim_{\mathrm{L}}[0,1] \cup \{2\} = \dim_{\mathrm{L}}^{\theta}[0,1] \cup \{2\} = 0 < 1 = \dim_{\mathrm{H}}[0,1] \cup \{2\}.$$

However, if $F$ is closed, then one can prove that the lower dimension is the 'smallest of all dimensions'. This was first proved by Larman in the 1960s [179, 180]. To prove this result we rely on the mass distribution principle, Lemma 3.4.2. This is a well-known result in fractal geometry (for example, see [70]) but we include the proof since it is so simple.

**Lemma 3.4.2** Let $F \subseteq \mathbb{R}^d$ be a Borel set and suppose that for some constants $C, s > 0$ there exists a non-trivial Borel measure on $F$ such that for all $x \in F$ and $r > 0$

$$\mu(B(x,r)) \leqslant Cr^s.$$

Then $\dim_{\mathrm{H}} F \geqslant s$.

*Proof* Let $\{U_i\}_i$ be an arbitrary $r$-cover of $F$. Then

$$0 < \mu(F) \leqslant \sum_i \mu(U_i) \leqslant \sum_i C|U_i|^s$$

and taking infimum over all $r$-covers we obtain $\mathcal{H}_r^s(F) \geqslant \mu(F)/C > 0$ and therefore $\mathcal{H}^s(F) \geqslant \mu(F)/C > 0$ and $\dim_{\mathrm{H}} F \geqslant s$ as required. □

**Theorem 3.4.3** If $F \subseteq \mathbb{R}^d$ is closed, then

$$\dim_{\mathrm{L}} F \ \leqslant\ \dim_{\mathrm{ML}} F \ \leqslant\ \dim_{\mathrm{H}} F.$$

*Proof* We first prove that $\dim_{\mathrm{L}} F \leqslant \dim_{\mathrm{H}} F$, and we may assume the lower dimension is strictly positive as otherwise there is nothing to prove. Fix $0 < s < \dim_{\mathrm{L}} F$. We claim that there exists $c \in (0,1)$ such that, given any $x \in F$ and $0 < \rho < |F|$, we can find a collection of $c^{-s}$ points $x_i \in F$ such that the closed balls $B(x_i, c\rho)$ are pairwise disjoint and $B(x_i, c\rho) \subseteq B(x, \rho)$. To see this, let $s < t < \dim_{\mathrm{L}} F$ and note that it follows directly from the definition of lower dimension that there exists $C > 0$ such that, for all $x \in F$ and $0 < r < R < |F|$, we can find a collection of $C(R/r)^t$ points $x_i \in F$ such that the closed balls $B(x_i, r)$ are pairwise disjoint and $B(x_i, r) \subseteq B(x, R)$. (This is taking advantage of the relationship between minimal covers and maximal packings.) The claim follows by choosing $R = \rho$ and $r = c\rho$ such that $c \in (0,1)$ is chosen small enough to guarantee $C(R/r)^t = Cc^{s-t}c^{-s} \geqslant c^{-s}$.

Let $B_0$ be any closed ball centred in $F$ with radius $\rho \in (0, |F|)$ and apply the claim to obtain a collection of $c^{-s}$ pairwise disjoint closed balls of radius $c\rho$ centred in $F$ which are all contained in $B_0$. Apply this procedure inductively on every ball present at each stage of the construction to produce a decreasing sequence of sets $\mathcal{B}_k$, (that is, $B_0 = \mathcal{B}_0 \supset \mathcal{B}_1 \supset \mathcal{B}_2 \supset \dots$), where each $\mathcal{B}_k$ consists of $c^{-sk}$ many pairwise disjoint closed balls of radius $c^k\rho$ centred in $F$. The intersection

$$F' = \bigcap_k \mathcal{B}_k \tag{3.4}$$

is a non-empty compact set and since $F$ is closed we necessarily have $F' \subseteq F$. Let $\mu$ be a Borel probability measure defined on $F'$ where each ball in $\mathcal{B}_k$ has mass $c^{ks}$. Given $r \in (0,1)$, let $k$ be the unique integer such that $c^k\rho \leqslant r < c^{(k-1)}\rho$. By a simple volume argument, for all $x \in F$ and $r > 0$, the ball $B(x,r)$ can intersect no more than a constant $C(d)$ many balls from $\mathcal{B}_k$, where $C(d)$ is a constant depending only on the ambient spatial dimension $d$. Therefore

$$\mu(B(x,r)) \leqslant C(d)c^{ks} \leqslant C(d)\rho^{-s}r^s$$

and it follows from the mass distribution principle Lemma 3.4.2 that $\dim_{\mathrm{H}} F \geqslant s$. Letting $s \to \dim_{\mathrm{L}} F$ yields $\dim_{\mathrm{H}} F \geqslant \dim_{\mathrm{L}} F$, as required.

Moreover, for any $E \subseteq F$,

$$\dim_{\mathrm{L}} E = \dim_{\mathrm{L}} \overline{E} \leqslant \dim_{\mathrm{H}} \overline{E} \leqslant \dim_{\mathrm{H}} \overline{F} = \dim_{\mathrm{H}} F$$

where the final inequality uses the fact that $F$ is closed. Therefore $\dim_{\mathrm{ML}} F \leqslant \dim_{\mathrm{H}} F$. □

The lower spectrum is *not* generally bounded above by the modified lower dimension or by the Hausdorff dimension, even for closed sets. This is perhaps a little counter-intuitive because the proof given above almost looks like it should work with the lower spectrum replacing the lower dimension. However, the problem is that when we iterate the information we get from the lower spectrum for a fixed $\theta$ the radii of the balls in the nested sequence $\mathcal{B}_k$ decrease super-exponentially, which is not sufficient to establish the mass distribution estimate. For a concrete example of a compact set $F$ with

$$\dim_{\mathrm{L}}^{\theta} F > \dim_{\mathrm{H}} F$$

for some $\theta \in (0,1)$, see the self-affine carpets considered in Section 8.

García and Hare proved that the *quasi*-lower dimension of a closed set *is* bounded above by its Hausdorff dimension, that is, for closed $F \subseteq \mathbb{R}^d$, $\dim_{\mathrm{qL}} F \leqslant \dim_{\mathrm{H}} F$, see [118, Proposition 10]. Therefore, one always has

$$\lim_{\theta \to 1} \dim_{\mathrm{L}}^{\theta} F \leqslant \dim_{\mathrm{H}} F$$

for closed $F$, even if the spectrum is not uniformly bounded above by the Hausdorff dimension. Slightly more information can be derived from [120, Proposition 2.12], but we omit the details.

Despite the upper box dimension being a lower bound for the Assouad spectrum, one can also bound the Assouad spectrum above in terms of the upper box dimension. The following result, first proved in [111, Proposition 3.1], will be used throughout the book.

**Lemma 3.4.4** For a non-empty bounded set $F$ and $\theta \in (0, 1)$,

$$\overline{\dim}_{\mathrm{B}} F \leqslant \dim_{\mathrm{A}}^{\theta} F \leqslant \min \left\{ \frac{\overline{\dim}_{\mathrm{B}} F}{1 - \theta}, \ \dim_{\mathrm{qA}} F \right\}.$$

*Proof* The lower bound is from Lemma 3.4.1. The fact that the quasi-Assouad dimension is an upper bound is trivial, and also from Lemma 3.4.1, and so it remains to prove the 'other' upper bound.

Fix $\theta \in (0, 1)$, $s > \overline{\dim}_{\mathrm{B}} F$ and let $x \in F$ and $r \in (0, 1)$. Appealing to the definition of the upper box dimension, there exists $C > 0$ depending only on $s$ such that

$$N_r(B(x, r^{\theta}) \cap F) \leqslant N_r(F) \leqslant C r^{-s} = C \left( \frac{r^{\theta}}{r} \right)^{s/(1-\theta)},$$

which implies $\dim_{\mathrm{A}}^{\theta} F \leqslant \frac{s}{1-\theta}$. Since $s > \overline{\dim}_{\mathrm{B}} F$ was arbitrary, the desired upper bound follows. □

The bounds in Lemma 3.4.4 are sharp. For a family of examples where the upper bound is obtained, see Theorem 3.4.7, and for a method for constructing examples where the lower bound is obtained, see Theorem 3.3.5 or Theorem 9.4.2.

The following counter-intuitive result follows from Lemma 3.4.4 and the fact that the Assouad spectrum approaches the quasi-Assouad dimension as $\theta \to 1$, see Corollary 3.3.7. We give an alternative direct proof here, following [118].

**Lemma 3.4.5** Let $F \subseteq \mathbb{R}^d$ be a non-empty bounded set. Then $\overline{\dim}_{\mathrm{B}} F = 0$ if and only if $\dim_{\mathrm{qA}} F = 0$.

*Proof* This result was first proved by García and Hare [118, Proposition 15], and we give their elegant proof, which does not rely on the Assouad spectrum. Since $\overline{\dim}_{\mathrm{B}} F \leqslant \dim_{\mathrm{qA}} F$ always holds, it suffices to prove that if $\overline{\dim}_{\mathrm{B}} F = 0$, then $\dim_{\mathrm{qA}} F = 0$. Let $\varepsilon, \delta > 0$ and note that $\overline{\dim}_{\mathrm{B}} F = 0$ means that there exists a constant $C > 0$ such that for all $r \in (0,1)$, $N_r(F) \leqslant C r^{-\varepsilon\delta/(1+\delta)}$. It follows that, for all $x \in F$ and $0 < r \leqslant R^{1+\delta} \leqslant 1$,

$$N_r(B(x,R) \cap F) \leqslant N_r(F) \leqslant C r^{-\varepsilon\delta/(1+\delta)} \leqslant C(R/r)^{\varepsilon}$$

which proves that $h_F(\delta) \leqslant \varepsilon$ and since $\varepsilon, \delta > 0$ are arbitrary the result follows. □

The 'only if' part of Lemma 3.4.5 is surprising since we generally have control in the opposite direction, that is, $\overline{\dim}_{\mathrm{B}} F \leqslant \dim_{\mathrm{qA}} F$, and, moreover, the *Assouad* dimension can take on any value in $[0, d]$ even in cases where the box dimension is 0. We are not aware of such a mutual dependency result holding for any other pair of dimensions.

**Corollary 3.4.6** Let $F \subseteq \mathbb{R}^d$ be a non-empty bounded set with $\overline{\dim}_{\mathrm{B}} F < \dim_{\mathrm{qA}} F$. Suppose

$$\rho = \inf\{\theta \in (0,1) : \dim_{\mathrm{A}}^{\theta} F = \dim_{\mathrm{qA}} F\} = 1 - \frac{\overline{\dim}_{\mathrm{B}} F}{\dim_{\mathrm{qA}} F}$$

noting that this is as small as the phase transition $\rho$ can be by Lemma 3.4.4. Then, for all $\theta \in (0,1)$,

$$\dim_{\mathrm{A}}^{\theta} F = \min\left\{ \frac{\overline{\dim}_{\mathrm{B}} F}{1-\theta}, \dim_{\mathrm{qA}} F \right\}.$$

*Proof* This follows from Corollary 3.3.4 and Lemma 3.4.4. □

To motivate Corollary 3.4.6, we briefly return to the example considered in Theorem 2.1.1. We take the example a little further by including an exponent $p$ as well as computing the Assouad spectrum and demonstrating sharpness of the upper bound from Lemma 3.4.4.

**Theorem 3.4.7** For $p > 0$ and $F_p = \{0\} \cup \{1/n^p : n \in \mathbb{N}\}$,

$$\dim_{\mathrm{A}}^{\theta} F_p = \min\left\{ \frac{1}{(1+p)(1-\theta)}, 1 \right\}.$$

Therefore,

$$\dim_{\mathrm{B}} F_p = \frac{1}{1+p},$$

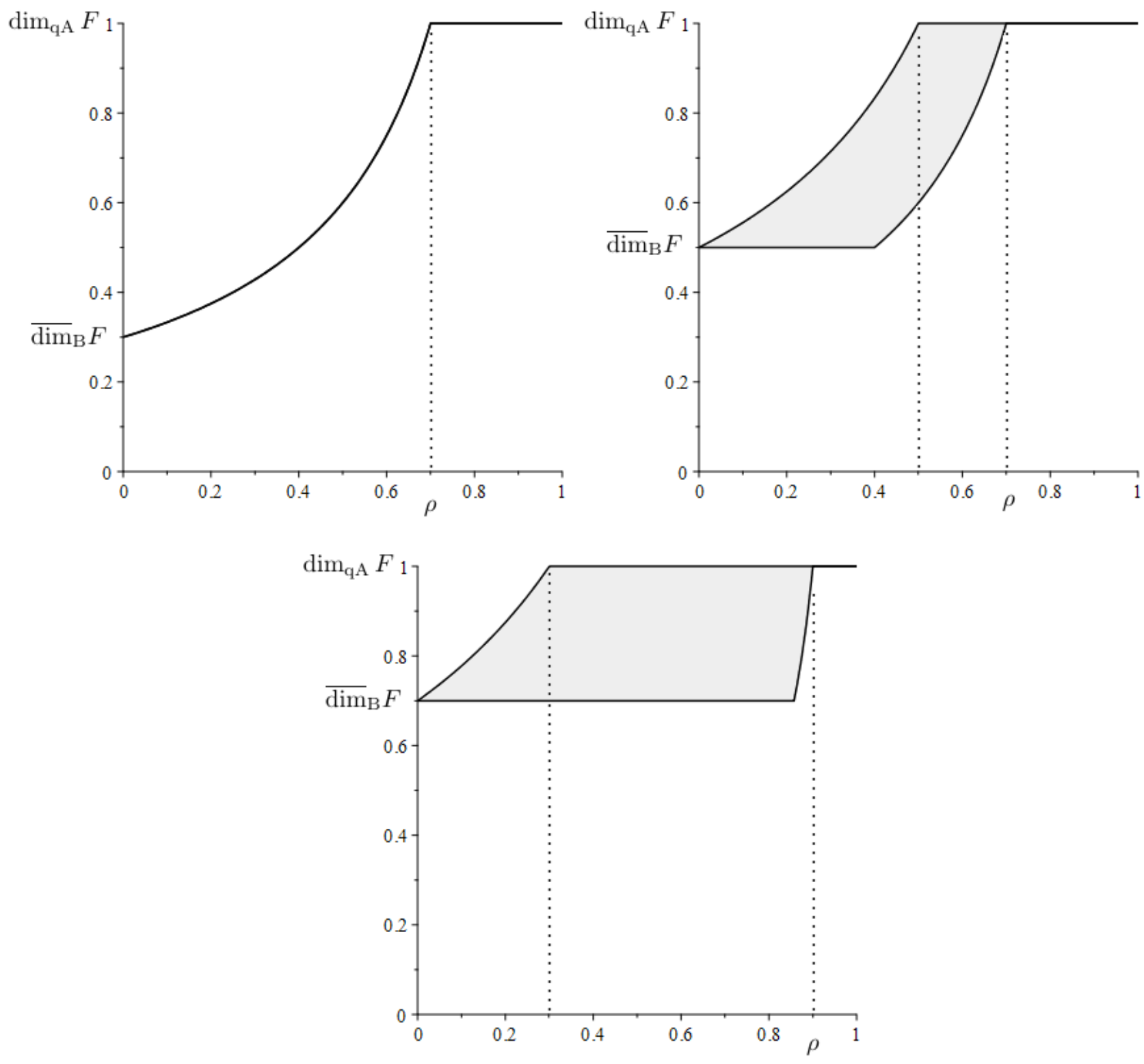

Figure 3.4 Plots of the bounds on the Assouad spectrum as functions of $\theta$ in terms of the upper box, quasi-Assouad dimension and the phase transition $\rho$ in different cases. The upper bounds are from Lemma 3.4.4 and the lower bounds are from Corollary 3.3.4. The gap in between the upper and lower bounds is shaded in grey. Vertical dotted lines are included to indicate $\rho$ and $1 - \frac{\overline{\dim}_{\mathrm{B}} F}{\dim_{\mathrm{qA}} F}$. Corollary 3.4.6 applies to the top left example and so the upper and lower bounds coincide.

and

$$\dim_{\mathrm{A}} F_p = \dim_{\mathrm{qA}} F_p = 1.$$

*Proof* In light of Lemma 3.4.4 and Corollary 3.4.6 it suffices to prove that $\dim_{\mathrm{B}} F_p = \frac{1}{1+p}$ and $\dim_{\mathrm{A}}^{p/(p+1)} F_p \geqslant 1$. Notice that both of these facts were covered in the case $p = 1$ in Theorem 2.1.1. There it was proved that $\dim_{\mathrm{A}} F_1 \geqslant 1$. However, the scales used to prove this satisfied $r = R^2$, which means that the stronger lower bound $\dim_{\mathrm{A}}^{1/2} F_1 \geqslant 1$ holds.

The extension of this argument to general $p$ is straightforward and left to the reader. □

Corollary 3.3.4 can be improved if the lower dimension is large.

**Theorem 3.4.8** For any set $F \subseteq \mathbb{R}^d$ and $0 < \theta_1 < \theta_2 < 1$,

$$\dim_{\mathrm{A}}^{\theta_1} F \geqslant \left(\frac{\theta_2 - \theta_1}{1-\theta_1}\right) \dim_{\mathrm{A}}^{\theta_1/\theta_2} F + \left(\frac{1-\theta_2}{1-\theta_1}\right) \dim_{\mathrm{L}}^{\theta_2} F.$$

In particular, if

$$\rho = \inf\{\theta \in (0,1) : \dim_{\mathrm{A}}^{\theta} F = \dim_{\mathrm{qA}} F\}$$

exists, $\rho \in (0,1)$ and $\dim_{\mathrm{L}} F = \overline{\dim}_{\mathrm{B}} F$, then for $\theta \in (0,\rho)$

$$\dim_{\mathrm{A}}^{\theta} F \ \geqslant\ \overline{\dim}_{\mathrm{B}} F + \frac{(1-\rho)\theta}{(1-\theta)\rho}(\dim_{\mathrm{qA}} F - \overline{\dim}_{\mathrm{B}} F).$$

*Proof* The proof strategy here is similar to Theorem 3.3.1. Given $r \in (0,1)$, first choose $x \in F$ which maximises

$$N_{r^{\theta_2}}(B(x, r^{\theta_1}) \cap F) \tag{3.5}$$

and take a maximal $r^{\theta_2}$-packing of $B(x, r^{\theta_1}) \cap F$ by balls of radius $r^{\theta_2}$. We may consider $r$-covers of each of these balls independently, losing at most a constant in our estimate for $N_r(B(x, r^{\theta_1}) \cap F)$. However, we can only bound the $r$-covering number for each of these balls from below by

$$\inf_{x \in F} N_r(B(x, r^{\theta_2}) \cap F). \tag{3.6}$$

Combining (3.5) and (3.6), for arbitrary $\varepsilon > 0$,

$$N_{r^{\theta_2}}(B(x, r^{\theta_1}) \cap F) \ \geqslant\ C \left(\frac{r^{\theta_1}}{r^{\theta_2}}\right)^{\dim_{\mathrm{A}}^{\theta_1/\theta_2} F - \varepsilon} \left(\frac{r^{\theta_2}}{r}\right)^{\dim_{\mathrm{L}}^{\theta_2} F - \varepsilon}$$

for a constant $C$ depending only on $\varepsilon$. This proves the result upon letting $\varepsilon \to 0$. □

The lower bound in Theorem 3.4.8 is harder to establish than the lower bound in Theorem 3.3.1 but — perhaps surprisingly — either estimate may be better, depending on how relatively large the lower dimension is. The lower bound from Theorem 3.4.8 turns out to be a natural candidate for the Assouad spectrum of several examples, both when the lower dimension is equal to the box dimension and when it is strictly smaller, see Section 17.7.

Somewhat surprisingly, there does not appear to be an analogue of Lemma 3.4.4 for the lower spectrum. In general, one cannot improve on $\dim_{\mathrm{L}} F \leqslant \dim_{\mathrm{L}}^{\theta} F \leqslant \underline{\dim}_{\mathrm{B}} F$, see [111, Proposition 3.9]. However, it was proved in [42, Theorem 1.1] that

$$\lim_{\theta \to 0} \dim_{\mathrm{L}}^{\theta} F \tag{3.7}$$

exists, and can take any value in the interval $[\dim_{\mathrm{qL}} F, \underline{\dim}_{\mathrm{B}} F]$. It would be interesting to see if (3.7) has any particular meaning outside the context of the lower spectrum. We will prove later that (3.7) always coincides with the lower box dimension if $F$ is the attractor of an IFS of bi-Lipschitz contractions, see Theorem 6.3.1.

### 3.4.2 Further basic properties

In this section we summarise some other basic properties satisfied by the variants of the Assouad dimension met so far. This discussion will culminate with Table 3.1, which should be compared with Table 2.1 on page 20. This discussion mirrors that in Section 2.4 and the reader should refer back for definitions of the basic properties we discuss, such as stable under closure, monotonicity etc. For clarity, we say a property holds for the Assouad or lower spectrum if it holds simultaneously for all $\theta \in (0,1)$. For example, the Assouad spectrum is stable under closure because, for all $F \subseteq \mathbb{R}^d$, $\dim_{\mathrm{A}}^{\theta} \overline{F} = \dim_{\mathrm{A}}^{\theta} F$ simulataneously for all $\theta \in (0,1)$.

The lower dimension is stable under closure, but the *modified* lower dimension is not. For example, let

$$F \;=\; \{(p/q, 1/q) \;:\; p, q \in \mathbb{N},\, p \leqslant q,\, \gcd(p,q) = 1\} \tag{3.8}$$

and observe that every point $x \in F$ is isolated and therefore any subset of $F$ has an isolated point, see Figure 3.5. This implies that $\dim_{\mathrm{ML}} F = 0$. However, $[0,1] \times \{0\} \subseteq \overline{F}$ and so $\dim_{\mathrm{ML}} \overline{F} = 1$.

On the other hand, the open set property is clearly satisfied for modified lower dimension. However, it fails for the lower dimension and lower spectrum. Let

$$U \;=\; \bigcup_{n=1}^{\infty} \left(1/n - 2^{-n},\ 1/n + 2^{-n}\right) \;\subseteq\; \mathbb{R}$$

which is an open set. Let $\theta \in (0,1)$ and $s, C > 0$. For $n \geqslant 1$, let $r(n) = 2/2^n$ and observe that, for all sufficiently large $n$,

$$r(n)^{\theta} < 1/n - 1/(n+1) - 2^{-(n+1)}$$

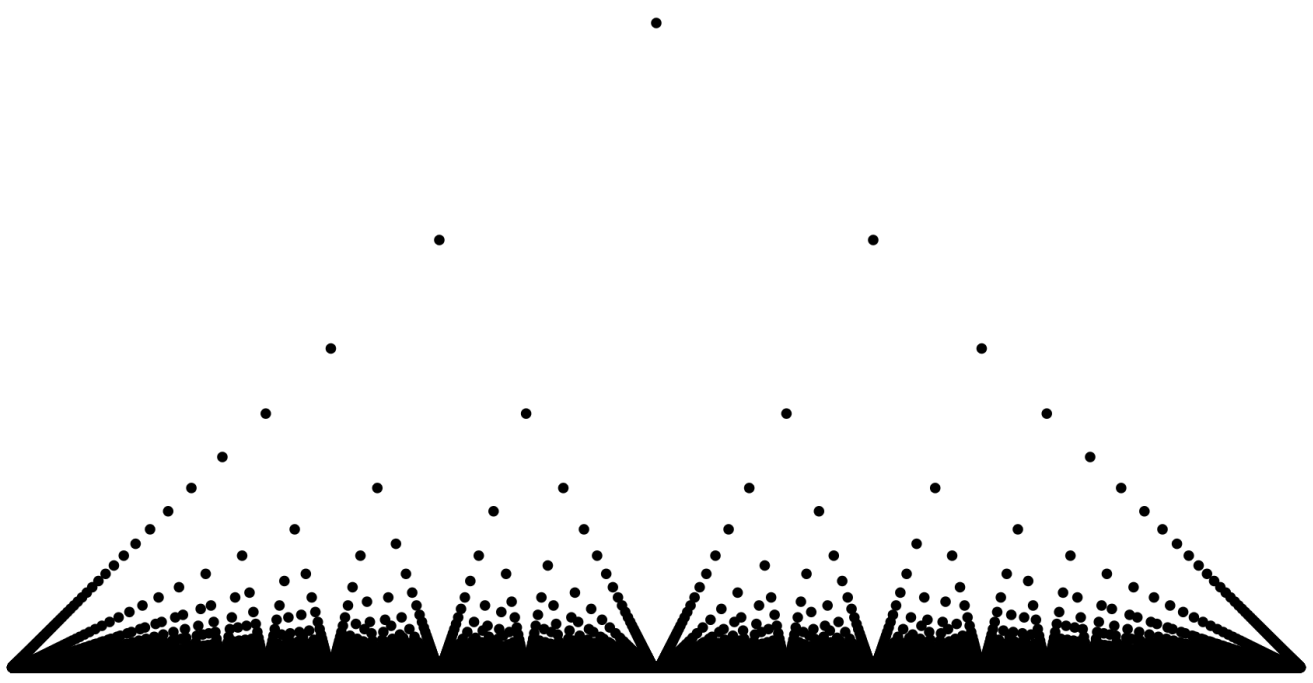

Figure 3.5 The set (3.8), which coincides with the graph of the 'popcorn function' restricted to the rationals. The popcorn function is defined to be identically 0 on the irrationals. This function has several notable properties, including being continuous at all irrational points and discontinuous at all rational points.

and therefore

$$N_{r(n)}\Big(B\big(1/n, r(n)^\theta\big) \cap U\Big) = 1$$

since $B\big(1/n, r(n)^\theta\big) \cap U = \big(1/n - 2^{-n},\ 1/n + 2^{-n}\big)$. Therefore we may choose $n$ large enough to ensure

$$N_{r(n)}\Big(B\big(1/n, r(n)^\theta\big) \cap U\Big) \ = \ 1 \ < \ C\left(\frac{r(n)^\theta}{r(n)}\right)^s$$

which implies $\dim_{\mathrm{L}} U \leqslant \dim_{\mathrm{L}}^\theta U \leqslant s$ and letting $s \to 0$ yields $\dim_{\mathrm{L}} U = \dim_{\mathrm{L}}^\theta U = 0$. We also note that the unbounded open set $\mathbb{R} \times (0,1) \subseteq \mathbb{R}^2$ has lower dimension 1, which can be seen by choosing $R \to \infty$ and $r = 1$.

The lower dimension and spectrum are not stable under finite unions: consider the union of $[0,1] \cup \{2\}$ and $\{0\} \cup [1,2]$. One can say more if the sets in the union are uniformly separated, but the behaviour is rather different from the other dimensions.

**Lemma 3.4.9** Let $E, F \subseteq \mathbb{R}^d$ be non-empty and satisfy

$$\inf_{x \in E, y \in F} |x - y| > 0.$$

Then,

$$\dim_{\mathrm{L}} E \cup F \ = \ \min\{\dim_{\mathrm{L}} E,\ \dim_{\mathrm{L}} F\},$$

and, for all $\theta \in (0,1)$,

$$\dim_{\mathrm{L}}^\theta E \cup F \ = \ \min\Big\{\dim_{\mathrm{L}}^\theta E,\ \dim_{\mathrm{L}}^\theta F\Big\}.$$

*Proof* The formulae for the lower dimension and lower spectrum are proved similarly. The inequalities $\min\{\dim_{\mathrm{L}} E, \dim_{\mathrm{L}} F\} \leqslant \dim_{\mathrm{L}}(E \cup F)$ and $\min\{\dim_{\mathrm{L}}^{\theta} E, \dim_{\mathrm{L}}^{\theta} F\} \leqslant \dim_{\mathrm{L}}^{\theta}(E \cup F)$ hold without the uniform separation assumption since, for $0 < r < R < 1$,

$$\inf_{x \in E \cup F} N_r\big(B(x,R) \cap (E \cup F)\big)$$
$$\geqslant \min\left\{ \inf_{x \in E} N_r\big(B(x,R) \cap E\big),\ \inf_{x \in F} N_r\big(B(x,R) \cap F\big) \right\}.$$

The reverse inequalities hold, assuming the separation assumption, since, for $R < \min\{\inf_{x \in E, y \in F} |x-y|,\ 1\}$, the above inequality can be replaced by equality. This is because if $x \in E$, then $B(x,R) \cap F = \varnothing$ and vice versa. □

The question of finite stability for the modified lower dimension is a little tricky. In [111, Proposition 4.3] it was proved that modified lower dimension is finitely stable for sets $E, F$ which are 'properly separated', in the sense of Lemma 3.4.9, with the general case left as an open question [111, Question 9.3]. Here we prove finite stability holds provided the sets in question are closed. The proof is based on a simple adaptation of the proof of Lemma 3.4.3.

**Lemma 3.4.10** If $E, F \subseteq \mathbb{R}^d$ are non-empty and closed, then

$$\dim_{\mathrm{ML}} E \cup F \;=\; \max\{\dim_{\mathrm{ML}} E,\ \dim_{\mathrm{ML}} F\}.$$

*Proof* The lower bound ($\geqslant$) follows immediately by monotonicity, but the upper bound ($\leqslant$) is more subtle. We may assume $0 < \dim_{\mathrm{ML}} E \cup F$, let $0 < s < \dim_{\mathrm{ML}} E \cup F$, and choose a set $X \subseteq E \cup F$ with $\dim_{\mathrm{L}} X \geqslant s$ by applying the definition of modified lower dimension. Since $E$ and $F$ are closed, the union $E \cup F$ is closed and therefore the closure $\overline{X} \subseteq E \cup F$. We now follow the proof of Theorem 3.4.3 with $F$ replaced by $\overline{X}$ to obtain a Cantor set

$$\overline{X}' = \bigcap_k \mathcal{B}_k \subseteq \overline{X} \subseteq E \cup F$$

as in (3.4) on page 40 which satisfies $\dim_{\mathrm{L}} \overline{X}' = \dim_{\mathrm{H}} \overline{X}' \geqslant s$. There are now two cases. If $\overline{X}' \subseteq E$, then

$$\dim_{\mathrm{ML}} E \geqslant \dim_{\mathrm{L}} \overline{X}' \geqslant s.$$

On the other hand, if $\overline{X}'$ is not contained in $E$, then, since $E$ is closed, there must exist $x \in \overline{X}'$ and $r > 0$ such that $B(x,r) \cap \overline{X}' = F \cap \overline{X}'$. Recall that the sets $\mathcal{B}_k$ are made up of $c^{-sk}$ many pairwise disjoint closed balls $\{B_k^i\}_i$ of radius $c^k \rho$. It follows that for some $k \geqslant 1$ and $i$, $B_k^i \cap \overline{X}' \subseteq F$. Therefore

$$\dim_{\mathrm{ML}} F \geqslant \dim_{\mathrm{L}}(B_k^i \cap \overline{X}') \geqslant s$$

where the final inequality follows due to the homoegeneity of the construction of $\overline{X}'$, see the proof of Theorem 3.4.3 for the details. We have thus established

$$\max\{\dim_{\mathrm{ML}} E,\ \dim_{\mathrm{ML}} F\} \geqslant s$$

proving the result. □

The proof of Lemma 2.4.2 can be suitably adapted to prove that the other dimensions and spectra we have met in this Chapter are stable under bi-Lipschitz maps and we leave the details as an exercise.

**Lemma 3.4.11** The lower, modified lower and quasi-Assouad dimensions, as well as the Assouad and lower spectra, are stable under bi-Lipschitz maps.

Despite being bi-Lipschitz stable, none of these notions are stable under Lipschitz maps. For the lower dimension and spectrum this is trivial, since the projection of the set $[0,1] \times \{0\} \cup \{(0,1)\} \subseteq \mathbb{R}^2$ onto the first coordinate is Lipschitz and easily seen to increase the lower dimension and spectrum from 0 to 1. For the modified lower dimension, the graph of the popcorn function again does the job, see Figure 3.5, where we again consider the projection onto the first coordinate.

Demonstrating that the Assouad spectrum can increase under Lipschitz maps is rather more complicated. Since the upper box dimension cannot increase under Lipschitz maps, using Lemma 3.4.4 we always have

$$\dim_{\mathrm{A}}^{\theta} T(F) \leqslant \min\left\{\frac{\overline{\dim}_{\mathrm{B}} T(F)}{1-\theta},\ d\right\} \leqslant \min\left\{\frac{\overline{\dim}_{\mathrm{B}} F}{1-\theta},\ d\right\} \tag{3.9}$$

for bounded sets $F \subseteq \mathbb{R}^d$ and Lipschitz maps $T$. This estimate turns out to be sharp and shows that there is some control on how badly the Assouad spectrum can increase under Lipschitz map, but the situation for the quasi-Assouad and Assouad dimensions are as wild as possible. Note that if $\overline{\dim}_{\mathrm{B}} F = 0$, then the Assouad spectrum and quasi-Assouad dimension cannot increase under Lipschitz maps by Corollary 3.4.5, but the Assouad dimension may still increase arbitrarily, see for example [103]. An example showing that the quasi-Assouad dimension can increase under Lipschitz maps can be found in [118, Proposition 12]. Here we show that the Assouad spectrum can also increase and that the upper bound (3.9) is sharp.

**Theorem 3.4.12** For given $s \in (0,1)$, there exists a compact set $F \subseteq \mathbb{R}^2$ such that for all $\theta \in (0,1)$

$$\dim_{\mathrm{A}}^{\theta} F = \overline{\dim}_{\mathrm{B}} F = \dim_{\mathrm{A}} F = s$$

and

$$\dim_{\mathrm{A}}^{\theta} \pi F = \min\left\{\frac{s}{1-\theta}, 1\right\}$$

where $\pi : \mathbb{R}^2 \to \mathbb{R}$ denotes projection onto the first coordinate. In particular, $\pi$ is Lipschitz and so the Assouad spectrum, quasi-Assouad dimension and Assouad dimension may increase under Lipschitz maps.

We delay the proof of Theorem 3.4.12 until Section 10.1 since we will rely on some material covered later in the book. We are now able to summarise the basic properties of the dimensions and spectra we have met in this Chapter, see Table 3.1.

| Property | $\dim_{\mathrm{L}}$ | $\dim_{\mathrm{ML}}$ | $\dim_{\mathrm{qL}}$ | $\dim_{\mathrm{L}}^{\theta}$ | $\dim_{\mathrm{A}}^{\theta}$ | $\dim_{\mathrm{qA}}$ |
|---|---|---|---|---|---|---|
| Monotone | × | ✓ | × | × | ✓ | ✓ |
| Finitely stable | × | ? | × | × | ✓ | ✓ |
| Countably stable | × | × | × | × | × | × |
| Lipschitz stable | × | × | × | × | × | × |
| Bi-Lipschitz stable | ✓ | ✓ | ✓ | ✓ | ✓ | ✓ |
| Stable under closure | ✓ | × | ✓ | ✓ | ✓ | ✓ |
| Open set property | × | ✓ | × | × | ✓ | ✓ |

Table 3.1 *This table summarises which properties are satisfied by which dimensions and spectra and should be compared with Table 2.1.*

### 3.4.3 Hölder distortion

Bi-Hölder maps are a natural generalisation of bi-Lipschitz maps, where more distortion is allowed. Let $X$ be a domain lying in $\mathbb{R}^d$. We say an injective map $f : X \to \mathbb{R}^d$ is $(\alpha, \beta)$-Hölder, or bi-Hölder, for $0 < \alpha \leqslant 1 \leqslant \beta < \infty$ if there exists a constant $C \geqslant 1$ such that, for all $x, y \in X$,

$$C^{-1}|x-y|^{\beta} \leqslant |f(x)-f(y)| \leqslant C|x-y|^{\alpha}. \tag{3.10}$$

Notice that $(1,1)$-Hölder is the same as bi-Lipschitz.

Dimensions are not generally preserved under bi-Hölder maps, but one can often control the distortion in terms of the Hölder exponents $\alpha$ and $\beta$. For example, if dim is the Hausdorff, packing, or upper or lower box dimension, and $f$ is $(\alpha, \beta)$-Hölder, then

$$\frac{\dim X}{\beta} \leqslant \dim f(X) \leqslant \frac{\dim X}{\alpha}, \tag{3.11}$$

see [70, Proposition 3.3]. Notably, the Assouad dimension does not satisfy such bounds, see [191, Proposition 1.2], Theorem 10.2.4, or the papers [103, 106].

The Assouad spectrum, which is inherently more regular than the Assouad dimension, *can* be controlled in this context but the control is more complicated than (3.11). The following lemma was first proved in [111, Proposition 4.7].

**Lemma 3.4.13** Suppose $f : X \to \mathbb{R}^d$ is $(\alpha, \beta)$-Hölder. Then, for all $\theta \in (0,1)$,

$$\frac{1-\beta\theta/\alpha}{\beta(1-\theta)} \dim_{\mathrm{A}}^{\beta\theta/\alpha} X \ \leqslant\ \dim_{\mathrm{A}}^{\theta} f(X) \ \leqslant\ \frac{1-\alpha\theta/\beta}{\alpha(1-\theta)} \dim_{\mathrm{A}}^{\alpha\theta/\beta} X$$

and

$$\frac{1-\beta\theta/\alpha}{\beta(1-\theta)} \dim_{\mathrm{L}}^{\beta\theta/\alpha} X \ \leqslant\ \dim_{\mathrm{L}}^{\theta} f(X) \ \leqslant\ \frac{1-\alpha\theta/\beta}{\alpha(1-\theta)} \dim_{\mathrm{L}}^{\alpha\theta/\beta} X$$

where $\dim_{\mathrm{A}}^{\beta\theta/\alpha} X$ and $\dim_{\mathrm{L}}^{\beta\theta/\alpha} X$ are taken to equal 0 if $\beta\theta/\alpha \geqslant 1$.

*Proof* Fix $\theta \in (0,1)$ and let $t > \dim_{\mathrm{A}}^{\alpha\theta/\beta} X$. Let $x \in X$ and $r \in (0,1)$. Therefore, by (3.10),

$$f(B(x, r^{\theta/\beta} C^{1/\beta})) \supseteq B(f(x), r^{\theta}),$$

and for every set $U \subseteq X$ with $|U| \leqslant (r/C)^{1/\alpha}$

$$|f(U)| \leqslant C|U|^{\alpha} \leqslant r.$$

This means that the image of an $(r/C)^{1/\alpha}$-cover of $B(x, r^{\theta/\beta}C^{1/\beta})) \cap X$ under $f$ yields an $r$-cover of $B(f(x), r^\theta) \cap f(X)$. Therefore, applying the definition of the Assouad spectrum of $X$, there is a constant $C_0$ depending only on $t$ and $C$ such that

$$N_r(B(f(x), r^\theta) \cap f(X)) \leqslant N_{(r/C)^{1/\alpha}}(B(x, r^{\theta/\beta}C^{1/\beta}) \cap X)$$

$$\leqslant C_0 \left( \frac{r^{\theta/\beta}C^{1/\beta}}{(r/C)^{1/\alpha}} \right)^t$$

$$\leqslant C_0 C^{t(1/\alpha+1/\beta)} \left(r^{\theta-1}\right)^{\frac{t(\theta/\beta-1/\alpha)}{\theta-1}} .$$

The dependence of $C_0$ on $C$ is just to account for the fact that the 'big' scale $r^{\theta/\beta}C^{1/\beta}$ is not quite the 'small' scale $(r/C)^{1/\alpha}$ to the power $\alpha\theta/\beta$. It follows that

$$\dim_{\mathrm{A}}^\theta f(X) \leqslant \frac{t(1/\alpha - \theta/\beta)}{1-\theta} = \frac{1-\theta\alpha/\beta}{\alpha(1-\theta)} t$$

and letting $t \to \dim_{\mathrm{A}}^{\alpha\theta/\beta} X$ yields the desired upper bound. The lower bound for $\dim_{\mathrm{A}}^\theta f(X)$ can be obtained by applying the upper bound to the map $f^{-1} : f(X) \to X$ with $\theta' = \beta\theta/\alpha$ for $\theta < \alpha/\beta$. This uses the fact that $f^{-1}$ is $(\beta^{-1}, \alpha^{-1})$-Hölder. The estimates for the lower spectrum are proved in a similar way and we leave the details to the reader. □

The following immediate corollary of Lemma 3.4.13 was first recorded in [111, Theorem 4.11] and allows one to recover some information about the behaviour of the Assouad dimension under Hölder images.

**Corollary 3.4.14** Suppose $f : X \to \mathbb{R}^d$ is $(\alpha, \beta)$-Hölder. Then

$$\dim_{\mathrm{A}} f(X) \geqslant \dim_{\mathrm{qA}} f(X) \geqslant \frac{(1-\rho)}{\beta - \alpha\rho} \dim_{\mathrm{qA}} X$$

where $\rho = \inf\{\theta \in (0,1) : \dim_{\mathrm{A}}^\theta X = \dim_{\mathrm{qA}} X\}$.

Another consequence of Lemma 3.4.13 is that the Assouad spectrum and quasi-Assouad dimension of compact sets are preserved under quasi-Lipschitz maps. Recall that the quasi-Assouad dimension of a compact set is preserved under quasi-Lipschitz maps. This was proved in [191], see Section 3.2. Straight from the definitions, one sees that if a map $f : X \to Y$ is quasi-Lipschitz for compact $X, Y$, then it is $(1-\varepsilon, 1+\varepsilon)$-Hölder for all $\varepsilon \in (0,1)$. This was also observed in [276].

**Corollary 3.4.15** Suppose $f : X \to \mathbb{R}^d$ is $(1-\varepsilon, 1+\varepsilon)$-Hölder for all $\varepsilon \in (0,1)$. Then, for all $\theta \in (0,1)$,

$$\dim_{\mathrm{A}}^{\theta} f(X) = \dim_{\mathrm{A}}^{\theta} X \quad \text{and} \quad \dim_{\mathrm{L}}^{\theta} f(X) = \dim_{\mathrm{L}}^{\theta} X$$

and also

$$\dim_{\mathrm{qA}} f(X) = \dim_{\mathrm{qA}} X \quad \text{and} \quad \dim_{\mathrm{qL}} f(X) = \dim_{\mathrm{qL}} X.$$

In particular, this holds if $f$ is quasi-Lipschitz and $X$ is compact.

*Proof* This follows from Lemma 3.4.13 by letting $\varepsilon \to 0$. The quasi-Assouad invariance then follows by letting $\theta \to 1$ and applying Corollary 3.3.7 and, for the quasi-lower dimension, letting $\varepsilon \to 0$. If $f$ is quasi-Lipschitz and $X$ is compact, then $Y$ is compact, $f$ is a bijection and

$$\frac{\log |f(x) - f(y)|}{\log |x - y|} \to 1$$

uniformly as $|x - y| \to 0$. Therefore, for all $\varepsilon \in (0,1)$, there exists $\delta \in (0,1)$ such that for $|x - y| \leqslant \delta$,

$$|x - y|^{1+\varepsilon} \leqslant |f(x) - f(y)| \leqslant |x - y|^{1-\varepsilon},$$

which does not quite prove that $f$ is $(1-\varepsilon, 1+\varepsilon)$-Hölder since we need an estimate like this for *all* $x, y \in X$, not just $x, y$ which are $\delta$-close. However, we can 'upgrade' the upper bound since if $|x - y| > \delta$, then

$$|f(x) - f(y)| \leqslant |f(X)| = \frac{|f(X)|}{\delta^{1-\varepsilon}} \delta^{1-\varepsilon} \leqslant \frac{|f(X)|}{\delta^{1-\varepsilon}} |x - y|^{1-\varepsilon}.$$

For the lower bound, suppose

$$\inf_{x,y \in X : |x-y| > \delta} \frac{|f(x) - f(y)|}{|x - y|^{1+\varepsilon}} = 0.$$

Therefore, since $X$ is compact, we can find $x, y \in X$ with $|x - y| \geqslant \delta$ such that $f(x) = f(y)$, which contradicts the assumptions that $f$ is a bijection. This completes the proof that $f$ is $(1-\varepsilon, 1+\varepsilon)$-Hölder for all $\varepsilon \in (0,1)$ and the corollary follows. □

### 3.4.4 An example with non-monotonic Assouad spectrum

The following unpublished example is due to Han Yu.[5] For $n \in \mathbb{N}$, let $c_n = 2^{-2^n}$ and $A_n$ be a finite set formed as follows, see Figure 3.6. First take the interval $[0, c_n]$ and inside place $c_n^{-1/2}$ many evenly spaced intervals of length $c_n^2$. Then, inside each of these intervals place $c_n^{-1}$ many evenly spaced points with distance $c_n^3$ between them. Note that $c_{n+1} = c_n^2$, which will be useful below. Finally, let

$$F = \bigcup_{n=1}^{\infty} \left(A_n + 2^{-n}\right) \subseteq \mathbb{R}.$$

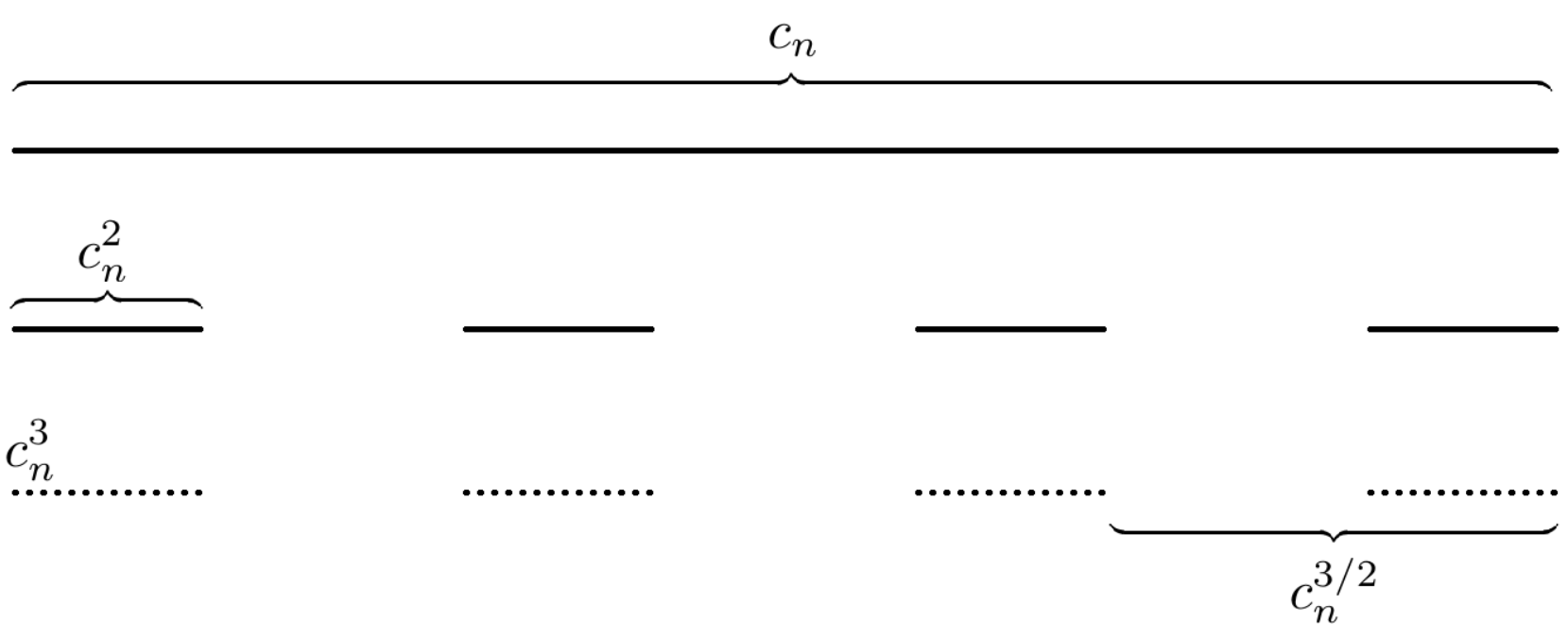

Figure 3.6 A figure showing the construction of $A_n$. For obvious reasons this is not to scale!

**Theorem 3.4.16** The Assouad spectrum of $F$ is not monotonic and is given by

$$\dim_{\mathrm{A}}^{\theta} F = \begin{cases} \frac{1/2}{1-\theta} & \theta \in (0, 1/3] \\ \frac{5/6-\theta}{1-\theta} & \theta \in (1/3, 1/2] \\ \frac{1/3}{1-\theta} & \theta \in (1/2, 2/3] \\ 1 & \theta \in (2/3, 1) \end{cases}$$

To prove Theorem 3.4.16 one should keep Theorem 3.3.1 and Corollary 3.3.4 in mind. In fact this example also serves to demonstrate the sharpness of these results and, moreover, by applying Theorem 3.3.1 it suffices to prove

[5] This example was originally intended for [111, 280] but was replaced by a more complicated example where the spectrum was not even (finitely) piecewise monotonic.

(i) $\overline{\dim}_{\mathrm{B}} F \leqslant 1/2$
(ii) $\dim_{\mathrm{A}}^{1/3} F \geqslant 3/4$
(iii) $\dim_{\mathrm{A}}^{1/2} F \leqslant 2/3$
(iv) $\dim_{\mathrm{A}}^{2/3} F \geqslant 1$,

which greatly simplifies the argument.

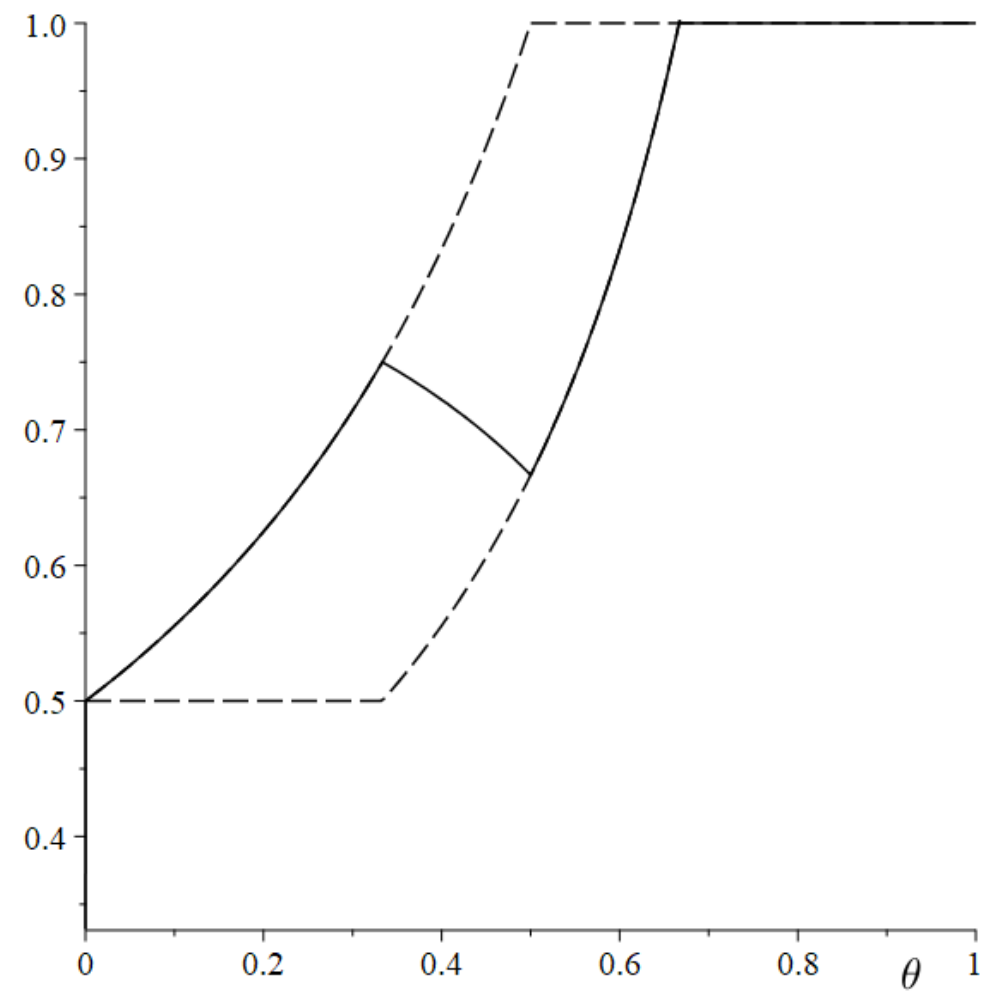

Figure 3.7 A plot of $\dim_{\mathrm{A}}^{\theta} F$ (solid). The general upper and lower bounds from Lemma 3.4.4 as well as the lower bound from Corollary 3.3.4 based on knowledge of $\rho = 2/3$ are shown as dashed lines.

*Proof* We consider each of the cases (i)–(iv). The awkward cases are (i) and (iii) since here we need to understand all locations and scales. Cases (ii) and (iv), on the other hand, are rather straightforward due to the construction of $F$, more specifically, $A_n$, making it easy to find locations and scales where the set looks big.

(i) Let $R \in (0, 1/2)$, which we may assume is the reciprocal of a power of 2, and let $k \geqslant 1$ be the unique integer satisfying

$$c_{k+1} \leqslant R < c_k.$$

For $n \geqslant k+1$ each set $A_n$ may be covered by a single interval of length $R$ and so

$$N_R\left(\bigcup_{n=k+1}^{\infty} (A_n + 2^{-n})\right) \leqslant N_R\left(\{2^{-n} : n \geqslant 1\}\right) \leqslant 3|\log_2 R|.$$

For the set $A_k$, if $c_k^{3/2} \leqslant R \leqslant c_k$, then

$$N_R(A_k) \leqslant \frac{c_k}{R} \leqslant \frac{c_k}{c_k^{3/2}} = c_k^{-1/2} \leqslant R^{-1/2}$$

and, if $c_k^2 = c_{k+1} \leqslant R < c_k^{3/2}$, then

$$N_R(A_k) \leqslant c_k^{-1/2} \leqslant R^{-1/2}.$$

For the set $A_{k-1}$, we have $c_{k-1}^4 \leqslant R < c_{k-1}^2$. If $c_{k-1}^3 \leqslant R \leqslant c_{k-1}^2$, then

$$N_R(A_{k-1}) \leqslant \frac{c_{k-1}^2}{R} c_{k-1}^{-1/2} = \frac{c_{k-1}^{3/2}}{R} \leqslant R^{-1/2}$$

and, if $c_{k-1}^4 \leqslant R < c_{k-1}^3$, then

$$N_R(A_{k-1}) \leqslant \# A_{k-1} = c_{k-1}^{-3/2} \leqslant R^{-1/2}.$$

Finally, for the sets $A_n$ with $n \leqslant k-2$,

$$N_R(A_n) \leqslant \# A_n = c_n^{-3/2}.$$

Moreover, since $R < c_k = c_{k-2}^4 \leqslant c_n^3$ we have $c_n^{-3/2} \leqslant R^{-1/2}$, and so

$$\sum_{n=1}^{k-2} N_R(A_n) \leqslant \sum_{n=1}^{k-2} c_n^{-3/2} \leqslant R^{-1/2} k \leqslant R^{-1/2} \log_2 |\log_2 R|.$$

Putting these estimates together proves $\overline{\dim}_{\mathrm{B}} F \leqslant 1/2$, as required.

(ii) Let $\theta = 1/3$, $R = c_n$ and $x = 2^{-n} \in (A_n + 2^{-n}) \subseteq F$. Since $R^3 = c_n^3$,

$$N_{R^3}(B(x,R) \cap F) \geqslant \# A_n = c_n^{-1/2} c_n^{-1} = \left(\frac{R}{R^3}\right)^{3/4}$$

and taking $n$ arbitrarily large we conclude $\dim_{\mathrm{A}}^{1/3} F \geqslant 3/4$, as required.

(iii) Let $\theta = 1/2$, $R \in (0, 1/2)$ which we again may assume is the reciprocal of a power of 2, and $x \in A_n \subseteq F$ for $n \in \mathbb{N}$ which we may assume is large. Suppose $B(x,R)$ intersects $A_{n'}$ for some $n' \neq n$. Therefore $R > 2^{-n}/2$ and $R^2 > 2^{-2n}/4 > c_n$ and so

$$\begin{aligned} N_{R^2}(B(x,R) \cap F) \leqslant N_{R^2}(B(0,5R) \cap F) &\leqslant N_{R^2}\left(\{2^{-n} : n \geqslant 1\}\right) \\ &\leqslant 6 |\log_2 R|. \end{aligned}$$

On the other hand, suppose $B(x,R)$ only intersects $A_n$, in which case we consider the following seven cases. In order to follow the subsequent covering arguments in each case, we recommend referring to Figure 3.6.

(a) If $R < c_n^3$, then $N_{R^2}(B(x,R) \cap F) = 1$.

(b) If $c_n^3 \leqslant R < c_n^2$, then $R^2 < c_n^4 < c_n^3$ and therefore

$$N_{R^2}(B(x,R) \cap F) \leqslant \frac{2R}{c_n^3} = 2\frac{R}{(c_n^4)^{3/4}} \leqslant 2\frac{R}{(R^2)^{3/4}} = 2\left(\frac{R}{R^2}\right)^{1/2}.$$

(c) If $c_n^2 \leqslant R < c_n^{3/2}$, then $R^2 < c_n^3$ and therefore

$$N_{R^2}(B(x,R) \cap F) \leqslant 3c_n^{-1} \leqslant 3\left(\frac{R}{R^2}\right)^{2/3}.$$

(d) If $c_n^{3/2} \leqslant R < c_n$, then $c_n^3 \leqslant R^2 < c_n^2$ and therefore

$$N_{R^2}(B(x,R) \cap F) \leqslant \frac{2R}{c_n^{3/2}}\frac{c_n^2}{R^2} \leqslant 2\left(\frac{R}{R^2}\right)^{2/3}.$$

(e) If $c_n \leqslant R < c_n^{3/4}$, then $c_n^2 \leqslant R^2 < c_n^{3/2}$ and therefore

$$N_{R^2}(B(x,R) \cap F) \leqslant c_n^{-1/2} \leqslant \left(\frac{R}{R^2}\right)^{2/3}.$$

(f) If $c_n^{3/4} \leqslant R < c_n^{1/2}$, then $c_n^{3/2} \leqslant R^2 < c_n$ and therefore

$$N_{R^2}(B(x,R) \cap F) \leqslant \frac{c_n}{R^2} \leqslant \left(\frac{R}{R^2}\right)^{2/3}.$$

(g) If $R \geqslant c_n^{1/2}$, then $R^2 \geqslant c_n$ and so

$$N_{R^2}(B(x,R) \cap F) \leqslant 1.$$

Putting together all of these cases we conclude $\dim_{\mathrm{A}}^{1/2} F \leqslant 2/3$, as required.

(iv) Let $\theta = 2/3$, $R = c_n$ and $x = 2^{-n} \in (A_n + 2^{-n}) \subseteq F$. Since $R^{3/2} = c_n^{3/2}$,

$$N_{R^{3/2}}(B(x,R) \cap F) \geqslant \frac{c_n}{c_n^{3/2}} = \left(\frac{R}{R^{3/2}}\right)$$

and taking $n$ arbitrarily large we conclude $\dim_{\mathrm{A}}^{2/3} F \geqslant 1$, as required.

□

# 4
# Dimensions of measures

A key aspect of dimension theory is how dimensions of sets relate to dimensions of measures. To this end, all the familiar notions of dimension for sets have natural counterparts in the context of measures, recall the discussion Section 1.2.1. The Assouad dimension and its many variants are no different and in this chapter we explore the Assouad type dimensions and spectra of measures.

## 4.1 Assouad and lower dimensions of measures

Throughout we assume $\mu$ is a locally finite Borel measure on $\mathbb{R}^d$, although much of the theory applies more generally. A measure $\mu$ on $\mathbb{R}^d$ is *locally finite* if $\mu(B(x,r)) < \infty$ for all $x \in \mathbb{R}^d$ and $r > 0$. The Assouad and lower dimensions of $\mu$ describe the optimal global control on the relative measure of concentric balls. These dimensions were introduced formally in [153, 152] motivated by previous work on the existence of doubling measures, see for example [193, 173, 40].[1]

As usual, we write $|E|$ for the diameter of a non-empty (possibly unbounded) set $E$ and $\mathrm{supp}(\mu)$ for the support of a measure $\mu$, which is

[1] When they were first introduced, the Assouad and lower dimensions of measures were referred to as the upper and lower regularity dimensions, respectively. This terminology is still popular, but not universally. In fact, when I worked with these notions in the past I used these names. However, I now favour the Assouad and lower dimension being used. This is for two reasons: firstly, it is important that the community agrees on a common set of terminology and notation and, secondly, using the same terminology for sets and measures serves to emphasise the connection and interplay between the notions.

necessarily closed, see (1.2). The *Assouad dimension* of $\mu$ is defined by

$$\dim_{\mathrm{A}} \mu = \inf \Bigg\{ s \geqslant 0 \;:\; \text{there exists } C > 0 \text{ such that, for all } 0 < r < R < |\mathrm{supp}(\mu)| \text{ and } x \in \mathrm{supp}(\mu), \; \frac{\mu(B(x,R))}{\mu(B(x,r))} \leqslant C \left(\frac{R}{r}\right)^s \Bigg\}$$

and, provided $|\mathrm{supp}(\mu)| > 0$, the *lower dimension* of $\mu$ is defined by

$$\dim_{\mathrm{L}} \mu = \sup \Bigg\{ s \geqslant 0 \;:\; \text{there exists } C > 0 \text{ such that, for all } 0 < r < R < |\mathrm{supp}(\mu)| \text{ and } x \in \mathrm{supp}(\mu), \; \frac{\mu(B(x,R))}{\mu(B(x,r))} \geqslant C \left(\frac{R}{r}\right)^s \Bigg\}$$

and otherwise it is 0. We adopt the convention that $\inf \varnothing = +\infty$.

There is a straightforward connection between the Assouad and lower dimensions of sets and measures, which also provides a powerful method to estimate the dimensions of sets. This is in some sense dual to the approach provided by weak tangents (discussed in Chapter 5) since the estimates are in the opposite direction. Recall that a measure $\mu$ is *doubling* (on its support) if there is a constant $c \geqslant 1$ such that, for all $x \in \mathrm{supp}(\mu)$ and $r > 0$,

$$\mu(B(x, 2r)) \leqslant c\mu(B(x, r)).$$

The following fundamental equivalence was first noted in [147, Lemma 3.2].

**Lemma 4.1.1** A measure $\mu$ is doubling if and only if $\dim_{\mathrm{A}} \mu < \infty$.

*Proof* The fact that finite Assouad dimension implies doubling is immediate from the definition. The other implication requires a little more work. Suppose $\mu$ is doubling with doubling constant $c \geqslant 1$. Let $x \in \mathrm{supp}(\mu)$ and $0 < r < R$. Let $n \geqslant 0$ be the unique integer satisfying

$$2^{-n}R \leqslant 2r < 2^{-n+1}R,$$

noting that $n \leqslant \log(R/r)/\log 2$. Then

$$\begin{aligned}\frac{\mu(B(x,R))}{\mu(B(x,r))} &= \frac{\mu(B(x,2^{-n}R))}{\mu(B(x,r))} \prod_{k=0}^{n-1} \frac{\mu(B(x,2^{-k}R))}{\mu(B(x,2^{-(k+1)}R))} \\ &\leqslant c^{n+1} \\ &\leqslant c\left(\frac{R}{r}\right)^{\frac{\log c}{\log 2}}\end{aligned}$$

which proves $\dim_A \mu \leqslant \frac{\log c}{\log 2} < \infty$ as required. □

It will be important for us that the Assouad and lower dimensions of a measure control the Assouad and lower dimensions of its support, see [153, Section 3].

**Lemma 4.1.2** If a measure $\mu$ is doubling and fully supported on a closed set $F$, then

$$\dim_L \mu \leqslant \dim_L F \leqslant \dim_A F \leqslant \dim_A \mu.$$

*Proof* We begin by proving the first inequality. Let $t < \dim_L \mu$ and consider $x \in F$ and $R > r > 0$. Let $\{U_i\}_{i=1}^N$ be a minimal $r$-cover of $B(x,R) \cap F$. Therefore

$$\mu(B(x,R)) \leqslant \sum_{i=1}^{N} \mu(U_i) \leqslant N \sup_{y \in B(x,R)\cap F} \mu(B(y,r)).$$

Moreover, since $\mu$ is assumed to be doubling, for all $y \in B(x,R) \cap F$, $\mu(B(x,R)) \geqslant c^{-1}\mu(B(x,2R)) \geqslant c^{-1}\mu(B(y,R))$, where $c$ is the doubling constant of $\mu$. Therefore, for all such $y$,

$$\frac{\mu(B(x,R))}{\mu(B(y,r))} \geqslant c^{-1}\frac{\mu(B(y,R))}{\mu(B(y,r))} \geqslant C\left(\frac{R}{r}\right)^t$$

for a constant $C$ which is independent of $x, y, r$ and $R$. It follows that

$$N_r(B(x,R)\cap F) = N \geqslant C\left(\frac{R}{r}\right)^t$$

which proves $\dim_L F \geqslant t$ and letting $t \to \dim_L \mu$ yields the desired inequality.

The middle inequality is trivial and so it remains to prove the final inequality. Let $s > \dim_A \mu$ and consider $x \in F$ and $R > r > 0$. Suppose $\{x_i\}_{i=1}^N$ is a collection of $N$ points in $F$ such that the balls $B(x_i, r)$ are pairwise disjoint and contained in $B(x, R)$. Therefore

$$\mu(B(x,R)) \geqslant \sum_{i=1}^{N} \mu(B(x_i,r)) \geqslant N \min_i \mu(B(x_i,r)) > 0.$$

Moreover, for all $i$, $B(x,R) \subseteq B(x_i, 2R)$ and therefore

$$\frac{\mu(B(x,R))}{\mu(B(x_i,r))} \leqslant \frac{\mu(B(x_i,2R))}{\mu(B(x_i,r))} \leqslant C\left(\frac{R}{r}\right)^s$$

for a constant $C$ which is independent of $x, i, r$ and $R$. It follows that $N \leqslant C\left(\frac{R}{r}\right)^s$ and, moreover, that

$$N_{6r}(B(x,R)\cap F) \leqslant C\left(\frac{R}{r}\right)^s$$

since if the collection $\{x_i\}$ is chosen such that $N$ is maximised, then the collection of balls $B(x_i, 3r)$ is a $6r$-cover of $B(x,R) \cap F$. This is enough to prove $\dim_{\mathrm{A}} F \leqslant s$ and letting $s \to \dim_{\mathrm{A}} \mu$ yields the desired inequality. □

A much deeper fact is that the previous lemma is sharp, allowing us to express the dimensions of sets in terms of the corresponding dimensions of measures supported on the sets. Recall, for example, that $\dim_{\mathrm{H}} F = \sup\{\dim_{\mathrm{H}} \mu : \mathrm{supp}(\mu) \subseteq F\}$, see (1.3).

**Theorem 4.1.3** For closed $F \subseteq \mathbb{R}^d$,

$$\dim_{\mathrm{A}} F = \inf\{\dim_{\mathrm{A}} \mu : \mathrm{supp}(\mu) = F\}$$

and

$$\dim_{\mathrm{L}} F = \sup\{\dim_{\mathrm{L}} \mu : \mathrm{supp}(\mu) = F\}.$$

We omit the proof of this result. This Assouad dimension result was proved by Luukkainen and Saksman [193] (in general) and Konyagin and Vol'berg [173] (for compact $F$) and the lower dimension result was proved by Bylund and Gudayol [40]. See [152] for further extensions and clarifications.

It should be clear that the measures $\mu$ considered in Theorem 4.1.3 must be fully supported (compare with (1.3) which does not require this). Given Theorem 4.1.3, it is natural to ask if there always exists a measure which *precisely* realises the Assouad and lower dimensions, that is, with $\dim_{\mathrm{A}} \mu = \dim_{\mathrm{A}} F$ or $\dim_{\mathrm{L}} \mu = \dim_{\mathrm{L}} F$. However, this is not always possible, at least for the Assouad dimension. Konyagin and Vol'berg [173, Theorem 4] provided an example of a compact set $F$ for which $\dim_{\mathrm{A}} \mu > \dim_{\mathrm{A}} F$ for all doubling measures $\mu$ supported on $F$. Käenmäki and Lehrbäck [152] observed that the example from [173] is a special case of an inhomogeneous self-similar set, and that this class of set exhibits this surprising behaviour quite generally. Inhomogeneous self-similar sets are a generalisation of the self-similar sets we will study in Chapter 7. They were introduced by Barnsley and Demko [24] (see also [128]) and their dimension theory has been studied in some detail, see [39, 87, 220, 255]. Käenmäki and Lehrbäck [152] also observed that inhomogeneous self-similar sets with suitable separation conditions cannot be used to exhibit examples where the *lower* dimension is not precisely realised, and it remains an interesting open question to determine whether or not the lower dimension of a set can always be precisely realised by a doubling measure, see Question 17.2.1.

Following on from Theorem 4.1.3, another natural question is: given a set $F$, what are the possible values of $\dim_{\mathrm{A}} \mu$ for measures supported on $F$? We know from Theorem 4.1.3 that this set of values is contained in $[\dim_{\mathrm{A}} F, \infty]$ and may or may not include the point $\dim_{\mathrm{A}} F$. Hare, Mendivil and Zuberman [126] proved that if $F \subseteq \mathbb{R}$ is a compact set with $\dim_{\mathrm{A}} F > 0$, then for all $s \in (\dim_{\mathrm{A}} F, \infty)$, there exists a measure $\mu$ on $F$ with $\dim_{\mathrm{A}} \mu = s$. It was also shown in [126] that if $F \subseteq \mathbb{R}$ is a compact set with $\dim_{\mathrm{L}} F > 0$, then for all $t \in [0, \dim_{\mathrm{L}} F)$, there exists a measure $\mu$ on $F$ with $\dim_{\mathrm{L}} \mu = t$. Interestingly, the assumption that $\dim_{\mathrm{A}} F > 0$ is necessary for the first result, since they provide an example of a set $E \subseteq \mathbb{R}$ with $\dim_{\mathrm{A}} E = 0$ such that

$$\{\dim_{\mathrm{A}} \mu : \mathrm{supp}(\mu) = E\} = \{0, \infty\}.$$

In fact, one can choose

$$E = \{2^{-2^n} : n \geqslant 1\} \cup \{0\}.$$

These results were extended to arbitrary complete metric spaces $F$ by Suomala [262].

## 4.2 Assouad spectrum and box dimensions of measures

By now it should be obvious how to define the Assouad and lower spectra of measures and one might even be able to guess some of the relationships between these spectra and the Assouad and lower dimensions. To make explicit, for $\theta \in (0,1)$, the *Assouad spectrum* of $\mu$ is defined by

$$\dim_{\mathrm{A}}^{\theta} \mu = \inf \Bigg\{ s \geqslant 0 \ : \ \text{there exists } C > 0 \text{ such that, for all } 0 < r < |\mathrm{supp}(\mu)| \text{ and } x \in \mathrm{supp}(\mu), \ \frac{\mu(B(x, r^{\theta}))}{\mu(B(x,r))} \leqslant C \left( \frac{r^{\theta}}{r} \right)^{s} \Bigg\}$$

and, provided $|\mathrm{supp}(\mu)| > 0$, the *lower spectrum* of $\mu$ is defined by

$$\dim_{\mathrm{L}}^{\theta} \mu = \sup \Bigg\{ s \geqslant 0 \ : \ \text{there exists } C > 0 \text{ such that, for all } 0 < r < |\mathrm{supp}(\mu)| \text{ and } x \in \mathrm{supp}(\mu), \ \frac{\mu(B(x, r^{\theta}))}{\mu(B(x,r))} \geqslant C \left( \frac{r^{\theta}}{r} \right)^{s} \Bigg\}$$

and otherwise it is 0. The following basic properties may be established following arguments similar to those we have already documented, and also appear in [75].

**Lemma 4.2.1** For any measure $\mu$,

$$\dim_{\mathrm{L}} \mu \leqslant \dim_{\mathrm{L}}^{\theta} \mu \leqslant \dim_{\mathrm{A}}^{\theta} \mu \leqslant \dim_{\mathrm{A}} \mu$$

and if a measure $\mu$ is fully supported on a closed set $F$, then

$$\dim_{\mathrm{L}}^{\theta} \mu \leqslant \dim_{\mathrm{L}}^{\theta} F \leqslant \dim_{\mathrm{A}}^{\theta} F \leqslant \dim_{\mathrm{A}}^{\theta} \mu.$$

An interesting feature of the Assouad spectrum of measures is that, in searching for an analogue of Lemma 3.4.4, we are led to the appropriate definition of the *upper box* dimension of a measure. Specifically, the upper box dimension of a measure $\mu$ is

$$\overline{\dim}_{\mathrm{B}}\mu = \inf\Big\{ \alpha : \text{ there exists a constant } C>0 \text{ such that, for all } 0<r<|\mathrm{supp}(\mu)| \text{ and } x\in\mathrm{supp}(\mu), \ \mu(B(x,r)) \geqslant Cr^{\alpha} \Big\}.$$

The upper box dimension for measures was introduced in [75], although it was referred to as the upper Minkowski dimension. A natural dual notion is the lower box dimension of $\mu$ defined by

$$\underline{\dim}_{\mathrm{B}}\mu = \inf\Big\{ \alpha : \text{ there exists a constant } C>0 \text{ such that, for all } r_0>0, \text{ there exists } 0<r<r_0 \text{ such that, for all } x\in\mathrm{supp}(\mu),\ \mu(B(x,r)) \geqslant Cr^{\alpha} \Big\}.$$

If $\overline{\dim}_{\mathrm{B}}\mu = \underline{\dim}_{\mathrm{B}}\mu$, then we refer to the common value as the box dimension of $\mu$, and denote it by $\dim_{\mathrm{B}}\mu$. Different definitions of the box dimensions of measures have been proposed in the past, for example, see [234, Section 7.1]. However, we favour the definitions above since they seem consistent with the theory of the Assouad spectrum which we are developing. The following analogue of Lemma 3.4.4 was proved in [75].

**Lemma 4.2.2** For a measure $\mu$ supported on a compact set $F$ and $\theta \in (0,1)$,

$$\overline{\dim}_{\mathrm{B}}\mu \leqslant \dim_{\mathrm{A}}^{\theta}\mu \leqslant \frac{\overline{\dim}_{\mathrm{B}}\mu}{1-\theta}.$$

*Proof* We start by proving the second inequality. Fix $0<\theta<1$ and let $s > \overline{\dim}_{\mathrm{B}}\mu/(1-\theta)$. This means that there exists a constant $C>0$ such that $\mu(B(x,r)) \geqslant Cr^{(1-\theta)s}$ for all $x \in F$ and $0<r<1$. Therefore,

$$\frac{\mu(B(x,r^{\theta}))}{\mu(B(x,r))} \leqslant \mu(F)C^{-1}r^{(\theta-1)s}$$

for all $x\in F$ and $0<r<1$, which yields $\dim_{\mathrm{A}}^{\theta}\mu \leqslant s$ as required.

The first inequality, similar to Lemma 3.4.4, is more awkward. Fix $0<\theta<1$ and let $\dim_{\mathrm{A}}^{\theta}\mu < s < t$. This means that there is $C>0$ such that

$$\frac{\mu(B(x,r^{\theta}))}{\mu(B(x,r))} \leqslant Cr^{(\theta-1)s}$$

for all $x \in F$ and $0 < r < 1$. Since $F$ is compact, we may assume by rescaling if necessary that $|F| < 1/2$. Fix $0 < r < 1$ and choose $k \geqslant 1$ such that $r^{\theta^{k-1}} < 1/2 \leqslant r^{\theta^k}$. Therefore

$$\begin{aligned} \mu(B(x,r)) &= \frac{\mu(B(x,r))}{\mu(B(x,r^\theta))} \frac{\mu(B(x,r^\theta))}{\mu(B(x,r^{\theta^2}))} \cdots \frac{\mu(B(x,r^{\theta^{k-1}}))}{\mu(B(x,r^{\theta^k}))} \mu(B(x,r^{\theta^k})) \\ &\geqslant C^{-k} r^{(1-\theta)s} r^{\theta(1-\theta)s} \cdots r^{\theta^{k-1}(1-\theta)s} \mu(B(x,\tfrac{1}{2})) \\ &\geqslant C^{-1} \left( \frac{-\log 2}{\log r} \right)^{\frac{\log C}{\log \theta}} r^{(1-\theta^k)s} \mu(F) \\ &\geqslant \mu(F) C^{-1} \left( \frac{\log 2}{-\log r} \right)^{\frac{\log C}{\log \theta}} \frac{1}{r^{-(s-t)}} r^t \end{aligned}$$

for all $x \in F$. It follows that there is a constant $C_0 > 0$ such that $\mu(B(x,r)) \geqslant C_0 r^t$ for all $x \in F$ and $0 < r < 1$ and hence, $\overline{\dim}_{\mathrm{B}} \mu \leqslant t$ completing the proof. □

We may express the upper box dimension of a set in terms of the upper box dimension of measures fully supported on the set. This is further motivation for the definition of upper box dimension for measures and should be compared with (1.3) for Hausdorff dimension and Theorem 4.1.3 for Assouad dimension. One direction follows from Lemma 4.2.1 upon letting $\theta \to 0$, and the other direction follows from Tricot's 'anti-Frostman lemma' [265, Lemma 4], although the proof is only sketched there (see also [44, Lemma 3.2] and [75]). We prove the upper bound here directly using an elegant argument due to Kenneth Falconer, see also [75].

**Theorem 4.2.3** For a non-empty compact set $F \subseteq \mathbb{R}^d$,

$$\overline{\dim}_{\mathrm{B}} F = \inf \left\{ \overline{\dim}_{\mathrm{B}} \mu : \ \mathrm{supp}(\mu) = F \right\}.$$

*Proof* The upper bound ($\leqslant$) follows from the analogous result for the Assouad spectrum given in Lemma 4.2.1 and it remains to prove the upper bound ($\geqslant$). Let $t > s > \overline{\dim}_{\mathrm{B}} F$ and for all $k \in \mathbb{N}$ let $\{x_{k,i}\}_{i=1}^{N_k}$ be a maximal $2^{-k}$-separated subset of $F$. We can choose $c > 0$ such that for all $k$

$$N_k \leqslant c 2^{ks}.$$

Let

$$\mu = \sum_{k \in \mathbb{N}} \sum_{i=1}^{N_k} 2^{-kt} \delta_{x_{k,i}}$$

noting that this is a finite measure since

$$\mu(\mathbb{R}^d) = \sum_{k \in \mathbb{N}} \sum_{i=1}^{N_k} 2^{-kt} \leqslant \sum_{k \in \mathbb{N}} c 2^{ks} 2^{-kt} < \infty$$

and, moreover, $\mathrm{supp}(\mu) = F$. Finally, given $x \in F$ and $r \in (0,1)$, choose the unique $k \in \mathbb{N}$ such that $2^{-k} < r \leqslant 2^{-k+1}$, and observe that there must exist $i \in \{1, \dots, N_k\}$ such that $x_{k,i} \in B(x,r)$. Therefore

$$\mu(B(x,r)) \geqslant 2^{-kt} \geqslant (1/2)^t r^t$$

proving $\overline{\dim}_{\mathrm{B}} \mu \leqslant t$, as required. □

The 'tale of two spectra' analysis for measures, including the introduction and discussion of the quasi-Assouad and quasi-lower dimensions of measures was conducted in [127].

# 5
# Weak tangents and microsets

One of the most effective ways to bound the Assouad dimension of a set from below is to use *weak tangents*; an approach pioneered by Mackay and Tyson [195] with ideas going back to Keith and Laakso [164]. Weak tangents are tools for capturing the extremal local structure of a set which should make their connection to Assouad and lower dimensions evident. In this chapter we explore this connection. We prove Mackay and Tyson's estimate, see Theorem 5.1.2, and a stronger result of Furstenberg which characterises the Assouad dimension entirely in terms of the Hausdorff dimension of weak tangents, see Theorem 5.1.3. This latter theorem has many useful applications, see Theorem 10.2.1 concerning orthogonal projections in Section 10.2 and Theorem 11.1.3 concerning distance sets in Section 11.1. We also explore the connections between weak tangents and the lower dimension in Section 5.2, the Assouad spectra in Section 5.3, and the dimensions of measures in Section 5.4.

## 5.1 Weak tangents and the Assouad dimension

Let $\mathcal{K}(\mathbb{R}^d)$ denote the set of all non-empty compact subsets of $\mathbb{R}^d$. This is a complete metric space when equipped with the *Hausdorff metric* $d_{\mathcal{H}}$, see [242, Chapter 2] or [205, Section 4.13], defined by

$$d_{\mathcal{H}}(A, B) = \inf\{\delta > 0 \ : \ A \subseteq B_\delta \text{ and } B \subseteq A_\delta\} \tag{5.1}$$

where, for any $C \in \mathcal{K}(\mathbb{R}^d)$,

$$C_\delta \ = \ \{x \in \mathbb{R}^d \ : \ |x - y| < \delta \text{ for some } y \in C\}$$

denotes the open $\delta$-neighbourhood of $C$. We will also consider the space $\mathcal{K}(X)$ for a fixed non-empty compact set $X \subseteq \mathbb{R}^d$. This is the set of all non-empty compact subsets of $X$ and, importantly, this is a compact subset of $\mathcal{K}(\mathbb{R}^d)$.

**Definition 5.1.1** Let $X \in \mathcal{K}(\mathbb{R}^d)$ be a fixed reference set (usually the closed unit ball or cube) and let $E, F \subseteq \mathbb{R}^d$ be closed sets with $E \subseteq X$. Suppose there exists a sequence of similarity maps $T_k : \mathbb{R}^d \to \mathbb{R}^d$ such that $d_{\mathcal{H}}(E, T_k(F) \cap X) \to 0$ as $k \to \infty$. Then $E$ is called a *weak tangent* to $F$.

Recall that a *similarity map* $T : \mathbb{R}^d \to \mathbb{R}^d$ is a map of the form $x \mapsto cO(x) + t$ where $c > 0$ is a scalar (called the *similarity ratio*), $O \in \mathcal{O}(\mathbb{R}, d)$ is a real orthogonal matrix and $t \in \mathbb{R}^d$ is a translation. For example, similarities mapping $\mathbb{R}$ to itself are maps of the form $T(x) = cx + t$ for $c \in \mathbb{R} \setminus \{0\}$ and $t \in \mathbb{R}$. We say $T$ is a *homothety* if $O$ is the identity. Observe that for all $x, y \in \mathbb{R}^d$,

$$|T(x) - T(y)| = c|x - y|.$$

The following is a fundamental tool in the study of Assouad dimension and was proved in [195, Proposition 6.1.5].

**Theorem 5.1.2** Let $F \subseteq \mathbb{R}^d$ be closed, $E \subseteq \mathbb{R}^d$ be compact, and suppose $E$ is a weak tangent to $F$. Then $\dim_A F \geqslant \dim_A E$.

*Proof* Let $s > \dim_A F$ and $c_k > 0$ denote the similarity ratio of $T_k$. It follows that there exists a constant $C > 0$ depending only on $s$ such that, for all $k \in \mathbb{N}$, $0 < r < R$ and $x \in T_k(F)$,

$$N_r\big(B(x, R) \cap T_k(F)\big) \;\leqslant\; C\left(\frac{R}{r}\right)^s.$$

Recall the proof of Lemma 2.4.2 where images of covers under bi-Lipschitz maps are considered. Fix $0 < r < R$ and $y \in E$. Choose $k \in \mathbb{N}$ such that $d_{\mathcal{H}}(E, T_k(F) \cap X) < r/2$. It follows that there exists $x \in T_k(F) \cap X$ such that $B(y, R) \cap E \subseteq B(x, 2R)$ and hence, given any $r/2$-cover of $B(x, 2R) \cap T_k(F)$, we may find an $r$-cover of $B(y, R) \cap E$ by the same number of sets (or fewer). Thus,

$$\begin{aligned} N_r\big(B(y, R) \cap E\big) \;&\leqslant\; N_{r/2}\big(B(x, 2R) \cap T_k(F)\big) \leqslant C\left(\frac{2R}{r/2}\right)^s \\ &= C\,4^s\left(\frac{R}{r}\right)^s \end{aligned}$$

which proves that $\dim_A E \leqslant s$ and therefore $\dim_A E \leqslant \dim_A F$, as required. □

Theorem 5.1.2 has a remarkable converse: one can always achieve the Assouad dimension as the Hausdorff dimension of a weak tangent. This beautiful fact was first explicitly stated by Käenmäki, Ojala and Rossi [154, Proposition 5.7], but has origins in the work of Furstenberg [114, 115, 116] and Bishop and Peres [34], see the discussion in Section 2.3. It has turned out to be a very useful technical tool and often allows one to pass a geometrical problem to the level of weak tangents. This approach provides a mechanism for finding 'Assouad dimension analogues' of known results concerning the Hausdorff dimension (and other dimensions) as well as proving stronger results, not known for other dimensions, see [89] for examples of this approach in action.

The argument used to prove [154, Proposition 5.7] is based on an argument of Bishop and Peres [34, Lemma 2.4.4], which in turn takes inspiration from Furstenberg's work. Bishop and Peres [34, Lemma 2.4.4] proved that, given a compact set $F \subseteq \mathbb{R}$, there exists a weak tangent with Hausdorff dimension equal to the *upper box* dimension of $F$. We present a modification of the proof of Käenmäki, Ojala and Rossi [154]. See also [155, Corollary 5.2] and [154, Propositions 5.8].

**Theorem 5.1.3** Let $F \subseteq \mathbb{R}^d$ be closed and non-empty with $\dim_A F = s \in [0, d]$. Then there exists a compact set $E \subseteq \mathbb{R}^d$ with $\mathcal{H}^s(E) > 0$ which is a weak tangent to $F$. In particular, $\dim_H E = \dim_A E = s$. Moreover, one may assume that all of the similarities used to define $E$ are homotheties.

*Proof* For integers $m$, let $\mathcal{Q}(m)$ denote the collection of closed dyadic $2^{-m}$-cubes imposed on $\mathbb{R}^d$ oriented with the coordinate axes. Moreover, given $Q \in \mathcal{Q}(m)$ and $n > m$ write

$$M_n(Q) = \# \{Q' \in \mathcal{Q}(n) \ : \ Q' \subseteq Q \text{ and } Q' \cap F \neq \varnothing\} .$$

Directly from the definition of the Assouad dimension we can find sequences $m_k, n_k \in \mathbb{N}$ and $Q_k \in \mathcal{Q}(m_k)$ such that

$$M_{n_k}(Q_k) \geqslant 2^{(n_k - m_k)(s-1/k)} \tag{5.2}$$

and $n_k - m_k \to \infty$ as $k \to \infty$. Let $\mu_k$ be a Borel probability measure supported on $Q_k \cap F$ such that

$$\mu_k(Q') = \frac{1}{M_{n_k}(Q_k)}$$

for all $Q' \in \mathcal{Q}(n_k)$ such that $Q' \subseteq Q_k$ and $Q' \cap F \neq \varnothing$.

*Claim:* For all $t < s$ and $l \in \mathbb{N}$ there exists $k(l,t)$ and $Q(l,t) \in \mathcal{Q}(L)$ for some $L \geqslant m_{k(l,t)} + l$ with $Q(l,t) \subseteq Q_{k(l,t)}$ such that

$$\frac{\mu_{k(l,t)}(Q')}{\mu_{k(l,t)}(Q(l,t))} \leqslant \frac{|Q'|^t}{|Q(l,t)|^t}$$

for all $Q' \in \mathcal{Q}(l')$ with $L \leqslant l' \leqslant L + l$.

*Proof of claim:* Suppose the claim does not hold, in which case we can find $t < s$ and $l$ such that, for all $k$ and $Q_l \in \mathcal{Q}(m_k + l)$,

$$\frac{\mu_k(Q_l^1)}{\mu_k(Q_l)} > \frac{|Q_l^1|^t}{|Q_l|^t}$$

for some $Q_l^1 \in \mathcal{Q}(l_1)$ with $m_k + l < l_1 \leqslant m_k + 2l$. Fix such a 'bad cube' $Q_l$, which, by the pigeonhole principle, we may assume satisfies

$$\mu_k(Q_l) \geqslant 2^{-ld}, \tag{5.3}$$

and $Q_l^1 \subseteq Q_l$. Consider $Q_l^1$ and apply the failure of the claim again to obtain

$$\frac{\mu_k(Q_l^2)}{\mu_k(Q_l^1)} > \frac{|Q_l^2|^t}{|Q_l^1|^t}$$

for some $Q_l^2 \in \mathcal{Q}(l_2)$ with $l_1 < l_2 \leqslant l_1 + l$. Iterating this procedure to generate a sequence of 'bad cubes' $Q_l \supset Q_l^1 \supset Q_l^2 \supset \cdots \supset Q_l^q$ where $q$ is chosen to be the largest integer satisfying $|Q_l^q| > 2^{-n_k}\sqrt{d}$. In particular,

$$n_k - l \leqslant l_q \leqslant n_k. \tag{5.4}$$

Telescoping, and applying (5.3), we obtain

$$\mu_k(Q_l^q) > \frac{|Q_l^q|^t}{|Q_l|^t}\mu_k(Q_l) > \frac{2^{-n_k t} d^{t/2}}{2^{-(m_k+l)t} d^{t/2}} 2^{-ld}$$

$$= 2^{-(n_k - m_k)t} 2^{-l(d-t)}. \tag{5.5}$$

Let $t < t' < s$ and observe that, for large enough $k$,

$$M_{n_k}(Q_k) \geqslant 2^{(n_k - m_k)t'}$$

by (5.2). Since $Q_l^q$ contains at most $2^{dl}$ cubes of level $n_k$ (by (5.4)) we have

$$\mu_k(Q_l^q) \leqslant \frac{2^{dl}}{M_{n_k}(Q_k)} \leqslant 2^{-(n_k - m_k)t'} 2^{dl}. \tag{5.6}$$

Therefore, combining (5.5) and (5.6), we obtain

$$2^{-(n_k-m_k)t}2^{-l(d-t)} \leqslant 2^{-(n_k-m_k)t'}2^{dl}$$

which is clearly false for sufficiently large $k$ since $n_k - m_k \to \infty$, proving the claim.

We are now ready to use the claim to build a suitable weak tangent. Let $t_l$ be a sequence such that $t_l < s$ and $t_l \to s$. For $l \in \mathbb{N}$, let $T_l$ be the unique homothetic similarity which sends the cube $Q(l, t_l)$ from the claim to $[0,1]^d$. By taking a subsequence if necessary, we may assume there is a compact set $E \subseteq [0,1]^d$ such that

$$T_l(F) \cap [0,1]^d \to E$$

in $d_{\mathcal{H}}$ and a Borel probability measure $\mu$ supported on $E$ such that

$$\mu^l := \frac{1}{\mu_{k(l,t_l)}(Q(l,t_l))}\mu_{k(l,t_l)} \circ T_l^{-1} \to \mu$$

in the weak-$*$ topology. Observe that $E$ is a weak tangent to $F$.

Let $r \in (0,1)$ and $Q \in \mathcal{Q}(n)$ where $n$ is the largest integer such that $2^{-n} \geqslant r$. Then, by the claim, for sufficiently large $l$,

$$\mu^l(Q) = \frac{\mu_{k(l,t_l)}(T_l^{-1}(Q))}{\mu_{k(l,t_l)}(Q(l,t_l))} \leqslant \frac{|T_l^{-1}(Q)|^{t_l}}{|Q(l,t_l)|^{t_l}} \leqslant 2^{-nt_l} \leqslant 2^s r^{t_l}.$$

Moreover, for $x \in E$, the ball $B(x,r)$ can intersect at most $4^d$ cubes $Q \in \mathcal{Q}(n)$. Therefore

$$\mu(B(x,r)) \leqslant \liminf_{l\to\infty} \mu^l(B(x,r)) \leqslant \liminf_{l\to\infty} 4^d 2^s r^{t_l} = 4^d 2^s r^s.$$

It follows from the mass distribution principle, Lemma 3.4.2, that $\mathcal{H}^s(E) > 0$, and therefore $\dim_{\mathrm{H}} E \geqslant s$. Finally, the fact that $\dim_{\mathrm{H}} E = \dim_{\mathrm{A}} E = s = \dim_{\mathrm{A}} F$ follows from Theorem 5.1.2. □

We call $E$ in Definition 5.1.1 a *strong tangent* to $F$ if $X = B(0,1)$ (the closed unit ball), there is a point $x_0 \in F$ such that the sequence of similarity maps $T_k$ satisfies $T_k(x_0) = 0$ for all $k$ (we zoom-in to the specific point $x_0 \in F$), and $c_k \nearrow \infty$ as $k \to \infty$ (we genuinely zoom-in). Strong tangents are not the right notion to study the Assouad dimension, as the following theorem shows. This was first demonstrated by Le Donne and Rajala [182, Example 2.20], see also [89].

**Theorem 5.1.4** There exists a compact set $F \subseteq [0,1]$ with $\dim_A F = 1$, such that all strong tangents to $F$ have Assouad dimension equal to 0.

*Proof* Let $F \subseteq [0,1]$ be given by

$$F = \{0\} \ \cup \ \bigcup_{k=1}^{\infty} \bigcup_{l=0}^{k} \{2^{-k} + l4^{-k}\}.$$

For each $k$ let $T_k$ be the similarity map defined by $T_k(x) = k^{-1}4^k(x - 2^{-k})$ and observe that

$$T_k(F) \cap [0,1] = \bigcup_{l=0}^{k} \{l/k\} \to [0,1]$$

in the Hausdorff metric $d_{\mathcal{H}}$ as $k \to \infty$. Therefore $[0,1]$ is a weak tangent to $F$ and $\dim_A F = 1$ by Theorem 5.1.2. Since $F$ has only one accumulation point (at 0), the only non-trivial strong tangents are obtained by zooming in at this point. As such, suppose $E$ is a strong tangent to $F$ at 0. Therefore, there exists a sequence of similarity maps $S_n$ $(n \geqslant 1)$ of the form $S_n = c_n x$ where $1 < c_n \nearrow \infty$ as $n \to \infty$ such that

$$S_n(F) \cap B(0,1) \to E$$

in $d_{\mathcal{H}}$ as $n \to \infty$. We have assumed without loss of generality that $c_n > 0$, which we may do by passing to a subsequence and introducing a reflection if necessary. Let

$$F_0 = \{0\} \ \cup \ \left\{2^{-k} : k \in \mathbb{N}\right\}.$$

It is a straightforward exercise to show that $\dim_A F_0 = 0$. For comparison, see Theorem 3.4.7 and note that $2^{-k}$ is 'faster' than $1/k^p$ for any $p$. We show that $E$ is also a strong tangent to $F_0$, which will complete the proof since any strong or weak tangent to $F_0$ necessarily has Assouad dimension 0 by Theorem 5.1.2. For each $n$, let

$$m(n) = \min\{\, k \geqslant 1 \ : \ c_n 2^{-k} \leqslant 1\}$$

and note that $m(n) \to \infty$ as $n \to \infty$. It follows that

$$S_n(F) \cap B(0,1) \ \subseteq \ \{0\} \cup \bigcup_{k=m(n)}^{\infty} \left[c_n 2^{-k},\, c_n 2^{-k} + c_n k 4^{-k}\right]$$

and

$$S_n(F_0) \cap B(0,1) \;=\; \{0\} \cup \bigcup_{k=m(n)}^{\infty} \{c_n 2^{-k}\}.$$

Therefore

$$\begin{aligned} d_{\mathcal{H}}\Big(S_n(F) \cap B(0,1),\, S_n(F_0) \cap B(0,1)\Big) &\leqslant \sup_{k \geqslant m(n)} c_n k 4^{-k} \\ &= c_n m(n) 4^{-m(n)} \\ &\leqslant m(n) 2^{-m(n)} \to 0 \end{aligned}$$

which yields $S_n(F_0) \cap B(0,1) \to E$ in $d_{\mathcal{H}}$, completing the proof. □

If $F \subseteq \mathbb{R}^d$ and $\dim_{\mathrm{A}} F = d$, that is, the Assouad dimension of $F$ is the ambient spatial dimension, then even more information can be gleaned about the weak tangents.

**Theorem 5.1.5** If $F \subseteq \mathbb{R}^d$ and $\dim_{\mathrm{A}} F = d$, then $[0,1]^d$ is a weak tangent to $F$.

*Proof* For a given $R > 0$, let $\mathcal{Q}(R)$ be the natural tiling of $\mathbb{R}^d$ consisting of all translations of $[0,R]^d$ by $t \in R\mathbb{Z}^d$. For $r > 0$ and a given $Q \in \mathcal{Q}(R)$, let

$$M_r(Q) \;=\; \#\{Q' \in \mathcal{Q}(r) : Q' \subseteq Q \text{ and } Q' \cap F \neq \varnothing\}.$$

Since $Q \in \mathcal{Q}(R)$ are both contained in a ball of radius comparable to $R$ and contain a ball of radius comparable to $R$, we may work with cubes instead of balls in the definition of Assouad dimension. It is also sufficient just to consider scales $R$ and $r$ which are dyadic, that is, numbers of the form $2^k$, for $k \in \mathbb{Z}$.

Let $F \subseteq \mathbb{R}^d$ be such that $\dim_{\mathrm{A}} F = d$. In order to reach a contradiction, assume that $[0,1]^d$ is *not* a weak tangent to $F$. This means that there exists a dyadic $\varepsilon \in (0,1)$ such that for all dyadic $R > 0$ and all $Q \in \mathcal{Q}(R)$, there exists a $Q' \in \mathcal{Q}(\varepsilon R)$ such that $Q' \subseteq Q$ and $F \cap Q' = \varnothing$. If $Q \in \mathcal{Q}(R)$ and $r > \varepsilon R$, then

$$M_r(Q) \;\leqslant\; \varepsilon^{-d}.$$

and so we may assume from now on that $r < \varepsilon R$. Therefore, for all

$Q \in \mathcal{Q}(R)$ and all dyadic $r < \varepsilon R$,

$$\begin{aligned} M_r(Q) &\leqslant \#\{Q'' \in \mathcal{Q}(r) : Q'' \subseteq Q\} \ - \ \#\{Q'' \in \mathcal{Q}(r) : Q'' \subseteq Q'\} \\ &= \left(\frac{R}{r}\right)^d - \left(\frac{\varepsilon R}{r}\right)^d \\ &= \left(\frac{R}{r}\right)^d (1-\varepsilon^d). \end{aligned} \tag{5.7}$$

This does not yet contradict the assumption that $\dim_A F = d$, and so we proceed iteratively by cutting out more cubes of increasingly smaller size. Provided $r < \varepsilon^2 R$, we may use (5.7) to 'cut out' one cube from $\mathcal{Q}(\varepsilon^2 R)$ within each of the cubes from $\mathcal{Q}(\varepsilon R)$ which do intersect $F$. Therefore

$$M_r(Q) \leqslant \left(\frac{R}{r}\right)^d \left(1-\varepsilon^d\right)^2.$$

We may continue this process of 'cutting and reducing' as long as $r < \varepsilon^k R$. Therefore, if we choose $m \in \mathbb{N}$ such that $\varepsilon^{m+1} R \leqslant r < \varepsilon^m R$, then we obtain

$$\begin{aligned} M_r(Q) &\leqslant \left(\frac{R}{r}\right)^d \left(1-\varepsilon^d\right)^m \\ &\leqslant \left(\frac{R}{r}\right)^d \left(1-\varepsilon^d\right)^{\log(r/R)/\log\varepsilon - 1} \\ &= \frac{1}{1-\varepsilon^d}\left(\frac{R}{r}\right)^{d-\frac{\log(1-\varepsilon^d)}{\log\varepsilon}}. \end{aligned}$$

This proves

$$\dim_A F \leqslant d - \frac{\log(1-\varepsilon^d)}{\log\varepsilon} < d,$$

a contradiction. □

Theorem 5.1.5 is related to the notion of *porosity*. We say a set $F \subseteq \mathbb{R}^d$ is *α-porous* if there exists $\alpha \in (0, 1/2)$ such that for all $x \in F$ and all $r > 0$ there exists $y \in B(x,r)$ such that $B(y, \alpha r) \cap F = \varnothing$. This should be interpreted as saying that there are uniformly large holes in $F$ at all locations and scales. We say $F$ is *porous* if it is *α-porous* for some $\alpha$. This notion goes back to Väisälä [273]. It is straightforward to show that a set $F \subseteq \mathbb{R}^d$ is porous if and only if $[0,1]^d$ is not a weak tangent to the closure of $F$. Therefore Theorem 5.1.5 shows that $F \subseteq \mathbb{R}^d$ is porous if and only if $\dim_A F < d$. This fact was proved by Luukkainen [192, Theorem 5.2], who in fact established the following result, following Salli's celebrated results bounding the Hausdorff dimension of porous sets [246].

**Theorem 5.1.6** There is a constant $\varepsilon(\alpha, d) > 0$ depending only on $\alpha$ and $d$ such that if $F \subseteq \mathbb{R}^d$ is $\alpha$-porous, then

$$\dim_A F \leqslant d - \varepsilon(\alpha, d).$$

## 5.2 Weak tangents for the lower dimension?

There are dual analogues of Theorems 5.1.2, 5.1.3 and 5.1.5 which hold for the *lower* dimension, although one has to be careful. Consider the set $F = [0,1]$, which has $\dim_L F = 1$. However, taking $X = [0,1]$ and $T_k$ defined by $T_k(x) = (x-1)/k$ we find $T_k(F) \cap X = \{0\}$ for all $k$. Therefore $\{0\}$ is a weak tangent to $F$ according to Definition 5.1.1 and we cannot hope for $\dim_L F \leqslant \dim_L E$ (or even $\dim_L F \leqslant \dim_A E$) for all weak tangents $E$, which would be the naïve analogue of Theorem 5.1.2.

It should be clear that the problem here is that we are restricting the set $T_k(F)$ in the 'wrong way'. This can be rectified by insisting that the weak tangents intersect the interior of $X$, see [99, Theorem 1.1] for the details and proof.

**Theorem 5.2.1** Let $F \subseteq \mathbb{R}^d$ be closed, $E \subseteq \mathbb{R}^d$ be compact and suppose that $E$ is a weak tangent to $F$ where $X = [0,1]^d$ in the definition of weak tangent and $E \cap (0,1)^d \neq \varnothing$. Then $\dim_L F \leqslant \dim_H E$.

Note that one cannot replace $\dim_H E$ with $\dim_L E$ in the Theorem 5.2.1 since weak tangents which intersect the interior can still have unwanted isolated points on the boundary. There is a way of modifying the result to force $\dim_L E$ to appear, see [88, Proposition 7.7]. However, this is unlikely to yield better applications than Theorem 5.2.1 since the useful weak tangents tend to have equal lower and Hausdorff dimension.

With this basic inequality established, it is natural to look for an analogue of Theorem 5.1.3. The following was established in [99, Theorem 1.1].

**Theorem 5.2.2** Let $F \subseteq \mathbb{R}^d$ be closed with $\dim_L F = s \in [0,d]$. Then there exists a compact set $E \subseteq \mathbb{R}^d$ such that $E$ is a weak tangent intersecting the interior of $X$ as in Theorem 5.2.1 and satisfies $\dim_H E = \overline{\dim}_B E = s$.

Despite the utility of Theorem 5.1.2, the dual result for lower dimension has not yet provided many applications, but it does paint an interesting picture of the tangent structure of a given set. For concreteness, given a closed set $F \subseteq \mathbb{R}^d$, let $\mathcal{W}(F)$ be the set consisting of all weak tangents $E$ to $F$ derived using the reference set $X = [0,1]^d$, satisfying $E \cap (0,1)^d \neq \varnothing$, and such that the similarity ratios of the maps $T_k$ monotonically increase without bound, see Definition 5.1.1. That is, $\mathcal{W}(F)$ is the set of all reasonable weak tangents to $F$. Note that $\mathcal{W}(F)$ is not quite a gallery as in Section 2.3 but is similar in spirit. Moreover, write

$$\Delta(F) = \{\dim_{\mathrm{H}} E : E \in \mathcal{W}(F)\}$$

for the set of Hausdorff dimensions obtained by reasonable weak tangents. Using this notation, we have already seen that

$$\dim_{\mathrm{A}} F = \max \Delta(F)$$

and

$$\dim_{\mathrm{L}} F = \min \Delta(F)$$

noting that the maximum and minimum are indeed attained, see Theorems 5.1.2, 5.1.3, 5.2.1 and 5.2.2. In particular, $\Delta(F)$ is a singleton if and only if $\dim_{\mathrm{L}} F = \dim_{\mathrm{A}} F$. Otherwise it contains at least two points, corresponding to the Assouad and lower dimensions. It was shown in [99, Theorem 1.3] that if $\Delta \subseteq [0,d]$ is an $\mathcal{F}_\sigma$ set which contains its maximum and minimum, then there exists a compact set $F \subseteq \mathbb{R}^d$ such that

$$\Delta(F) = \Delta.$$

Recall that an $\mathcal{F}_\sigma$ set is a set which can be expressed as a countable union of closed sets. Finally, the analogue of Theorem 5.1.5, proved in [99, Theorem 1.4], is that a closed set $F \subseteq \mathbb{R}^d$ has $\dim_{\mathrm{L}} F = 0$ if and only if there is a set $E \in \mathcal{W}(F)$ which is a singleton.

## 5.3 Weak tangents for spectra?

Given how readily weak tangents relate to the Assouad dimension, it is at first sight a little surprising that they do not directly connect with the Assouad spectrum, or quasi-Assouad dimension. For example, the polynomial sequence $F_p = \{1/n^p : n \geqslant 1\}$ has the unit interval $[0,1]$ as a weak tangent (even a 'strong tangent' at 0), but yet the Assouad spectrum is strictly smaller than 1 for all $\theta < p/(p+1)$. A more striking example occurs later, where we will see that the limit set of Mandelbrot percolation in $\mathbb{R}^d$ has $[0,1]^d$ as a weak tangent, and so has Assouad dimension $d$, but has quasi-Assouad dimension strictly less than $d$, see Section 9.4. With these examples in mind, it might seem hopeless to connect weak tangents to the Assouad spectrum and quasi-Assouad dimension. However, García and Hare [118] observed that if one knows the speed of convergence to a given weak tangent, then one can get stronger lower bounds on dimension; that is, one can bound the Assouad spectrum, not just the Assouad dimension.

**Theorem 5.3.1** Let $F \subseteq \mathbb{R}^d$ be closed, $E \subseteq \mathbb{R}^d$ be compact, and suppose $E$ is a weak tangent to $F$ such that there is a constant $A > 0$ such that

$$d_{\mathcal{H}}(E, T_k(F) \cap X) \leqslant A c_k^{(\theta-1)/\theta} \tag{5.8}$$

where $c_k > 1$ is the similarity ratio of $T_k$. Then $\dim_{\mathrm{A}}^{\theta} F \geqslant \underline{\dim}_{\mathrm{B}} E$.

*Proof* We may assume that $c_k \to \infty$ since, if not, then $E \subseteq S(F)$ for some similarity $S$ and the result is trivial. Moreover, we may assume $\underline{\dim}_{\mathrm{B}} E > 0$ since otherwise there is nothing to prove. Let $0 < s < \underline{\dim}_{\mathrm{B}} E$ and, without loss of generality, we further assume $|E| = 1$.

For $k \geqslant 1$, let $r_k = A c_k^{(\theta-1)/\theta}$ and choose $y \in E$. For $k$ sufficiently large, $r_k \leqslant 1$ and we may find a collection $\{x_i\}_i$ of $N \geqslant r_k^{-s}$ many points which all lie in $E$ and are $3r_k$ separated. The rate of convergence condition (5.8) implies that for each $i$ we can find $y_i \in B(x_i, r_k) \cap T_k(F)$. Notice that the points $y_i$ are necessarily $r_k$ separated. For each $i$, let $z_i = T_k^{-1}(y_i)$ which gives us a collection of $N$ points $z_i \in F$ which are

$$r_k / c_k = A c_k^{-1/\theta}$$

separated. Moreover, each $z_i$ lies in the ball $B(z_1, 3c_k^{-1})$. It follows that

$$N_{A c_k^{-1/\theta}}(B(z_1, 3c_k^{-1}) \cap F) \geqslant N \geqslant r_k^{-s} = A^s \left( \frac{c_k^{-1}}{c_k^{-1/\theta}} \right)^s$$

which is enough to show $\dim_{\mathrm{A}}^{\theta} F \geqslant s$ and letting $s \to \underline{\dim}_{\mathrm{B}} E$ proves the result. □

In the above proof we only needed the extra rate of convergence condition to hold for infinitely many $k$, not all $k$ as stated. However, one can always pass to a subsequence of the $T_k$ such that the convergence condition holds for all $k$, and so this yields no advantage. Several variants on Theorem 5.3.1 are considered in [118, Theorem 6]. For example, if one assumes the $c_k$ grow at most exponentially, then the result may be strengthened, and the similarity condition on $T_k$ may be weakened somewhat. García and Hare also consider analogous results for the lower spectrum, where they had to be careful which weak tangents to consider, as seen in the previous section. We can formulate the result in terms of the quasi-Assouad dimension, rather than the Assouad spectrum, in which case we only need the convergence condition with

$$d_{\mathcal{H}}(E, T_k(F) \cap X) \leqslant A c_k^{-\varepsilon}$$

for some $\varepsilon > 0$, see [118]. Finally, we note we can recover Theorem 3.4.7 using this approach, which also goes some way to proving sharpness of Theorem 5.3.1. Recall $F_p = \{1/n^p : n \geqslant 1\}$ for fixed $p > 0$. Let $T_k$ be defined by $T_k(x) = k^p x$ and observe that

$$d_{\mathcal{H}}([0,1], T_k(F_p) \cap [0,1]) = 1 - k^p(1/(k+1)^p) \leqslant k^{-1} = (k^p)^{-1/p}.$$

Using Theorem 5.3.1 this implies $\dim_{\mathrm{A}}^{\theta} F_p \geqslant 1$ for $\theta$ satisfying $(\theta - 1)/\theta = -1/p$ and therefore

$$\dim_{\mathrm{A}}^{\theta} F_p = 1$$

for all $\theta \geqslant p/(p+1)$ and this, combined with Corollary 3.3.4, yields

$$\dim_{\mathrm{A}}^{\theta} F_p = \frac{1}{(1+p)(1-\theta)}$$

for $0 < \theta < p/(p+1)$.

## 5.4 Weak tangents for measures?

The ‘weak tangent’ approach to the Assouad dimension has proved rather fruitful in tackling several different problems. Indeed, it will show up throughout this book. The weak tangent approach for measures is somewhat less developed, partly because it does not work without suitable modification. Even once the theory is in place, however, it may still be less applicable than the theory for sets since understanding tangent measures is rather harder than understanding tangent sets.

We say that a probability measure $\nu$ on $\mathbb{R}^d$ is a *weak tangent measure* to a measure $\mu$ on $\mathbb{R}^d$ if there exists a sequence of similarity maps $T_k$ on

$\mathbb{R}^d$ and a sequence of weights $p_k > 0$ such that

$$p_k \mu \circ T_k^{-1} \rightharpoonup \nu$$

where $\rightharpoonup$ denotes weak convergence. Note the key difference between this and the notion of weak tangents for sets is that we do not restrict to a fixed compact set $X$. This means that weak tangent measures typically have unbounded support. The following analogous result to Theorem 5.1.2 was proved in [98, Theorem 2.2].

**Theorem 5.4.1** Suppose $\nu$ is a weak tangent measure to a measure $\mu$ on $\mathbb{R}^d$. Then $\dim_{\mathrm{A}} \nu \leqslant \dim_{\mathrm{A}} \mu$.

Combining this with Lemma 4.1.1 we see that all weak tangent measures of doubling measures are doubling.

**Corollary 5.4.2** If $\nu$ is a weak tangent measure to a doubling measure $\mu$ on $\mathbb{R}^d$, then $\nu$ is doubling.

It turns out that, for the conclusions of Theorem 5.4.1 and Corollary 5.4.2 to hold, it is necessary not to restrict the support of the measure. We demonstrate this using the following example, adapted from [98]. We construct a doubling measure on the plane which upon restriction to the unit ball becomes non-doubling, even with infinite box dimension! This shows that if one restricts the measures $p_k \mu \circ T_k^{-1}$ in the definition of weak tangent measures to the unit ball, then the analogue of Theorem 5.4.1 fails.

We begin with the square shown on the left of Figure 5.1. Let $r_i$ be a sequence of radii which decay exponentially to 0 and $x_i$ be a sequence of centres moving clockwise on the unit circle $\mathbb{S}^1$ which converge polynomially to a point $z \in \mathbb{S}^1$ as shown. Consider the balls $B_i = B(x_i, r_i)$ and remove from the original square the points in the interior of $B_i \cap B(0,1)$ which are at distance strictly greater than $r_i^i$ from $\mathbb{S}^1$ — see the grey region on the right in Figure 5.1. Label the remaining set, that is, the square with the grey regions removed, by $F$.

Let $\mu$ be the restriction of 2-dimensional Lebesgue measure to $F$. It is straightforward to see that $\mu$ is doubling, and even $\dim_{\mathrm{L}} \mu = \dim_{\mathrm{A}} \mu = 2$. However, $\nu = \mu|_{B(0,1)}$ is not doubling and moreover, $\dim_{\mathrm{B}} \mu = \dim_{\mathrm{A}} \nu = \infty$. For large enough $i$ we have

$$\nu(B(x_i, 2r_i)) \geqslant (\pi(2r_i)^2 - \pi r_i^i)/3 \geqslant \pi r_i^2$$

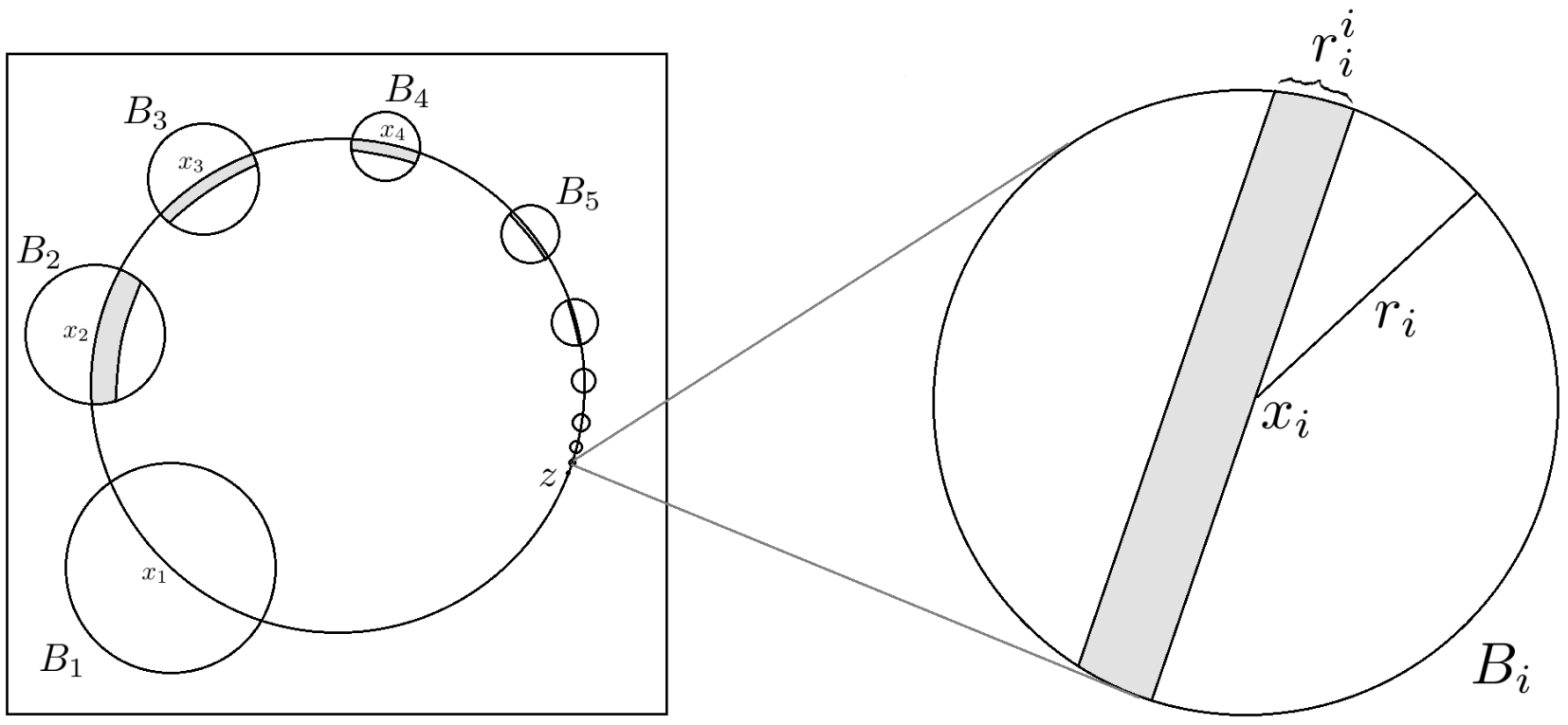

Figure 5.1 A doubling measure whose restriction to the unit ball is not doubling.

and $\nu(B(x_i, r_i)) \leqslant 3r_i r_i^i = 3r_i^{i+1}$. Therefore

$$\frac{\nu(B(x_i, 2r_i))}{\nu(B(x_i, r_i))} \geqslant \frac{\pi r_i^2}{3r_i^{i+1}} = \frac{\pi}{3r_i^{i-1}} \to \infty$$

which is enough to see that $\nu$ is not doubling directly from the definition. However, since for large enough $i$ we have $\nu(B(x_i, r_i)) \leqslant 3r_i^{i+1}$, we can conclude the stronger result that $\dim_{\mathrm{B}} \nu = \infty$.

The above example motivates the general question of when the restriction of a doubling measure is doubling and this turns out to be a subtle problem, see [216]. One also may ask questions about how much the box dimension or Assouad spectrum may increase upon restriction to a certain family of sets. The example above shows that in general they can jump from finite to infinite.

# PART TWO

# EXAMPLES

# 6
# Iterated function systems

One of the most important and widely used methods for constructing fractal sets is via iterated function systems. These were introduced in a celebrated article of Hutchinson [145] and we refer the reader to [70, Chapter 9] for a thorough overview. We consider the special case of self-similar sets in Chapter 7, self-affine sets in Chapter 8, and self-conformal sets in Section 9.1, but in this chapter we introduce the theory of iterated function systems and their attractors in full generality. We prove that the lower spectrum of such attractors is surprisingly well-behaved, see Theorem 6.3.1, and by considering two forms of quasi-self-similarity we will establish many dimension estimates by extending the implicit theorems of Falconer and McLaughlin to include the Assouad and lower dimensions, see Theorem 6.4.3 and Corollary 6.4.4.

## 6.1 IFS attractors and symbolic representation

An *iterated function system (IFS)* is a finite collection of contractions mapping a given compact domain $X \subseteq \mathbb{R}^d$ into itself. A *contraction* on $X$ is a map $S : X \mapsto X$ such that there is a constant $c \in (0, 1)$ such that

$$|S(x) - S(y)| \leqslant c|x - y| \tag{6.1}$$

for all $x, y \in X$. The constant $c$ is often referred to as the *contraction ratio.* The prominence of IFSs in the literature dates back to Hutchinson's seminal 1981 paper [145], which contained the following theorem.

**Theorem 6.1.1** For a given IFS $\{S_i\}_{i \in \mathcal{I}}$ there exists a unique non-empty compact set $F$, called the *attractor*, which satisfies

$$F = \bigcup_{i \in \mathcal{I}} S_i(F).$$

*Proof* This can be proved by an elegant application of *Banach's contraction mapping theorem*. Define the Hutchinson operator $\Phi : \mathcal{K}(X) \to \mathcal{K}(X)$ by

$$\Phi(K) = \bigcup_{i \in \mathcal{I}} S_i(K),$$

recalling that $\mathcal{K}(X)$ is the set of non-empty compact subsets of $X$ equipped with the Hausdorff metric, see (5.1). It follows from the fact that each of the maps $S_i$ is a contraction on $X$ that $\Phi$ is a contraction on $\mathcal{K}(X)$, with the constant appearing in (6.1) bounded above by the largest of the constants associated with the maps $S_i$. Since $\mathcal{K}(X)$ is complete, it follows from Banach's contraction mapping theorem that $\Phi$ has a unique fixed point $F \in \mathcal{K}(X)$ satisfying $\Phi(F) = F$, which proves the theorem. □

Attractors of IFSs are often complicated fractal sets but the invariance under the Hutchinson operator forces each such attractor to enjoy an approximate form of self-similarity, that is, the attractor is made up of scaled down copies of itself. Despite this apparently simple description, two major difficulties can arise when trying to understand the geometry of the attractor. The first is that hidden in the phrase 'scaled down' could be a significant amount of distortion, which is then amplified on smaller scales (upon iteration). The second is that the scaled down copies may overlap with each other and these overlaps may be subtle and hard to understand. With these two difficulties in mind, one often makes further assumptions on either the regularity of the maps $S_i$ (that is, controlling the distortion) or the separation of the sets $S_i(F)$ (that is, controlling the amount of overlap). Common conditions on the regularity of the maps include:

(i) Each map $S_i$ is a *similarity*, in which case the attractor $F$ is called *self-similar*. Here, for each $S_i$, distances are contracted by the same amount at all locations and in all directions.

(ii) Each map $S_i$ is *affine*, in which case the attractor $F$ is called *self-affine*. Here, for each $S_i$, distances are contracted according to a linear map and therefore different amounts of contraction can occur in different directions.

(iii) Each map $S_i$ is *conformal*, in which case the attractor $F$ is called *self-conformal.* Here, for each $S_i$, distances are contracted by different amounts at different locations globally but the contraction is uniform when considered locally. More precisely, at every point the Jacobian derivative of $S_i$ is a similarity, which is allowed to depend on the point.

These special families of IFSs and attractors will be discussed in more detail in Chapter 7, Chapter 8, and Section 9.1, respectively. Moving our attention to the problem of overlaps, the following separation conditions are fundamental in the theory of IFSs.

**Definition 6.1.2** An IFS, $\{S_i\}_{i\in\mathcal{I}}$, with attractor $F$ satisfies the *strong separation condition (SSC)*, if the sets $\{S_i(F)\}_{i\in\mathcal{I}}$ are pairwise disjoint. It satisfies the *open set condition (OSC)* if there exists a non-empty open set $U$ such that

$$\bigcup_{i\in\mathcal{I}} S_i(U) \subseteq U$$

with the sets $\{S_i(U)\}_{i\in\mathcal{I}}$ pairwise disjoint.

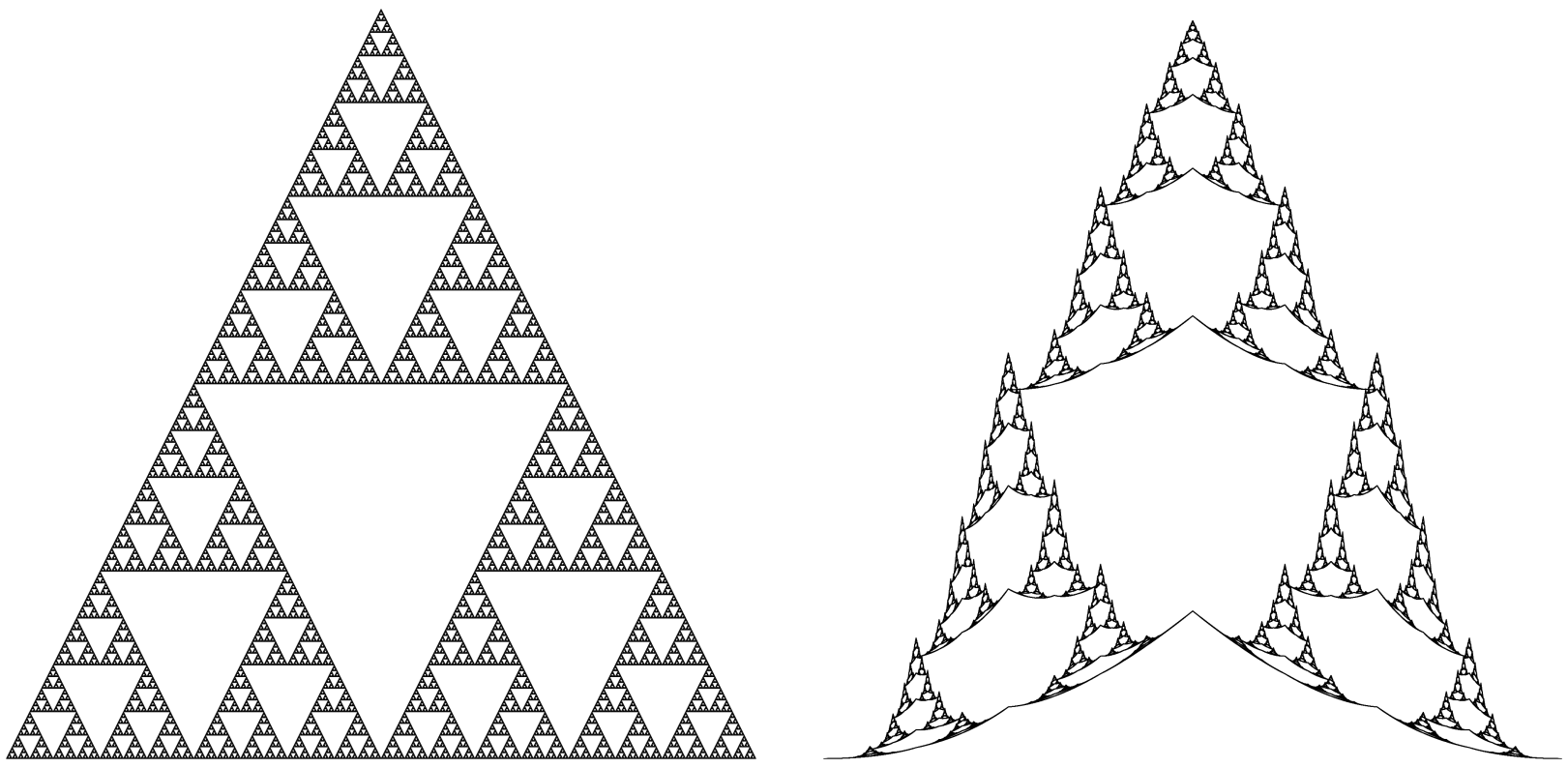

Figure 6.1 Two IFS attractors. Left: the Sierpiński Triangle. Right: an attractor of a nonlinear non-conformal IFS. Both attractors are made up of three scaled down copies of themselves. The mappings used on the left are strict similarities, making the Sierpiński Triangle self-similar, whereas the mappings used on the right are more complicated nonlinear contractions. The IFS on the left satisfies the OSC (with $U$ the interior of the convex hull of the attractor) but not the SSC.

One often uses a symbolic space built from the index set $\mathcal{I}$ to model the attractor of an IFS. It can be more convenient to work with the geometry of the symbolic space and there is a straightforward and canonical method to transfer information from the symbolic space to the fractal. We briefly describe this technique and fix some notation which will be used in later chapters. Let $\mathcal{I}^* = \bigcup_{k\geqslant 1}\mathcal{I}^k$ denote the set of all finite sequences with entries in $\mathcal{I}$ and for

$$\boldsymbol{i} = \big(i_1, i_2, \dots, i_k\big) \in \mathcal{I}^k$$

write

$$S_{\boldsymbol{i}} = S_{i_1} \circ S_{i_2} \circ \cdots \circ S_{i_k}.$$

Write $\mathcal{I}^{\mathbb{N}}$ to denote the set of all infinite $\mathcal{I}$-valued sequences and for $\boldsymbol{i} \in \mathcal{I}^{\mathbb{N}}$ or $\mathcal{I}^l$ with $l \geqslant k$ write $\boldsymbol{i}|_k \in \mathcal{I}^k$ to denote the restriction of $\boldsymbol{i}$ to its first $k$ entries. Let $\Pi : \mathcal{I}^{\mathbb{N}} \to F$ be the *symbolic coding map* from the symbolic space $\mathcal{I}^{\mathbb{N}}$ to $F$ defined by

$$\Pi(\boldsymbol{i}) = \bigcap_{k\in\mathbb{N}} S_{\boldsymbol{i}|_k}(X).$$

Then $\Pi(\mathcal{I}^{\mathbb{N}}) = F$, the attractor. For $\boldsymbol{i}, \boldsymbol{j} \in \mathcal{I}^*$, we write $\boldsymbol{i} \prec \boldsymbol{j}$ if $\boldsymbol{j}|_k = \boldsymbol{i}$ for some $k \leqslant |\boldsymbol{j}|$, where $|\boldsymbol{j}|$ is the length of the sequence $\boldsymbol{j}$. That is, $\boldsymbol{i} \prec \boldsymbol{j}$ if $\boldsymbol{i}$ is a *prefix* of $\boldsymbol{j}$. For

$$\boldsymbol{i} = (i_1, i_2, \dots, i_{k-1}, i_k) \in \mathcal{I}^*$$

let

$$\boldsymbol{i}^- = (i_1, i_2, \dots, i_{k-1}) \in \mathcal{I}^* \cup \{\omega\},$$

where $\omega$ is the empty word. For notational convenience the map $S_\omega$ is taken to be the identity. The above establishes a symbolic representation of $F$, which is also useful for visualisation purposes. For example, for all $k \geqslant 1$,

$$F = \bigcup_{\boldsymbol{i}\in\mathcal{I}^k} S_{\boldsymbol{i}}(F)$$

and the 'pieces' $S_{\boldsymbol{i}}(F)$ thus provide a natural cover of $F$ at a scale which tends to zero as $k \to \infty$. If $F$ is self-similar or self-conformal and satisfies the OSC, this approach leads to essentially optimal covers. However, for self-affine sets this is far from true. For $\boldsymbol{i} \in \mathcal{I}^*$, we write $[\boldsymbol{i}] = \{\boldsymbol{j} \in \mathcal{I}^{\mathbb{N}} : \boldsymbol{i} \prec \boldsymbol{j}\}$, so that $\Pi([\boldsymbol{i}]) = S_{\boldsymbol{i}}(F)$.

## 6.2 Invariant measures

One is often interested in measures supported on attractors of IFSs which enjoy analogous invariance properties, that is, they are not just arbitrary measures on the attractor, but measures which are invariant under the IFS. One natural approach for constructing such measures is to project Bernoulli measures from the symbolic space onto the attractor. Given an IFS $\{S_i\}_{i\in\mathcal{I}}$, with attractor $F$, associate a probability vector $\{p_i\}_{i\in\mathcal{I}}$ where each $p_i \in (0,1)$ and $\sum_{i\in\mathcal{I}} p_i = 1$. Then, analogously to Theorem 6.1.1, there is a unique Borel probability measure $\mu$ satisfying

$$\mu = \sum_{i\in\mathcal{I}} p_i \mu \circ S_i^{-1} \tag{6.2}$$

where $\mu \circ S_i^{-1}$ is the pushforward of $\mu$ under $S_i$, defined by $\mu \circ S_i^{-1}(E) = \mu(S_i^{-1}(E))$ for Borel sets $E$. This can be proved in a similar way to Theorem 6.1.1 by demonstrating that the associated operator is a contraction on the space of Borel probability measures on $X$ equipped with the Wasserstein distance, or another suitable metric.

We generally refer to measures $\mu$ given by (6.2) as *Bernoulli measures*, but when the attractor is self-similar, we say $\mu$ is self-similar, etc. A question of interest to us is, given a notion of dimension dim and an IFS attractor $F$, when does there exist a Bernoulli measure $\mu$ such that $\dim \mu = \dim F$, or at least a sequence of Bernoulli measures $\mu_k$ such that $\dim \mu_k \to \dim F$? The answer to this question is trivial if one simply requires $\mu$ to be supported on $F$, see (1.3), Theorems 4.1.3, and 4.2.3. Therefore, the potential interest in the question concerns the existence of *invariant* measures realising the dimension of the attractor. Surprisingly, the answer to this question is not always positive. For example, Das and Simmons [49] provided an example of a self-affine set in $\mathbb{R}^3$ such that the supremum of the Hausdorff dimensions of all invariant measures is strictly smaller than the Hausdorff dimension of the set. This answered an explicit question of Kenyon and Peres [168].

This interplay between dynamically invariant sets and measures is central to the dimension theory of dynamical systems. IFSs give rise to one form of dynamical invariance, although the natural dynamically invariant measures may not be Bernoulli. The natural measures might be projections of more general shift invariant measures on the symbolic space, such as Markov measures, or Gibbs measures, or they could simply be projections of shift ergodic measures. Moreover, other forms of invariance are possible, such as measures preserved under a dynamical system, or a group action. The question of interest becomes: when does a given invariant set support an invariant measure which *realises* its dimension? This question can then take on different flavours depending on the type of invariance and the dimensions involved. Concerning Hausdorff dimension, *realising the dimension* means that there is an invariant measure of maximal dimension, or at least that the supremum of the Hausdorff dimensions of invariant measures equals the Hausdorff dimension of the set, recall (1.3). For the Assouad/lower dimension, *realising the dimension* means that there is an invariant measure with minimal/maximal dimension, or at least that the infimum/supremum of the Assouad/lower dimension of invariant measures equals the Assouad/lower dimension of the set. It is particularly interesting to us whether or not an invariant measure can *simultaneously* realise the lower, Hausdorff, and Assouad dimensions, in a situation where they are distinct. We refer to this as the *simultaneous realisation problem*, which is deliberately vague at this point, and will revisit it several times later in the book as our family of examples grows.

## 6.3 Dimensions of IFS attractors

As we alluded to above, it is very difficult to say anything about the dimension theory of IFS attractors in general. One simple but important fact is that we may estimate the Hausdorff and box dimensions in terms of the Lipschitz constants of the defining maps. Specifically, if each $S_i$ is bi-Lipschitz, with constants $0 < b_i \leqslant c_i < 1$ such that for all $x, y \in X$

$$b_i|x-y| \leqslant |S_i(x) - S_i(y)| \leqslant c_i|x-y|,$$

and the IFS satisfies the SSC, then

$$t \leqslant \dim_{\mathrm{H}} F \leqslant \overline{\dim}_{\mathrm{B}} F \leqslant s \tag{6.3}$$

where $s$ and $t$ are defined by

$$\sum_{i \in \mathcal{I}} b_i^t = \sum_{i \in \mathcal{I}} c_i^s = 1.$$

See [70, Propositions 9.6-9.7]. The upper bound holds without any separation conditions and the lower bound, somewhat surprisingly, can fail even if the OSC is assumed. For example, consider the self-affine set defined by the IFS consisting of the maps $(x, y) \mapsto (x/2, x/2 + y/2)$ and $(x, y) \mapsto (x/2 + y/2, y/2)$. This IFS satisfies the OSC with the open set given by the open unit square $(0,1)^2$ but the attractor is the singleton $\{(0,0)\}$.

It is natural to ask if these bounds also apply to the Assouad and lower dimensions but they do not. We give a simple example of a self-affine set exhibiting this failure in Section 8.5.

It was proved in [278] that the lower dimension is strictly positive for any self-affine set (irrespective of overlaps), and it is natural to ask if this extends to more general IFSs. However, it was also shown in [278] that this is not true if the defining maps are only assumed to be invertible. Moreover, it can fail trivially if the maps are not invertible since, for example, one can force an attractor to have an isolated point by adding a constant map $x \mapsto c$ for some $c$ not already in the attractor.

If the maps are assumed to be bi-Lipschitz, then the lower spectrum is strictly positive, for at least some range of $\theta$. Moreover, $\dim_{\mathrm{L}}^{\theta} F \to \underline{\dim}_{\mathrm{B}} F$ as $\theta \to 0$. Recall that Lemma 3.4.4 shows that $\dim_{\mathrm{A}}^{\theta} F \to \overline{\dim}_{\mathrm{B}} F$ as $\theta \to 0$ for any set $F$, but it is not true in general that $\dim_{\mathrm{L}}^{\theta} F \to \underline{\dim}_{\mathrm{B}} F$.

**Theorem 6.3.1** Let $F$ be the attractor of an IFS consisting of bi-Lipschitz contractions and assume $F$ is not a singleton. Then there exists $\theta_0 \in (0,1)$ such that

$$\dim_{\mathrm{L}}^{\theta} F > 0$$

for all $\theta \in (0, \theta_0)$ and $\dim_{\mathrm{L}}^{\theta} F \to \underline{\dim}_{\mathrm{B}} F$ as $\theta \to 0$.

*Proof* First note that, even without the SSC, (6.3) guarantees that the lower box dimension of $F$ is strictly positive. Since we assume $F$ is not a singleton, we can find distinct maps $S_1$ and $S_2$ in the IFS with distinct fixed points. Therefore, for $n$ sufficiently large the IFS consisting of the two maps $\{S_1^n, S_2^n\}$ satisfies the SSC and therefore the attractor, which is a subset of $F$, has strictly positive lower box dimension by (6.3). Since the lower box dimension is positive, for any $0 < s < \underline{\dim}_{\mathrm{B}} F$, there exists $C > 0$ such that $M_r(F) \geqslant C r^{-s}$ for all $r \in (0,1)$, where $M_r(F)$ is the maximal size of an $r$-separated subset of $F$.

Write $c_{\max} = \max_{i \in \mathcal{I}} c_i < 1$ and $b_{\min} = \min_{i \in \mathcal{I}} b_i > 0$. Let $r$ be such that $0 < r^{\theta} < |F|$ and $x \in F$. Choose $\boldsymbol{i} \in \mathcal{I}^*$ such that $x \in S_{\boldsymbol{i}}(F) \subseteq B(x, r^{\theta})$ but $S_{\boldsymbol{i}^-}(F) \nsubseteq B(x, r^{\theta})$. This forces

$$c_{\max}^{|\boldsymbol{i}|-1} |F| \geqslant |S_{\boldsymbol{i}^-}(F)| \geqslant r^{\theta}$$

and so $|\boldsymbol{i}| \leqslant \log(r^{\theta}/|F|)/\log c_{\max} + 1$. Let $Y$ be a maximal $2r/b_{\min}^{|\boldsymbol{i}|}$-separated subset of $F$ and consider the set $S_{\boldsymbol{i}}(Y)$ which is necessarily a

$2r$-separated subset of $B(x, r^\theta) \cap F$. It follows that

$$N_r(B(x,r^\theta)\cap F) \;\geqslant\; \#Y \geqslant C(2r/b_{\min}^{|i|})^{-s}$$

$$\geqslant C b_{\min}^s 2^{-s}|F|^{-s\frac{\log b_{\min}}{\log c_{\max}}}\, r^{\theta s \frac{\log b_{\min}}{\log c_{\max}}-s}$$

which, since $s$ can be chosen arbitrarily close to $\underline{\dim}_{\mathrm{B}} F$, shows that

$$\dim_{\mathrm{L}}^\theta F \geqslant \underline{\dim}_{\mathrm{B}} F \frac{1-\theta\frac{\log b_{\min}}{\log c_{\max}}}{1-\theta}.$$

Therefore, $\dim_{\mathrm{L}}^\theta F > 0$ for $0 < \theta < \frac{\log c_{\max}}{\log b_{\min}}$ and $\dim_{\mathrm{L}}^\theta F \to \underline{\dim}_{\mathrm{B}} F$ as $\theta \to 0$. □

Theorem 6.3.1 combined with the fact that $\dim_{\mathrm{qL}} F \leqslant \dim_{\mathrm{H}} F$ for any closed set $F$, see [118, Proposition 10], shows that the lower spectrum of an attractor of a bi-Lipschitz IFS is non-constant whenever the Hausdorff and lower box dimension are distinct.

We can adapt the example from [278] to show that the previous theorem is not true in general if the bi-Lipschitz condition is relaxed to invertible. Consider the IFS consisting of two maps acting on $[0, 1]$ given by

$$S_1(x) = x^{1/x}/3$$

for $x > 0$ and $S_1(0) = 0$ and

$$S_2(x) = (x + 2)/3$$

and let $F$ be the attractor, see Figure 6.3. This IFS satisfies the SSC and the defining maps are invertible but $S_1$ is not bi-Lipschitz.

Fix $\theta \in (0, 1)$, a large integer $n \geqslant 1$, and consider $x = 0$ and $R = S_1^n(1/2)$. Then

$$B(0, R) \cap F = S_1^{n+1}(F)$$

but for $n$ sufficiently large $|S_1^{n+1}(F)| \leqslant R^{1/\theta}$ and so

$$N_{R^{1/\theta}}(B(0, R) \cap F) = 1$$

which proves that

$$\dim_{\mathrm{L}}^\theta F = 0 \tag{6.4}$$

for all $\theta \in (0, 1)$.

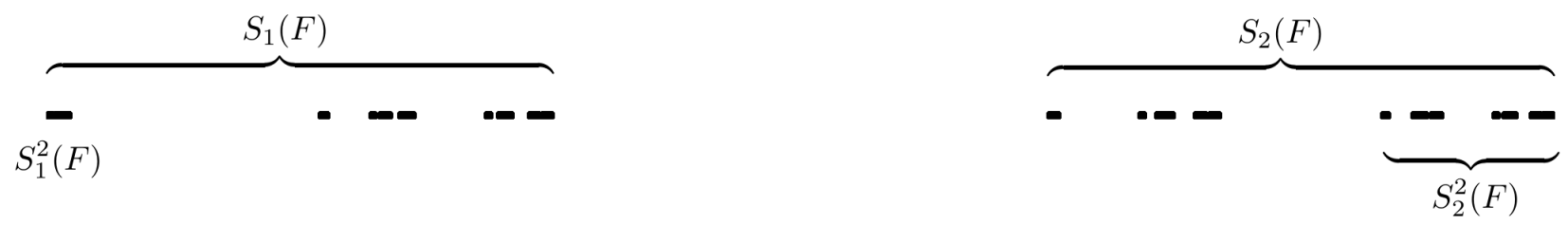

Figure 6.2 The IFS attractor defined above using invertible maps and satisfying the SSC, but with $\dim_{\mathrm{L}}^{\theta} F = 0$ for all $\theta \in (0,1)$. Notice how quickly the images $S_1^n(F)$ get smaller as $n$ increases.

## 6.4 Ahlfors regularity and quasi-self-similarity

In this section we take a short detour to investigate the Assouad and lower dimensions of sets displaying some degree of self-similarity, such as attractors of IFSs. This will allow us to determine certain relationships which hold in general between the dimensions and spectra we consider. This is a worthwhile pursuit before we consider specific examples because it will focus our attention on the interesting questions. Recall that a compact set $F \subseteq \mathbb{R}^d$ is called *Ahlfors regular* if there exists a constant $c > 0$ such that

$$\tfrac{1}{c}\, r^{\dim_{\mathrm{H}} F} \ \leqslant\ \mathcal{H}^{\dim_{\mathrm{H}} F}\big(B(x,r) \cap F\big) \ \leqslant\ c\, r^{\dim_{\mathrm{H}} F}$$

for all $x \in F$ and all $0 < r < |F|$, see [130, Chapter 8]. In a certain sense Ahlfors regular spaces are the most homogeneous spaces. This is reflected in the following proposition. Recall that the lower dimension and Assouad dimension are extremal and so establishing that these two dimensions are equal also establishes equality with the other notions of dimension and spectra.

**Theorem 6.4.1** If $F$ is compact and Ahlfors regular, then

$$\dim_{\mathrm{L}} F = \dim_{\mathrm{H}} F = \dim_{\mathrm{A}} F.$$

*Proof* Write $s = \dim_{\mathrm{H}} F$, $\mu = \mathcal{H}^s|_F$ for the $s$-dimensional Hausdorff measure restricted to $F$ and normalised such that $\mu(F) = 1$. Note that since $F$ is bounded, Ahlfors regularity guarantees $0 < \mathcal{H}^s(F) < \infty$ and so this normalisation is well-defined. Let $0 < r < R < |F|$ and $x \in F$ and note

$$c^{-2}\left(\frac{R}{r}\right)^s \leqslant \frac{\mu(B(x,R))}{\mu(B(x,r))} \leqslant c^2\left(\frac{R}{r}\right)^s.$$

It follows from Lemma 4.1.2 that

$$s \leqslant \dim_{\mathrm{L}} \mu \leqslant \dim_{\mathrm{L}} F \leqslant \dim_{\mathrm{A}} F \leqslant \dim_{\mathrm{A}} \mu \leqslant s$$

proving the result. □

The notion of *quasi-self-similarity* was introduced by Falconer [66] and McLaughlin [209] and serves to weaken the notion of Ahlfors regularity whilst maintaining some of the consequences. These results allow one to deduce facts about the dimensions and measures of a set without having to calculate them explicitly, which leads to them often being referred to as *implicit theorems.* This is done by assuming that either large parts of the set can be mapped to small parts without too much distortion or vice versa.

**Definition 6.4.2** (Quasi-self-similarity) Let $F \subseteq \mathbb{R}^d$ be a non-empty compact set.

(i) We say $F$ is *quasi-self-similar of type (1)* if there exists $a > 0$ and $r_0 > 0$ such that for every set $U$ that intersects $F$ with $|U| \leqslant r_0$, there is a function $g : F \cap U \to F$ satisfying

$$a\,|U|^{-1}\,|x - y| \;\leqslant\; |g(x) - g(y)| \tag{6.5}$$

for all $x, y \in F \cap U$.

(ii) We say $F$ is *quasi-self-similar of type (2)* if there exists $a > 0$ and $r_0 > 0$ such that for every closed ball $B$ with centre in $F$ and radius $R \leqslant r_0$, there is a function $g : F \to F \cap B$ satisfying

$$a\,R\,|x - y| \;\leqslant\; |g(x) - g(y)| \tag{6.6}$$

for all $x, y \in F$.

Temporarily write $s = \dim_{\mathrm{H}} F$. It was shown in [209, 66] that $F$ being quasi-self-similar of type (1) guarantees that $\mathcal{H}^s(F) \geqslant a^s > 0$ and $\underline{\dim}_{\mathrm{B}}F = \overline{\dim}_{\mathrm{B}}F = s$. On the other hand, it was shown in [66] that $F$ being quasi-self-similar of type (2) guarantees $\mathcal{H}^s(F) \leqslant 4^s\,a^{-s} < \infty$ and $\underline{\dim}_{\mathrm{B}}F = \overline{\dim}_{\mathrm{B}}F = s$, also see [69, Chapter 3]. These results were extended to include the Assouad and lower dimensions in [88]. By considering the relationships between the various dimensions and spectra discussed in this book, one may also glean implicit theorems for the other notions by sandwiching them in between the dimensions considered below.

**Theorem 6.4.3** Let $F \subseteq \mathbb{R}^d$ be compact and non-empty.

(i) If $F$ is quasi-self-similar of type (1), then

$$\dim_L F \leqslant \dim_H F = \dim_B F = \dim_A F.$$

(ii) If $F$ is quasi-self-similar of type (2), then

$$\dim_L F = \dim_H F = \dim_B F \leqslant \dim_A F.$$

(iii) If $F$ is quasi-self-similar of type (1) and (2), then

$$\dim_L F = \dim_H F = \dim_B F = \dim_A F$$

and, moreover, $F$ is Ahlfors regular.

*Proof* Throughout the proof write $s = \dim_H F = \underline{\dim}_B F = \overline{\dim}_B F$. It follows immediately from the definition of box dimension that for all $\varepsilon > 0$ there exists a constant $C \geqslant 1$ such that, for all $r > 0$,

$$C^{-1} r^{-s+\varepsilon} \leqslant N_r(F) \leqslant C r^{-s-\varepsilon}. \tag{6.7}$$

*Proof of (i).* Let $0 < r < R \leqslant r_0/2$ and $x \in F$. Since $F$ is quasi-self-similar of type (1), there exists $g_1 : B(x,R) \cap F \to F$ satisfying (6.5) with $U = B(x,R)$. If $\{U_i\}_i$ is an $(ar/2R)$-cover of $g_1(B(x,R) \cap F)$, then $\{g_1^{-1}(U_i)\}_i$ is an $r$-cover of $B(x,R) \cap F$. Then, applying (6.7),

$$\begin{aligned} N_r\big(B(x,R) \cap F\big) &\leqslant N_{ar/2R}\big(g_1(B(x,R) \cap F)\big) \leqslant N_{ar/2R}(F) \\ &\leqslant C(2/a)^{s+\varepsilon}\left(\frac{R}{r}\right)^{s+\varepsilon} \end{aligned}$$

which shows $\dim_A F \leqslant s+\varepsilon$ and letting $\varepsilon \to 0$ completes the proof. □

*Proof of (ii).* Let $0 < r < R \leqslant r_0/2$ and $x \in F$. Since $F$ is quasi-self-similar of type (2), there exists $g_2 : F \to B(x,R) \cap F$ satisfying (6.6) with $B = B(x,R)$. If $\{U_i\}_i$ is an $r$-cover of $g_2(F)$, then $\{g_2^{-1}(U_i)\}_i$ is an $(r/aR)$-cover of $F$. Then, applying (6.7),

$$N_r\big(B(x,R) \cap F\big) \geqslant N_r\big(g_2(F)\big) \geqslant N_{r/aR}(F) \geqslant C^{-1} a^{s-\varepsilon}\left(\frac{R}{r}\right)^{s-\varepsilon}$$

which shows $\dim_L F \geqslant s-\varepsilon$ and letting $\varepsilon \to 0$ completes the proof. □

*Proof of (iii).* It remains to prove that $F$ is Ahlfors regular. Let $r \in$

$(0, r_0/2)$, $x \in F$ and $g_1$ and $g_2$ be maps satisfying (6.5)-(6.6) with $U = B = B(x,r)$. The results from [66] mentioned above guarantee

$$a^s \leqslant \mathcal{H}^s(F) \leqslant 4^s a^{-s}.$$

Therefore, using the *scaling property* for Hausdorff measure, see [70, Proposition 2.2],

$$\mathcal{H}^s(B(x,r) \cap F) \leqslant \mathcal{H}^s\big(g_1^{-1}(F)\big) \leqslant (2r/a)^s \mathcal{H}^s(F) \leqslant 8^s a^{-2s} r^s.$$

and

$$\mathcal{H}^s(B(x,r) \cap F) \geqslant \mathcal{H}^s\big(g_2(F)\big) \geqslant (ar)^s \mathcal{H}^s(F) \geqslant a^{2s} r^s$$

completing the proof. □

The implicit theorems allow us to deduce information about the dimensions and spectra of various classes of IFS attractor as well as some variants, including graph-directed IFSs and sub-self-similar sets. See [69] for more information on these examples.

**Corollary 6.4.4**
(a) The following classes of sets are Ahlfors regular and, in particular, have equal Assouad and lower dimension:

(i) self-similar sets satisfying the OSC,
(ii) self-conformal sets satisfying the OSC,
(iii) graph-directed self-similar or self-conformal sets satisfying the graph-directed analogue of the OSC.

(b) The following classes of sets have equal Assouad dimension, quasi-Assouad and Hausdorff dimension, as well as having constant Assouad spectrum:

(i) sub-self-similar sets satisfying the OSC.

(c) The following classes of sets have equal lower dimension, quasi-lower dimension and Hausdorff dimension, as well as having constant lower spectrum, regardless of separation conditions:

(i) self-similar sets,
(ii) self-conformal sets,
(iii) graph-directed self-similar or self-conformal sets.

*Proof* This follows immediately from Theorem 6.4.3, see [66, 68, 88] for the details. □

Theorem 6.4.3 is sharp, in that the inequalities in parts (i) and (ii) cannot be replaced with equalities in general. To see this, note that the inequality in part (i) is sharp because the set $[0, 1] \cup \{2\}$ is quasi-self-similar of type (1), but has lower dimension strictly less than Hausdorff dimension. Also, the inequality in part (ii) is sharp because self-similar sets can have Assouad dimension strictly larger than Hausdorff dimension and are quasi-self-similar of type (2), see Theorem 7.2.1.

# 7
# Self-similar sets

Self-similar sets are a special case of IFS attractors and have been studied extensively in the literature. In this chapter we consider self-similar sets and measures in detail. Corollary 6.4.4 already implies that self-similar sets have equal lower, Hausdorff and box dimension, but we will see in Theorem 7.2.1 that the Assouad dimension may be strictly larger.

## 7.1 Self-similar sets and the Hutchinson-Moran formula

The key reason why self-similar sets are easier to deal with than self-affine sets and more general attractors is that the images of the set under compositions of maps from the IFS form natural covers for the original set. This leads to simple formulae for the dimensions of a self-similar set. Recall, a map $S : \mathbb{R}^d \to \mathbb{R}^d$ is called a *similarity* if there exists a constant $c > 0$ such that, for all $x, y \in \mathbb{R}^d$,

$$|S(x) - S(y)| = c|x - y|. \tag{7.1}$$

The constant $c$ is called the *similarity ratio* of $S$. If $c \in (0, 1)$ then $S$ is also a contraction, and we use the term *contracting similarity*. This should be compared with the definition of a contraction (6.1), noting that similarities yield a precise formula for the distortion. Moreover, similarities are bi-Lipschitz with equal upper and lower Lipschitz constants and so (6.3) immediately yields a precise formula for the dimensions of self-similar sets when the SSC is satisfied. However, this observation can be improved.

Given an IFS, $\{S_i\}_{i \in \mathcal{I}}$, consisting of contracting similarities with similarity ratios $\{c_i\}_{i \in \mathcal{I}}$, the *similarity dimension* is defined to be the unique solution to the *Hutchinson-Moran formula*

$$\sum_{i \in \mathcal{I}} c_i^s = 1. \tag{7.2}$$

This formula first appeared in [213] (see also [61, Chapter 13]) and later in [145]. The similarity dimension is always an upper bound for the upper box dimension of the attractor (but not the Assouad dimension, see Theorem 7.2.1). However, if the IFS satisfies the OSC, then all of the dimensions discussed in this book are equal to the similarity dimension, see [70, Section 9.3] and [88]. In fact, the attractor is Ahlfors regular by Corollary 6.4.4.

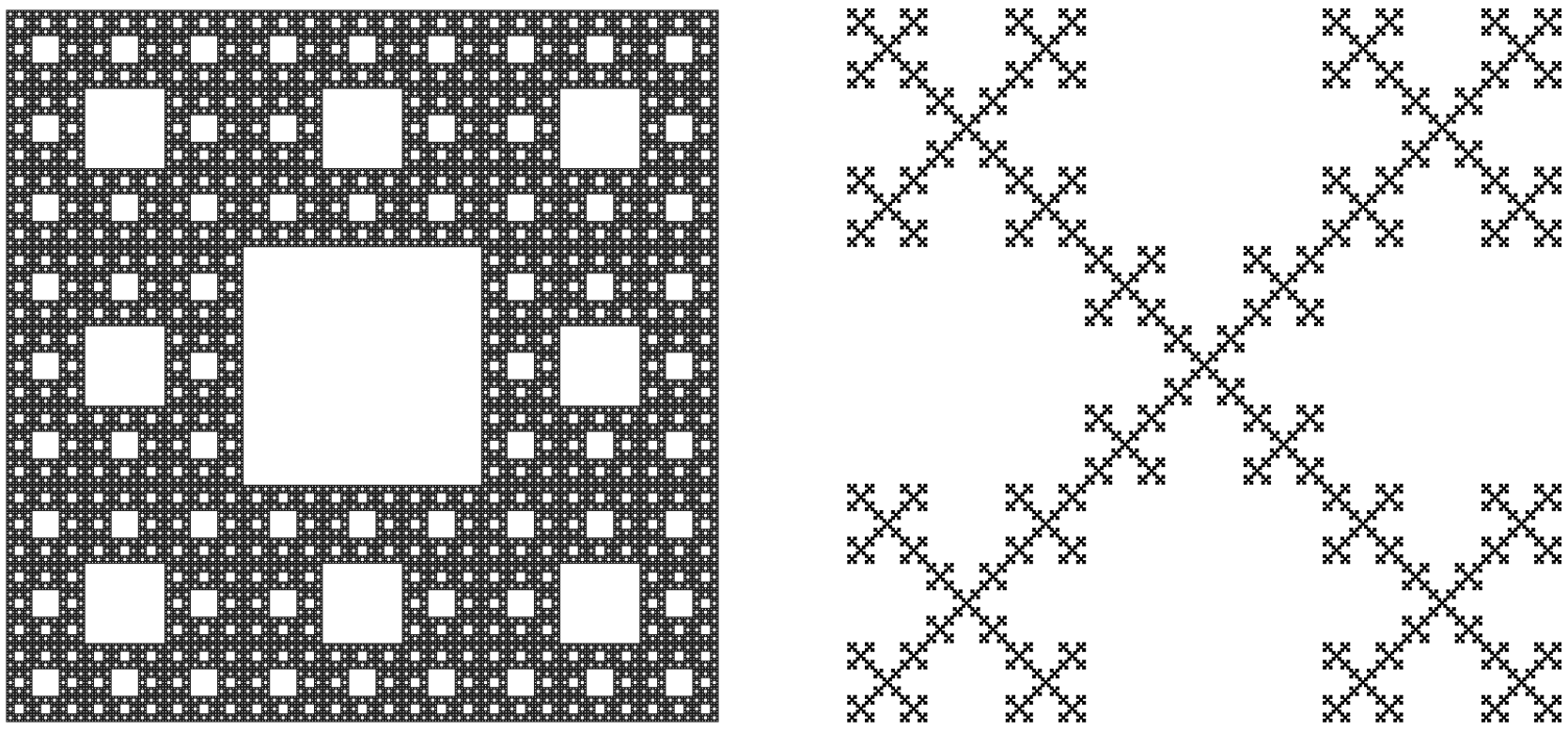

Figure 7.1 Two self-similar sets.

The case when the OSC is not satisfied is not fully understood. In $\mathbb{R}^d$, a 'dimension drop' (Hausdorff dimension strictly smaller than the minimum of the similarity dimension, $s$, and the ambient spatial dimension, $d$) can occur if different iterates of maps in the IFS overlap exactly, but it is a major open problem to decide if this is the only way the dimension can drop, see, for example, [233]. An important step towards solving this conjecture was made by Hochman [132]. For example, he verifies that the only cause for a dimension drop is exact overlaps, provided we are working in $\mathbb{R}$ and the defining parameters for the IFS are algebraic. Recall that a real number is called *algebraic* if it is the root of a non-zero polynomial with integer coefficients. Results in higher dimensions are also available, but there the situation is more complicated [133].

Even in the overlapping case, one can still say a lot about the dimensions of a self-similar set. It follows from Corollary 6.4.4 that for all self-similar sets $F$ one has

$$\dim_{\mathrm{L}} F = \dim_{\mathrm{qL}} F = \dim_{\mathrm{L}}^{\theta} F = \dim_{\mathrm{H}} F = \dim_{\mathrm{P}} F = \underline{\dim}_{\mathrm{B}} F = \overline{\dim}_{\mathrm{B}} F$$

for all $\theta \in (0,1)$. However, in the overlapping case the Assouad dimension can be strictly larger than the box dimension. Whether or not this was possible was asked by Olsen [219, Question 1.3] and answered in [88, Section 3.1]. We discuss the Assouad dimension in detail in the following section.

## 7.2 The Assouad dimension of self-similar sets

We start with an example from [88], demonstrating that the Assouad dimension can strictly exceed the box dimension. We also use this example to show that Assouad dimension can increase under Lipschitz maps.

**Theorem 7.2.1** There exists a self-similar set $F \subseteq \mathbb{R}$ such that $\dim_{\mathrm{B}} F < \dim_{\mathrm{A}} F$.

*Proof* Let $0 < \gamma < \alpha < \beta < 1$ be such that $\log \beta / \log \alpha \notin \mathbb{Q}$ and define similarity maps $S_1, S_2, S_3$ on $[0,1]$ by

$$S_1(x) = \alpha x, \qquad S_2(x) = \beta x \qquad \text{and} \qquad S_3(x) = \gamma x + (1-\gamma).$$

Let $F$ be the self-similar attractor of $\{S_1, S_2, S_3\}$, see Figure 7.2. We show that $[0,1]$ is a weak tangent to $F$. For each integer $k \geqslant 1$ let $T_k$ be defined by $T_k(x) = \beta^{-k} x$ and let

$$E_k = \{0\} \cup \left\{\alpha^m \beta^n : m \in \mathbb{N}, n \in \mathbb{Z} \text{ with } n \geqslant -k\right\} \cap [0,1],$$

noting that $E_k \subseteq T_k(F) \cap [0,1]$ for each $k$. Therefore, it suffices to show that $E_k \to [0,1]$ in the Hausdorff metric $d_{\mathcal{H}}$. Using the fact that the sets $E_k$ are nested,

$$\begin{aligned} E_k &\to \overline{\bigcup_{k \in \mathbb{N}} E_k} \cap [0,1] \\ &= \overline{\{\alpha^m \beta^n : m \in \mathbb{N}, n \in \mathbb{Z}\}} \cap [0,1] \\ &= [0,1], \end{aligned}$$

where the final equality is a consequence of Dirichlet's theorem[1] as follows. We argue that

$$\{m \log \alpha + n \log \beta : m \in \mathbb{N}, n \in \mathbb{Z}\}$$

is dense in $(-\infty, 0)$, which is sufficient. We have

$$m \log \alpha + n \log \beta = n \log \alpha \left( \frac{m}{n} + \frac{\log \beta}{\log \alpha} \right)$$

and Dirichlet's theorem gives that there exists infinitely many $n$ such that

$$\left| \frac{m}{n} + \frac{\log \beta}{\log \alpha} \right| < 1/n^2$$

for some $m$. Since $\log \beta / \log \alpha$ is irrational, we may choose $m, n$ such that

$$0 < |m \log \alpha + n \log \beta| < \frac{|\log \alpha|}{n}$$

with $n$ arbitrarily large, and thus $m \log \alpha + n \log \beta$ arbitrarily small (and negative). We can therefore find infinite arithmetic progressions

$$\{-\varepsilon k : k \in \mathbb{N}\} \subseteq \{m \log \alpha + n \log \beta : m \in \mathbb{N}, n \in \mathbb{Z}\}$$

for arbitrarily small

$$\varepsilon = |m \log \alpha + n \log \beta| > 0,$$

which completes the proof that $[0, 1]$ is a weak tangent to $F$. It follows that $\dim_{\mathrm{A}} F = 1$ by applying Theorem 5.1.2.

The box dimension of $F$ exists and is bounded above by the similarity dimension $s$, which is the unique solution of the Hutchinson-Moran formula (7.2)

$$\alpha^s + \beta^s + \gamma^s = 1.$$

If we choose $\alpha, \beta, \gamma$ such that $s < 1$ (whilst maintaining the irrationality assumption), then

$$\dim_{\mathrm{B}} F \leqslant s < 1 = \dim_{\mathrm{A}} F$$

as required. For example, it is sufficient to choose $\alpha = 1/3$, $\beta = 1/2$ and $\gamma = 1/8$, see Figure 7.2. □

We can use this example to prove that Assouad dimension can increase under Lipschitz maps. This counter-intuitive phenomenon was first observed in [192, Example A.6 2], via a different example.

[1] See Section 14.2 for more on Dirichlet's theorem.

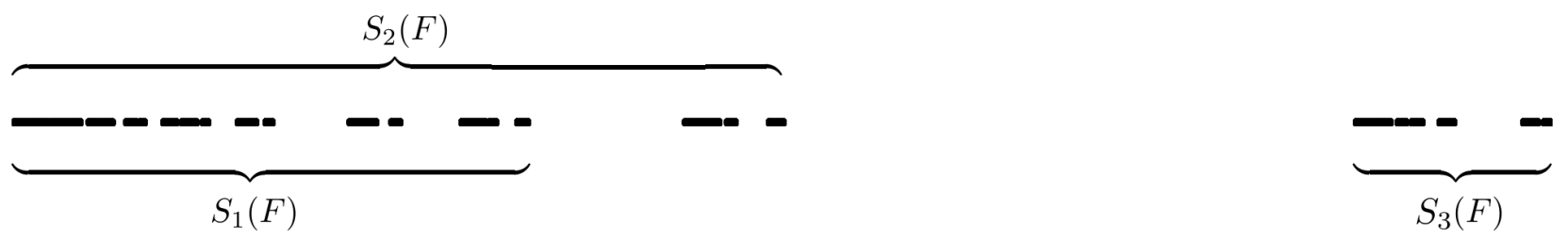

Figure 7.2 The overlapping self-similar set with $\alpha = 1/3$, $\beta = 1/2$ and $\gamma = 1/8$. Stretching our imagination slightly, we can see the tangent structure emerging around the origin. The box dimension of the self-similar set is bounded above by the similarity dimension $s \approx 0.9582 < 1$ but the Assouad dimension is 1.

**Corollary 7.2.2** The Assouad dimension can increase under Lipschitz maps. For example, there exists a compact set $E \subseteq \mathbb{R}^2$ such that $\dim_A \pi(E) = 1 > \dim_A E$ where $\pi$ denotes projection onto the first coordinate (which is a Lipschitz map).

*Proof* Let $\alpha, \beta, \gamma \in (0, 1)$ be chosen as in the proof of Theorem 7.2.1 and consider the similarity maps $T_1, T_2, T_3$ on $[0, 1]^2$ given by

$$T_1(x, y) = (\alpha x, \alpha y),$$

$$T_2(x, y) = (\beta x, \beta y) + (0, 1 - \beta)$$

and

$$T_3(x) = (\gamma x, \gamma y) + (1 - \gamma, 0)$$

and let $E$ be the attractor of $\{T_1, T_2, T_3\}$, see Figure 7.3. If $\alpha, \beta, \gamma$ are chosen such that $\alpha + \beta, \beta + \gamma, \alpha + \gamma \leqslant 1$ and such that the similarity dimension $s < 1$, then $\{T_1, T_2, T_3\}$ satisfies the OSC and therefore the Assouad dimension of $E$ is equal to $s$. However, the set $F$ from Theorem 7.2.1 is the projection of $E$ onto the first coordinate and therefore $\dim_A \pi(E) = 1 > \dim_A E$. □

The Assouad dimension of self-similar sets in $\mathbb{R}$ is well-understood. Fraser, Henderson, Olson and Robinson [96] established a precise dichotomy: either $\dim_A F = \dim_H F$ or $\dim_A F = 1$. Moreover, the dichotomy is precisely determined by the *weak separation property* of Lau and Ngai [181], and Zerner [282]. We write $I$ for the identity map on $\mathbb{R}^d$ and define

$$\mathcal{E} = \{S_{\boldsymbol{i}}^{-1} \circ S_{\boldsymbol{j}} \ : \ \boldsymbol{i}, \boldsymbol{j} \in \mathcal{I}^*, \ \boldsymbol{i} \neq \boldsymbol{j}\}$$

which is a subset of the space of all similarities mapping $\mathbb{R}^d$ to itself.

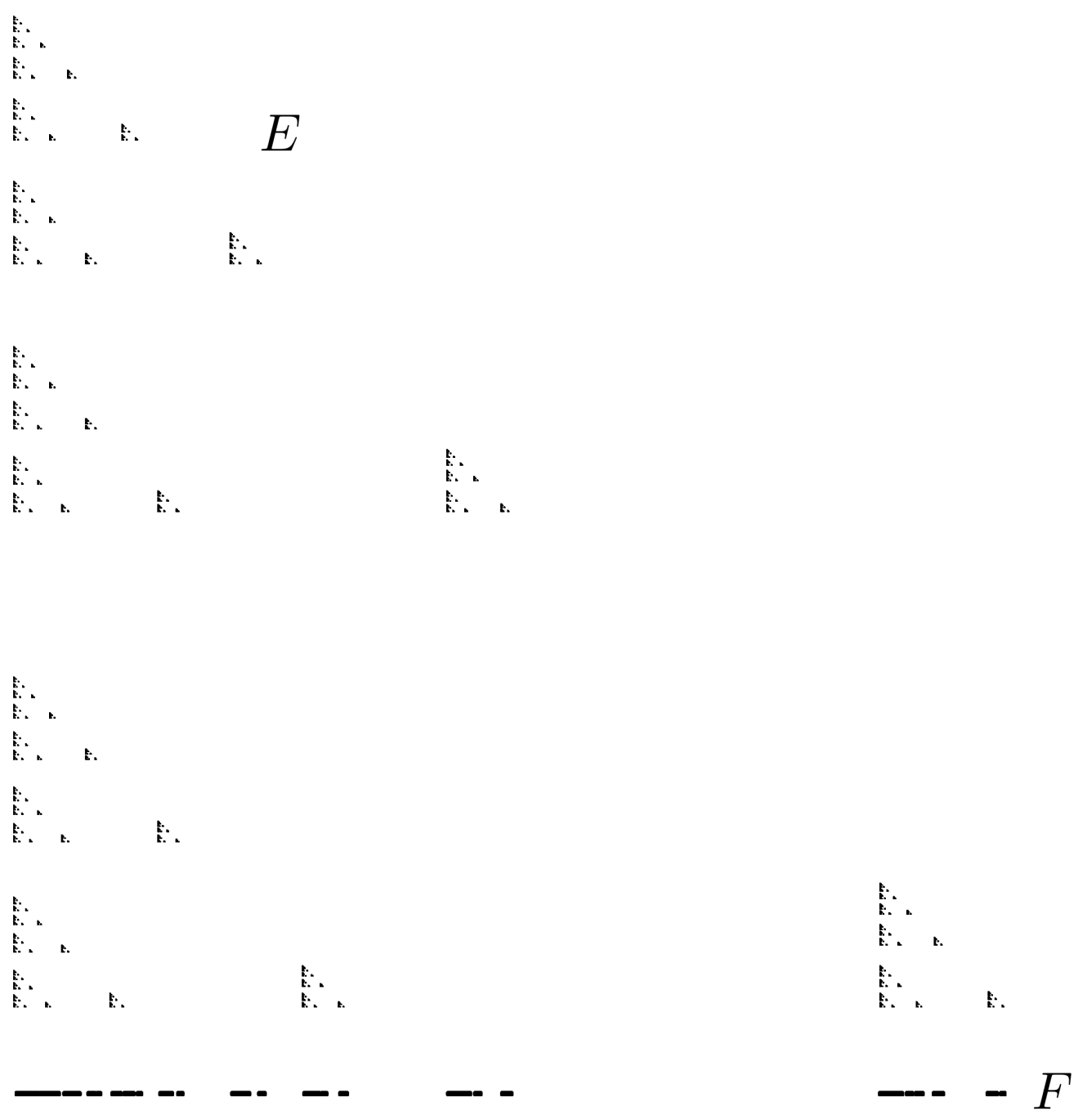

Figure 7.3 The set $E$ and its projection $\pi(E) = F$ for $\alpha = 1/3$, $\beta = 1/2$ and $\gamma = 1/8$. Here $\dim_A E \approx 0.9582 < \dim_A F = 1$.

We equip this space with the metric induced by the supremum norm restricted to the unit ball. More precisely, we define

$$d(f, g) = \|(f - g)|_{B(0,1)}\|_\infty$$

for similarities $f, g$.

**Definition 7.2.3** An IFS of similarities satisfies the *weak separation property (WSP)* if

$$I \notin \overline{\mathcal{E} \setminus \{I\}}.$$

Zerner [282] proved that the OSC is equivalent to $I \notin \overline{\mathcal{E}}$ which makes the OSC strictly stronger than the WSP.

**Theorem 7.2.4** Let $F \subseteq \mathbb{R}$ be a self-similar set containing at least two points. If the defining IFS satisfies the WSP, then $\dim_A F = \dim_H F$ and otherwise $\dim_A F = 1$.

Using this result, and some work concerning the Hausdorff measure of self-similar sets, Farkas and Fraser [81, Corollary 3.2] proved the following set of equivalences. A neat consequence of this is that the WSP is intrinsic to the set $F$, provided the Hausdorff dimension is not full: either all IFSs defining $F$ satisfy the WSP or none do. This does not hold for the OSC or SSC, for example.

**Theorem 7.2.5** Let $F \subseteq \mathbb{R}$ be a self-similar set containing at least two points with $s = \dim_{\mathrm{H}} F < 1$. Then the following are equivalent:

(i) the defining IFS satisfies the WSP
(ii) $\mathcal{H}^s(F) > 0$
(iii) $0 < \mathcal{H}^s(F) < \infty$
(iv) $F$ is Ahlfors regular
(v) $\dim_{\mathrm{A}} F = \dim_{\mathrm{H}} F$.

*Proof* The fact that (i) implies (ii) was proved by Zerner [282, Corollary of Proposition 2]. It follows from [66] that $\mathcal{H}^s(F) < \infty$ and therefore (ii) and (iii) are equivalent. The fact that (ii) implies (iv) follows from [81, Theorem 2.1 and Lemma 1.1] which gives that for any $r > 0$, $\mathcal{H}^s(E) = \mathcal{H}^s_r(E)$ for any $\mathcal{H}^s$-measurable subset $E \subseteq F$. Therefore, for any $x \in F$ and $0 < r \leqslant |F|$,

$$\mathcal{H}^s(F \cap B(x,r)) = \mathcal{H}^s_{2r}(F \cap B(x,r)) \leqslant 2^s r^s$$

where the final inequality comes by considering the trivial cover $\{B(x,r)\}$. Moreover, the reverse inequality is straightforward since $F \cap B(x,r)$ contains $S_{\boldsymbol{i}}(F)$ for some $\boldsymbol{i} \in \mathcal{I}^*$ such that the similarity ratio of $S_{\boldsymbol{i}}$ is at least $cr$ for a uniform constant $c > 0$ depending only on the IFS. Therefore, using the scaling property for Hausdorff measure (see [70, Proposition 2.2]),

$$\mathcal{H}^s(F \cap B(x,r)) \geqslant \mathcal{H}^s(S_{\boldsymbol{i}}(F)) = c^s r^s \mathcal{H}^s(F) > 0$$

which proves that $F$ is Ahlfors regular. Theorem 6.4.1 gives that (iv) implies (v). Finally, Theorem 7.2.4 and the fact that $s < 1$ shows that (v) implies (i), completing the proof. □

The self-similar set constructed in Theorem 7.2.1 clearly fails the WSP. This constitutes one of two well-known methods for constructing such examples. The other is due to Bandt and Graf [13] and works as follows. One advantage of the Bandt-Graf approach is that all of the similarity ratios can be the same. Let $c \in (0,1)$ be the common similarity ratio and consider the IFS consisting of the three maps

$$S_1(x) = cx, \qquad S_2(x) = c(x+t), \qquad S_3(x) = cx + (1-c) \tag{7.3}$$

acting on $[0,1]$ where $t \in (0,(1-c)/c)$ is

$$t = t(c) = \frac{1-c}{c} \sum_{k=0}^{\infty} c^{2^k}. \tag{7.4}$$

The elements $S_{\boldsymbol{i}}^{-1} \circ S_{\boldsymbol{j}}$ where $|\boldsymbol{i}| = |\boldsymbol{j}| = n$ are of the form

$$x \mapsto x + \sum_{k=1}^{n} c^{-k} a_k \tag{7.5}$$

where

$$a_k \in \{0, \pm ct, \pm(1-c), \pm(ct-(1-c))\}.$$

Consider $n = 2^m$ for some large integer $m \geqslant 1$ and choose the coefficients $\{a_k\}_{k=1}^n$ as

$$a_k = \begin{cases} -(1-c) & \text{if } k = 2^m - 2^l \text{ for some } l = 0, \ldots, m \\ ct & \text{if } k = n \\ 0 & \text{otherwise} \end{cases}$$

in which case (7.5) becomes

$$x \mapsto x + c^{-2^m}(ct) - (1-c)\sum_{l=0}^{m} c^{2^l - 2^m} = x + (1-c) \sum_{k=m+1}^{\infty} c^{2^k - 2^m}$$

using (7.4). Therefore, choosing $m$ large we can approximate the identity arbitrarily well by non-identity elements, proving the IFS fails the WSP.

Theorem 7.2.4 says that for self-similar sets in $\mathbb{R}$ the Assouad dimension is either equal to the Hausdorff dimension or equal to the ambient spatial dimension (which in this case is 1). We can use the example of Bandt and Graf to show that this does not hold in higher dimensions, following [96]. Consider the IFS acting on $[0,1]^2$ consisting of the maps

$$S_1(x) = x/5, \qquad S_2(x) = x/5 + (t/5, 0),$$

$$S_3(x) = x/5 + (4/5, 0), \qquad S_4(x) = x/5 + (0, 4/5)$$

where $t = t(1/5)$, as in (7.4). Let $E$ denote the attractor, see Figure 7.4. Note that $E$ is not contained in any line (and thus the true ambient dimension is 2) and also satisfies $\dim_{\mathrm{H}} E \leqslant \log 4/\log 5 < 1$ by the Hutchinson-Moran formula (7.2). Also, note that $E$ contains the attractor of the IFS (7.3) with $c = 1/5$, which we denote by $F$, and therefore $\dim_A E \geqslant \dim_A F = 1 > \dim_H E$. It remains to show that the Assouad dimension of $E$ is not 2, but this follows easily by noting $E$ is contained in the product of its projection onto the horizontal axis, which is $F$, with its projection onto the vertical axis, which is a self-similar set satisfying the SSC which we denote by $E'$. Therefore, applying known upper bounds for the Assouad dimension of products, see Corollary 10.1.2, and the Hutchinson-Moran formula (7.2),

$$\dim_A E \leqslant \dim_A F \times E' \leqslant \dim_A F + \dim_A E' = 1 + \log 2/\log 5 < 2 \quad (7.6)$$

as required.

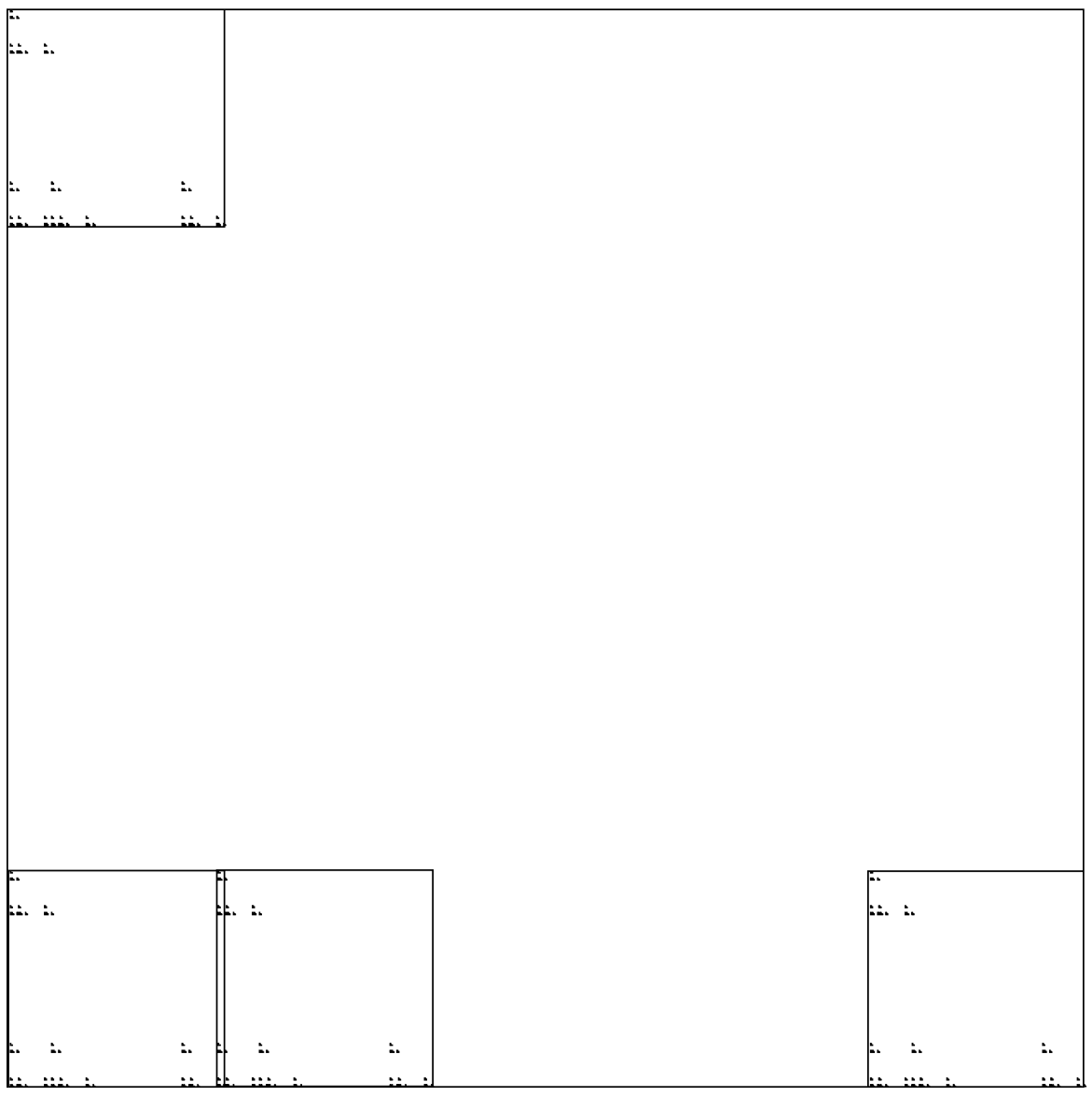

Figure 7.4 A self-similar set $E$ in the plane with $\dim_H E < \dim_A E < 2$. We included the images of the unit square under the four maps in the defining IFS.

The Assouad dimension of self-similar sets in $\mathbb{R}^d$ was considered by García [117] and comprehensive results were obtained. The situation is necessarily rather more complicated than that of Theorem 7.2.4 but it turns out that, provided the WSP fails, the Assouad dimension of $F$ can be expressed as the dimension of the subspace spanned by the 'overlapping directions' plus the Assouad dimension of the projection of $F$ onto the orthogonal complement of this subspace. We refer the reader to [117, Definition 1.1 and Theorem 1.4] for the details. However, we note that for the example $E$ from Figure 7.4 the subspace spanned by the overlapping directions is the horizontal axis and so García's theorem yields that in fact

$$\dim_{\mathrm{A}} E = 1 + \log 2 / \log 5,$$

improving on (7.6). As a final example where García's theorem applies, consider the overlapping variant of the Sierpiński triangle where we include the additional similarity $x \mapsto x/3$, see Figure 7.2.[2] The Hausdorff dimension of the attractor is bounded above by the similarity dimension $s \approx 1.7999$ (in fact this upper bound is strict since it is easy to find exact overlaps) and García's theorem implies that the Assouad dimension is equal to 2.

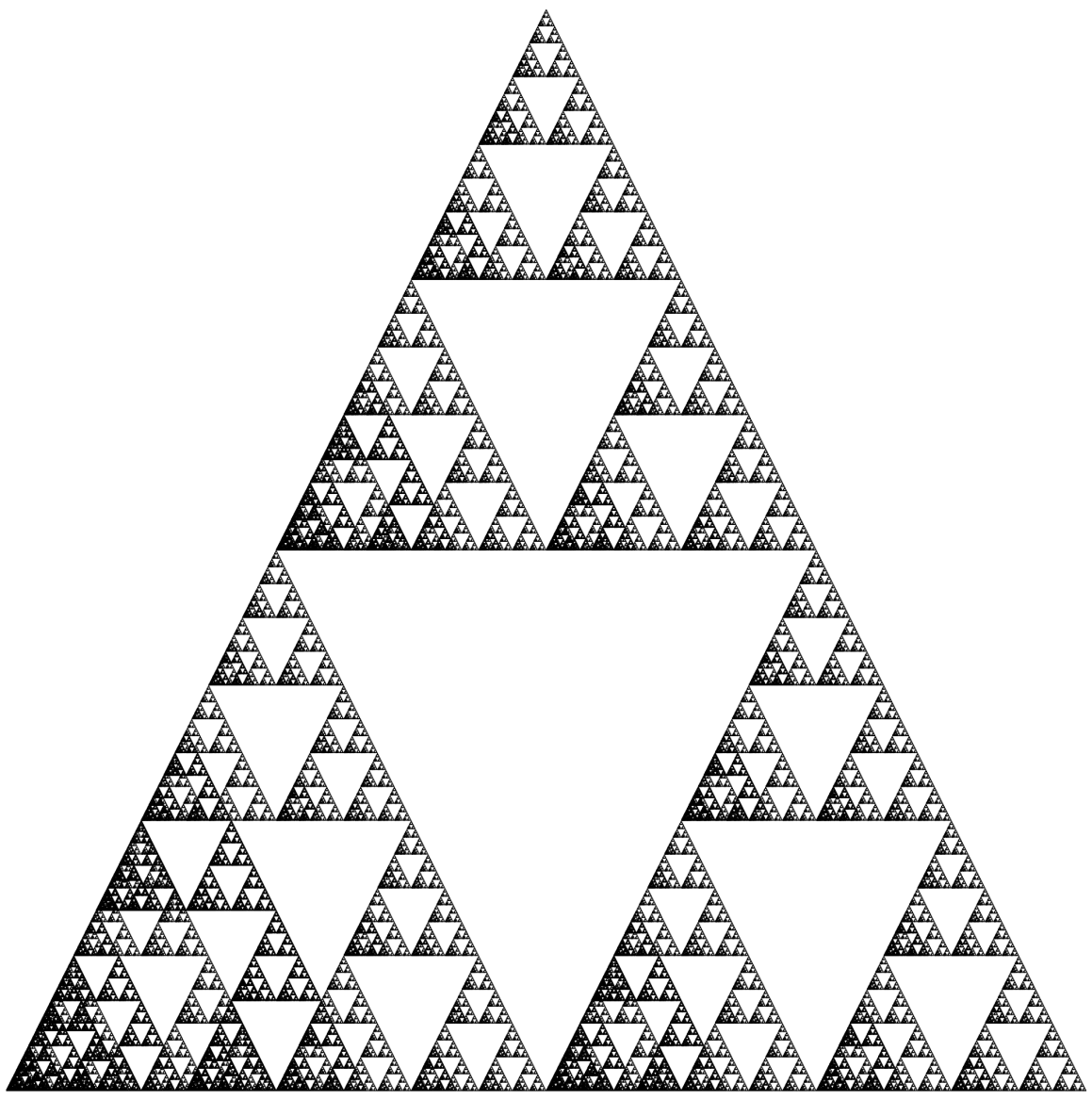

Figure 7.5 An overlapping variant of the Sierpiński triangle with $\dim_{\mathrm{H}} F < 1.7999 < \dim_{\mathrm{A}} F = 2$.

The following elegant result holds in the 'full rotations' case, see [117, Corollary 1.6].

[2] Eric Olson specifically discussed this example during a Warwick ETDS Seminar in February 2014, conjecturing that the Assouad dimension is equal to 2.

**Theorem 7.2.6** Let $F \subseteq \mathbb{R}^d$ be a self-similar set containing at least two points. If the defining IFS does not satisfy the WSP and the group generated by the orthogonal components of the defining similarities is dense in $SO(d)$ or $O(d)$, then

$$\dim_{\mathrm{A}} F = d.$$

Essentially the rotational components of the defining similarities plus failure of the WSP guarantee that the span of the overlapping directions is the whole of $\mathbb{R}^d$.

## 7.3 The Assouad spectrum of self-similar sets

Compared to the Assouad dimension, the Assouad spectrum is less well-understood for self-similar sets, even in the line. Some estimates based on the $L^q$-spectrum of an associated self-similar measure are provided in [112] and there is some evidence to suggest that if $F \subseteq \mathbb{R}^d$ is self-similar, then

$$\dim_{\mathrm{B}} F = \dim_{\mathrm{A}}^{\theta} F = \dim_{\mathrm{qA}} F$$

for all $\theta \in (0, 1)$. Relying on a powerful result of Shmerkin and work of Hochman, this can be verified in some cases.

**Theorem 7.3.1** If $F \subseteq \mathbb{R}$ is self-similar such that there are no exact overlaps in the construction and the defining parameters for the IFS are algebraic, then

$$\dim_{\mathrm{B}} F = \dim_{\mathrm{A}}^{\theta} F = \dim_{\mathrm{qA}} F$$

for all $\theta \in (0, 1)$.

*Proof* In the setting of the theorem, Hochman [132] proved that the box dimension is given by the similarity dimension, that is $\dim_{\mathrm{B}} F = s$, and we may assume $s < 1$ since otherwise there is nothing to prove. Let $\mu$ be the self-similar measure supported by $F$ associated to the Bernoulli weights $\{c_i^s\}$, see (7.2). Shmerkin [254, Theorem 6.6] strengthened Hochman's result by proving that the $L^q$-spectrum of $\mu$ is affine, that is,

$$\tau_\mu(q) = s(1 - q)$$

for $q > 0$. This is a strengthening of Hochman's result since $\tau_\mu(0)$ returns the box dimension of $F$. If the OSC is satisfied, then these results are straightforward, and go back much further, and so the novelty of Shmerkin's result (and Theorem 7.3.1) is in the overlapping case. Recall the (upper) $L^q$-*spectrum* is defined by

$$\tau_\mu(q) = \limsup_{r\to 0} \frac{\log M_r^q(\mu)}{-\log r}$$

where

$$M_r^q(\mu) = \sup\Big\{ \sum_i \mu(U_i)^q \Big\}$$

with the supremum taken over all centred packings $\{U_i\}_i$ of $F$ by balls of radius $r$. Let $\varepsilon \in (0,1)$ and choose $q > 1$ large enough so that $-(s-\varepsilon)q > s(1-q) = \tau_\mu(q)$. It follows that there exists a constant $C \geqslant 1$ such that

$$M_r^q(\mu) \ \leqslant \ Cr^{(s-\varepsilon)q}$$

for all $r \in (0,1)$. Then, since $\{B(x,r)\}$ is an $r$-packing for any particular choice of $x$ and $r \in (0,1)$, we have

$$\mu\left(B(x,r)\right)^q \ \leqslant \ M_r^q(\mu) \ \leqslant \ Cr^{(s-\varepsilon)q}$$

and therefore

$$\mu\left(B(x,r)\right) \ \leqslant \ C^{1/q} r^{(s-\varepsilon)} \tag{7.7}$$

for all $x \in F$ and $r \in (0,1)$.

For $\boldsymbol{i} = (i_1, i_2, \dots, i_k) \in \mathcal{I}^*$, we write $c_{\boldsymbol{i}} = c_{i_1} c_{i_2} \cdots c_{i_k}$ for the similarity ratio of $S_{\boldsymbol{i}}$. Let $x \in F$, $r > 0$ be small, and $\theta \in (0,1)$. For $\delta \in (0,1)$, let $\mathcal{I}(\delta) \subseteq \mathcal{I}^*$ be the $\delta$-*stopping*, defined by

$$\mathcal{I}(\delta) \ = \ \{\boldsymbol{i} \in \mathcal{I}^* \ : \ c_{\boldsymbol{i}} \leqslant \delta < c_{\boldsymbol{i}^-}\}$$

and

$$M(x, r^\theta) \ = \# \left\{\boldsymbol{i} \in \mathcal{I}(r^\theta) \ : \ S_{\boldsymbol{i}}(F) \cap B(x, r^\theta) \neq \varnothing \right\}.$$

All of the sets $S_{\boldsymbol{i}}(F)$ contributing to $M(x, r^\theta)$ lie completely inside the ball $B\left(x, ar^\theta\right)$ for some fixed constant $a > 1$ which is independent of $r^\theta$ and $x$. Moreover, using (7.7) and the definition of $\mu$,

$$M(x, r^\theta)\, c_{\min}^s r^{\theta s} \leqslant \mu\left(B\left(x, ar^\theta\right)\right) \ \leqslant \ C^{1/q} a^{(s-\varepsilon)} r^{\theta(s-\varepsilon)}$$

where $c_{\min} = \min_{i \in \mathcal{I}}\{c_i\}$, and so

$$M(x, r^\theta) \ \leqslant \ c_{\min}^{-s} C^{1/q} a^{(s-\varepsilon)} r^{-\varepsilon\theta}. \tag{7.8}$$

We wish to cover $B(x, r^\theta)$ by sets of diameter less than or equal to $r$

and to do this we iterate the construction of $F$ inside each of the sets $S_{\boldsymbol{i}}(F)$ contributing to $M(x, r^\theta)$. More precisely, the sets

$$\left\{S_{\boldsymbol{ij}}(F) \ : \ \boldsymbol{i} \in \mathcal{I}(r^\theta) \text{ with } S_{\boldsymbol{i}}(F) \cap B(x, r^\theta) \neq \varnothing, \ \boldsymbol{j} \in \mathcal{I}\left(r^{1-\theta}\right)\right\}$$

all have diameters at most $r$ and form a cover of $B(x, r^\theta)$. Therefore

$$N\left(B(x, r^\theta) \cap F, \ r\right) \leqslant M(x, r^\theta)\left(\#\mathcal{I}\left(r^{1-\theta}\right)\right) \tag{7.9}$$

and so it remains to estimate $\#\mathcal{I}\left(r^{1-\theta}\right)$. We have

$$\begin{aligned}\left(\#\mathcal{I}\left(r^{1-\theta}\right)\right) r^{(s+\varepsilon)(1-\theta)} \ &\leqslant \ c_{\min}^{-2} \sum_{\boldsymbol{j} \in \mathcal{I}(r^{1-\theta})} c_{\boldsymbol{j}}^{s+\varepsilon} \leqslant c_{\min}^{-2} \sum_{k=1}^{\infty} \sum_{\boldsymbol{j} \in \mathcal{I}^k} c_{\boldsymbol{j}}^{s+\varepsilon} \\ &= c_{\min}^{-2} \sum_{k=1}^{\infty} \left(\sum_{j \in \mathcal{I}} c_j^{s+\varepsilon}\right)^k \\ &=: b < \infty\end{aligned}$$

using the Hutchinson-Moran formula (7.2) and therefore

$$\#\mathcal{I}\left(r^{1-\theta}\right) \ \leqslant \ b r^{(s+\varepsilon)(\theta-1)}. \tag{7.10}$$

Thus, by (7.8), (7.9) and (7.10),

$$\begin{aligned}N\left(B(x, r^\theta) \cap F, \ r\right) &\leqslant M(x, r^\theta)\left(\#\mathcal{I}\left(r^{1-\theta}\right)\right) \\ &\leqslant c_{\min}^{-s} C^{1/q} a^{(s-\varepsilon)} r^{-\varepsilon\theta} \ b r^{(s+\varepsilon)(\theta-1)} \\ &\leqslant c_{\min}^{-s} C^{1/q} a^{(s-\varepsilon)} b \left(r^{\theta-1}\right)^{s+\frac{\varepsilon}{1-\theta}}\end{aligned}$$

which proves

$$\dim_{\mathrm{A}}^\theta F \leqslant s + \frac{\varepsilon}{1-\theta}$$

and letting $\varepsilon \to 0$ proves the result. □

Following [252], the requirement in Theorem 7.3.1 that there are no exact overlaps and the defining parameters for the IFS are algebraic can be replaced with a weaker but more technical assumption that there is no 'super-exponential concentration of cylinders'. The super-exponential concentration condition was first considered in [132] where it was proved that if the defining parameters are algebraic then super-exponential concentration is equivalent to the presence of exact overlaps. In fact, until 2019 the *only* known mechanism for generating super-exponential concentration was via exact overlaps. However, papers by Baker [11] and Bárány-Käenmäki [20] proved that super-exponential concentration without exact overlaps is possible. For example, [11] demonstrates this by constructing an explicit non-algebraic IFS using continued fraction expansions.

## 7.4 Dimensions of self-similar measures

The Assouad and box dimensions of self-similar measures satisfying the SSC were computed in [98, Theorem 2.4] and the lower dimension in [127, Proposition 2.3]. The fact that such measures are doubling goes back much further, for example to [217].

**Theorem 7.4.1** Let $\mu$ be a self-similar measure satisfying the SSC. Then $\mu$ is doubling and

$$\dim_{\mathrm{A}} \mu = \dim_{\mathrm{B}} \mu = \max_{i \in \mathcal{I}} \frac{\log p_i}{\log c_i}$$

and

$$\dim_{\mathrm{L}} \mu = \min_{i \in \mathcal{I}} \frac{\log p_i}{\log c_i}.$$

*Proof* Let $\delta$ be the minimal distance between distinct sets $S_i(F)$ and $S_j(F)$, that is,

$$\delta = \min_{i \neq j} \inf_{\substack{x \in S_i(F) \\ y \in S_j(F)}} |x - y| > 0$$

which is strictly positive since the SSC holds. Given $x \in F$ with

$$\Pi(i_1, i_2, \dots) = x$$

for $(i_1, i_2, \dots) \in \mathcal{I}^{\mathbb{N}}$, which is uniquely defined since the SSC holds, and $r > 0$ small, let $n(x,r)$ be the largest integer such that $2r/\delta \leqslant c_{i_1} c_{i_2} \cdots c_{i_{n(x,r)}}$ and $m(x,r)$ be the smallest non-negative integer such that

$$\Pi([i_1, \dots, i_{m(x,r)}]) \subseteq B(x, r/2).$$

Using these definitions,

$$p_{i_1} \cdots p_{i_{m(x,r)}} \leqslant \mu(B(x,r))$$

and

$$\Pi([i_1, \dots, i_{m(x,r)-1}]) \not\subseteq B(x, r/2)$$

and so $c_{i_1} \cdots c_{i_{m(x,r)-1}} \geqslant r/2$. Moreover,

$$\mu(B(x,r)) \leqslant p_{i_1} \cdots p_{i_{n(x,r)}}$$

since $B(x,r) \cap F \subseteq \Pi([i_1, \dots, i_{n(x,r)}])$.

Let $i \in \mathcal{I}$ be such that $t := \frac{\log p_i}{\log c_i}$ is maximised and $x = \Pi(\bar{i})$ where $\bar{i} = (i, i, \dots) \in \mathcal{I}^{\mathbb{N}}$. It follows that $n(x,r) > \log(2r/\delta)/\log c_i - 1$ and therefore

$$\mu(B(x,r)) \leqslant p_i^{n(x,r)} \leqslant p_i^{-1}(2r/\delta)^{\log p_i/\log c_i} = p_i^{-1}(2/\delta)^t r^t$$

which shows $\dim_{\mathrm{A}} \mu \geqslant \overline{\dim}_{\mathrm{B}} \mu \geqslant \underline{\dim}_{\mathrm{B}} \mu \geqslant t$.

For the upper bound, let $x \in F$, $0 < r < R$ and assume that $R/r$ is large enough to guarantee $n(x,R) < m(x,r)$. Hence

$$\begin{aligned} \frac{\mu(B(x,R))}{\mu(B(x,r))} &\leqslant \frac{p_{i_1} \cdots p_{i_{n(x,R)}}}{p_{i_1} \cdots p_{i_{m(x,r)}}} \\ &= \prod_{l=n(x,R)+1}^{m(x,r)} \frac{1}{p_{i_l}} \\ &= \prod_{l=n(x,R)+1}^{m(x,r)} \left(\frac{1}{c_{i_l}}\right)^{\log p_{i_l}/\log c_{i_l}} \\ &\leqslant \prod_{l=n(x,R)+1}^{m(x,r)} \left(\frac{1}{c_{i_l}}\right)^{t} \\ &= \left(\frac{1}{c_{i_{n(x,R)+1}} c_{i_{m(x,r)}}}\right)^t \left(\frac{c_{i_1} c_{i_2} \cdots c_{i_{n(x,R)+1}}}{c_{i_1} c_{i_2} \cdots c_{i_{m(x,r)-1}}}\right)^t \\ &\leqslant \left(\frac{4}{\delta \min_{i \in \mathcal{I}} c_i^2}\right)^t \left(\frac{R}{r}\right)^t \end{aligned} \tag{7.11}$$

proving $\dim_{\mathrm{A}} \mu \leqslant t$, as required. The lower dimension calculation is similar and omitted. □

For comparative purposes, the Hausdorff dimension of a self-similar measure satisfying the OSC is given by

$$\dim_{\mathrm{H}} \mu = \frac{\sum_{i \in \mathcal{I}} p_i \log p_i}{\sum_{i \in \mathcal{I}} p_i \log c_i}, \tag{7.12}$$

see [70, 217]. This formula can be usefully interpreted as the *entropy* of the measure divided by the *Lyapunov exponent* of the measure. .

Consider the following simple example. Let $F$ be the middle third

Cantor set, see Figure 7.6, which is the self-similar set generated by the maps $x \mapsto x/3$ and $x \mapsto x/3 + 2/3$. Since there are only two maps, the family of self-similar measures supported by $F$ is the 1-parameter family $\mu_p$ corresponding to the probability vector $\{p, 1-p\}$ as $p$ varies in $(0,1)$. The dimensions of $\mu_p$ as functions of $p$ are

$$\dim_{\mathrm{L}} \mu_p = \min\left\{-\frac{\log p}{\log 3}, -\frac{\log(1-p)}{\log 3}\right\},$$

$$\dim_{\mathrm{H}} \mu_p = \frac{-p\log p - (1-p)\log(1-p)}{\log 3}$$

and

$$\dim_{\mathrm{A}} \mu_p = \dim_{\mathrm{B}} \mu_p = \max\left\{-\frac{\log p}{\log 3}, -\frac{\log(1-p)}{\log 3}\right\}.$$

Figure 7.6 The middle third Cantor set.

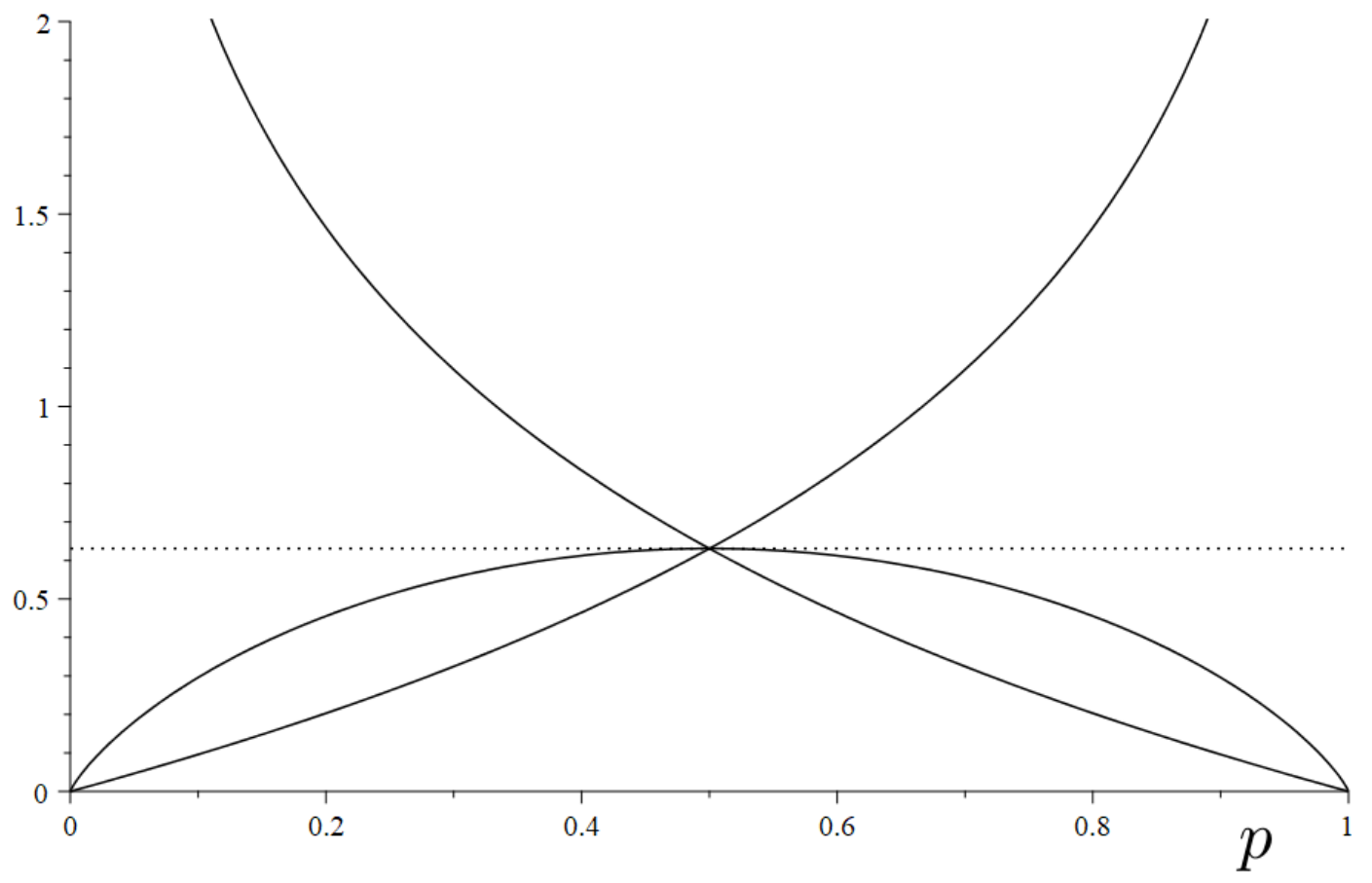

Figure 7.7 Plots of $\dim_{\mathrm{L}} \mu_p \leqslant \dim_{\mathrm{H}} \mu_p \leqslant \dim_{\mathrm{B}} \mu_p = \dim_{\mathrm{A}} \mu_p$ as functions of $p \in (0,1)$. The lower, Hausdorff, box and Assouad dimensions of the middle third Cantor set are all $\log 2/\log 3$ (dotted line) and are all realised by the measure corresponding to $p = 1/2$.

The assumption that the SSC is satisfied seems rather strong in Theorem 7.4.1. However, it is necessary since there are very straightforward examples of self-similar measures satisfying only the OSC which fail to even be doubling. Note how the coefficient of $(R/r)^t$ at the end of the proof of Theorem 7.4.1 blows up as $\delta \to 0$, see (7.11). This gives some indication that the SSC is playing a key role. For example, consider the IFS consisting of the maps $S_1 : x \mapsto x/2$ and $S_2 : x \mapsto (x+1)/2$ and the self-similar measure $\mu$ associated with $\{p, 1-p\}$ for $p > 1/2$. The support of $\mu$ is the attractor which in this case is $[0,1]$. Let $k \geqslant 1$ be large and $x = 1/2 - 2^{-k} \in [0,1]$ such that

$$B(x, 2^{-k}) = \Pi([1, 2, 2, \dots, 2])$$

where 2 has been chosen $k-2$ times. This means

$$B(x, 3 \cdot 2^{-k}) \supset \Pi([2, 1, 1, \dots, 1])$$

where this time 1 has been chosen $k-2$ times. Therefore

$$\frac{\mu(B(x, 3 \cdot 2^{-k}))}{\mu(B(x, 2^{-k}))} \geqslant \frac{(1-p)p^{k-2}}{p(1-p)^{k-2}} = \left(\frac{p}{1-p}\right)^{k-3} \to \infty$$

proving that $\mu$ is not doubling and hence $\dim_{\mathrm{A}} \mu = \infty$. In fact, this argument can be adapted to show that $\dim_{\mathrm{qA}} \mu = \infty$ and so $\mu$ is not even 'quasi-doubling', see [125, Example 3.2]. The paper [125] contains a thorough investigation of the quasi-Assouad dimension of self-similar measures which fail the SSC and the paper [127] considers the lower dimension. For example, [127, Corollary 3.5] shows that the lower dimension of any self-similar measure which is not an atom is positive. Recall that a measure is called an *atom* if it is supported on a single point. This is perhaps surprising since the analogous statement for Assouad dimension fails dramatically, since the Assouad dimension can be infinite very easily.

**Theorem 7.4.2** Let $\mu$ be a self-similar measure which is not an atom. Then

$$\dim_{\mathrm{L}} \mu > 0.$$

# 8
# Self-affine sets

Self-affine sets are another special case of IFS attractors. Since the defining maps may contract distance by different amounts in different directions, the dimension theory of self-affine sets and measures is more complicated, and much richer, than that of self-similar sets. In this chapter we study self-affine sets in detail, paying particular attention to Bedford-McMullen carpets where our theory can be developed explicitly. Formulae for the dimensions of Bedford-McMullen carpets are given in Theorem 8.3.1 and for the Assouad and lower spectra in Theorem 8.3.3. The dimension theory of self-affine measures is considered in Section 8.6, focusing on several explicit examples.

## 8.1 Self-affine sets and two strands of research

An *affine* map is a map consisting of two parts: a linear map and a translation. An attractor of an IFS is a *self-affine set* if all the defining maps are affine. We note that in $\mathbb{R}^d$ self-similar sets are a special case of self-affine sets, since similarity maps are necessarily affine. Even in the presence of separation properties, self-affine sets are notoriously difficult to handle in comparison with self-similar sets and there are still many fascinating open problems in the area. For example, does the box dimension of a self-affine set necessarily exist? (Recall that the box dimension of an arbitrary self-*similar* set always exists by the implicit theorems, see Corollary 6.4.4.) In the 1980s there were two important results in the dimension theory of self-affine sets, which began two interconnected but complementary strands of research. The survey paper [71] describes both strands of research as well as other areas connected with the dimension theory of self-affine sets.

The first strand of research was pioneered by Falconer, beginning with the seminal papers [65, 67] from 1988 and 1992 respectively. Falconer considers general self-affine sets where analysis of the dimension is a subtle problem but, instead of considering specific examples, Falconer established a formula for the Hausdorff and box dimension which holds 'generically' in an appropriate sense. Falconer's formula for the dimension is known as the *affinity dimension*, and is the self-affine analogue of the similarity dimension, see (7.2). The other strand of research began with the work of Bedford [26] and McMullen [210] in 1984, where a special case was considered — examples now known as Bedford-McMullen carpets. The simplicity of this model facilitated exact calculation of the dimensions, whilst maintaining many of the interesting features of self-affine sets. This strategy led to various different classes of self-affine carpet being introduced with increasing levels of generality, see [122, 15, 83, 86, 172, 168, 102]. Traditionally, a 'carpet' is a planar self-affine set generated by diagonal matrices. However, many of the key features of a carpet are manifest in more general models. For example, the sets considered in [86] also allow anti-diagonal matrices and the examples in [172] allow triangular matrices. Also, higher dimensional analogues have been considered, where the term 'carpet' is often replaced by 'sponge' [168, 49, 102].

Since 2010 there has been substantial progress on the dimension theory of self-affine sets, see [18, 79, 135]. The emerging phenomenon is that the affinity dimension gives the correct value for the Hausdorff and box dimensions in all but a few particular cases, much like the similarity dimension in the self-similar setting. The particular cases where the affinity dimension does *not* give the correct answer include many self-affine carpets, where dimension drop can occur due to excessive alignment of cylinders on small scales. Therefore, developments in the general theory have brought these two strands of research back together again: one strand studies the exceptions to the other.

## 8.2 Falconer's formula and the affinity dimension

The *singular values* of an invertible linear map, $A : \mathbb{R}^d \to \mathbb{R}^d$, are the positive square roots of the eigenvalues of $A^T A$. Viewed geometrically, these numbers are the lengths of the semi-axes of the image of the unit ball under $A$. Thus, roughly speaking, the singular values correspond to how much the map contracts (or expands) in different directions. For

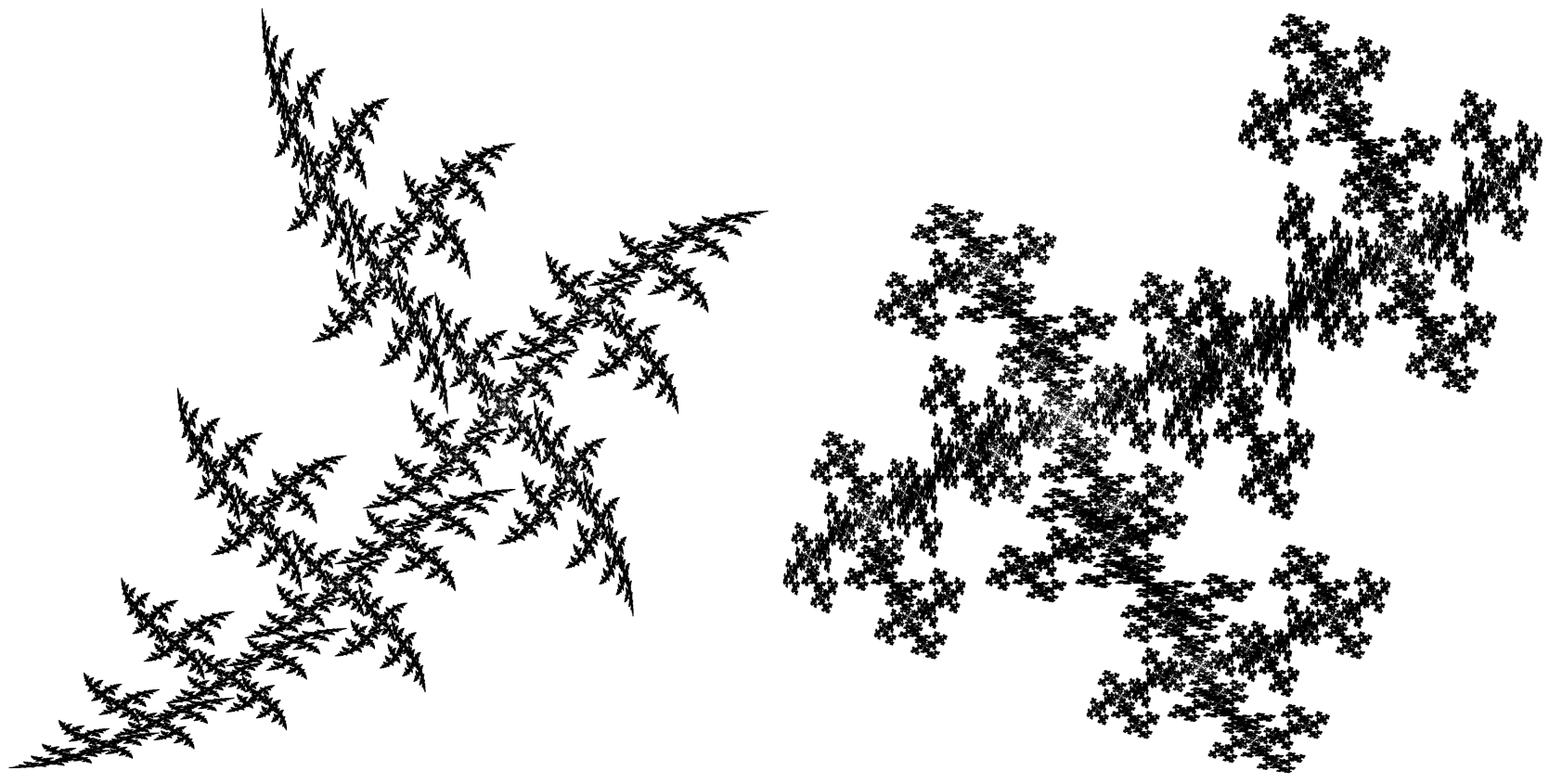

Figure 8.1 Two self-affine sets.

$s \in [0, d]$, the *singular value function* $\phi^s(A)$ is defined by

$$\phi^s(A) = \alpha_1 \alpha_2 \dots \alpha_{\lceil s \rceil - 1} \alpha_{\lceil s \rceil}^{s - \lceil s \rceil + 1} \tag{8.1}$$

where $\alpha_1 \geqslant \dots \geqslant \alpha_d > 0$ are the singular values of $A$. This function has played a vital role in the study of self-affine sets over the past 30 years. Let $\{A_i : i \in \mathcal{I}\}$ be a finite collection of linear contractions on $\mathbb{R}^d$, and write $m = \#\mathcal{I}$. The *affinity dimension* is defined by

$$s(A_i : i \in \mathcal{I}) = \inf \left\{ s : \sum_{k=1}^{\infty} \sum_{i \in \mathcal{I}^k} \phi^s(A_i) < \infty \right\}. \tag{8.2}$$

Falconer proved the following celebrated result in 1988 [65]. We write $\mathcal{L}^{md}$ to denote the $m$-fold product of $d$-dimensional Lebesgue measure, supported on the space $\times_{i \in \mathcal{I}} \mathbb{R}^d$.

**Theorem 8.2.1** Let $\{A_i\}_{i \in \mathcal{I}}$ be a finite collection of invertible linear contractions on $\mathbb{R}^d$. For $\mathbf{t} = (t_1, \dots, t_m) \in \times_{i \in \mathcal{I}} \mathbb{R}^d$, write $F_{\mathbf{t}}$ for the self-affine set generated by the IFS $\{A_i + t_i\}_{i \in \mathcal{I}}$. Then, for all $\mathbf{t}$,

$$\dim_{\mathrm{H}} F_{\mathbf{t}} \leqslant \overline{\dim}_{\mathrm{B}} F_{\mathbf{t}} \leqslant \min\{s(A_i : i \in \mathcal{I}), d\}.$$

If, in addition, each of the linear maps has singular values strictly less than 1/2, then, for $\mathcal{L}^{md}$-almost all $\mathbf{t}$,

$$\dim_{\mathrm{H}} F_{\mathbf{t}} = \dim_{\mathrm{B}} F_{\mathbf{t}} = \min\{s(A_i : i \in \mathcal{I}), d\}.$$

In fact, the initial proof in [65] required that the singular values be strictly less than 1/3 but this was relaxed to 1/2 by Solomyak [256], who also noted that 1/2 is the optimal constant, based on an example of Przytycki and Urbański [237]. See [150] for a similar result where measures other than Lebesgue measure $\mathcal{L}^{md}$ are used in determining the 'generic result'.

Theorem 8.2.1 shows that the affinity dimension is always an upper bound for the box and Hausdorff dimensions of a self-affine set. However, the affinity dimension is *not* necessarily an upper bound for the Assouad dimension, see [194, 88], or quasi-Assouad dimension, see [111]. Despite the failure of the sure upper bound, one might still hope that the Assouad dimension is almost surely given by the affinity dimension, or is at least almost surely constant. However, Fraser and Orponen [106] proved that this is not true based on an example of Peres, Simon and Solomyak [232], work of Farkas and Fraser [81] and of Fraser, Henderson, Olson and Robinson [96]. Recall that self-similar sets are a special case of self-affine sets.

**Theorem 8.2.2** Fix $c \in (1/5, 1/3)$ and let $s = -\log 3/\log c < 1$. For $\mathbf{t} = (t_1, t_2, t_3) \in \mathbb{R}^3$, let $F_{\mathbf{t}}$ denote the self-similar attractor of the IFS $\{x \mapsto cx + t_i\}_{i=1}^3$. Then there exist two non-empty disjoint open sets $U, V \subseteq \mathbb{R}^3$ such that:

(i) $\dim_{\mathrm{A}} F_{\mathbf{t}} = s$ for all $\mathbf{t} \in U$
(ii) $\dim_{\mathrm{A}} F_{\mathbf{t}} = 1$ for almost all $\mathbf{t} \in V$.

In particular, $\mathbf{t} \mapsto \dim_{\mathrm{A}} F_{\mathbf{t}}$ is not an almost surely constant function.

*Proof* The open set $U$ is easy to find: choose $\mathbf{t}' \in \mathbb{R}^3$ such that the SSC is satisfied, which can be done since $c < 1/3$, and observe that the SSC is satisfied for all $\mathbf{t}$ in some open neighbourhood, $U$, of $\mathbf{t}'$ in $\mathbb{R}^3$.

Finding the second set $V$ is more subtle. Observe that for all $\mathbf{t} \in \mathbb{R}^3$ the associated self-similar set $F_{\mathbf{t}}$ is just an affine image of the projection $\pi E$ of the self-similar set $E$ in the plane generated by the maps $x \mapsto cx$, $x \mapsto cx + (1-c, 0)$, $x \mapsto cx + (0, 1-c)$ for some $\pi$. These projections were considered by Peres, Simon and Solomyak [232] and they proved that there exists a non-empty open set $J$ of projections $\pi$ such that for almost all $\pi \in J$, $\mathcal{H}^s(\pi F) = 0$. This implies that there is a non-empty open set $V \subseteq \mathbb{R}^3$ such that for almost all $\mathbf{t} \in V$, $\mathcal{H}^s(F_{\mathbf{t}}) = 0$. Therefore, Theorems 7.2.4 and 7.2.5 imply that, for almost all $\mathbf{t} \in V$, we have $\dim_{\mathrm{A}} F_{\mathbf{t}} = 1$. □

We will discuss the projection results of Peres, Simon and Solomyak [232] in more detail in Section 10.2.

## 8.3 Self-affine carpets

We first recall the construction introduced by Bedford and McMullen. Take the unit square, $[0,1]^2$, and divide it up into an $m \times n$ grid for some fixed positive integers $n > m > 1$. Select a subset of the rectangles formed by the grid and consider the IFS consisting of the orientation preserving affine maps which map $[0,1]^2$ onto each chosen rectangle. The attractor of this IFS is the self-affine *Bedford-McMullen carpet*, see Figure 8.2.

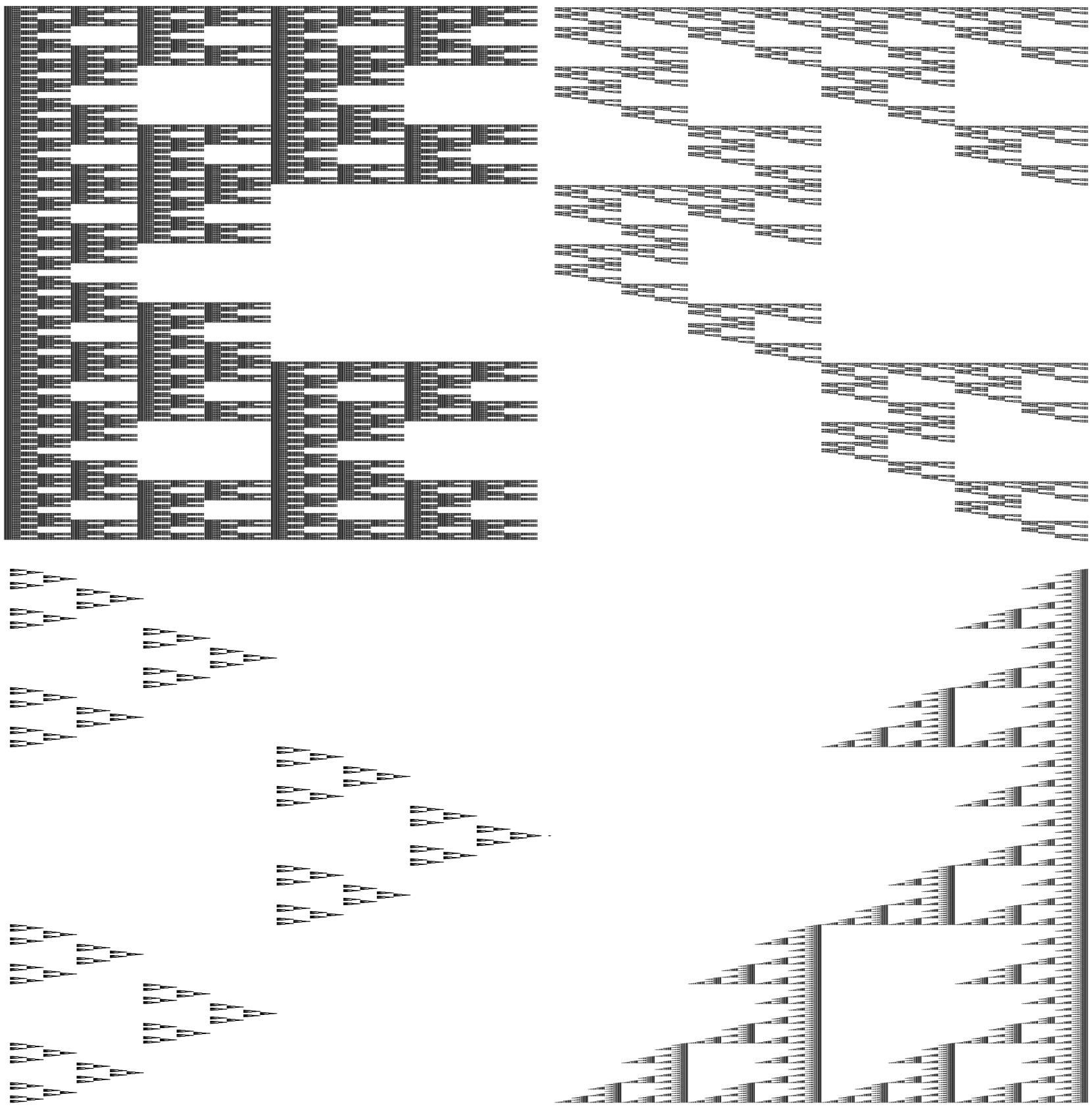

Figure 8.2 Four self-affine Bedford-McMullen carpets based on the $2 \times 3$ grid.

Bedford [26] and McMullen [210] independently obtained explicit formulae for the box and Hausdorff dimensions of the attractor and, more recently, Mackay [194] computed the Assouad dimension, and Fraser computed the lower dimension [88]. In order to state the dimension formulae, we require more notation. Let $N$ be the number of chosen rectangles, $N_0$ be the number of columns containing at least one chosen rectangle, and $N_i > 0$ be the number of rectangles chosen in the $i$th non-empty column.

**Theorem 8.3.1** Let $F$ be a Bedford-McMullen carpet. Then

$$\dim_{\mathrm{A}} F = \frac{\log N_0}{\log m} + \max_i \frac{\log N_i}{\log n},$$

$$\dim_{\mathrm{B}} F = \frac{\log N_0}{\log m} + \frac{\log(N/N_0)}{\log n},$$

$$\dim_{\mathrm{H}} F = \frac{\log \sum_{i=1}^{N_0} N_i^{\log m/\log n}}{\log m}$$

and

$$\dim_{\mathrm{L}} F = \frac{\log N_0}{\log m} + \min_i \frac{\log N_i}{\log n}.$$

Formulae for the Assouad and lower dimensions of more general self-affine carpets are available in [88, 101]. When computing the Assouad dimension given above, Mackay used weak tangents to obtain the lower bound, recall Theorem 5.1.2. We give a different proof below, based on a direct covering argument, but it is instructive to appreciate the weak tangents Mackay was able to construct. He proved that

$$\pi F \times E$$

is a weak tangent to $F$ where $\pi F$ is the projection of $F$ onto the first coordinate and $E$ is another self-similar set constructed using an IFS associated with the column maximising $N_i$. More precisely, $E$ is generated by an IFS consisting of $N_i$ contractions each with similarity ratio $1/n$. The set $E$ can be thought of as the 'maximal vertical slice' of $F$, which shows up at the level of weak tangents because of the affine distortion in the construction of $F$. Weak tangents were also used in [88] to bound the lower dimension from above.

Note that the Assouad, box, Hausdorff, and lower dimensions are equal if $N_i$ is constant whenever it is non-zero (the 'uniform fibres' case) and otherwise they are all distinct (the 'non-uniform fibres' case).

**Corollary 8.3.2** Let $F$ be a Bedford-McMullen carpet. Then either

$$\dim_{\mathrm{L}} F < \dim_{\mathrm{H}} F < \dim_{\mathrm{B}} F < \dim_{\mathrm{A}} F$$

or

$$\dim_{\mathrm{L}} F = \dim_{\mathrm{H}} F = \dim_{\mathrm{B}} F = \dim_{\mathrm{A}} F.$$

It is natural to ask if this dichotomy holds for more general self-affine carpets. However, the self-affine carpet depicted in Figure 8.3 satisfies $\dim_{\mathrm{L}} F = 1 < 3/2 = \dim_{\mathrm{H}} F = \dim_{\mathrm{B}} F = \dim_{\mathrm{A}} F$, see [88] for the details and also Table 17.1. This set is the attractor of the IFS acting on $[0,1]^2$ consisting of the maps

$$S_1((x,y)) = (x/5, y/4) + (0, 1/4), \quad S_2((x,y)) = (x/5, y/4) + (0, 1/2),$$

$$S_3((x,y)) = (4x/5, y/4) + (1/5, 0), \quad S_4((x,y)) = (4x/5, y/4) + (1/5, 3/4).$$

This is an example of the type of self-affine carpet introduced and studied by Barański [15]. The key feature which distinguishes it from the Bedford-McMullen setting is that the strongest contraction is not always in the same direction. The Assouad and lower dimensions of Barański carpets were computed in [88].

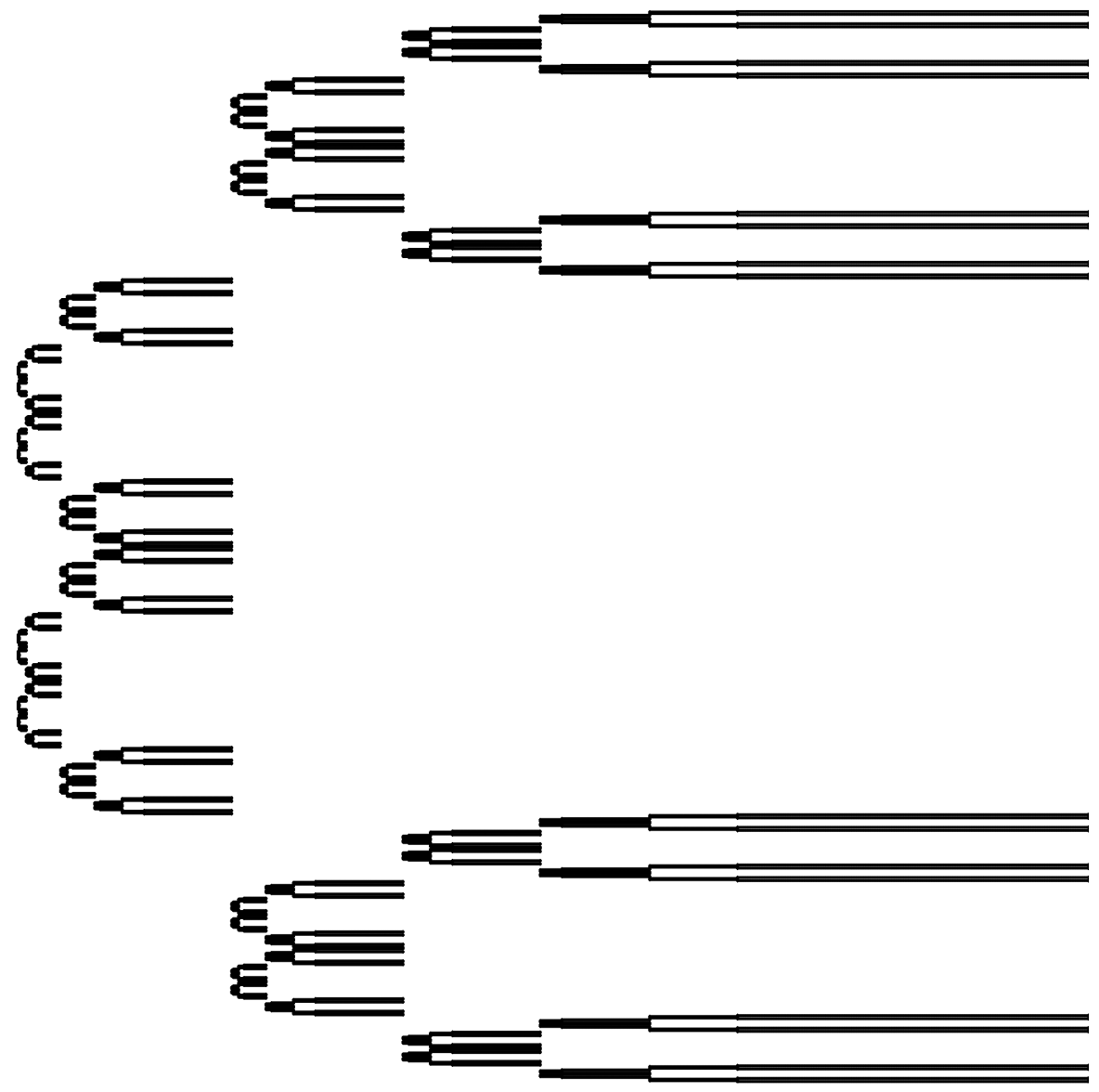

Figure 8.3 A self-affine set with $\dim_{\mathrm{L}} F < \dim_{\mathrm{H}} F = \dim_{\mathrm{B}} F = \dim_{\mathrm{A}} F$.

Finally we note that, despite how apparently easy it is to have lower dimension equal to zero, it is easy to see from Theorem 8.3.1 that the lower dimension of a self-affine carpet is always strictly positive. It turns out that the lower dimension of *any* self-affine set which is not a single point is strictly positive. This follows from [278], which proved that non-trivial self-affine sets are uniformly perfect, which is equivalent to strictly positive lower dimension, see Theorem 13.1.2. This result is a manifestation of the mild homogeneity possessed by IFS attractors.

We omit the proof of Theorem 8.3.1 but note that the formulae for the Assouad, box, and lower dimensions follow from the next theorem, which we will prove in part. We write $N_{\min} = \min_i N_i$ and $N_{\max} = \max_i N_i$, noting that with this notation

$$\dim_{\mathrm{A}} F \;=\; \frac{\log N_0}{\log m} + \frac{\log N_{\max}}{\log n},$$

and

$$\dim_{\mathrm{L}} F \;=\; \frac{\log N_0}{\log m} + \frac{\log N_{\min}}{\log n}.$$

The Assouad and lower spectra of Bedford-McMullen carpets were computed by Fraser and Yu [112].

**Theorem 8.3.3** Let $F$ be a Bedford-McMullen carpet. Then, for all $0 < \theta \leqslant \log m/\log n$,

$$\dim_{\mathrm{A}}^{\theta} F \;=\; \frac{\dim_{\mathrm{B}} F \;-\; \theta\left(\frac{\log(N/N_{\max})}{\log m} + \frac{\log N_{\max}}{\log n}\right)}{1-\theta}$$

and

$$\dim_{\mathrm{L}}^{\theta} F \;=\; \frac{\dim_{\mathrm{B}} F \;-\; \theta\left(\frac{\log(N/N_{\min})}{\log m} + \frac{\log N_{\min}}{\log n}\right)}{1-\theta}$$

and, for all $\log m/\log n \leqslant \theta < 1$,

$$\dim_{\mathrm{A}}^{\theta} F \;=\; \frac{\log N_0}{\log m} + \frac{\log N_{\max}}{\log n}$$

and

$$\dim_{\mathrm{L}}^{\theta} F \;=\; \frac{\log N_0}{\log m} + \frac{\log N_{\min}}{\log n}.$$

An immediate consequence of Theorem 8.3.3 is that for any Bedford-McMullen carpet with non-uniform fibres the Assouad spectrum is strictly smaller than the general upper bound given by Lemma 3.4.4 for all $\theta < \log m/\log n$.

Theorem 8.3.3 also shows that, in the non-uniform fibres case,

$$\dim_{\mathrm{L}}^{\theta} F \to \dim_{\mathrm{B}} F > \dim_{\mathrm{H}} F$$

as $\theta \to 0$. This provides the first concrete example of a compact set where the lower spectrum is not uniformly bounded above by the Hausdorff dimension. Recall that Theorem 6.3.1 shows this will happen for any self-affine set with distinct Hausdorff and lower box dimensions.

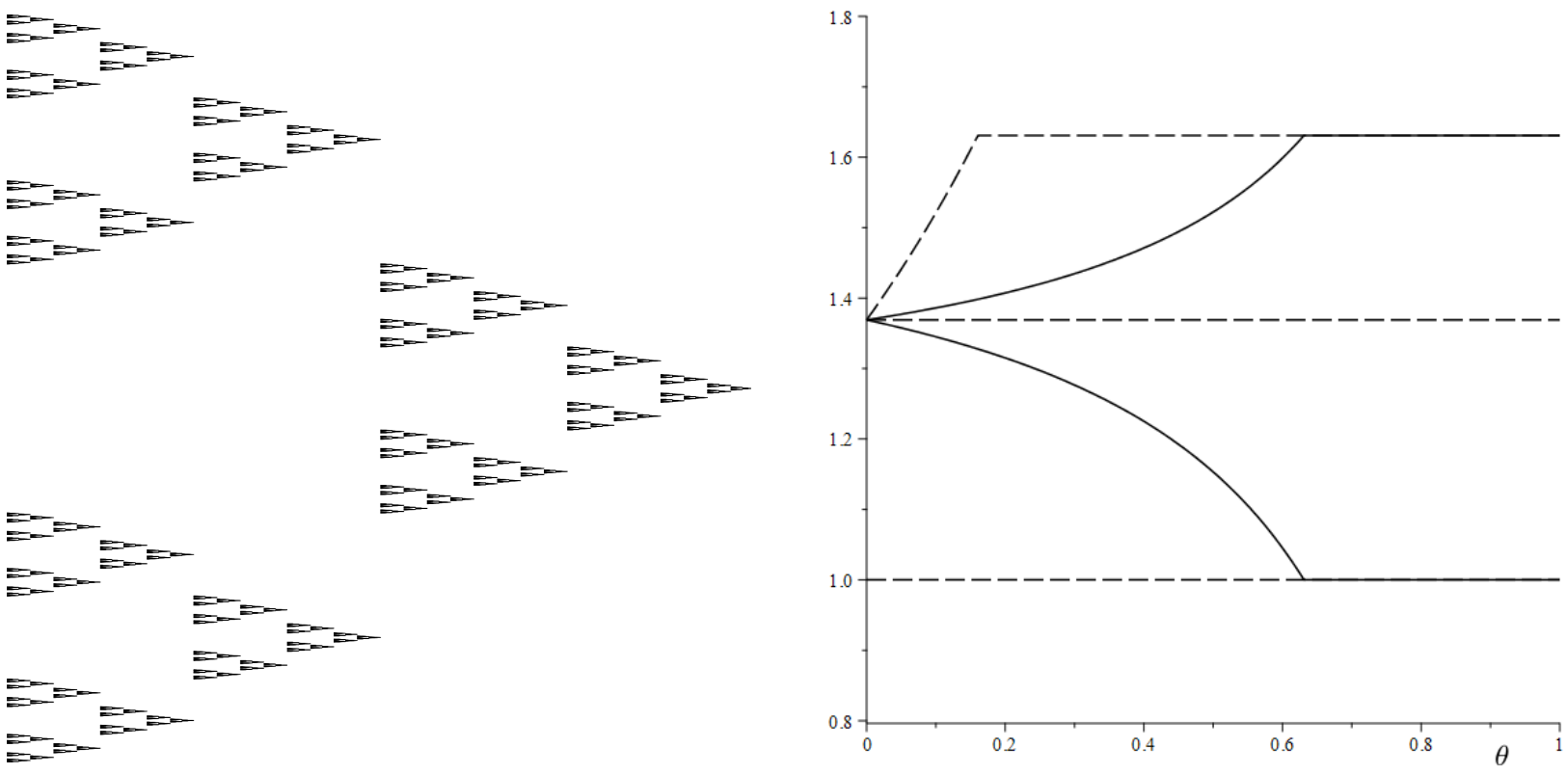

Figure 8.4 Left: a self-affine Bedford McMullen carpet, where $m = 2$, $n = 3$ and we have chosen two rectangles from the first column and one from the second. Right: plots of the Assouad and lower spectra with the basic bounds from Section 3.4.1 plotted as dashed lines for comparison.

Before proving Theorem 8.3.3 we need some more notation, which is standard in the study of self-affine carpets. Let $\mathcal{I} = \{0, \dots, m-1\}$ and $\mathcal{J} = \{0, \dots, n-1\}$ and label the $mn$ rectangles in the regular grid by $\mathcal{I} \times \mathcal{J}$ counting from bottom left to top right. Let $\mathcal{D} \subseteq \mathcal{I} \times \mathcal{J}$ correspond to the rectangles chosen when defining the Bedford-McMullen carpet. For each $d = (i,j) \in \mathcal{D}$, let $S_d : [0,1]^2 \to [0,1]^2$ be the contraction defined by

$$S_d(x,y) = (x/m + i/m,\ y/n + j/n)$$

which is the affine map which takes the unit square to the rectangle labelled by $d$. Therefore $F$ is the attractor of the IFS $\{S_d\}_{d \in \mathcal{D}}$. We are abusing notation slightly here since $\mathcal{D}$ is now the index set for the IFS, rather than $\mathcal{I}$. Instead we are using $\mathcal{I}$ to index the projection onto the first coordinate, which is a self-similar set. Let

$$\mathcal{D}^\infty = \{\mathbf{d} = (d_1, d_2, \dots) : d_l = (i_l, j_l) \in \mathcal{D}\}$$

be the set of infinite words over $\mathcal{D}$ and let $\Pi : \mathcal{D}^\infty \to [0,1]^2$ be the symbolic coding map, recalling that $F = \Pi(\mathcal{D}^\infty)$ is the Bedford-McMullen carpet. Let $\mathbf{d} \in \mathcal{D}^\infty$ and $r > 0$ be a small radius. Define $l_1(r), l_2(r)$ to be the unique natural numbers satisfying

$$m^{-l_1(r)} \leqslant r < m^{-l_1(r)+1} \tag{8.3}$$

and

$$n^{-l_2(r)} \leqslant r < n^{-l_2(r)+1}. \tag{8.4}$$

The *approximate square* centred at $\mathbf{d} = ((i_1, j_1), (i_2, j_2), \dots) \in \mathcal{D}^\infty$ with radius $r > 0$ is defined by

$$Q(\mathbf{d}, r) = \Big\{\mathbf{d}' = ((i'_1, j'_1), (i'_2, j'_2), \dots) \in \mathcal{D}^\infty \quad : \\ i'_l = i_l \ \text{ for all } l \leqslant l_1(r) \text{ and } j'_l = j_l \text{ for all } l \leqslant l_2(r)\Big\}.$$

The image of this symbolic object on the carpet, $\Pi\big(Q(\mathbf{d}, r)\big)$, is a subset of $F$ which contains the point $\Pi(\mathbf{d})$ and naturally sits inside a rectangle which is 'approximately a square' of radius $r$ in that it has width $m^{-l_1(r)} \in (r/m, r]$ and height $n^{-l_2(r)} \in (r/n, r]$. One obtains equivalent definitions of the Assouad and lower spectra if one replaces $B(x, R)$ by $\Pi\big(Q(\mathbf{d}, R)\big)$, that is, we may use approximate squares instead of balls.

The following lemma is the first indication that there will be a phase transition in the spectra at $\theta = \log m / \log n$.

**Lemma 8.3.4**

(i) If $\theta \leqslant \log m / \log n$ and $r \in (0,1)$ is small enough, then

$$l_2(r^\theta) \leqslant l_1(r^\theta) \leqslant l_2(r) \leqslant l_1(r).$$

(ii) If $\theta \geqslant \log m / \log n$ and $r \in (0,1)$ is small enough, then

$$l_2(r^\theta) \leqslant l_2(r) \leqslant l_1(r^\theta) \leqslant l_1(r).$$

*Proof* This follows immediately from the definitions of $l_1(\cdot)$ and $l_2(\cdot)$. □

*Proof of Theorem 8.3.3:* We will only prove the formula for the Assouad spectra and leave the lower spectra calculation (which is similar) as an exercise. First consider $0 < \theta < \log m / \log n$. Let

$$\mathbf{d} = ((i_1, j_1), (i_2, j_2), \dots) \in \mathcal{D}^\infty$$

and $r \in (0, 1)$ be small. Consider the approximate square $\Pi\big(Q(\mathbf{d}, r^\theta)\big)$ and observe that it is made up of several horizontal strips with width $m^{-l_1(r^\theta)}$, which is the same as that of the approximate square, and height $n^{-l_1(r^\theta)}$, which is smaller, but still larger than $r/n$ since $l_1(r^\theta) \leqslant l_2(r)$ by Lemma 8.3.4. These strips contain images of $F$ under $l_1(r^\theta)$-fold compositions of maps from $\{S_d\}_{d\in\mathcal{D}}$ and so to keep track of how many there are, we count rectangles in the appropriate columns for each map in the composition. The total number of such strips is therefore

$$\prod_{l=l_2(r^\theta)+1}^{l_1(r^\theta)} N_{i_l}.$$

We need to cover these strips by sets of diameter $r$ and so we iterate the construction inside each horizontal strip until the heights of the construction rectangles are $n^{-l_2(r)}$, which is approximately $r$. This takes $l_2(r) - l_1(r^\theta)$ iterations and, this time, for every iteration we pick up $N$ smaller rectangles (one for every map in the IFS), rather than just those inside a particular column. We are left with a collection of rectangles with height approximately $r$ and width $m^{-l_2(r)}$, which is larger than $r$. We continue iterating inside each of these rectangles until the widths of the construction rectangles are $m^{-l_1(r)}$, which is roughly $r$. This takes a further $l_1(r) - l_2(r)$ iterations. We can then cover the resulting collection of rectangles by small sets of diameter $r$, observing that sets formed by this last stage of iteration can be efficiently covered simultaneously, provided they are in the same column. Therefore, for each of the last iterations we only require a factor of $N_0$ more covering sets (one for each non-empty column).

Putting all these estimates together, yields a good estimate for

$$N_r\big(\Pi\big(Q(\mathbf{d}, r^\theta)\big)\big)$$

which is accurate up to constants which we can safely ignore. More precisely, there is a constant $C > 0$, which may depend on $\theta$ but not on

**d** or $r$, such that

$$
\begin{aligned}
&N_r\big(\Pi\big(Q(\mathbf{d}, r^\theta)\big)\big) \\
&\leqslant C \left( \prod_{l=l_2(r^\theta)+1}^{l_1(r^\theta)} N_{i_l} \right) \left( N^{l_2(r)-l_1(r^\theta)} \right) \left( N_0^{l_1(r)-l_2(r)} \right) \qquad (8.5) \\
&\leqslant C \left( N_{\max} \right)^{l_1(r^\theta)-l_2(r^\theta)} \left( N^{l_2(r)-l_1(r^\theta)} \right) \left( N_0^{l_1(r)-l_2(r)} \right) \\
&\leqslant C N^3 \left(N_{\max}\right)^{\log r^\theta/\log n - \log r^\theta/\log m} \\
&\qquad \cdot N^{\log r^\theta/\log m - \log r/\log n} \; N_0^{\log r/\log n - \log r/\log m} \\
&= C N^3 r^{\theta \log N_{\max}/\log n + \theta \log(N/N_{\max})/\log m + \log(N_0/N)/\log n \, - \, \log N_0/\log m}.
\end{aligned}
$$

This proves

$$
\dim_{\mathrm{A}}^\theta F \leqslant \frac{\left(\frac{\log N_0}{\log m} + \frac{\log(N/N_0)}{\log n}\right) - \theta\left(\frac{\log(N/N_{\max})}{\log m} + \frac{\log N_{\max}}{\log n}\right)}{1-\theta},
$$

as required. The lower bound follows similarly by choosing $\mathbf{d} \in \mathcal{D}^\infty$, such that $N_{i_l} = N_{\max}$ for all $l$. For such **d** the upper bound established above is also a lower bound, up to a constant.

Our work so far has shown that for $\theta = \log m/\log n$

$$
\dim_{\mathrm{A}}^\theta F = \frac{\log N_0}{\log m} + \frac{\log N_{\max}}{\log n}
$$

which is the Assouad dimension of $F$. It therefore follows from Corollary 3.3.3 and Theorem 8.3.1 and that for $\theta \in (\log m/\log n, 1)$

$$
\dim_{\mathrm{A}}^\theta F = \dim_{\mathrm{qA}} F = \dim_{\mathrm{A}} F.
$$

In fact the above argument requires only a minor modification to prove Mackay's result (see Theorem 8.3.1)

$$
\dim_{\mathrm{A}} F = \frac{\log N_0}{\log m} + \frac{\log N_{\max}}{\log n} \qquad (8.6)
$$

which would complete the proof of Theorem 8.3.3 without relying on *a priori* information about the Assouad dimension. Since $\dim_{\mathrm{A}} F \geqslant \dim_{\mathrm{A}}^\theta F$, all that remains is to prove the reverse inequality. Consider $\mathbf{d} \in \mathcal{D}^\infty$, chosen such that $N_{i_l} = N_{\max}$ for all $l$, and $0 < r < R < 1$. If $r \leqslant R^{1/\theta}$, then the argument above already proves

$$N_r\big(\Pi\big(Q(\mathbf{d}, R)\big)\big) \leqslant c \left(\frac{R}{r}\right)^{\frac{\log N_0}{\log m} + \frac{\log N_{\max}}{\log n}}$$

for a constant $c > 0$ and, therefore, we may assume $r > R^{1/\theta}$. Following the 'decomposing and covering' approach from above, with $R$ replacing $r^\theta$, we obtain a family of horizontal strips with height approximately $r$ after $l_2(r)$ steps. This happens before the widths become smaller than $R$, since in this case $l_2(r) \leqslant l_1(R)$, see Lemma 8.3.4. The effect is that the 'middle term' in (8.5) (concerning powers of $N$) disappears. This gives, for a constant $c' > 0$,

$$\begin{aligned} &N_r\big(\Pi\big(Q(\mathbf{d}, R)\big)\big) \\ &\leqslant c' \left( \prod_{l=l_2(R)+1}^{l_2(r)} N_{i_l} \right) \left( N_0^{l_1(r)-l_1(R)} \right) \\ &= c' \left( N_{\max} \right)^{l_2(r)-l_2(R)} \quad N_0^{l_1(r)-l_1(R)} \\ &\leqslant c' N^{-2} \left(N_{\max}\right)^{\log R/\log n - \log r/\log n} \quad N_0^{\log R/\log m - \log r/\log m} \\ &= c' N^{-2} \left(\frac{R}{r}\right)^{\frac{\log N_0}{\log m} + \frac{\log N_{\max}}{\log n}} \end{aligned}$$

completing the proof of (8.6) and Theorem 8.3.3. □

We can deduce the quasi-Assouad and quasi-lower dimensions immediately from Theorem 8.3.3.

**Corollary 8.3.5** Let $F$ be a Bedford-McMullen carpet. Then

$$\dim_{\mathrm{qA}} F = \dim_{\mathrm{A}} F = \frac{\log N_0}{\log m} + \frac{\log N_{\max}}{\log n}$$

and

$$\dim_{\mathrm{qL}} F = \dim_{\mathrm{L}} F = \frac{\log N_0}{\log m} + \frac{\log N_{\min}}{\log n}.$$

We can also compute the modified lower dimension from Theorem

8.3.3, although it requires one more ingredient. Surprisingly it turns out to be equal to the *Hausdorff* dimension, not the lower dimension.

**Corollary 8.3.6** Let $F$ be a Bedford-McMullen carpet. Then

$$\dim_{\mathrm{ML}} F = \dim_{\mathrm{H}} F = \frac{\log \sum_{i=1}^{N_0} N_i^{\log m/\log n}}{\log m}.$$

Corollary 8.3.6 follows immediately from Theorem 8.3.3 and the following lemma due to Ferguson, Jordan and Shmerkin [84, Lemma 4.3]. This shows how to construct subsets of carpets which have lower dimension arbitrarily close to the Hausdorff dimension of the original carpet. This proves the lower bound and, since the modified lower dimension of a compact set is always bounded above by the Hausdorff dimension by Theorem 3.4.3, the result follows. We include the proof of the result of Ferguson, Jordan and Shmerkin, which introduced a useful technique based on 'approximating systems from within' by choosing a subsystem according to digit frequencies and then applying Stirling's formula.

**Lemma 8.3.7** Given a Bedford-McMullen carpet $F$ and $\varepsilon > 0$, there exists a subset $E \subseteq F$ that is also a Bedford-McMullen carpet, has uniform fibres, and satisfies

$$\dim_{\mathrm{H}} E = \dim_{\mathrm{L}} E \geqslant \dim_{\mathrm{H}} F - \varepsilon.$$

*Proof* Let $s = \dim_{\mathrm{H}} F$ and, for $d = (i, j) \in \mathcal{D}$, let

$$p_d = N_i^{(\log m/\log n)-1}/m^s$$

noting that $\sum_{d\in\mathcal{D}} p_d = 1$. The probability weights $p_d$ are those defining the *McMullen measure*, that is, the unique self-affine measure of maximal Hausdorff dimension. For $k \geqslant 1$ a large integer, let

$$l(k) = \sum_{d\in\mathcal{D}} \lfloor kp_d \rfloor \tag{8.7}$$

and

$$\mathcal{D}(k) = \left\{ \mathbf{d} = (d_1, \dots, d_{l(k)}) \in \mathcal{D}^{l(k)} \ : \ \text{for all } d \in \mathcal{D},\ \#\{t : d_t = d\} = \lfloor kp_d \rfloor \right\}.$$

In particular, $\mathcal{D}(k)$ selects the subset of $l(k)$ level words whose digits appear with the 'correct' frequency from the point of view of Hausdorff dimension.

Consider the IFS $\{S_{\mathbf{d}} : \mathbf{d} \in \mathcal{D}(k)\}$ and denote its attractor by $E(k)$. Clearly $E(k)$ is a Bedford-McMullen carpet and, since each element of $\mathcal{D}(k)$ is composed of the same collection of maps (just in a different order), $E$ also has uniform fibres. It follows from Theorem 8.3.1 that

$$\dim_{\mathrm{L}} E(k) = \dim_{\mathrm{H}} E(k) = \frac{\log N_0(k)}{\log m(k)} + \frac{\log(N(k)/N_0(k))}{\log n(k)}$$

where $m(k) = m^{l(k)}$ and $n(k) = n^{l(k)}$ are the integer parameters defining the grid associated with $E(k)$ and $N_0(k)$ and $N(k)$ are the number of non-empty columns and the number of rectangles in the construction of $E(k)$ respectively. These latter parameters can be computed explicitly by a simple counting argument and we find

$$N(k) = \frac{l(k)!}{\prod_{d \in \mathcal{D}} \lfloor k p_d \rfloor !}$$

and

$$N_0(k) = \frac{l(k)!}{\prod_{i=1}^{N_0} \left( \sum_{j=1}^{N_i} \lfloor k p_{(i,j)} \rfloor \right)!}.$$

We can therefore estimate the dimension of $E(k)$ from below using Stirling's approximation

$$a \log a - a \ \leqslant \ \log a! \ \leqslant \ a \log a - a + \log a \qquad (a \geqslant 7) \tag{8.8}$$

and (8.7) to get

$$\begin{aligned} \dim_{\mathrm{L}} E(k) &= \frac{\log l(k)! - \sum_{i=1}^{N_0} \log \left( \sum_{j=1}^{N_i} \lfloor k p_{(i,j)} \rfloor \right)!}{l(k) \log m} \\ &\qquad + \frac{\sum_{i=1}^{N_0} \log \left( \sum_{j=1}^{N_i} \lfloor k p_{(i,j)} \rfloor \right)! - \sum_{d \in \mathcal{D}} \log \lfloor k p_d \rfloor !}{l(k) \log n} \\ &\geqslant s - \varepsilon(k) \end{aligned}$$

where $\varepsilon(k) \to 0$ as $k \to \infty$. This proves the theorem, since we can choose $k$ arbitrarily large. □

## 8.4 Self-affine sets with a comb structure

A family of planar non-carpet self-affine sets was considered by Bárány, Käenmäki and Rossi [21]. The phenomenon that the Assouad dimension is given by the dimension of a projection plus the largest dimension of an orthogonal slice is again manifest although, as one moves away from carpets, it becomes less clear which projections to consider and the structure of both the projection and maximal slice are more complicated. The paper deals with the extra complexity of the projections by assuming that all of the projections are intervals. The question then becomes: over which collection of slices should the dimension be maximised? This is a subtle problem in general. Notice that in the case of Bedford-McMullen carpets one always considers slices in the direction of the strongest contraction, which is always the vertical direction. When the matrices associated with the defining maps are non-diagonal matrices, we must consider the slices over a larger family of directions, but the philosophy is the same.

An IFS of affine maps acting on the plane is said to be *reducible* if there is a line which is invariant under the defining matrices (for example, Bedford-McMullen carpets are generated by reducible IFSs, since the matrices preserve both the horizontal and vertical axes), *irreducible* if there is no such invariant line, and *strongly irreducible* if there is no finite union of lines which is invariant under the defining matrices. Moreover, it is said to be *proximal* if the semigroup generated by the defining matrices contains an element which has two real and distinct eigenvalues.

For example, if the defining matrices are all upper-triangular, then the IFS is reducible. Moreover, a (proximal) IFS defined by a collection of diagonal and anti-diagonal matrices is irreducible but not strongly irreducible since the union of the vertical and horizontal axes is invariant but neither the vertical nor horizontal axes are invariant on their own. In fact, the case where there are both diagonal and anti-diagonal matrices is the only one which is irreducible but not strongly irreducible, see [21, Lemma 2.2]. The proximal condition should be interpreted as guaranteeing that the attractor of the IFS is 'genuinely self-affine'. For example, it rules out self-similar sets. We write $G(1,2)$ for the collection of 1-dimensional subspaces of $\mathbb{R}^2$ and $\pi F$ for the projection of $F$ onto $\pi \in G(1,2)$. The space $G(1,2)$ is often referred to as the Grassmannian manifold, see Chapter 10.2. We also write $\pi^\perp \in G(1,2)$ for the 1-dimensional subspace perpendicular to $\pi$.

The following theorem was proved in [21, Theorem 3.2], although their result is more general than we state.

**Theorem 8.4.1** Let $F \subseteq \mathbb{R}^2$ be a self-affine set satisfying the SSC and such that $\pi F$ is an interval for all $\pi \in G(1,2)$. Furthermore, assume that the defining IFS is either: proximal and reducible; or strongly irreducible and all the defining matrices have strictly positive entries. Then there exists $\pi \in G(1,2)$ such that

$$\dim_A F = 1 + \max_{x \in F} \dim_H \big( (\pi^\perp + x) \cap F \big).$$

The projection $\pi$ appearing in Theorem 8.4.1 is made more explicit in [21], but we omit discussion of the details. It turns out that one needs to consider all *Furstenberg projections* and then choose $\pi$ such that $\max_{x \in F} \dim_H (\pi^\perp + x) \cap F$ is maximised.

The proof of Theorem 8.4.1 relies on many more recent advances in the theory of self-affine sets. The tangent structure of the self-affine sets considered is crucial and it turns out that typically the weak tangents of self-affine sets have a 'fibred structure' given by the product of an interval and a Cantor set. (Here we restrict ourselves to weak tangents associated to sequences of similarities whose similarity ratios tend to infinity, that is, we genuinely 'zoom-in'.) The 'interval' here comes from the projection condition and the Cantor set comes from a slice of $F$. The tangent structure of self-affine sets has been developed in, for example, [14, 151]. Other more recent developments connecting the ergodic theory of the action of the collection of matrices in the defining IFS on projective space also play a vital role, see [17, 18, 19, 79].

One simple and well-known family of sets to which Theorem 8.4.1 applies is the family considered by Przytycki-Urbański [237]. These self-affine carpets fall into the more general family studied by Feng and Wang [83] and are related to Bernoulli convolutions. Fix $\lambda \in (1/2, 1)$ and consider the IFS consisting of $S_1$ and $S_2$ acting on $[0,1]^2$ defined by

$$S_1(x,y) = (\lambda x, y/2) \qquad \text{and} \qquad S_2(x,y) = (\lambda x + (1-\lambda), (y+1)/2)$$

writing $F$ for the attractor. This IFS is clearly proximal and reducible and so the formula from Theorem 8.4.1 applies and, moreover, there is only one 'Furstenberg direction' and so the family of slices one has to consider are simply the vertical slices, yielding

$$\dim_A F = 1 + \max_{x \in [0,1]} \dim_H V_x \cap F$$

where $V_x = \{(x,y) : y \in \mathbb{R}\}$. This family of self-affine sets was also studied by Fraser and Jordan [101] and the Assouad dimension was computed via a different approach. There it was proved, see [101, Corollary 2.3], that

$$\dim_A F = 2 + s \frac{\log \lambda}{\log 2}$$

where $s$ is the infimum over all local dimensions of the associated Bernoulli convolution. Combining the two approaches, one can relate the Hausdorff dimension of slices to the local dimensions of Bernoulli convolutions, two quantities which *a priori* seem unrelated to the Assouad dimension.

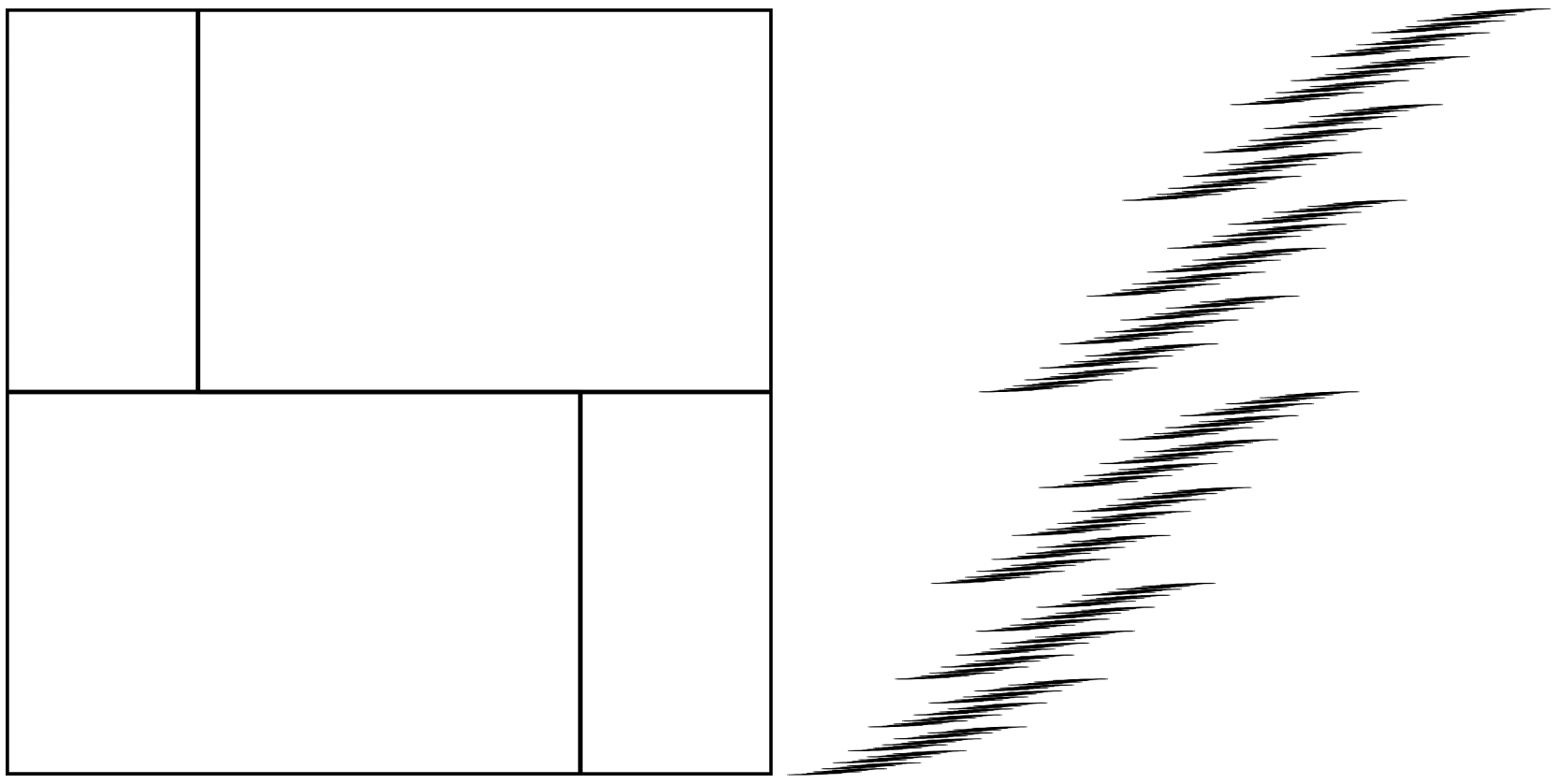

Figure 8.5 Rectangles depicting the IFS (left) and the associated self-affine Przytycki-Urbański set (right) for $\lambda = 3/4$.

## 8.5 A family of worked examples

In this section we present a simple family of examples, based on an example given in [97]. This family witnesses a certain 'discontinuity' in the Assouad and lower dimensions, not exhibited for the Hausdorff and box dimensions. It also serves to show that the bounds from (6.3) do not hold for the Assouad and lower dimension or spectra.

Fix $\lambda \in (0, 1/3]$ and consider the IFS $\{S_1, S_2, S_3\}$ acting on the square where

$$S_1(x, y) = (x/3, \lambda y),$$

$$S_2(x, y) = (x/3 + 2/3, \lambda y),$$

and

$$S_3(x, y) = (x/3, \lambda y + 1 - \lambda).$$

For $\lambda \in (0, 1/3)$, the attractor of this IFS is a self-affine carpet of the type considered by Lalley and Gatzouras [122] and for $\lambda = 1/3$ it is a self-similar set, sometimes known as the right-angled 1-dimensional

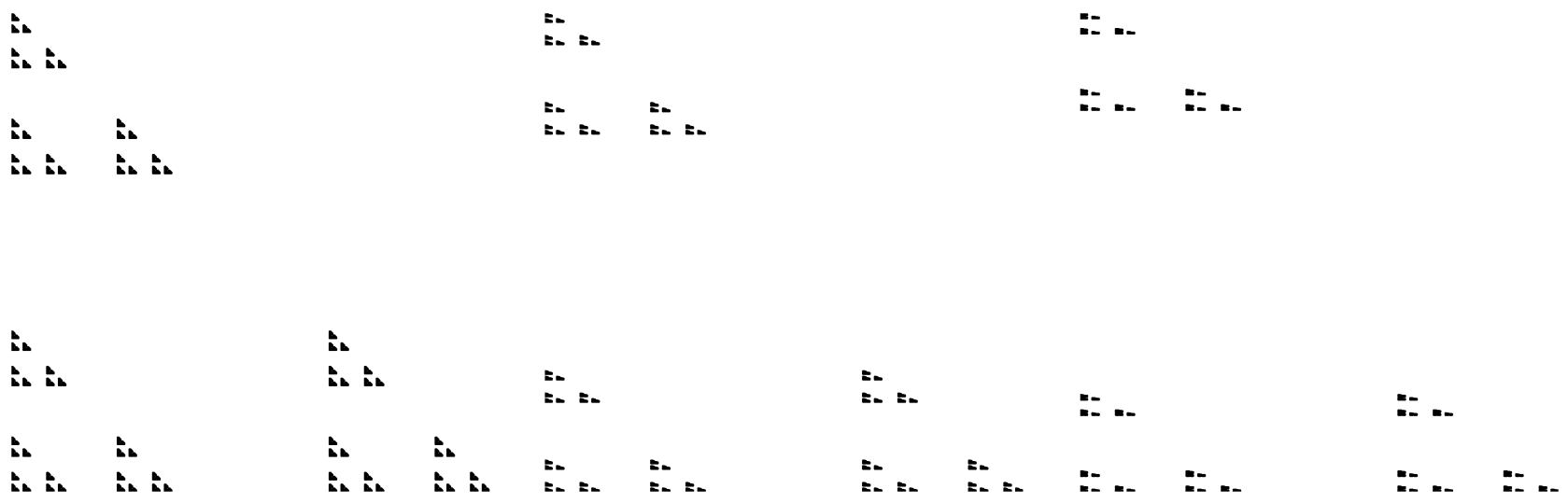

Figure 8.6 The attractor $F_\lambda$ for $\lambda$ equal to 1/3, 1/4 and 1/5, respectively.

Sierpiński triangle. We denote the attractor by $F_\lambda$. The dimension theory of self-affine sets in the Lalley-Gatzouras family is similar to the Bedford-McMullen family. For $\lambda \in (0, 1/3)$,

$$\dim_L F_\lambda = \frac{\log 2}{\log 3}, \qquad \dim_H F_\lambda = \frac{\log\left(2^{-\log 3/\log\lambda} + 1\right)}{\log 3}$$

$$\dim_B F_\lambda = \frac{\log 2}{\log 3} + \frac{\log(3/2)}{-\log\lambda}, \qquad \dim_A F_\lambda = \frac{\log 2}{\log 3} + \frac{\log 2}{-\log\lambda}$$

and for $\lambda = 1/3$

$$\dim_A F_\lambda = \dim_L F_\lambda = \dim_H F_\lambda = \dim_B F_\lambda = 1.$$

The function $\lambda \mapsto \dim F_\lambda$ is continuous for the box and Hausdorff dimensions but discontinuous at $\lambda = 1/3$ for Assouad and lower dimensions.

The situation for the Assouad and lower spectra is more subtle. Since each map in the defining IFS has the same linear part, a simple adaptation of the proof of Theorem 8.3.3 shows that for $\lambda \in (0, 1/3)$ and $0 < \theta \leqslant \frac{\log 3}{-\log\lambda}$ we have

$$\dim_A^\theta F_\lambda = \frac{\frac{\log 2}{\log 3} + \frac{\log(3/2)}{-\log\lambda} - \theta\left(\frac{\log(3/2)}{\log 3} + \frac{\log 2}{-\log\lambda}\right)}{1-\theta}$$

and

$$\dim_L^\theta F_\lambda = \frac{\frac{\log 2}{\log 3} + \frac{\log(3/2)}{-\log\lambda} - \theta}{1-\theta}$$

and, for $\theta > \frac{\log 3}{-\log\lambda}$, the Assouad and lower spectra are constantly equal to the Assouad and lower dimensions, respectively. See Theorem 8.3.3 for comparison. Therefore the discontinuity at $\lambda = 1/3$ is observed uniformly, but not point-wise. That is, for fixed $\theta$,

$$\dim_{\mathrm{A}}^{\theta} F_\lambda \to 1 = \dim_{\mathrm{A}}^{\theta} F_{1/3}$$

as $\lambda \to 1/3$, however,

$$\sup_{\theta \in (0,1)} \left| \dim_{\mathrm{A}}^{\theta} F_\lambda - \dim_{\mathrm{A}}^{\theta} F_{1/3} \right| \to 2\frac{\log 2}{\log 3} - 1 > 0$$

as $\lambda \to 1/3$. The situation for the lower spectrum is similar, and discontinuity is observed for the quasi-Assouad and quasi-lower dimensions.

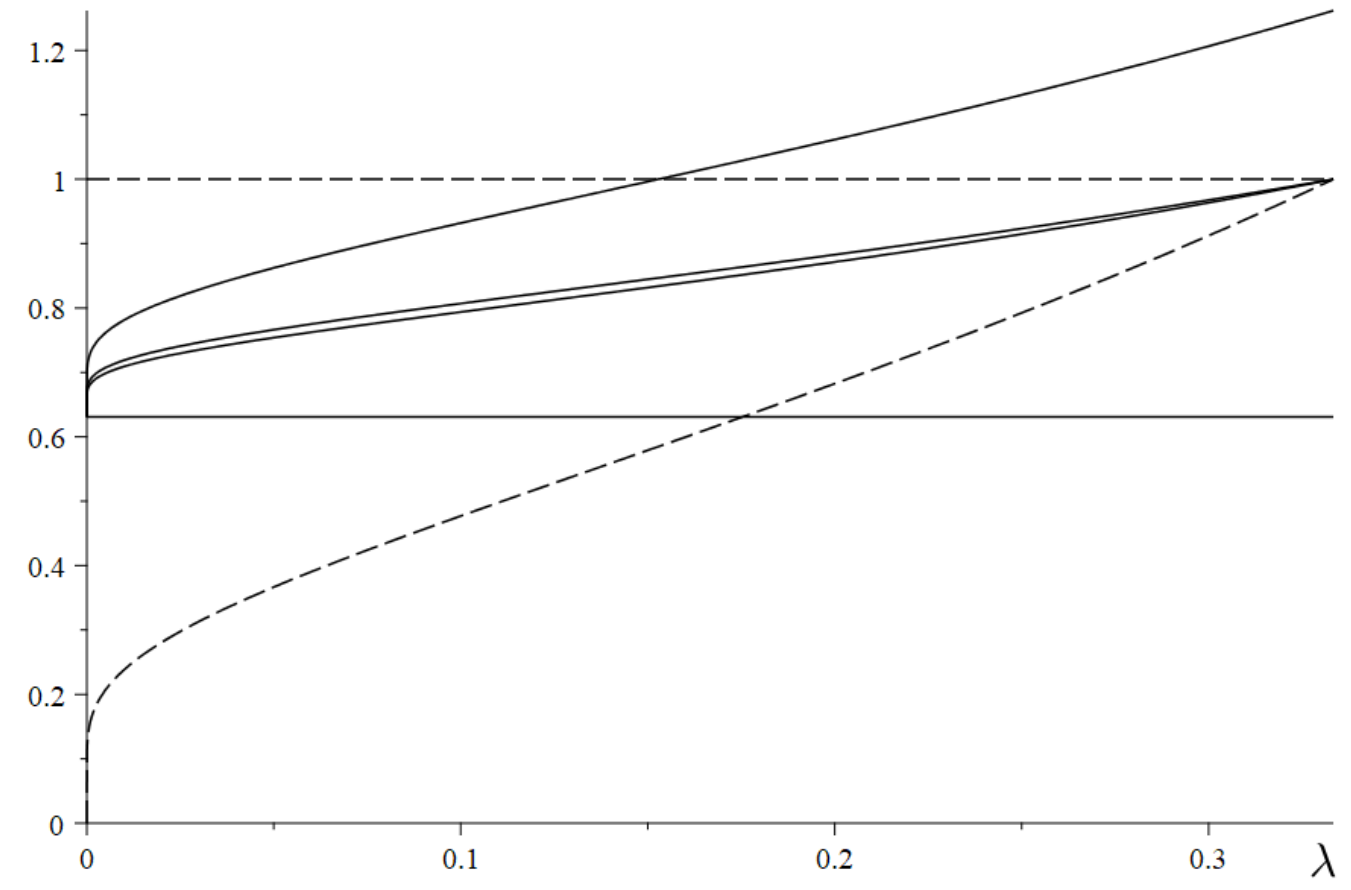

Figure 8.7 A plot of $\dim_{\mathrm{L}} F_\lambda$, $\dim_{\mathrm{H}} F_\lambda$, $\dim_{\mathrm{B}} F_\lambda$ and $\dim_{\mathrm{A}} F_\lambda$ for $\lambda \in (0, 1/3]$ with the upper and lower bounds from (6.3) on page 88 shown as dashed lines. Note that the Assouad and lower dimensions are outside the range given by these bounds for large enough $\lambda$.

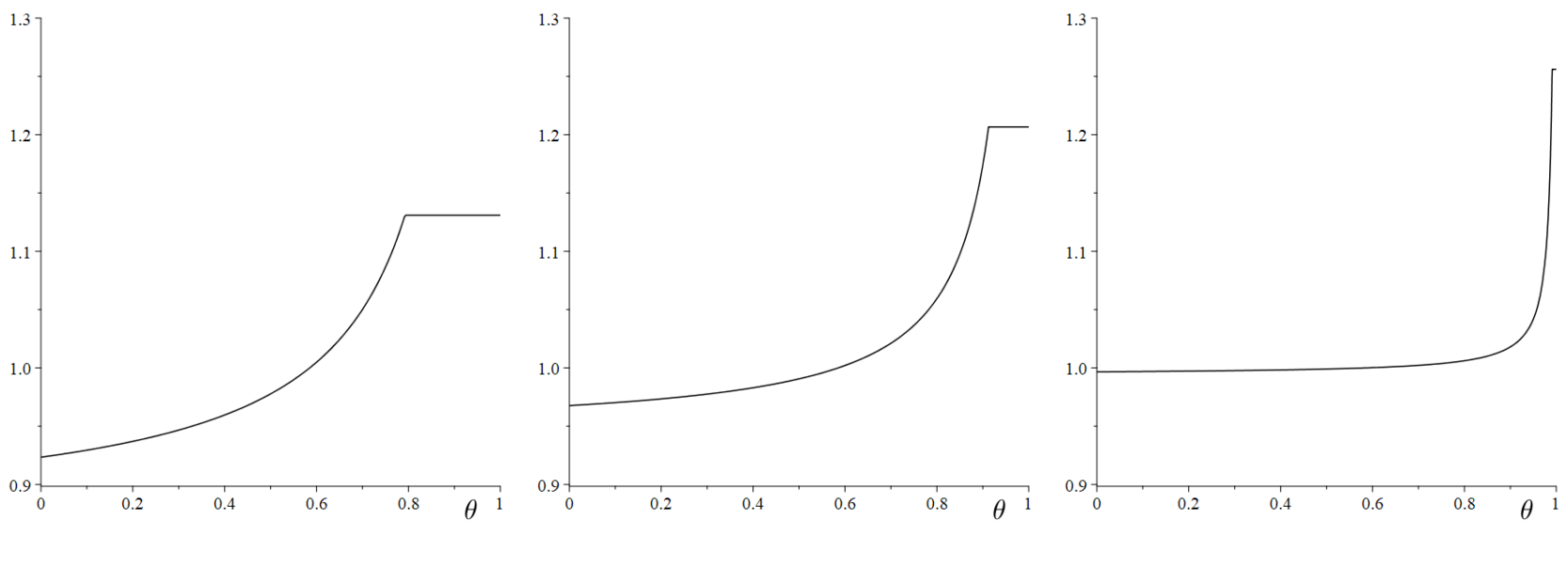

Figure 8.8 Plots of $\dim_{\mathrm{A}}^{\theta} F_\lambda$ for $\lambda = 1/4,\ 0.3,\ 0.33$, moving from left to right.

## 8.6 Dimensions of self-affine measures

For self-similar sets satisfying the OSC, the self-similar measure associated with Bernoulli weights $\{c_i^s\}$ (where $c_i$ is the $i$th similarity ratio and $s$ is the similarity dimension) is Ahlfors regular and therefore certainly doubling. In stark contrast to the self-similar case, there exist self-affine sets satisfying the OSC such that *all* associated self-affine measures fail to be doubling. This was first established by Li, Wei and Wen [190]. In fact, Bedford-McMullen carpets provide such examples. Let $m = 2$, $n = 4$ and choose three rectangles corresponding to the black rectangles in Figure 8.9. This carpet clearly satisfies the OSC.

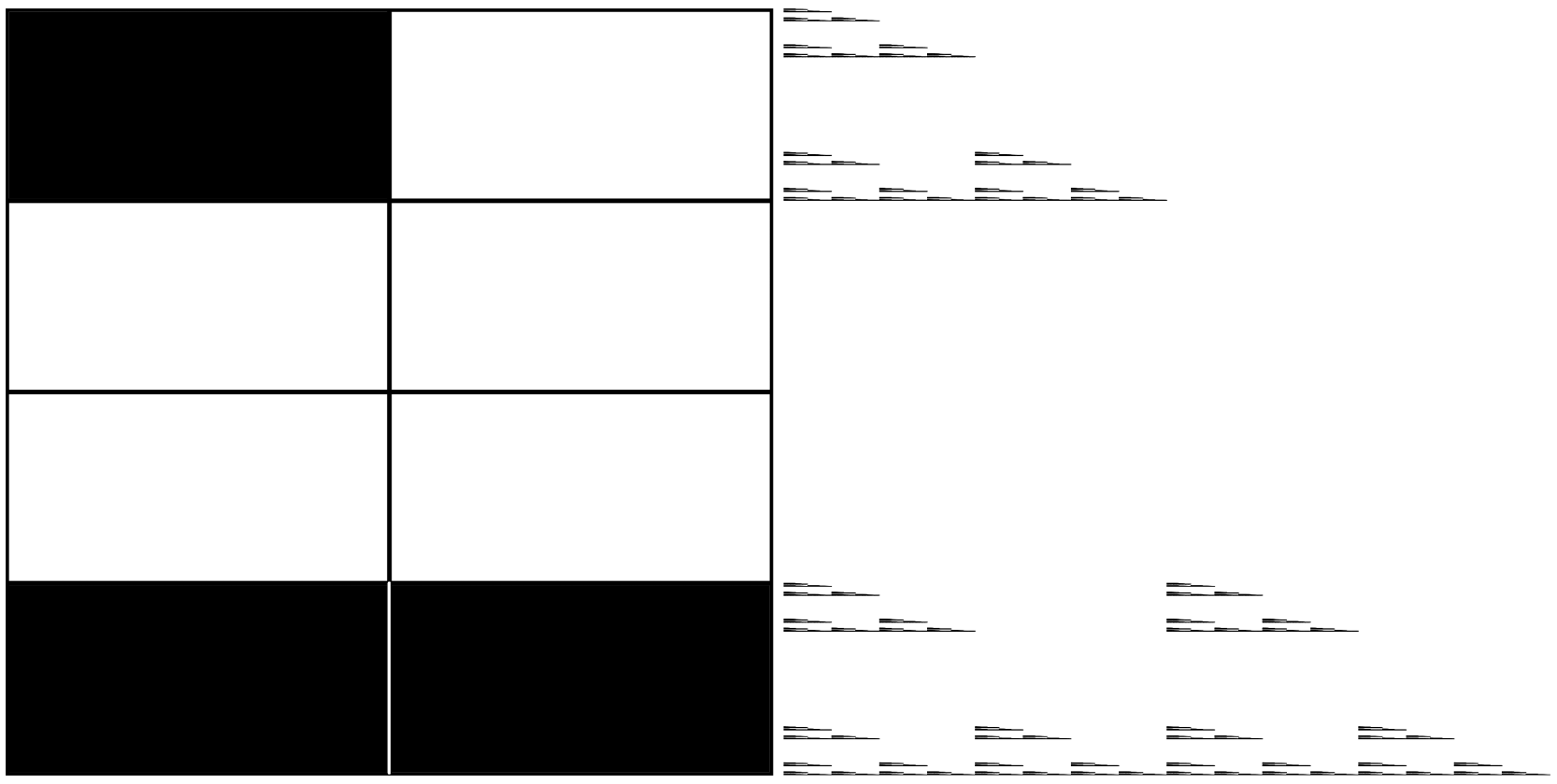

Figure 8.9 The defining pattern (left) and associated self-affine Bedford-McMullen carpet which satisfies the OSC but does not support a doubling self-affine measure.

Fix a strictly positive probability vector $(p_{(0,0)}, p_{(1,0)}, p_{(0,3)})$ and let $\mu$ be the associated self-affine measure.

Let $k \geqslant 1$ be a large integer, $R = 4^{-k}$ noting that $l_1(R) = 2k$ and $l_2(R) = k$, see (8.3)-(8.4). Firstly, consider $\mathbf{d}, \mathbf{d}' \in \mathcal{D}$ satisfying

$$\mathbf{d} = \left( \underbrace{(1,0), \dots, (1,0)}_{k \text{ times}}, (1,0), \underbrace{(0,0), \dots, (0,0)}_{k-1 \text{ times}}, \dots \right)$$

and

$$\mathbf{d}' = \left( \underbrace{(1,0), \dots, (1,0)}_{k \text{ times}}, (0,0), \underbrace{(1,0), \dots, (1,0)}_{k-1 \text{ times}}, \dots \right).$$

Crucially, the approximate squares $\Pi(Q(\mathbf{d},R))$ and $\Pi(Q(\mathbf{d}',R))$ are adjacent in that they share a common edge. Moreover,

$$\frac{\mu(\Pi(Q(\mathbf{d},R)))}{\mu(\Pi(Q(\mathbf{d}',R)))} = \frac{p_{(1,0)}^{k+1}(p_{(0,0)}+p_{(0,3)})^{k-1}}{p_{(1,0)}^{2k-1}(p_{(0,0)}+p_{(0,3)})} = \left(\frac{p_{(0,0)}+p_{(0,3)}}{p_{(1,0)}}\right)^{k-2}.$$

If $\mu$ is doubling, this ratio must remain bounded away from zero and infinity, which is only possible if

$$p_{(1,0)} = p_{(0,0)} + p_{(0,3)}, \tag{8.9}$$

which we assume from now on. Secondly, consider $\mathbf{d},\mathbf{d}' \in \mathcal{D}$ satisfying

$$\mathbf{d} = \left((1,0), \underbrace{(0,0),\ldots,(0,0)}_{2k-1 \text{ times}}, \ldots\right)$$

and

$$\mathbf{d}' = \left((0,0), \underbrace{(1,0),\ldots,(1,0)}_{2k-1 \text{ times}}, \ldots\right).$$

Again, the approximate squares $\Pi(Q(\mathbf{d},R))$ and $\Pi(Q(\mathbf{d}',R))$ are adjacent and this time, using (8.9),

$$\frac{\mu(\Pi(Q(\mathbf{d},R)))}{\mu(\Pi(Q(\mathbf{d}',R)))} = \frac{p_{(1,0)}(p_{(0,0)}+p_{(0,3)})^{2k-1}\left(\frac{p_{(0,0)}}{p_{(0,0)}+p_{(0,3)}}\right)^{k-1}}{(p_{(0,0)}+p_{(0,3)})p_{(1,0)}^{2k-1}\left(\frac{p_{(0,0)}}{p_{(0,0)}+p_{(0,3)}}\right)}$$
$$= \left(\frac{p_{(0,0)}}{p_{(0,0)}+p_{(0,3)}}\right)^{k-2}.$$

If $\mu$ is doubling, this ratio must also remain bounded away from zero and infinity, which is only possible if $p_{(0,0)} = p_{(0,0)} + p_{(0,3)}$. However, this is not possible since $p_{(0,3)} > 0$ and therefore $\mu$ cannot be doubling.

This time in parallel with the self-similar setting, the lower dimension of self-affine measures is rather better behaved than the Assouad dimension (recall Theorem 7.4.2). For example, [127, Theorem 3.9] shows that any non-trivial self-affine measure has positive lower dimension.

**Theorem 8.6.1** Let $\mu$ be a self-affine measure which is not an atom. Then

$$\dim_{\mathrm{L}} \mu > 0.$$

In order to guarantee that self-affine measures supported on Bedford-McMullen carpets are doubling, one needs to impose a stronger separation condition than the SSC. This condition, introduced by King [170], is known as the *very strong separation condition* (VSSC) and stipulates than the SSC is satisfied and, in addition, rectangles cannot be chosen from adjacent columns. Notably, this rules out the carpet from Figure 8.9. With this condition in place, formulae for the Assouad, lower and box dimensions of self-affine measures on carpets can be derived. The Assouad dimension was considered in [98, Theorem 2.6] and the lower dimension in [127, Theorem 3.11]. The box dimensions are derivable from [218, Proposition 3.3.1] (see the discussion in [98]) although they were not referred to as box dimensions at this time. The results in [98] and [218] also apply to Bedford-McMullen sponges, the higher dimensional analogue of Bedford-McMullen carpets.

In order to state the results, it is convenient to write, for $d = (i, j) \in \mathcal{D}$,

$$P(d) = \sum_{d'=(i,j')\in\mathcal{D}} p_{d'}$$

which is the sum of the probabilities associated to rectangles in the same column as $d = (i, j)$, that is, with $i$ fixed and $j'$ varying.

**Theorem 8.6.2** Let $\mu$ be a self-affine measure supported on a Bedford-McMullen carpet satisfying the VSSC. Then

$$\dim_A \mu = \max_{d\in\mathcal{D}} \frac{\log(1/P(d))}{\log m} + \max_{d\in\mathcal{D}} \frac{\log(P(d)/p_d)}{\log n},$$

$$\dim_L \mu = \min_{d\in\mathcal{D}} \frac{\log(1/P(d))}{\log m} + \min_{d\in\mathcal{D}} \frac{\log(P(d)/p_d)}{\log n}$$

and

$$\dim_B \mu = \max_{d\in\mathcal{D}} \log(1/P(d)) \left( \frac{1}{\log m} - \frac{1}{\log n} \right) + \max_{d\in\mathcal{D}} \frac{\log(1/p_d)}{\log n}.$$

Comparing this theorem to Theorem 8.3.1 we see that there always exists a self-affine measure which (simultaneously) realises the Assouad and lower dimensions of a Bedford-McMullen carpet.

**Corollary 8.6.3** Let $\mu$ be a self-affine measure supported on a Bedford-McMullen carpet $F$ with non-uniform fibres. Then

$$\dim_A F = \dim_A \mu$$

if and only if $P(d) = 1/N_0$ and $P(d)/p_d \leqslant N_{\max}$ for all $d \in \mathcal{D}$, and

$$\dim_{\mathrm{L}} F = \dim_{\mathrm{L}} \mu$$

if and only if $P(d) = 1/N_0$ and $P(d)/p_d \geqslant N_{\min}$ for all $d \in \mathcal{D}$. Moreover, there is always at least one self-affine measure $\mu$ satisfying

$$\dim_{\mathrm{A}} F = \dim_{\mathrm{A}} \mu > \dim_{\mathrm{L}} F = \dim_{\mathrm{L}} \mu.$$

The natural measure $\mu$ which simultaneously realises the Assouad and lower dimensions was introduced in [98] and called the *coordinate uniform measure* and can be constructed by first distributing weight equally among all used columns and then, within each column, distributing weight evenly between the rectangles in the column. There may be other measures which simultaneously realise the Assouad and lower dimensions, provided there exists a column with $N_{\min} < N_i < N_{\max}$. In this case there is room to distribute the weight not exactly equally between the rectangles in the intermediate column.

Consider the *simultaneous realisation problem*, discussed in Section 6.2, which asks if there exists a measure which simultaneously realises different notions of dimension (in situations when they are distinct). In the positive direction, the coordinate uniform measure solves this problem for Bedford-McMullen carpets with non-uniform fibres and satisfying the VSSC in the setting of Assouad and lower dimension. However, if we throw the Hausdorff dimension into the mix, then the situation is not as straightforward. In fact there does not exist a self-affine measure which simultaneously realises all three notions of dimension. In the case where the carpet does not have uniform fibres, the Assouad, Hausdorff, and lower dimensions of Bedford-McMullen carpets are necessarily distinct by Theorem 8.3.1. Moreover, it is well-known that there exists a unique invariant probability measure of maximal Hausdorff dimension, namely, the *McMullen measure*. The McMullen measure is the self-affine measure associated with

$$p_d = N_i^{(\log m/\log n)-1}/m^{\dim_{\mathrm{H}} F}$$

for $d = (i, j) \in \mathcal{D}$. Therefore, $P(d) = N_i^{\log m/\log n}/m^{\dim_{\mathrm{H}} F}$ varies with $i$ and it follows from Corollary 8.6.3 that the Assouad and lower dimensions of the McMullen measure are always distinct from the Assouad and lower dimensions of the carpet, provided the VSSC is satisfied and the carpet does not have uniform fibres.

In general the Hausdorff dimension of a self-affine measure on a Bedford-

McMullen carpet is given by

$$\dim_{\mathrm{H}} \mu = \frac{h(\mu)}{\log n} + \frac{\log n - \log m}{\log n} \dim_{\mathrm{H}} \pi\mu \tag{8.10}$$

where

$$h(\mu) = -\sum_{d \in \mathcal{D}} p_d \log p_d$$

is the entropy of $\mu$ and $\dim_{\mathrm{H}} \pi\mu$ is the Hausdorff dimension of the projection of $\mu$ onto the first coordinate. This projection is a self-similar measure satisfying the OSC and therefore the dimension can be computed easily by applying (7.12). Here $\log n$ and $\log m$ are the Lyapunov exponents of $\mu$. The formula (8.10) is known as a *Ledrappier-Young formula*, following [183, 184]. It was proved by Bárány and Käenmäki [19] that every planar self-affine measure satisfies the appropriate analogue of the Ledrappier-Young formula, which essentially says the Hausdorff dimension can be expressed in terms of entropy, Lyapunov exponents (which are more complicated for general self-affine measures) and dimensions of projections. See [17, 82] for further work in this direction.

The situation for box dimension is even worse: there does not exist a self-affine measure which realises the box dimension. One might have guessed that the coordinate uniform measure would work due to the connection between the upper box and (quasi-)Assouad dimensions via the Assouad spectrum. However, for the coordinate uniform measure $\mu$

$$\dim_{\mathrm{B}} \mu = \frac{\log N_0}{\log m} + \frac{\log N_{\max}}{\log n} = \dim_{\mathrm{A}} \mu = \dim_{\mathrm{A}} F.$$

Another reasonable candidate is the uniform measure which gives weight $1/N$ to each symbol. However, for this measure

$$\dim_{\mathrm{B}} F < \dim_{\mathrm{B}} \mu = \frac{\log(N/N_{\min})}{\log m} + \frac{\log N_{\min}}{\log n}.$$

For any self-affine measure

$$\max_{d \in \mathcal{D}} 1/P(d) \geqslant N_0$$

and

$$\max_{d \in \mathcal{D}} 1/p_d \geqslant N$$

by the pigeonhole principle. If both of these lower bounds could be simultaneously achieved then the box dimension of $F$ would be realised. However, this is clearly impossible since one cannot uniformly distribute across the individual rectangles and across the columns simultaneously. Observe that the formula for $\dim_{\mathrm{B}} \mu$ is continuous in $\{p_d\}_{d \in \mathcal{D}}$ and $\dim_{\mathrm{B}} \mu \to \infty$ if $\min p_{d \in \mathcal{D}} \to 0$. These observations, combined with the fact that $\dim_{\mathrm{B}} F$ is not realised, demonstrate that in fact there is a constant $\delta > 0$ such that $\dim_{\mathrm{B}} \mu \geqslant \dim_{\mathrm{B}} F + \delta$ for all self-affine $\mu$.

### 8.6.1 Another worked example

Let $F$ be the Bedford-McMullen carpet with $m = 3$ and $n = 5$ generated by four maps, as shown in Figure 8.10, and $\mu_\varepsilon$ be a 1-parameter family of measures corresponding to the probabilities indicated in Figure 8.10, where $\varepsilon$ varies in $(0, 1)$. Although the collection of all self-affine measures on $F$ is a 3-parameter family, when trying to realise the dimensions of $F$ it is more efficient to give all rectangles in a given column equal weight. This restriction reduces the family of appropriate measures to the 1-parameter family $\mu_\varepsilon$. The VSSC is satisfied for this IFS.

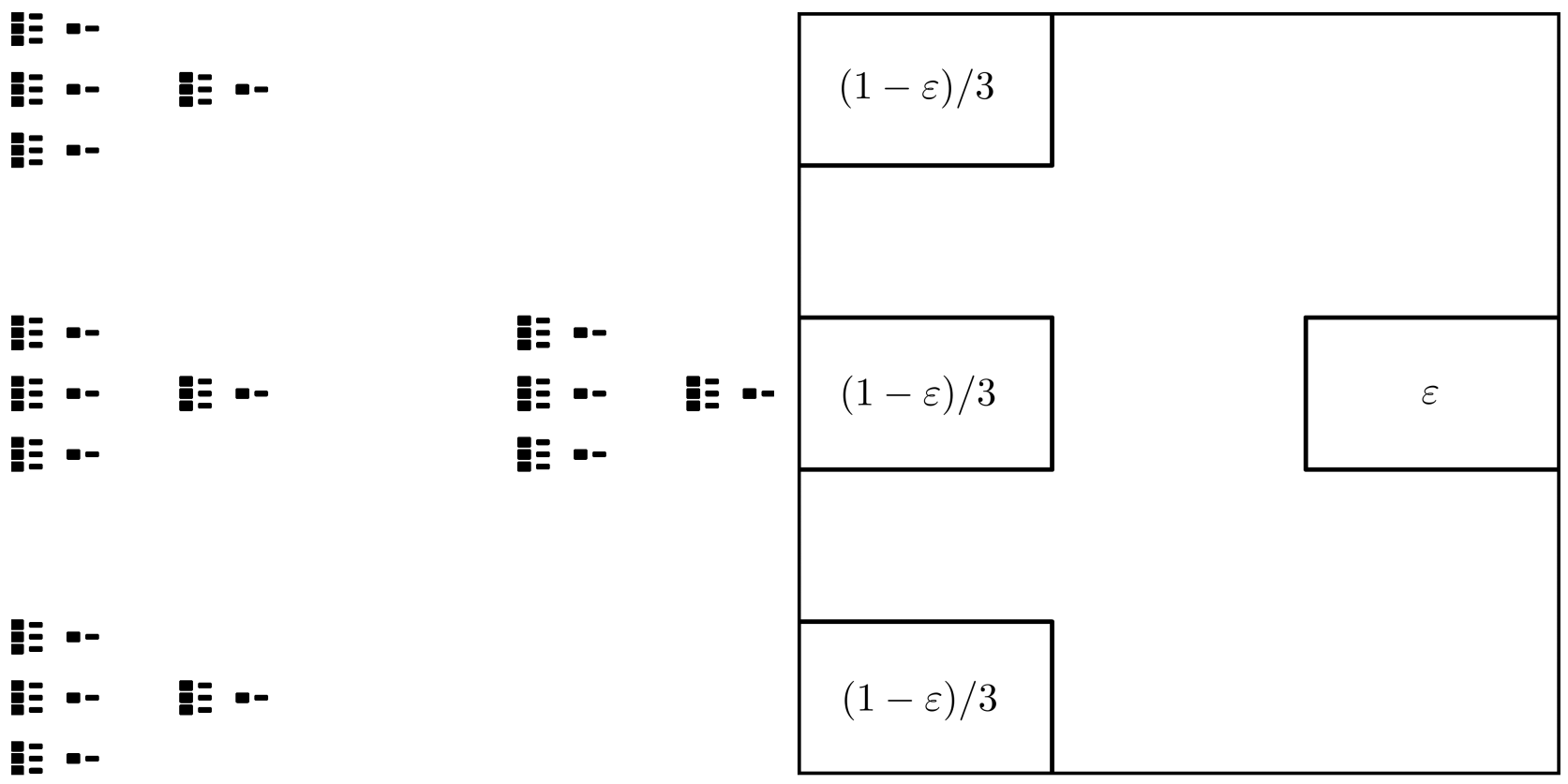

Figure 8.10 The Bedford-McMullen carpet (left) and the associated parameterised probabilities defining a self-affine measure supported on the carpet (right).

Applying the formulae from Theorem 8.6.2 and (8.10), we get

$$\dim_{\mathrm{L}} \mu_\varepsilon = \frac{\log\left(\min\left\{\frac{1}{\varepsilon}, \frac{1}{1-\varepsilon}\right\}\right)}{\log 3},$$

$$\dim_{\mathrm{H}} \mu_\varepsilon = \frac{h(\mu_\varepsilon)}{\log 5} - \frac{\log 5 - \log 3}{\log 5} \dim_{\mathrm{H}} \pi\mu_\varepsilon,$$

$$\dim_{\mathrm{B}} \mu_\varepsilon = \log\left(\max\left\{\frac{1}{\varepsilon}, \frac{1}{1-\varepsilon}\right\}\right)\left(\frac{1}{\log 3} - \frac{1}{\log 5}\right) + \frac{\log\left(\max\left\{\frac{1}{\varepsilon}, \frac{3}{1-\varepsilon}\right\}\right)}{\log 5}$$

and

$$\dim_{\mathrm{A}} \mu_\varepsilon = \frac{\log\left(\max\left\{\frac{1}{\varepsilon}, \frac{1}{1-\varepsilon}\right\}\right)}{\log 3} + \frac{\log 3}{\log 5},$$

see Figure 8.11, and also

$$\dim_{\mathrm{L}} F = \frac{\log 2}{\log 3} \approx 0.6309,$$

$$\dim_{\mathrm{H}} F = \frac{\log\left(1 + 3^{\frac{\log 3}{\log 5}}\right)}{\log 3} \approx 1.0347,$$

$$\dim_{\mathrm{B}} F = \frac{\log 2}{\log 3} + \frac{\log 2}{\log 5} \approx 1.0616$$

and

$$\dim_{\mathrm{A}} F = \frac{\log 2}{\log 3} + \frac{\log 3}{\log 5} \approx 1.3135.$$

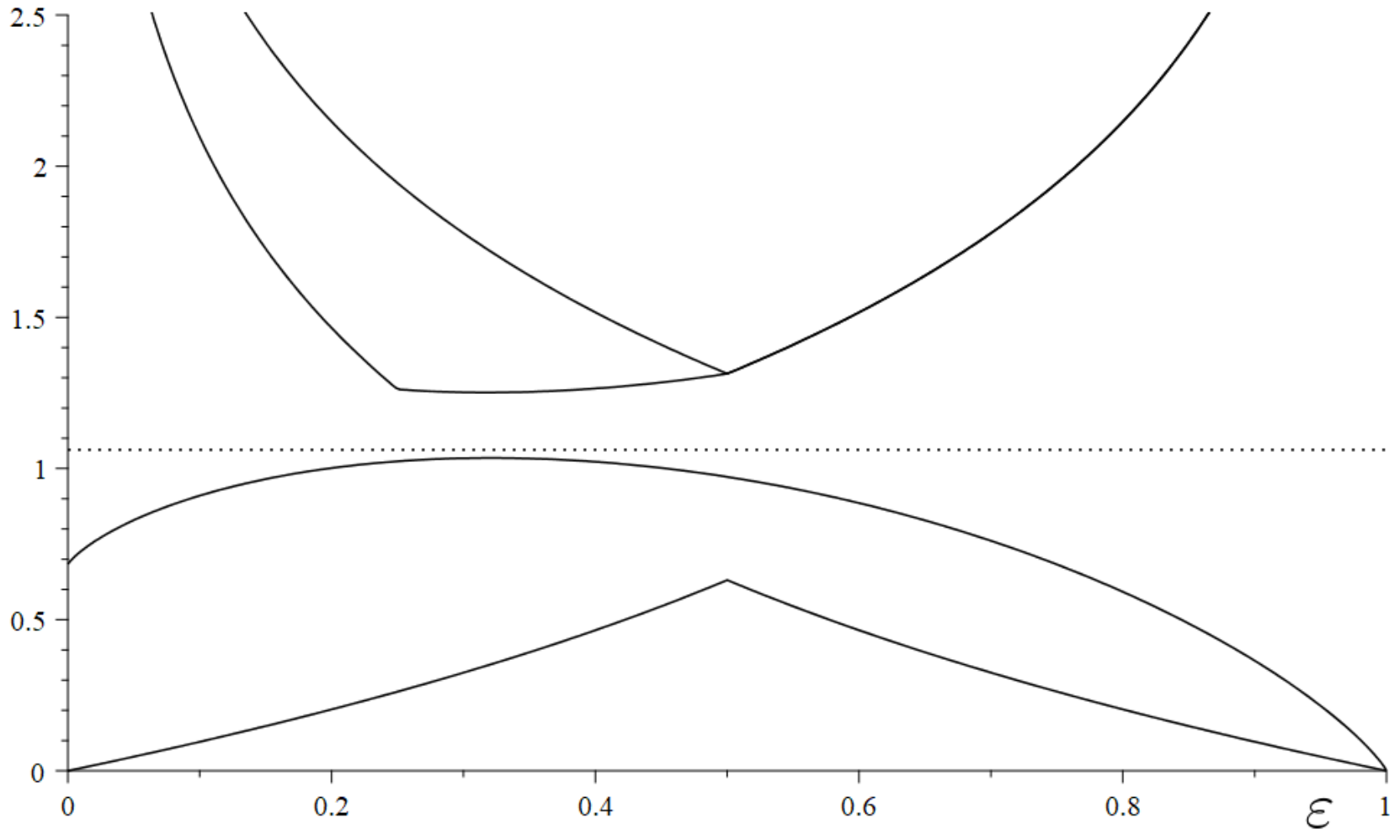

Figure 8.11 Plots of $\dim_L \mu_\varepsilon \leqslant \dim_H \mu_\varepsilon \leqslant \dim_B \mu_\varepsilon \leqslant \dim_A \mu_\varepsilon$ as functions of $\varepsilon \in (0,1)$. The lower and Assouad dimensions of $\mu_\varepsilon$ realise the lower and Assouad dimensions of the support for $\varepsilon = 1/2$ (the coordinate uniform measure) and the Hausdorff dimension of $\mu_\varepsilon$ realises the Hausdorff dimension of the support at $\varepsilon = 1/(1+3^{\log 3/\log 5}) \approx 0.321$ (the McMullen measure). The box dimension of the support (dotted line) is not realised. The uniform measure corresponds to $\varepsilon = 1/4$ and the box dimension of $\mu_\varepsilon$ is minimised at $\varepsilon = \log(5/3)/\log 5 \approx 0.317$.

# 9
# Further examples: attractors and limit sets

In this chapter we collect together discussion of several further families of fractal sets. In Section 9.1 we consider self-conformal sets, which are another special case of IFS attractor. In Section 9.2 we consider sets invariant under parabolic interval maps. While these sets are similar in spirit to IFS attractors, they are not necessarily attractors of IFSs since the inverse branches of the associated dynamical system fail to be strict contractions. The resulting parabolic behaviour greatly influences the Assouad dimension of such sets, see Theorem 9.2.1. In Section 9.3 we consider limit sets of Kleinian groups which are invariant sets for the group action on the boundary of hyperbolic space. Parabolic behaviour is again the main issue here and the dimension theory is developed in detail, see Theorem 9.3.3. In Section 9.4 we consider the random limit sets resulting from Mandelbrot percolation. Here we find a stark difference between the typical behaviour of the Assouad and quasi-Assouad dimension. We discuss a result of Troscheit which brings the $\phi$-Assouad dimensions, discussed in Section 3.3.3, into consideration in order to understand the gap between the quasi-Assouad and Assouad dimensions, see Theorem 9.4.4.

## 9.1 Self-conformal sets

Self-conformal sets are natural generalisations of self-similar sets. The key property that the cylinder sets $S_i(F)$ give rise to efficient covers of $F$ is preserved (this is not true for self-affine sets), but the defining maps need not be similarities. At first sight it seems quite difficult to construct examples of such IFSs and in fact we must pay careful attention to the derivatives of the defining maps. For this reason it takes a little more time to set up the theory. Given an open domain $U \subseteq \mathbb{R}^d$ ($d \geqslant 2$), a continuously differentiable ($C^1$) map $S : U \to \mathbb{R}^d$ is *conformal* if for all $x \in U$ the Jacobian derivative $S'(x) : \mathbb{R}^d \to \mathbb{R}^d$ is an invertible similarity. Our intuition that such maps are hard to construct is somewhat supported by the fact that, at least in high dimensions, there are not many ways to construct conformal maps. Liouville's theorem says that conformal maps on $\mathbb{R}^3$ are either similarities or compositions of reflections in spheres or lines.[1] In $\mathbb{R}^2$ (identified with $\mathbb{C}$) conformal maps correspond precisely to holomorphic functions with non-zero derivative on $U$. In $\mathbb{R}$ the situation is slightly different. Requiring that the derivatives are invertible similarities is not enough to ensure that the cylinders give rise to good covers and more regularity is necessary. (Indeed, *all* non-zero derivatives in $\mathbb{R}$ are invertible similarities, since all affine maps are similarities.) The appropriate condition is to impose more regularity on the derivative. Abusing terminology slightly, given an open interval $U$, we say $S : U \to \mathbb{R}$ is *conformal* if $S$ is $C^{1+\alpha}$ for some $\alpha \in (0, 1]$ with non-zero derivative. This means that $S$ is continuously differentiable with $|S'(x)| > 0$ for all $x \in U$ and the derivative is Hölder continuous, that is, there exist constants $C > 0$ and $\alpha \in (0, 1]$ such that

$$|S'(x) - S'(y)| \leqslant C|x - y|^\alpha$$

for all $x, y \in U$. We say a set $F \subseteq \mathbb{R}^d$ is *self-conformal* if it is the attractor of an IFS consisting of contractions on a compact set $X \subseteq \mathbb{R}^d$ which all extend to injective conformal maps on a common bounded convex open set $U \supset X$. Self-conformal sets in the line are often called 'cookie-cutters', see [69, Chapter 4], where they are sometimes more conveniently expressed as repellers of uniformly expanding interval maps. Here the contractions in the defining IFS correspond to the inverse branches of the expanding map. This is also the case in higher dimensions. Key examples in the plane include Julia sets for certain hyperbolic rational functions on $\mathbb{C}$, such as $z \mapsto z^2 + c$ where $c \in \mathbb{C}$ with $|c| \geqslant 2.475$, or limit sets of certain Kleinian groups.

Angelevska, Käenmäki and Troscheit [5] conducted a thorough investigation of self-conformal sets, mostly in the line, in an effort to develop conformal analogues of results from [96, 81] — specifically those presented here as Theorems 7.2.4 and 7.2.5. There are numerous complications here, not least with the definition of the weak separation property (WSP). Recall our definition of the WSP for IFSs of similarities (Defi-

[1] See the connection to Möbius maps in Section 9.3.

nition 7.2.3) which said the identity should not be in the closure of

$$\{S_{\boldsymbol{i}}^{-1} \circ S_{\boldsymbol{j}} : \boldsymbol{i}, \boldsymbol{j} \in \mathcal{I}^*, \boldsymbol{i} \neq \boldsymbol{j}\} \setminus \{I\}$$

where $I$ is the identity map. There is a problem applying this definition to the conformal case since the domain of $S_{\boldsymbol{i}}^{-1}$ generally will not be the range of $S_{\boldsymbol{j}}$ and so we cannot compose the maps. The reason this works so well in the self-similar case is that a similarity on $\mathbb{R}^d$ is determined by its action on $d+1$ points and so one can blindly extend the defining maps to the whole of $\mathbb{R}^d$. Following [5] we say a conformal IFS satisfies the WSP if

$$\inf \left\{ \|S_{\boldsymbol{i}}'\|^{-1} \sup_{x \in F} |S_{\boldsymbol{i}}(x) - S_{\boldsymbol{j}}(x)| \ : \ \boldsymbol{i}, \boldsymbol{j} \in \mathcal{I}^* \text{ and } S_{\boldsymbol{i}}|_F \neq S_{\boldsymbol{j}}|_F \right\} \ > \ 0.$$

Notice how inverses are not used and instead one normalises by the operator norm of the derivative $S_{\boldsymbol{i}}'$. Moreover, one must restrict attention to the action of the maps on $F$. Also notice that if the maps are similarities, then this is equivalent to the definition of WSP used previously. Actually this was not the definition of the WSP in [5] but it was proved to be equivalent to the WSP when $F$ is not a singleton, see [5, Theorem 3.2]. With this definition in place, the following analogues of Theorems 7.2.4 and 7.2.5 were obtained in [5, Theorem 4.1 and Corollary 4.2].

**Theorem 9.1.1** Let $F \subseteq \mathbb{R}$ be a self-conformal set containing at least two points. If the defining IFS satisfies the WSP, then $\dim_{\mathrm{A}} F = \dim_{\mathrm{H}} F$ and otherwise $\dim_{\mathrm{A}} F = 1$.

**Theorem 9.1.2** Let $F \subseteq \mathbb{R}$ be a self-conformal set containing at least two points with $s = \dim_{\mathrm{H}} F < 1$. Then the following are equivalent:

(i) the defining IFS satisfies the WSP
(ii) $\mathcal{H}^s(F) > 0$
(iii) $0 < \mathcal{H}^s(F) < \infty$
(iv) $F$ is Ahlfors regular
(v) $\dim_{\mathrm{A}} F = \dim_{\mathrm{H}} F$.

The 'best guess' for the Hausdorff (and box, packing and lower) dimension of a self-conformal set is given by the unique solution to *Bowen's formula*, which is the conformal analogue of the Hutchinson-Moran formula (7.2). The *pressure* $P : [0, \infty) \to \mathbb{R}$ is defined by

$$P(s) = \lim_{k \to \infty} \frac{1}{k} \log \sum_{\boldsymbol{i} \in \mathcal{I}^k} \|S_{\boldsymbol{i}}\|^s. \tag{9.1}$$

This limit exists by a standard subadditivity argument, recall *Fekete's lemma*. Moreover, $P$ is decreasing, convex and has a unique root. As in the self-similar case, we expect Bowen's formula to give the Hausdorff dimension of the attractor unless there is a 'good reason' for this not to hold, such as exact overlaps. That is, for self-conformal $F \subseteq \mathbb{R}^d$, it is natural to conjecture that $\dim_H F = \min\{s, d\}$ (where $s$ is given by $P(s) = 0$) provided there are no exact overlaps. We say the IFS has *exact overlaps* if $S_{\boldsymbol{i}}|_F = S_{\boldsymbol{j}}|_F$ for some $\boldsymbol{i}, \boldsymbol{j} \in \mathcal{I}^*$ with $\boldsymbol{i} \neq \boldsymbol{j}$. Note that exact overlaps clearly force $\dim_H F$ to be smaller than the solution to Bowen's formula $P(s) = 0$, and so the interesting case is when there are no exact overlaps. As a direct consequence of Theorem 9.1.2, [5, Theorem 4.3] obtains the following contribution towards the 'overlaps conjecture' in the self-conformal setting.

**Theorem 9.1.3** Let $F \subseteq \mathbb{R}$ be a self-conformal set with $\dim_H F = s \in (0, 1)$ and $\mathcal{H}^s(F) > 0$. Then $P(\dim_H F) = 0$ if and only if there are no exact overlaps.

This result follows from the previous theorem together with two fundamental results regarding self-conformal sets. The first is that the OSC is equivalent to $\mathcal{H}^s(F) > 0$ where $P(s) = 0$, which is a result of Peres, Rams, Simon and Solomyak [231]. The second is that if an IFS satisfies the WSP but does not have exact overlaps, then it in fact satisfies the OSC. This is a result of Deng and Ngai [52, Theorem 1.3].

## 9.2 Invariant sets for parabolic interval maps

The 'cookie cutter' sets mentioned in the previous setting can be viewed as invariant sets for uniformly expanding full branched interval maps modelled by the full shift. This is natural in the setting of IFSs (where the complexity comes from distortion or overlaps, rather than from the underlying dynamical system) but is quite restrictive if one is interested in more general dynamically invariant sets. See [69, 234, 283] for more examples of invariant sets arising from dynamical systems. In this section we consider parabolic (or non-uniformly expanding) interval maps and essentially show that the Assouad dimension of any invariant set is 1, which is maximal.

Let $S : [0, 1] \to [0, 1]$ be a contractive interval map with a parabolic fixed point. By *contractive* we mean that, for all $x, y \in [0, 1]$,

$$|S(x) - S(y)| < |x - y|$$

and by *parabolic fixed point* we mean a point $p \in [0,1]$ such that $S(p) = p$ and $S$ is differentiable at $p$ with $|S'(p)| = 1$. If $p = 0$ or $p = 1$, then we use one-sided derivatives. Note that $S$ is not a strict contraction and we require no regularity assumptions on $S$ other than differentiability at $p$, although it clearly has to be continuous. We are interested in sets which map into themselves under $S$. This is a weak form of dynamical invariance. Recall that $\{x, S(x), S^2(x), \dots\}$ is the *forward orbit* of the point $x$ under $S$.

**Theorem 9.2.1** Let $S$ be a contractive interval map with a parabolic fixed point $p$. If $F \subseteq [0,1]$ contains a point whose forward orbit avoids $p$ and $S(F) \subseteq F$, then $\dim_{\mathrm{A}} F = 1$.

*Proof* Let $z \in F$ be such that the forward orbit of $z$ under $S$ avoids $p$. Since $S$ is contractive, $S^n(z) \to p$ as $n \to \infty$. This follows from a mild generalisation of Banach's contraction mapping theorem which holds for contractive maps on compact spaces (rather than contractions on complete spaces). Moreover, $S^n(z) \neq p$ for all $n \geqslant 1$ and, by the invariance assumption, $F \supseteq \{S^n(z) : n \geqslant 1\}$. Assume that $S'(p) = 1$, noting that if $S'(p) = -1$, then we may replace $S$ with $S \circ S$. The sequence $S^n(z)$ is eventually monotonic and we assume that it is eventually decreasing without loss of generality. Moreover, the gaps $S^n(z) - S^{n+1}(z) > 0$ are decreasing in $n$ since $S$ is contractive and

$$\frac{S^n(z) - S^{n+1}(z)}{S^n(z) - p} = 1 - \frac{S\big(S^n(z)\big) - S(p)}{S^n(z) - p} \to 0$$

since $S'(p) = 1$. It follows that for similarities $T_n : \mathbb{R} \to \mathbb{R}$ defined by

$$T_n(x) = \frac{x - p}{S^n(z) - p}$$

we have

$$d_{\mathcal{H}}(T_n(F) \cap [0,1],\, [0,1]) \leqslant \frac{S^n(z) - S^{n+1}(z)}{S^n(z) - p} \to 0$$

where $d_{\mathcal{H}}$ is the Hausdorff metric. Therefore $[0,1]$ is a weak tangent to $F$ and $\dim_{\mathrm{A}} F = 1$, as required. □

This theorem is rather abstract but it applies to many dynamically invariant sets. For example, $S$ could be a member of a finite or infinite parabolic iterated function system with $F$ as the attractor, see Urbański's *parabolic Cantor sets* [272] (these are not to be confused with the standard IFSs where all the maps are strict contractions). Alternatively, $F$ could be the limit set of a Fuchsian group with at least one parabolic element. This particular example is discussed in more detail in the following section. A third example where Theorem 9.2.1 applies is if $S$ is an inverse branch of a non-uniformly hyperbolic expansive (parabolic) map $f$ where again this map can have finitely or infinitely many branches. In this second case, $F$ could be a dynamical repeller or attractor for $f$ or a survivor set. A useful concrete example to keep in mind is the *Manneville-Pomeau map* $f_\alpha : [0,1] \to [0,1]$ which, for a fixed $\alpha \in (0,1)$, is defined by

$$f_\alpha(x) = x + x^{1+\alpha} \pmod 1.$$

This is a well-studied and important example of a parabolic interval map. For any $k \geqslant 1$, the left-most inverse branch of $f_\alpha^k$ is a contractive interval map with parabolic fixed point 0.

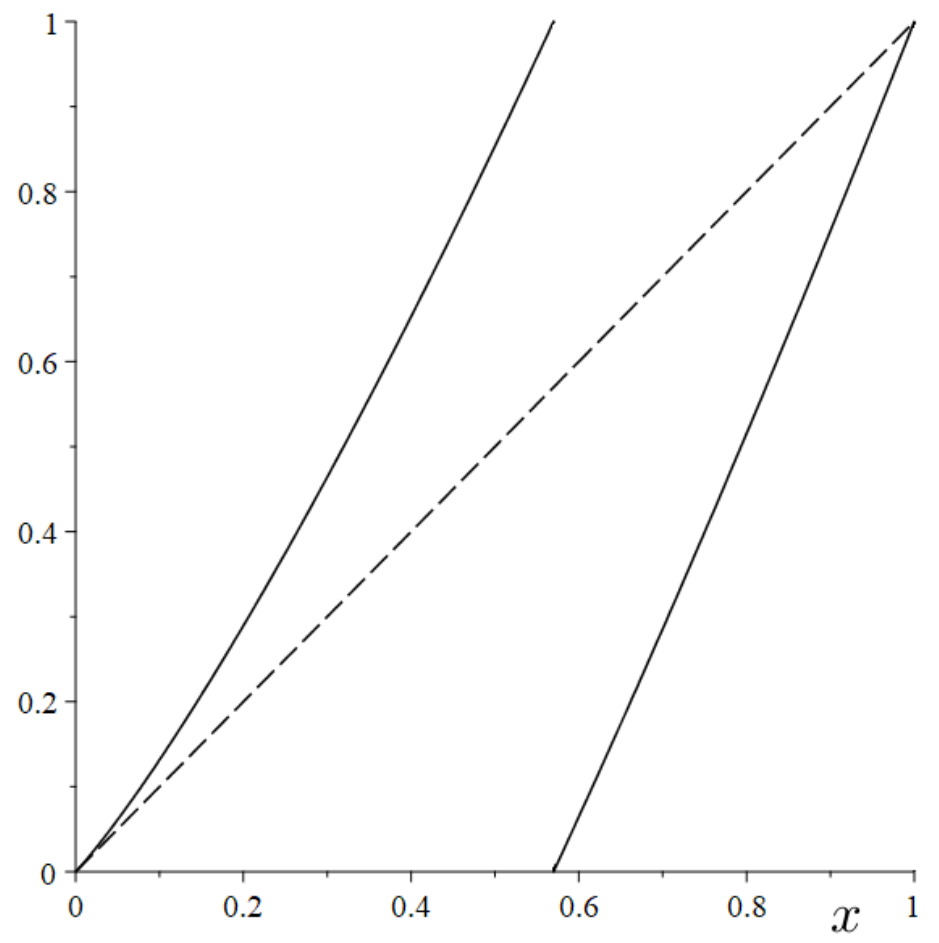

Figure 9.1 Manneville-Pomeau map for $\alpha = 1/2$ with parabolic fixed point at 0.

There is a vast body of literature studying the Hausdorff (and other notions of) dimension of invariant sets for interval maps. Usually strong regularity assumptions are required on the map in order to say anything concrete, see the 'cookie cutters' from [69]. Since the presence of one parabolic point (and one non-trivial orbit) is sufficient to force the Assouad dimension to be maximal we are able to state Theorem 9.2.1 without any such regularity assumptions on $f$. Indeed, we cannot say anything concrete about the Hausdorff dimension of the sets $F$ from Theorem 9.2.1.

Another concrete example to which Theorem 9.2.1 applies is the *Rényi map* or *backward continued fraction map* $\mathcal{R} : [0,1] \to [0,1]$ defined by

$$\mathcal{R}(x) \;=\; \left\{ \frac{1}{1-x} \right\}$$

where $\{y\}$ denotes the fractional part of $y \geqslant 0$. Every irrational $x \in [0,1]$ has a unique non-terminating backwards continued fraction expansion $x = [a_1, a_2, \dots]_{\mathcal{B}}$ where $a_i$ are integers greater than or equal to 2 and

$$x \;=\; [a_1, a_2, \dots]_{\mathcal{B}} \;:=\; 1 - \cfrac{1}{a_1 - \cfrac{1}{a_2 - \cfrac{1}{a_3 - \dots}}}.$$

The Rényi map acts as the left shift on backwards continued fraction expansions in that $\mathcal{R}(x) = [a_2, a_3, \dots]_{\mathcal{B}}$. The Rényi map is the dual to the more recognisable Gauss map, which is naturally associated with the forward continued fraction expansion (or just continued fraction expansion). The Gauss map $\mathcal{G} : [0,1] \to [0,1]$ is defined by

$$\mathcal{G}(x) \;=\; \left\{ \frac{1}{x} \right\}.$$

Similar to above, every irrational $x \in [0,1]$ has a unique non-terminating continued fraction expansion $x = [a_1, a_2, \dots]_{\mathcal{C}}$ where $a_i$ are integers greater than or equal to 1 and

$$x \;=\; [a_1, a_2, \dots]_{\mathcal{C}} \;:=\; \cfrac{1}{a_1 + \cfrac{1}{a_2 + \cfrac{1}{a_3 + \dots}}}.$$

The Gauss map acts as the left shift on continued fraction expansions in that $\mathcal{G}(x) = [a_2, a_3, \dots]_{\mathcal{C}}$.

**Corollary 9.2.2** Given a finite set of positive integers $D \subseteq \mathbb{N}$ with $\# D \geqslant 2$, let

$$\mathcal{C}_D \;=\; \{x \in [0,1] \;:\; x = [a_1, a_2, \dots]_{\mathcal{C}} \text{ with } a_i \in D \text{ for all } i\}$$

and

$$\mathcal{B}_D \;=\; \{x \in [0,1] \;:\; x = [a_1, a_2, \dots]_{\mathcal{B}} \text{ with } a_i \in D \text{ for all } i\}$$

be the sets of irrational numbers in $[0,1]$ whose continued fraction expansion and backwards continued fraction expansion respectively consist only of digits from $D$. Then

$$0 < \dim_{\mathrm{H}} \mathcal{C}_D \;=\; \dim_{\mathrm{A}} \mathcal{C}_D < 1$$

and, if $2 \notin D$, then

$$0 < \dim_{\mathrm{H}} \mathcal{B}_D \ = \ \dim_{\mathrm{A}} \mathcal{B}_D < 1$$

and, if $2 \in D$, then

$$0 < \dim_{\mathrm{H}} \mathcal{B}_D \ < \ \dim_{\mathrm{A}} \mathcal{B}_D = 1.$$

*Proof* For any choice of $D$, the closure of $\mathcal{C}_D$ is a self-conformal set satisfying the OSC and so the Assouad and Hausdorff dimensions coincide and are given by Bowen's formula, that is, the unique zero of the topological pressure (9.1). Moreover, this solution is readily seen to be strictly between 0 and 1. The dimensions of $\mathcal{C}_D$ and its closure coincide here: the Assouad dimension is stable under taking closure, and the Hausdorff dimension is not affected by taking the closure in this case because it only adds countably many (rational points) to the set. The conformal contractions in the defining IFS are the inverse branches of the Gauss map corresponding to the chosen digits. More precisely, the contraction associated with digit $d \geqslant 1$ is the inverse of $\mathcal{G}$ restricted to the interval $[1/(d+1), 1/d]$, see Figure 9.2.

If $2 \notin D$, then the same is true of $\mathcal{B}_D$, although the dimensions of $\mathcal{C}_D$ and $\mathcal{B}_D$ need not coincide. That is, the closure of $\mathcal{B}_D$ is a self-conformal set satisfying the OSC. Here the contraction associated with digit $d \geqslant 2$ is the inverse of $\mathcal{R}$ restricted to the interval $[(d-2)/(d-1), (d-1)/d]$, see Figure 9.2.

However, if $2 \in D$, then $\mathcal{B}_D$ is mapped into itself by the left most inverse branch of the Rényi map, that is, $S(x) = x/(1+x)$. This is a contractive interval map with parabolic fixed point $p = 0$ and Theorem 9.2.1 applies to give $\dim_{\mathrm{A}} \mathcal{B}_D = 1$. The fact that $0 < \dim_{\mathrm{H}} \mathcal{B}_D < 1$ follows by applying ideas from the thermodynamic formalism, see [272]. □

Both the Gauss map and the Rényi map have a point where the derivative has modulus 1. However, for the Gauss map this point is not fixed but for the Rényi map it is fixed, which leads to rather different behaviour in the context of Assouad dimension.

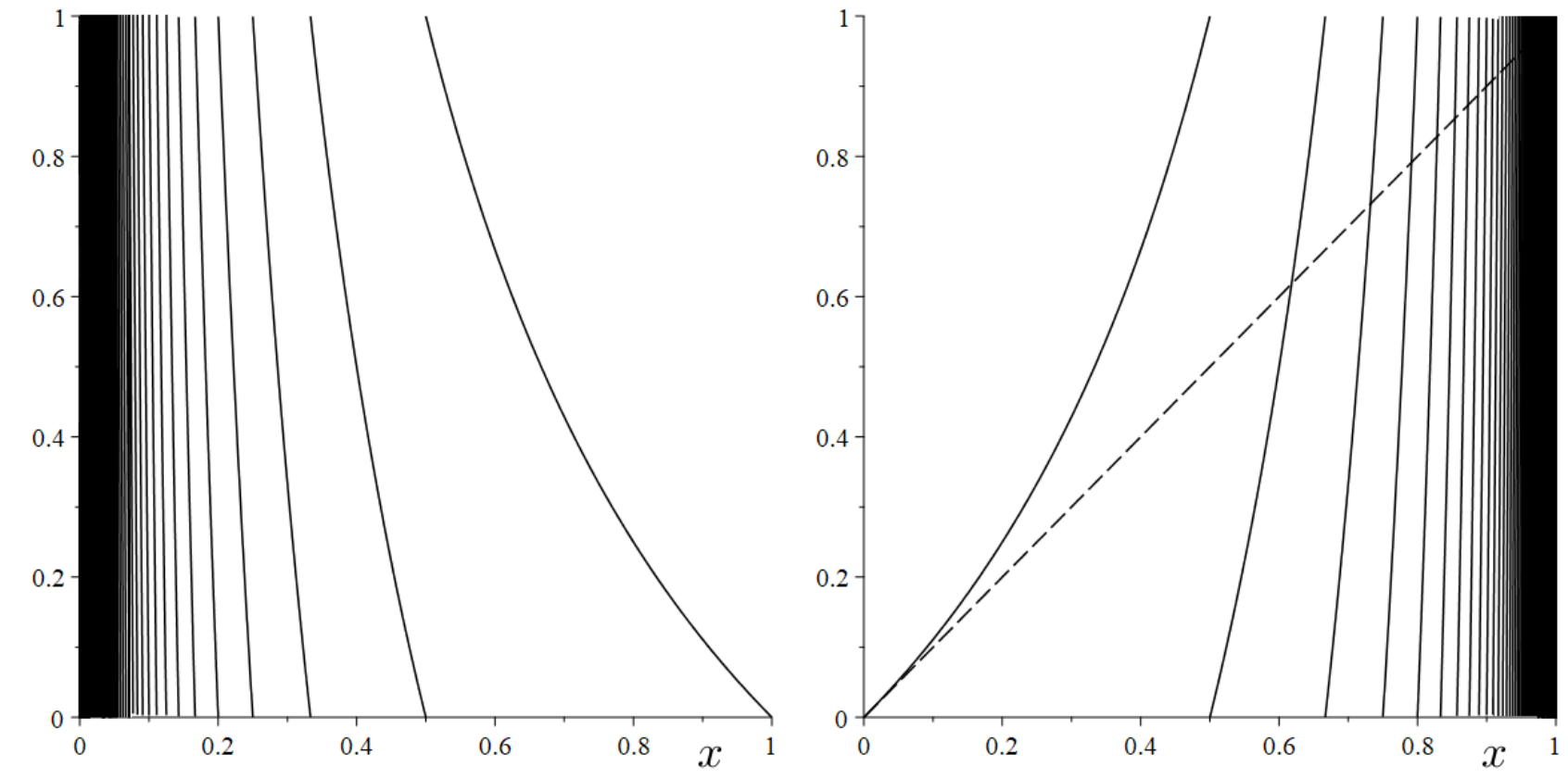

Figure 9.2 Left: the Gauss map. Right: the Rényi map.

## 9.3 Limit sets of Kleinian groups

The key property enjoyed by the examples in the previous section is non-uniform hyperbolicity, or *parabolicity.* This phenomenon occurs in many interesting contexts across dynamics, often leading to complicated behaviour. In this section we focus on another example, Kleinian limits sets, where we can take the theory further.

Beautiful fractal sets arise as limit sets of Kleinian groups acting on hyperbolic space. Here, parabolic elements in the group lead to parabolic dynamics on the limit set. For integer $d \geqslant 1$, the Poincaré ball

$$\mathbb{D}^{d+1} = \left\{ z \in \mathbb{R}^{d+1} \ : \ |z| < 1 \right\}$$

equipped with the hyperbolic metric $d$ defined by

$$ds = \frac{2|dz|}{1 - |z|^2}$$

is a popular model of $(d+1)$-dimensional hyperbolic space. The important geometric feature of this metric space is that distance blows up near the boundary, leading to a much 'larger' *boundary at infinity* than that possessed by Euclidean space. The boundary is $\mathbb{S}^d = \left\{ z \in \mathbb{R}^{d+1} \ : \ |z| = 1 \right\}$. This 'extra space' at infinity leads to a richer group of isometries, which also has a particularly elegant description: it is the stabiliser of $\mathbb{D}^{d+1}$ in the Möbius group acting on $\overline{\mathbb{R}^{d+1}}$. The *Möbius group* is the group of all maps which can be expressed as compositions of reflections in spheres. We will mostly be concerned with orientation preserving isometries of $(\mathbb{D}^{d+1}, d)$ and denote this group by $\mathrm{Con}(d)$. For a more detailed discussion of hyperbolic geometry, see [25, 203, 157].

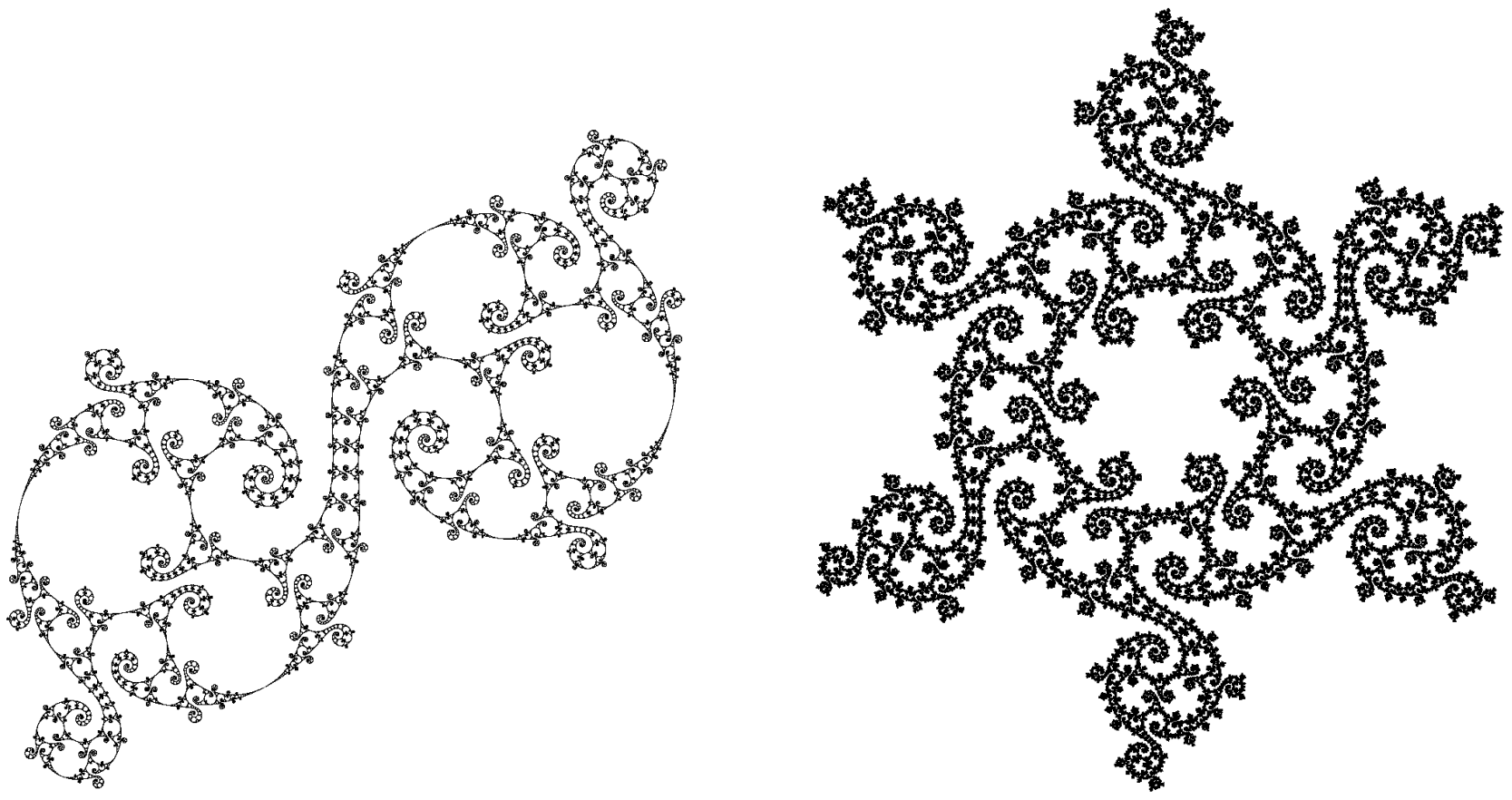

Figure 9.3 Two Kleinian limit sets.

A *Kleinian group* is a discrete subgroup of $\mathrm{Con}(d)$. We write $H \leqslant G$ for groups $H, G$ to mean that $H$ is a subgroup of $G$. Kleinian groups act *properly discontinuously* on $\mathbb{D}^{d+1}$. That is, all orbits are locally finite, meaning they intersect every compact set in at most finitely many points. The induced action of a Kleinian group on the boundary $\mathbb{S}^d$ of $\mathbb{D}^{d+1}$ is much more complicated and its associated limit set captures this complexity. Formally, writing $\mathbf{0} = (0, \ldots, 0) \in \mathbb{D}^{d+1}$, the *limit set* of a Kleinian group $\Gamma$ is defined by

$$L(\Gamma) = \overline{\Gamma(\mathbf{0})} \setminus \Gamma(\mathbf{0})$$

where $\Gamma(\mathbf{0})$ is the orbit of $\mathbf{0}$ under $\Gamma$ and the closure is the Euclidean closure. This is a compact subset of $\mathbb{S}^d$, which makes the limit set an inherently Euclidean object, despite its genesis in non-Euclidean geometry. We can choose to metrise the boundary with the inherited metric from $\mathbb{R}^{d+1}$, but sometimes it is more convenient to work with the spherical metric. However, these metrics are bi-Lipschitz equivalent and therefore the choice makes no difference to the dimension theory of the limit set.

**Lemma 9.3.1** Limit sets are closed subsets of $\mathbb{S}^d$ with respect to the Euclidean metric, and are thus compact.

*Proof* The fact that $L(\Gamma) \subseteq \mathbb{S}^d$ follows immediately since $\Gamma(\mathbf{0})$ is a discrete subset of $\mathbb{D}^{d+1}$ and so $\Gamma(\mathbf{0})$ cannot accumulate anywhere except on the boundary. Turning our attention to closedness, let $z_n \in L(\Gamma)$ be a sequence of points such that $z_n \to z$ in the Euclidean metric. For each $n$, let $g_{n,m} \in \Gamma$ be such that $g_{n,m}(\mathbf{0}) \to z_n$ in the Euclidean metric as $m \to \infty$. We can find such $g_{n,m}$ since $z_n \in L(\Gamma)$. For a given $n$, choose $m_n$ sufficiently large to ensure that

$$|g_{n,m_n}(\mathbf{0}) - z_n| \leqslant 1/n.$$

Then

$$|g_{n,m_n}(\mathbf{0}) - z| \leqslant |g_{n,m_n}(\mathbf{0}) - z_n| + |z_n - z| \leqslant 1/n + |z_n - z| \to 0$$

as $n \to \infty$ and therefore $z \in L(\Gamma)$. □

Limit sets provide another important example of dynamical invariance, where this time the dynamics comes from a group action.

**Lemma 9.3.2** Limit sets are (strongly) invariant under the action of the associated Kleinian group, that is, if $\Gamma$ is a Kleinian group, then, for all $g \in \Gamma$,

$$g(L(\Gamma)) = L(\Gamma).$$

*Proof* Let $z \in L(\Gamma)$, $g \in \Gamma$, and let $g_n \in \Gamma$ be such that $g_n(\mathbf{0}) \to z$. Since $g$ is a continuous map on the open Euclidean unit ball, it follows that

$$g(g_n(\mathbf{0})) \to g(z)$$

which implies $g(z) \in L(\Gamma)$. Thus we have proved that $g(L(\Gamma)) \subseteq L(\Gamma)$. The opposite inclusion follows by replacing $g$ with $g^{-1}$ and applying $g$. □

If the limit set is empty or consists only of one or two points, then the Kleinian group is called *elementary* and otherwise it is *non-elementary*, in which case it is necessarily uncountable. Roughly speaking, a Kleinian group is called *geometrically finite* if it has a fundamental domain with finitely many sides and this is the most widely studied class of Kleinian group. The dimension theory of the limit sets of geometrically finite Kleinian groups is particularly well-developed and elegant, whereas the geometrically infinite setting is rather more complicated. The *Poincaré exponent* of a Kleinian group $\Gamma$ is defined by

$$\delta(\Gamma) = \inf\left\{ s > 0 \ : \ \sum_{g \in \Gamma} \exp(-s\, d(\mathbf{0}, g(\mathbf{0}))) < \infty \right\}$$

and plays a central role in the geometry and dimension theory of $\Gamma$. In particular, the limit set of a non-elementary geometrically finite Kleinian group has Hausdorff dimension equal to $\delta(\Gamma)$. This seminal result goes back to the influential papers of Patterson (for Fuchsian groups with some assumptions on parabolic elements) [230, Theorems 4.1 and 5.1] and Sullivan (for the general higher dimensional case) [261, Theorem 1]. When $d = 1$, discrete subgroups of $\text{Con}(1)$ are called *Fuchsian*, instead of Kleinian. Around 20 years later it was shown that the packing and box dimensions of the limit set are also given by $\delta(\Gamma)$. This latter result is due independently to Bishop and Jones [33, Corollary 1.5] and Stratmann and Urbański [258, Theorem 3].

A particularly well-known example of a Kleinian limit set is the *Apollonian gasket* or *Apollonian circle packing*, see Figure 13.5. It can be constructed by starting with four mutually tangent circles lying in $\mathbb{C}$, one containing the other three, and then inductively adding in circles of the largest possible radius which lie tangent to three previously added circles, see [236, 229]. It is well-known that given any two circle packings formed in this way there is a Möbius transformation $g \in \text{PSL}(2, \mathbb{C}) \cong \text{Con}(2)$ taking one to the other, that is, there is a unique circle packing up to Möbius images. Therefore we may talk about *the* Apollonian circle packing. Moreover, the closure of the Apollonian circle packing is the limit set of a geometrically finite Kleinian group acting on the upper half-space model of hyperbolic geometry, where the boundary is $\mathbb{C} \cup \{\infty\}$ instead of $\mathbb{S}^2$. In the $\mathbb{D}^3$ model of hyperbolic space, which we focus on, the Apollonian circle packing would appear pasted onto the side of $\mathbb{S}^2$.

Limit sets of non-elementary geometrically finite Kleinian groups are also known to support a particularly nice measure, known as the *Patterson-Sullivan measure*. We write $\mu$ for the Patterson-Sullivan measure for the rest of this section. The Patterson-Sullivan measure is an atomless Borel probability measure. Moreover, it is conformal, meaning it respects the group action in a certain sense, and ergodic under the dynamics of the group action. It is also of Hausdorff dimension $\delta(\Gamma)$. Stratmann and Velani's *global measure formula*, see [259, Theorem 2], gives a formula for the Patterson-Sullivan measure of any ball up to uniform constants and is useful in deriving the Assouad and lower dimensions of the measure, which in turn have implications for the Assouad and lower dimensions of the limit set itself. In order to state the global measure formula, one needs to introduce several technical concepts from hyperbolic geometry and so we omit the precise statement. However, we provide a heuristic description as well as the formulae for the Assouad and lower dimen-sions which were calculated in [92], referring the interested reader to [92, 259, 258, 33].

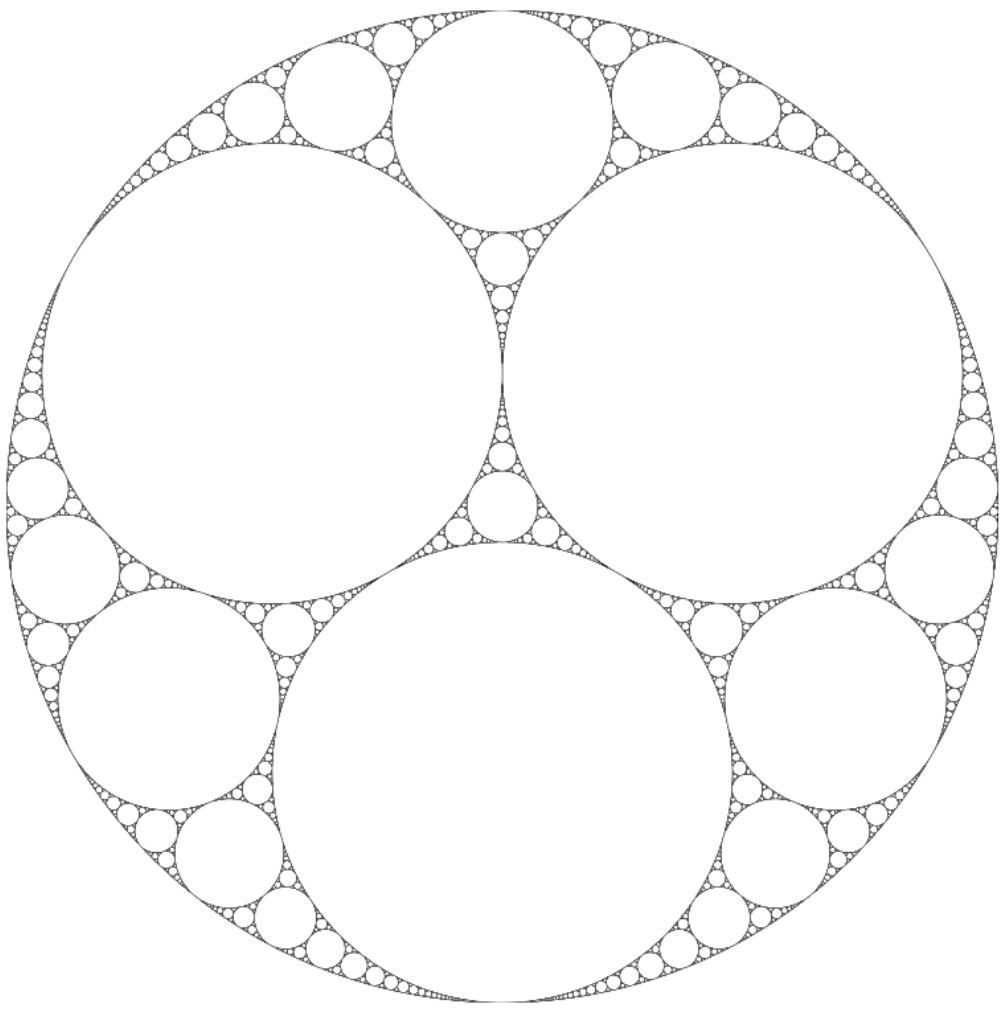

Figure 9.4 The Apollonian circle packing.

Fix a non-elementary geometrically finite Kleinian group $\Gamma$ and let $P \subseteq L(\Gamma)$ denote the countable set of all parabolic fixed points, that is, points fixed by parabolic elements of $\Gamma$. Isometries of $\mathbb{D}^{d+1}$ can be classified as hyperbolic, loxodromic, elliptic, or parabolic, depending on their fixed points. The parabolic elements are the ones with precisely one fixed point on the boundary $\mathbb{S}^d$ and are the elements which lead to interesting dynamical behaviour. The stabiliser of a parabolic fixed point $p$ cannot contain hyperbolic or loxodromic elements since if a subgroup of $\mathrm{Con}(d)$ contains a parabolic and a hyperbolic/loxodromic element which fix the same point, then the group is not discrete. Therefore, the parabolic elements in the stabiliser of $p$ in $\Gamma$ generate a free abelian group of finite index (as a subgroup of the stabiliser). Let $k(p) \in \{1, \dots, d\}$ be the maximal rank of a free abelian subgroup of the stabiliser of $p$ in $\Gamma$, which is necessarily generated by $k(p)$ parabolic elements all fixing $p$. The numbers $k(p)$ play an important role in the dimension theory of Kleinian limit sets. For example, it is well-known that $\delta(\Gamma) > k(p)/2$ for all $p \in P$.

Given $z \in L(\Gamma)$ and $R \in (0, 1)$, the global measure formula states that

$\mu(B(z,R))$ is given, up to uniform constants independent of $z$ and $R$, by

$$R^{\delta(\Gamma)} R^{\rho(z,R)(\delta(\Gamma)-k(z,R))} \tag{9.2}$$

where, roughly speaking, $k(z,R) = k(p)$ if $z$ is 'near' the parabolic fixed point $p \in P$ and $0 \leqslant \rho(z,R) \leqslant 1$ quantifies 'near' in terms of a fixed standard set of *horoballs*.[2] A horoball is a Euclidean sphere whose interior lies in $\mathbb{D}^{d+1}$ and is tangent to $\mathbb{S}^d$ at a parabolic fixed point. For more details on standard horoballs and ranks of parabolic elements, see [258]. A simple initial observation to make here is that if there are no parabolic elements, then the global measure formula (9.2) simplifies to

$$R^{\delta(\Gamma)}. \tag{9.3}$$

Therefore, applying Lemma 4.1.2, for non-elementary geometrically finite Kleinian groups with no parabolic fixed points,

$$\dim_{\mathrm{A}} L(\Gamma) = \dim_{\mathrm{L}} L(\Gamma) = \dim_{\mathrm{A}} \mu = \dim_{\mathrm{L}} \mu = \delta(\Gamma).$$

In light of this result, we assume from now on that $\Gamma$ contains a parabolic element and write $1 \leqslant k_{\min} \leqslant k_{\max} \leqslant d$ to denote the minimal and maximal ranks of parabolic fixed points, respectively. The main results from [92] are stated below.

**Theorem 9.3.3** Let $\Gamma \leqslant \mathrm{Con}(d)$ be a non-elementary geometrically finite Kleinian group with a parabolic element. Then

$$\dim_{\mathrm{A}} \mu = \max\{k_{\max},\, 2\delta(\Gamma) - k_{\min}\},$$

$$\dim_{\mathrm{L}} \mu = \min\{k_{\min},\, 2\delta(\Gamma) - k_{\max}\},$$

$$\dim_{\mathrm{A}} L(\Gamma) = \max\{k_{\max},\, \delta(\Gamma)\}$$

and

$$\dim_{\mathrm{L}} L(\Gamma) = \min\{k_{\min},\, \delta(\Gamma)\}.$$

Recall that, for any non-elementary geometrically finite Kleinian group $\Gamma$,

$$\dim_{\mathrm{H}} L(\Gamma) = \dim_{\mathrm{P}} L(\Gamma) = \dim_{\mathrm{B}} L(\Gamma) = \dim_{\mathrm{H}} \mu = \delta(\Gamma).$$

For $d = 2$ we can 'upgrade' the previous theorem to include finitely generated geometrically infinite Kleinian groups thanks to a result of Bishop and Jones [33]. They proved that if $\Gamma \leqslant \mathrm{Con}(2)$ is a finitely generated geometrically infinite Kleinian group, then $\dim_{\mathrm{H}} L(\Gamma) = 2$. In the geometrically infinite case, we still have $\delta(\Gamma) \leqslant \dim_{\mathrm{H}} L(\Gamma)$, but this inequality can be strict.

[2] Note that this is not how the global measure formula is usually presented, for example in [259, 92]. It is more common for the radius of the ball to be written as $R = e^{-t}$ for $t > 0$. We chose to express it in this way to make the notation consistent with the rest of the book. Finally, to simplify exposition, our function $\rho$ has been 'normalised' by dividing the usual $\rho$ by $|\log R|$.

**Corollary 9.3.4** Let $\Gamma \leqslant \mathrm{Con}(2)$ be a finitely generated non-elementary Kleinian group. Then

$$\dim_{\mathrm{A}} L(\Gamma) = \max\{k_{\max}, \dim_{\mathrm{H}} L(\Gamma)\}.$$

The full proofs of Theorem 9.3.3 are technical, but to see where the formulae come from is fairly straightforward. The following discussion is deliberately very rough and we refer the reader to [92] for the details. First, consider the Patterson-Sullivan measure $\mu$. If we choose $z = p \in P$, then $k(z,R) = k(p)$ for all $R > 0$ sufficiently small and $\rho(z,R) \to 1$ as $R \to 0$. Therefore, applying the global measure formula (9.2), for all $\varepsilon > 0$ there is a constant $C \geqslant 1$ such that

$$C^{-1}\left(\frac{R}{r}\right)^{2\delta(\Gamma)-k(p)-\varepsilon} \leqslant \frac{\mu(B(z,R))}{\mu(B(z,r))} \leqslant C\left(\frac{R}{r}\right)^{2\delta(\Gamma)-k(p)+\varepsilon}$$

which immediately gives

$$\dim_{\mathrm{L}} \mu \leqslant 2\delta(\Gamma) - k_{\max} \leqslant 2\delta(\Gamma) - k_{\min} \leqslant \dim_{\mathrm{A}} \mu. \tag{9.4}$$

Observing the behaviour

$$\left(\frac{R}{r}\right)^{k(p)}$$

is a little harder, and requires a more careful choice of $z$. If we choose $z$ 'close to $p$ at scale $R$' but 'far away from $p$ at scale $r$', then we can achieve

$$k(z,R) = k(p), \qquad k(z,r) = 0$$

$$\rho(z,R) = \frac{\log(R/r)}{|\log R|}, \qquad \rho(z,r) = 0$$

which, applying (9.2), gives $\frac{\mu(B(z,R))}{\mu(B(z,r))}$, up to constants, as

$$\frac{R^{\delta(\Gamma)} R^{\frac{\log(R/r)}{|\log R|}(\delta(\Gamma)-k(p))}}{r^{\delta(\Gamma)}} = \left(\frac{R}{r}\right)^{k(p)}$$

which yields

$$\dim_{\mathrm{L}} \mu \leqslant k_{\min} \leqslant k_{\max} \leqslant \dim_{\mathrm{A}} \mu. \tag{9.5}$$

Combining (9.4) and (9.5) yields the lower bound for $\dim_A \mu$ and the upper bound for $\dim_L \mu$ from Theorem 9.3.3. The remaining inequalities are proved by demonstrating that the scaling behaviour witnessed by the choices of $z, R$ and $r$ above are extremal. This requires more delicate understanding of the underlying hyperbolic geometry and we omit the details.

Concerning the Assouad dimension of the limit set $L(\Gamma)$, we have the lower bound $\dim_A L(\Gamma) \geqslant \dim_H L(\Gamma) = \delta(\Gamma)$ and the lower bound $\dim_A L(\Gamma) \geqslant k_{\max}$ is also fairly straightforward. In fact $[-1,1]^{k(p)}$ (embedded in $\mathbb{R}^{d+1}$) is a (weak) tangent to $L(\Gamma)$ at $p \in P$. This can be seen by taking the orbit of a point $z \in L(\Gamma)$ with $z \neq p$ under the free abelian group of rank $k(p)$ lying inside the stabiliser of $p$, which we denote by $\Gamma'$. Using the fact that $\Gamma$ is discrete, we find a $\mathbb{Z}^{k(p)}$ action with $p$ playing the role of infinity. This leads to $\Gamma'(z)$ being an inverted $\mathbb{Z}^{k(p)}$-lattice, accumulating at $p$, and letting $T_l$ $(l \geqslant 1)$ be a sequence of similarities which send $p \to 0 \in \mathbb{R}^{d+1}$ and have similarity ratio equal to $l$, we find (an embedding of) $[-1,1]^{k(p)}$ is contained in a weak tangent to $L(\Gamma)$ as required.

Similarly, the upper bound $\dim_L L(\Gamma) \leqslant \delta(\Gamma)$ is straightforward. The upper bound $\dim_L L(\Gamma) \leqslant k_{\min}$ seems like it should follow easily too, but there is an additional complication which does not occur in the Assouad case since we need to ensure that the weak tangent is precisely $[-1,1]^{k(p)}$ (not just a set containing $[-1,1]^{k(p)}$). However, a deep result of Bowditch [36] guarantees that the limit set lies within a bounded neighbourhood of a $k(p)$-dimensional affine plane containing $p$. This allows us to conclude that the weak tangent constructed above is precisely $[-1,1]^{k(p)}$.

The upper bound for the Assouad dimension and lower bound for the lower dimension are rather more difficult to prove and we omit the details. However, note that if $\delta(\Gamma) \leqslant (k_{\min} + k_{\max})/2$, then

$$\dim_A L(\Gamma) \leqslant \dim_A \mu = \max\{k_{\max},\, 2\delta(\Gamma) - k_{\min}\} = k_{\max}$$

and if $\delta(\Gamma) \geqslant (k_{\min} + k_{\max})/2$, then

$$\dim_L L(\Gamma) \geqslant \dim_L \mu = \min\{k_{\min},\, 2\delta(\Gamma) - k_{\max}\} = k_{\min}$$

and so the difficult bounds can be obtained 'cheaply' in these ranges. Interestingly, one of the tough bounds is always 'cheap'. Moreover, this highlights the situation when $\delta(\Gamma) = (k_{\min} + k_{\max})/2$ as particularly interesting (or straightforward). In this case

$$\dim_A L(\Gamma) = \dim_A \mu = k_{\max} \geqslant k_{\min} = \dim_L \mu = \dim_L L(\Gamma).$$

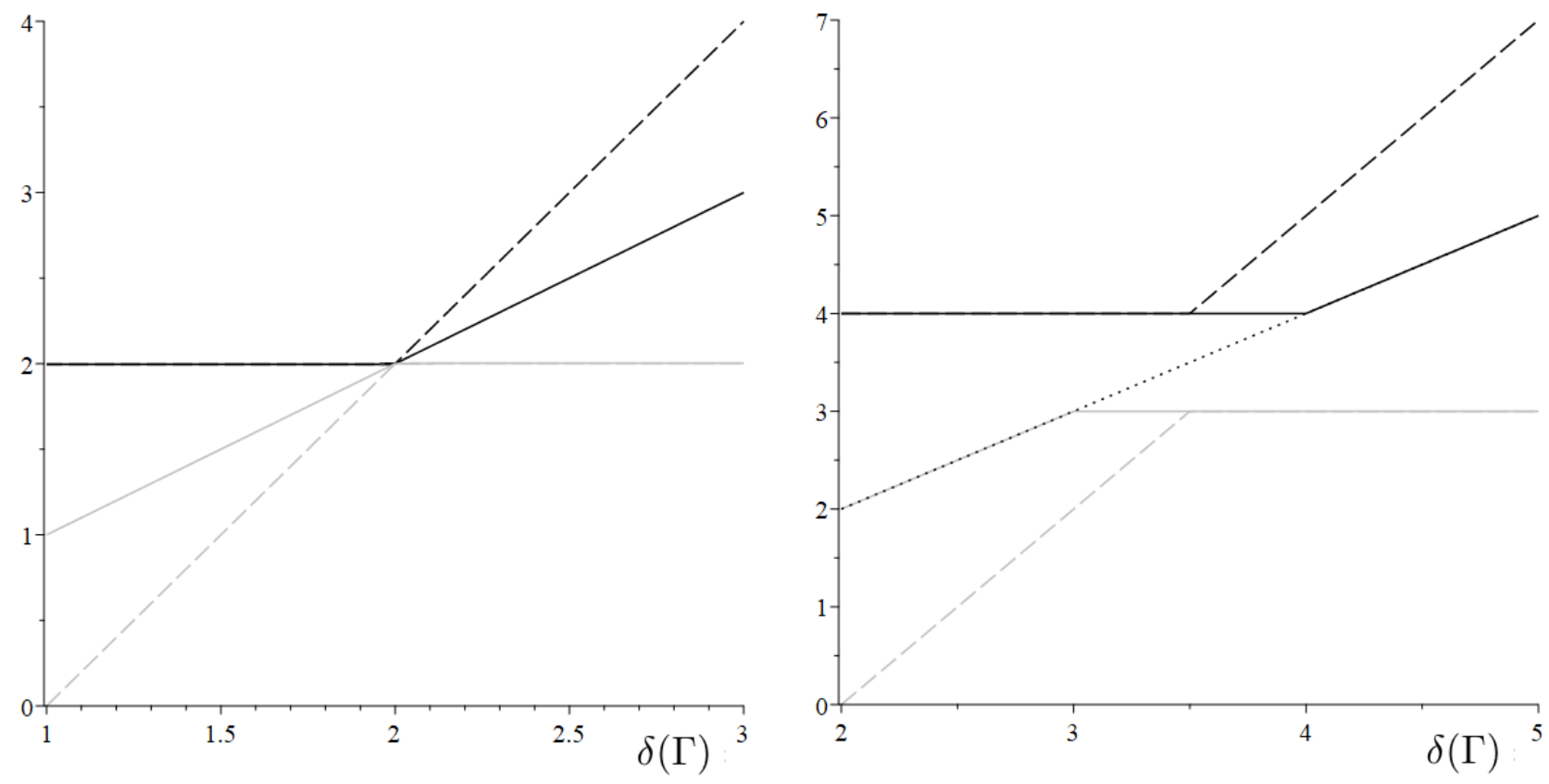

Figure 9.5 Two plots showing the dimensions from Theorem 9.3.3 as functions of $\delta(\Gamma) \in (k_{\max}/2, d]$ in two different settings. The Assouad and lower dimensions of $\mu$ are plotted with dashed lines, the Assouad and lower dimensions of $L(\Gamma)$ are plotted with solid lines, and $\delta(\Gamma) = \dim_H L(\Gamma) = \dim_H \mu$ is plotted with a dotted line for reference. The Assouad dimensions are plotted in black and the lower dimensions are plotted in grey. On the left $k_{\min} = k_{\max} = 2$ and $d = 3$, and on the right $k_{\min} = 3 < 4 = k_{\max}$ and $d = 5$.

**Example 9.3.5** Consider the Apollonian circle packing. The parabolic fixed points are the points of mutual tangency between two circles in the packing and it is straightforward to see that the rank of each of these points is 1, and therefore $k_{\min} = k_{\max} = 1$. Estimating the Poincaré exponent for this group is difficult, but has received a lot of attention in the literature and very good bounds are now available. In particular, $\delta(\Gamma) \approx 1.305\ldots$, see [211, 37, 229]. Therefore

$$\dim_A \mu = 2\delta(\Gamma) - 1 \approx 1.61\ldots$$

$$\dim_A L(\Gamma) = \delta(\Gamma) \approx 1.305\ldots$$

and

$$\dim_L \mu = \dim_L L(\Gamma) = 1.$$

Let us again revisit the *simultaneous realisation problem* from Section 6.2, which asks if there exists an invariant measure which simultaneously realises different notions of dimension. Corollary 9.3.6 below shows that the Patterson-Sullivan measure solves this problem for the Assouad, lower and Hausdorff dimensions in one very specific setting. This should be compared with the situation for Bedford-McMullen carpets, where there are no self-affine measures which solve the simultaneous realisation problem in the non-uniform fibres case. Note that $\mu$ *always* realises the Hausdorff dimension of $L(\Gamma)$, but only sometimes realises the Assouad and lower dimensions.

**Corollary 9.3.6** Let $\Gamma \leqslant \mathrm{Con}(d)$ be a non-elementary geometrically finite Kleinian group. Then

$$\dim_{\mathrm{L}} \mu = \dim_{\mathrm{L}} L(\Gamma) < \dim_{\mathrm{H}} \mu = \dim_{\mathrm{H}} L(\Gamma) < \dim_{\mathrm{A}} \mu = \dim_{\mathrm{A}} L(\Gamma)$$

if and only if $k_{\min} < k_{\max}$ and $\delta(\Gamma) = (k_{\min} + k_{\max})/2$.

Another common consideration in dimension theory is to identify possible relationships between dimensions in particular contexts, see Figure 17.1. Despite the fact that $\dim_{\mathrm{H}} L(\Gamma) = \dim_{\mathrm{B}} L(\Gamma)$ always holds for geometrically finite Kleinian limit sets, the behaviour regarding the Assouad and lower dimensions is quite rich. We note that $\dim_{\mathrm{H}} L(\Gamma) = \dim_{\mathrm{B}} L(\Gamma)$ does not always hold in the geometrically infinite case. In dimension 1, that is, when $\Gamma \leqslant \mathrm{Con}(1)$ is a Fuchsian group, we have a dichotomy similar to that observed for self-similar (and self-conformal) sets in $\mathbb{R}$, see Theorems 7.2.4 and 9.1.1.

**Corollary 9.3.7** Let $\Gamma \leqslant \mathrm{Con}(1)$ be a non-elementary geometrically finite Fuchsian group with $L(\Gamma) \neq \mathbb{S}^1$.

(i) If $\Gamma$ does not contain any parabolic elements, then

$$0 < \dim_{\mathrm{L}} L(\Gamma) = \dim_{\mathrm{H}} L(\Gamma) = \dim_{\mathrm{A}} L(\Gamma) = \delta(\Gamma) < 1.$$

(ii) If $\Gamma$ contains a parabolic element, then

$$1/2 < \dim_{\mathrm{L}} L(\Gamma) = \dim_{\mathrm{H}} L(\Gamma) = \delta(\Gamma) < \dim_{\mathrm{A}} L(\Gamma) = 1.$$

In higher dimensions, more relationships are possible. In fact we can achieve

$$\dim_{\mathrm{L}} L(\Gamma) < \dim_{\mathrm{H}} L(\Gamma) < \dim_{\mathrm{A}} L(\Gamma),$$

$$\dim_{\mathrm{L}} L(\Gamma) = \dim_{\mathrm{H}} L(\Gamma) < \dim_{\mathrm{A}} L(\Gamma)$$

or

$$\dim_{\mathrm{L}} L(\Gamma) < \dim_{\mathrm{H}} L(\Gamma) = \dim_{\mathrm{A}} L(\Gamma).$$

This should again be compared with the situation for Bedford-McMullen carpets, where we either have the lower, Hausdorff and Assouad dimensions equal to a common value or to three distinct values, see Corollary 8.3.2.

We have seen in this section that the lower dimension of the limit set of a non-elementary geometrically finite Kleinian group is positive. In fact this phenomenon has been observed in another context. Järvi and Vuorinen [148, Theorem 5.5] proved that limit sets of finitely generated Kleinian groups are *uniformly perfect*, which turns out to be equivalent, see Theorem 13.1.2. Uniform perfectness has been investigated extensively in the context of Kleinian limit sets, see [260] and the references therein.

## 9.4 Mandelbrot percolation

*Mandelbrot percolation* is a random process giving rise to fractals which are *statistically* self-similar, see [70, Section 15.2]. The model is simple, but it captures many key features of randomness and has attracted a lot of attention in the fractal geometry literature, see for example [238, 239].

We recall the construction, which begins with the unit cube $M_0 = [0,1]^d$, a fixed integer $m \geqslant 2$, and a probability $p \in (0,1)$. At the first step we divide $M_0$ into $m^d$ (closed) cubes of side length $m^{-1}$ and, for each cube, we independently choose to keep it with probability $p$, or remove it with probability $(1-p)$. Let $M_1$ be the union of cubes which are kept and then repeat the first step inside each kept cube independently. Iterating this procedure leads to a decreasing sequence of closed sets $M_k$ defined to be the union of cubes of sidelength $m^{-k}$ which are kept at stage $k$. The limit set is the random set $M = \cap_k M_k$, see Figure 9.6. It is well-known that if $p > m^{-d}$, then $M$ is non-empty with positive probability. We assume from now on that this is the case. To get a feel for why this is the threshold, note that inside a given cube the expected number of kept cubes at the next level is $pm^d$. If we condition on $M$ being non-empty, then almost surely

$$\dim_{\mathrm{H}} M = \dim_{\mathrm{B}} M = d + \frac{\log p}{\log m} \in (0, d).$$

The Assouad dimension behaves rather differently, as seen in the following result of Fraser, Miao and Troscheit [104]. This result can also be derived from earlier work of Berlinkov and Järvenpää on porosity [29]. An immediate and somewhat counter-intuitive consequence of this result is that almost surely one cannot embed $M$ in $\mathbb{R}^{d-1}$ via a bi-Lipschitz map, even if the Hausdorff dimension is almost surely much smaller than $d-1$ (even close to 0).

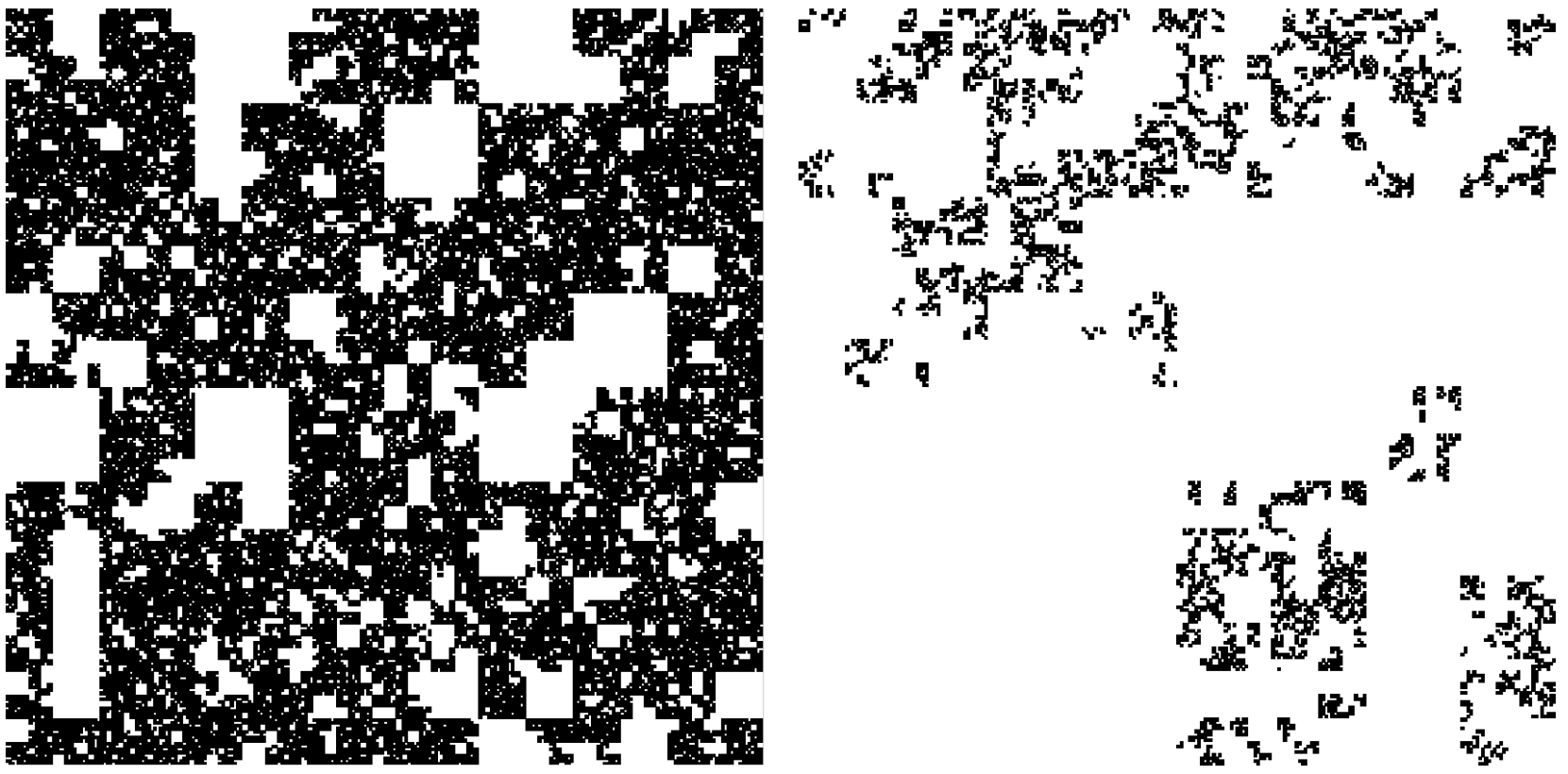

Figure 9.6 Two examples of limit sets for Mandelbrot percolation. In both cases $d = m = 2$. Can you tell what $p$ is? (A feature of finite resolution images is that we cannot rule out any $p$ with total certainty just by looking at the picture.)

**Theorem 9.4.1** Conditioned on $M$ being non-empty, almost surely

$$\dim_{\mathrm{A}} M = d.$$

*Proof* Since we condition on $M$ being non-empty we may assume that, for all levels $k \geqslant 1$, there exists at least one kept cube $Q_k \subseteq M_k$. For a given kept cube $Q_k \subseteq M_k$ and an integer $i \geqslant 0$, let $\mathcal{Q}_k^i$ denote the collection of all level $k+i$ cubes contained in $Q_k$ (not necessarily intersecting $M$). There are $m^{di}$ of these. Let $q$ be the probability that the interior of some $Q'$ intersects $M$, given that $Q'$ is a kept cube at some level. Due to the independence in the construction, $q$ does not depend on $Q'$ and $q > 0$ since we assume $p > m^{-d}$.

*Claim:* Almost surely there exists an increasing sequence of natural numbers $\{k_i\}_i$ and an associated sequence of kept cubes $Q_{k_i}$ at level $k_i$ such that all cubes in $\mathcal{Q}_{k_i}^i$ are kept and, moreover, intersect $M$.

*Proof of claim:* Let $r, s \in \mathbb{N}$ be given and fix a kept cube $Q_r$ at level $r$. First we estimate the probability $p_s^*$ that all cubes in $\mathcal{Q}_r^s$ are kept and intersect $M$. By independence this probability is independent of $r$. In order for all cubes to be kept we must first have chosen to keep all cubes in $\mathcal{Q}_r^1$, and then all cubes in $\mathcal{Q}_r^2$, and so on until finally keeping all cubes in $\mathcal{Q}_r^s$. This requires us to choose to keep a total of

$$L \leqslant \sum_{a=1}^{s} m^{da} = \frac{m^{d(s+1)} - m^d}{m^d - 1}$$

cubes. Since these choices are independent, the probability of keeping them all is $p^L$. Just because we have kept these cubes at level $(s+r)$, they may go on to become extinct later in the construction. However, the probability that the interiors of all of the cubes in $\mathcal{Q}_r^s$ will go on to intersect $M$ in the limit is at least $q^{m^{ds}}$. Therefore, the probability we are looking for is

$$p_s^* \geqslant p^L q^{m^{ds}}.$$

Define $l(s)$ inductively by $l(s+1) = l(s) + k(s)$ with $l(1) = 1$ and

$$k(s) = s \left\lceil \frac{-\log 2}{\log(1 - p_s^*)} \right\rceil. \tag{9.6}$$

Let $\mathcal{E}_s$ be the event that for at least one $j \in \{0, s, 2s, \ldots, k(s) - s\}$ there exists a kept cube $Q_{l(s)+j}$ such that all cubes in $\mathcal{Q}_{l(s)+j}^s$ survive and intersect $M$. Applying (9.6),

$$\mathbb{P}(\mathcal{E}_s) \geqslant 1 - (1 - p_s^*)^{k(s)/s} \geqslant 1/2.$$

The claim follows immediately by the Borel-Cantelli Lemma, noting that the events $\mathcal{E}_s$ are independent. □

To complete the proof we show that, conditioned on $M$ being non-empty, almost surely $X = [0,1]^d$ is a weak tangent to $M$. The result then follows from Theorem 5.1.2. Let $T_i$ be the natural orientation preserving similarity that maps the cube $Q_{k_i}$ to $X$. By the claim we have that, almost surely, each of the level $k_i + i$ cubes inside $Q_{k_i}$ intersect $M$. Therefore

$$d_{\mathcal{H}}(T_i(M) \cap X, X) \leqslant \sqrt{d}\, m^{-i}$$

and so $d_{\mathcal{H}}(T_i(M) \cap X, X) \to 0$ as $i \to \infty$ as required. □

The Assouad spectrum was considered in [112, 266, 280] and the outcome is rather different. The following result was first proved in [112] but an alternative and more direct proof was given by Troscheit [266]. We follow Troscheit's proof below.

**Theorem 9.4.2** Conditioned on $M$ being non-empty, almost surely

$$\dim_{\mathrm{A}}^{\theta} M = \dim_{\mathrm{B}} M$$

for all $\theta \in (0,1)$.

Proving Theorem 9.4.2, following [266], hinges on the following lemma, which appears to be fundamental in the study of the Assouad spectrum and quasi-Assouad dimension in random settings, see [266, 268]. It says that the probability of the number of cubes at level $k$ being much larger than the expectation is very small (super-exponentially small in $k$), which allows us to control the number of descendants of a given node uniformly and with high probability.

**Lemma 9.4.3** There exists $\tau > 0$ such that for all $\varepsilon \in (0,1)$ and all $k \geqslant 1$

$$\mathbb{P}(Z_k > m^{sk(1+\varepsilon)}) \leqslant \exp(-\tau m^{sk\varepsilon}).$$

This lemma constitutes a 'large deviations' estimate for Mandelbrot percolation. It applies more generally to supercritical *Galton-Watson processes* with finitely supported offspring distribution. For us, $Z_k$ is a Galton-Watson process because the construction proceeds independently inside each kept cube and the supercritical condition translates to $p > m^{-d}$. Moreover, the finitely supported offspring condition is satisfied since $m$ is fixed. The more general result is proved in [9], see also [268].

Write $Z_k$ to denote the number of kept cubes comprising $M_k$ and note that

$$\mathbb{E}(Z_k) = (pm^d)^k = m^{sk}$$

where $s$ is the almost sure Hausdorff dimension of the limit set $M$.

*Proof of Theorem 9.4.2.* Fix $\theta \in (0,1)$ and $\varepsilon > 0$. Throughout we condition on non-extinction. In order to prove that almost surely $\dim_{\mathrm{A}}^{\theta} M \leqslant s + \varepsilon$, it is enough to prove that almost surely for all sufficiently large integers $k \geqslant 1$ and $l = \lceil k/\theta \rceil$ that

$$\sup_{Q \subseteq M_k} Z_k^*(Q) \leqslant m^{(l-k)(s+\varepsilon)}$$

where $Z_k^*(Q)$ is the number of level $l$ descendants of a given level $k$ cube $Q \subseteq M_k$. Note the difference between $Z_k$ and $Z_k^*(Q)$.

Since $\mathbb{E}(Z_k) = (pm^d)^k = m^{sk}$, almost surely there is a random constant $C_1 > 0$ such that

$$Z_k \leqslant C_1 m^{sk(1+\varepsilon)} \tag{9.7}$$

for all $k \geqslant 1$. Moreover, using the statistical self-similarity of the process and applying Lemma 9.4.3, for all level $k$ cubes $Q$,

$$\mathbb{P}(Z_k^*(Q) > m^{s(l-k)(1+\varepsilon)}) = \mathbb{P}(Z_{l-k} > m^{s(l-k)(1+\varepsilon)})$$

$$\leqslant C \exp(-\tau m^{s(l-k)\varepsilon}). \tag{9.8}$$

Therefore, combining (9.7) and (9.8),

$$\mathbb{P}\left(\sup_{Q \subseteq M_k} Z_k^*(Q) > m^{s(l-k)(1+\varepsilon)}\right)$$

$$\leqslant C_1 m^{sk(1+\varepsilon)} C \exp(-\tau m^{s(l-k)\varepsilon})$$

$$\leqslant C_1 C \exp\left(k \log m^{s(1+\varepsilon)} - \tau m^{sk(1/\theta - 1)\varepsilon}\right)$$

$$\leqslant C_1 C \exp\left(-(\tau/2) m^{sk(1/\theta-1)\varepsilon}\right) \tag{9.9}$$

for large enough $k$. Therefore

$$\sum_{k \geqslant 1} \mathbb{P}\left(\sup_{Q \subseteq M_k} Z_k^*(Q) > m^{(l-k)(s+\varepsilon)}\right) < \infty$$

and by the Borel-Cantelli lemma

$$\sup_{Q \subseteq M_k} Z_k^*(Q) \leqslant m^{(l-k)(s+\varepsilon)}$$

for all but finitely many $k$ almost surely, completing the proof. □

Theorems 9.4.1 and 9.4.2 show that genuine interpolation between the box dimension and Assouad dimension is *not* achieved by the Assouad spectrum for Mandelbrot percolation. However, using the finer analysis provided by the $\phi$-Assouad dimensions $\dim_{\mathrm{A}}^{\phi} M$, the interpolation can be recovered, recall Section 3.3.3. Troscheit proved the following dichotomy in [268].

**Theorem 9.4.4** If

$$\frac{\log(R/\phi(R))}{\log|\log R|} \to 0$$

as $R \to 0$, then, conditioned on $M$ being non-empty, almost surely

$$\dim_{\mathrm{A}}^{\phi} M = d = \dim_{\mathrm{A}} M.$$

On the other hand, if

$$\frac{\log(R/\phi(R))}{\log|\log R|} \to \infty$$

as $R \to 0$, then, conditioned on $M$ being non-empty, almost surely

$$\dim_{\mathrm{A}}^{\phi} M = \dim_{\mathrm{B}} M = d + \frac{\log p}{\log m}.$$

The proof follows a similar strategy to the proof of Theorem 9.4.2 — especially the second part. Indeed, an inspection of the proof reveals that (9.9) is a very crude estimate and would also be possible for $\phi$ with $\phi(R)$ much larger than $R^{1/\theta}$. Also, note that Theorem 9.4.4 implies Theorem 9.4.2 by considering $\phi(R) = R^{1/\theta}$ and Theorem 9.4.1 by considering $\phi(R) = R$. A similar dichotomy, with the same threshold on $\phi$, was obtained in a different random setting based on 'complementary intervals' in [121], see also [119]. The Assouad spectrum of random self-affine carpets was considered in [108].

# 10
# Geometric constructions

In this chapter we consider how the dimension theory we have developed behaves in the context of three classical constructions in geometric measure theory: products, projections, and slices. In Section 10.1 we examine the relationship between the dimensions of sets $E$, $F$ with the product set $E \times F$. In Section 10.2 we consider how the dimensions of a set behave under orthogonal projection onto linear subspaces. The typical behaviour of the Hausdorff dimension under orthogonal projection is described in the celebrated Marstrand-Mattila projection theorems (10.3) and we find that there is no such result describing the Assouad dimension, see Theorem 10.2.4. However, one can obtain a one-sided Marstrand-Mattila type result, see Theorem 10.2.1, which turns out to have a very small set of exceptions due to a result of Orponen, see Theorem 10.2.2. In Section 10.2.3 we provide an application of the Assouad dimension and spectra to the *box* dimensions of orthogonal projections, see Corollary 10.2.8. Finally, in Section 10.3 we consider how the dimensions of the slices of a set relate to the dimensions of the set itself. We construct an example of a set with maximal box (and Assouad) dimension which has very small slices, see Theorem 10.3.2.

## 10.1 Products

Given sets $E \subseteq \mathbb{R}^{d_1}$ and $F \subseteq \mathbb{R}^{d_2}$, it is a natural question to consider how the dimensions of $E$ and $F$ relate to the dimensions of the product set

$$E \times F = \{(x, y) : x \in E,\ y \in F\} \subseteq \mathbb{R}^{d_1} \times \mathbb{R}^{d_2} \equiv \mathbb{R}^{d_1+d_2}.$$

Consider, for example, $E = F = [0, 1]$, in which case $E \times F = [0, 1]^2$ and we have

$$\dim(E \times F) = 2 = \dim E + \dim F.$$

Perhaps we can expect this additive behaviour more generally? This turns out to depend crucially on the notion of dimension used and the structure of $E$ and $F$. The theory is well-understood and the general phenomenon is that dimensions behave in pairs. It is a good exercise to prove that, for bounded sets $E, F$,

$$\underline{\dim}_{\mathrm{B}} E + \underline{\dim}_{\mathrm{B}} F \leqslant \underline{\dim}_{\mathrm{B}}(E \times F) \leqslant \underline{\dim}_{\mathrm{B}} E + \overline{\dim}_{\mathrm{B}} F$$

and

$$\underline{\dim}_{\mathrm{B}} E + \overline{\dim}_{\mathrm{B}} F \leqslant \overline{\dim}_{\mathrm{B}}(E \times F) \leqslant \overline{\dim}_{\mathrm{B}} E + \overline{\dim}_{\mathrm{B}} F.$$

Here upper and lower box dimensions behave as a pair and this pattern repeats for other dimension pairs. Marstrand proved that $\dim_{\mathrm{H}} E + \dim_{\mathrm{H}} F \leqslant \dim_{\mathrm{H}}(E \times F)$ for Borel sets $E, F$ [202] and Howroyd proved that the Hausdorff and packing dimension are a natural pair in this setting [142]. That is, for Borel sets $E, F$,

$$\dim_{\mathrm{H}} E + \dim_{\mathrm{H}} F \leqslant \dim_{\mathrm{H}}(E \times F) \leqslant \dim_{\mathrm{H}} E + \dim_{\mathrm{P}} F$$

and

$$\dim_{\mathrm{H}} E + \dim_{\mathrm{P}} F \leqslant \dim_{\mathrm{P}}(E \times F) \leqslant \dim_{\mathrm{P}} E + \dim_{\mathrm{P}} F.$$

The lower and Assouad dimensions also behave as a pair, see [88]. However, it is more natural to use the *modified* lower dimension when bounding the Assouad dimension of a product set from below, see Corollary 10.1.2. We prove the analogous statement for measures, from which the statement for sets will (at least partially) follow. This is an example of how results for sets are often obtained as corollaries to stronger results for measures. Given locally finite measures $\mu$ supported on $E$ and $\nu$ supported on $F$, the *product measure* $\mu \times \nu$ is defined to be the unique Borel measure satisfying

$$(\mu \times \nu)(E' \times F') = \mu(E')\nu(F')$$

for all Borel sets $E' \subseteq E$ and $F' \subseteq F$.

**Theorem 10.1.1** Let $\mu$ and $\nu$ be locally finite Borel measures on $\mathbb{R}^{d_1}$ and $\mathbb{R}^{d_2}$, respectively. Then

$$\dim_{\mathrm{L}} \mu + \dim_{\mathrm{A}} \nu \leqslant \dim_{\mathrm{A}}(\mu \times \nu) \leqslant \dim_{\mathrm{A}} \mu + \dim_{\mathrm{A}} \nu$$

and, if $\mu$ and $\nu$ are either both compactly supported or both have unbounded support, then

$$\dim_{\mathrm{L}} \mu + \dim_{\mathrm{L}} \nu \leqslant \dim_{\mathrm{L}}(\mu \times \nu) \leqslant \dim_{\mathrm{L}} \mu + \dim_{\mathrm{A}} \nu.$$

*Proof* We prove the result for $\dim_{\mathrm{A}}(\mu \times \nu)$. The lower dimension analogue is proved similarly and left as an exercise. To see why the product of a compactly supported measure and a measure with unbounded support behaves differently, see Example 10.1.3.

The upper bound is straightforward. Let $\varepsilon > 0$. For all $(x, y) \in E \times F$ and $0 < r < R < 1$,

$$\begin{aligned}\frac{\mu \times \nu(B((x,y),R))}{\mu \times \nu(B((x,y),r))} &\leqslant \frac{\mu \times \nu(B(x,R) \times B(y,R))}{\mu \times \nu(B(x,r/\sqrt{2}) \times B(y,r/\sqrt{2}))} \\ &= \frac{\mu(B(x,R))\nu(B(y,R))}{\mu(B(x,r/\sqrt{2}))\nu(B(y,r/\sqrt{2}))} \\ &\leqslant C\left(\frac{R}{r}\right)^{\dim_{\mathrm{A}} \mu + \dim_{\mathrm{A}} \nu + \varepsilon}\end{aligned}$$

for a constant $C \geqslant 1$ independent of $x, y, r$ and $R$. This proves that

$$\dim_{\mathrm{A}}(\mu \times \nu) \leqslant \dim_{\mathrm{A}} \mu + \dim_{\mathrm{A}} \nu + \varepsilon$$

and letting $\varepsilon \to 0$ gives the desired upper bound.

Turning our attention to the lower bound, let $\varepsilon > 0$ and choose sequences $y_k \in F$, $r_k, R_k \in (0,1)$ satisfying

$$\frac{\nu(B(y_k, R_k/\sqrt{2}))}{\nu(B(y_k, r_k))} \geqslant \left(\frac{R_k}{r_k}\right)^{\dim_{\mathrm{A}} \nu - \varepsilon} \tag{10.1}$$

for all $k$ and $R_k/r_k \to \infty$. Such sequences exist by the definition of

$\dim_{\mathrm{A}} \nu$. For $x \in E$ arbitrary,

$$\frac{\mu \times \nu(B((x, y_k), R_k))}{\mu \times \nu(B((x, y_k), r_k))} \geqslant \frac{\mu \times \nu(B(x, R_k/\sqrt{2}) \times B(y_k, R_k/\sqrt{2}))}{\mu \times \nu(B(x, r_k) \times B(y_k, r_k))}$$
$$= \frac{\mu(B(x, R_k/\sqrt{2}))\nu(B(y_k, R_k/\sqrt{2}))}{\mu(B(x, r_k))\nu(B(y_k, r_k))}$$
$$\geqslant C' \left(\frac{R_k}{r_k}\right)^{\dim_{\mathrm{L}} \mu - \varepsilon + \dim_{\mathrm{A}} \nu - \varepsilon}$$

for a constant $C' > 0$ independent of $x$ and $k$. The final inequality above follows from (10.1) and the definition of $\dim_{\mathrm{L}} \mu$. This proves that

$$\dim_{\mathrm{A}}(\mu \times \nu) \geqslant \dim_{\mathrm{L}} \mu + \dim_{\mathrm{A}} \nu - 2\varepsilon$$

and letting $\varepsilon \to 0$ gives the desired lower bound. □

We obtain the analogous result for sets as a corollary. The fact that the Assouad dimension is subadditive on products, that is, $\dim_{\mathrm{A}}(E \times F) \leqslant \dim_{\mathrm{A}} E + \dim_{\mathrm{A}} F$ was observed in [192, Theorem A.5], see also [240, Lemma 9] and [88, Proposition 2.1]. Sharpness of these product formulae was considered in [224, 223].

**Corollary 10.1.2** For non-empty sets $E, F$

$$\dim_{\mathrm{ML}} E + \dim_{\mathrm{A}} F \leqslant \dim_{\mathrm{A}}(E \times F) \leqslant \dim_{\mathrm{A}} E + \dim_{\mathrm{A}} F$$

and, if $E$ and $F$ are either both bounded or both unbounded, then

$$\dim_{\mathrm{L}} E + \dim_{\mathrm{L}} F \leqslant \dim_{\mathrm{L}}(E \times F) \leqslant \dim_{\mathrm{L}} E + \dim_{\mathrm{A}} F.$$

*Proof* Again, we only consider the statement for $\dim_{\mathrm{A}}(E \times F)$. The upper bound follows immediately from Theorem 10.1.1 and Theorem 4.1.3 by choosing a sequence of measures $\mu$ and $\nu$ supported on $E$ and $F$, respectively, whose Assouad dimensions approach the Assouad dimensions of $E$ and $F$. For the lower bound, we first obtain the statement with $\dim_{\mathrm{ML}} E$ replaced by $\dim_{\mathrm{L}} E$ by a similar approach to that used in the proof of the lower bound in Theorem 10.1.1. Let $\varepsilon > 0$ and choose sequences $y_k \in F$, $r_k, R_k \in (0, 1)$ satisfying

$$M_{r_k}(B(y_k, R_k/\sqrt{2}) \cap F) \geqslant \left(\frac{R_k}{r_k}\right)^{\dim_{\mathrm{A}} F - \varepsilon} \tag{10.2}$$

for all $k$ and $R_k/r_k \to \infty$, where $M_{r_k}(\cdot)$ is the maximum cardinality of an $r_k$-separated subset of a given set. Such sequences exist by the definition of $\dim_{\mathrm{A}} F$. For $x \in E$ arbitrary,

$$\begin{aligned} &M_{r_k}(B((x, y_k), R_k) \cap (E \times F)) \\ &\qquad \geqslant M_{r_k}\left( \left(B(x, R_k/\sqrt{2}) \cap E\right) \times \left(B(y_k, R_k/\sqrt{2}) \cap F\right) \right) \\ &\qquad \geqslant M_{r_k}(B(x, R_k/\sqrt{2}) \cap E) M_{r_k}(B(y_k, R_k/\sqrt{2}) \cap F) \\ &\qquad \geqslant C'' \left(\frac{R_k}{r_k}\right)^{\dim_{\mathrm{L}} E - \varepsilon + \dim_{\mathrm{A}} F - \varepsilon} \end{aligned}$$

for a constant $C'' > 0$ independent of $x$ and $k$, using (10.2) and the definition of $\dim_{\mathrm{L}} E$. This proves that $\dim_{\mathrm{A}}(E \times F) \geqslant \dim_{\mathrm{L}} E + \dim_{\mathrm{A}} F - 2\varepsilon$ and letting $\varepsilon \to 0$ gives the desired lower bound. Finally, we can upgrade this lower bound by replacing $\dim_{\mathrm{L}} E$ by $\dim_{\mathrm{ML}} E$ by passing to a subset of $E$ (and therefore a subset of $E \times F$) with lower dimension close to $\dim_{\mathrm{ML}} E$. □

The following example shows that if we take the product of a bounded set with an unbounded set, then the bounds for lower dimension from Theorem 10.1.1 and Corollary 10.1.2 do not hold. The only inequality which can go wrong is the first one.

**Example 10.1.3** There exist sets $E, F \subseteq \mathbb{R}$ such that

$$\dim_{\mathrm{L}} E + \dim_{\mathrm{L}} F > \dim_{\mathrm{L}}(E \times F)$$

and measures $\mu, \nu$ on $\mathbb{R}$ such that

$$\dim_{\mathrm{L}} \mu + \dim_{\mathrm{L}} \nu > \dim_{\mathrm{L}}(\mu \times \nu).$$

For example, let $E = [0, 1]$ and $F = \mathbb{R}$, which both have lower dimension 1. However, considering $r = 1$ and $R \to \infty$ to obtain the upper bound,

$$\dim_{\mathrm{L}}(E \times F) = 1.$$

For an example with measures, let $\mu$ and $\nu$ be Lebesgue measure on $E$ and $F$ respectively.

The product formulae presented above are useful in many contexts, particularly when constructing examples. We give an application of this here, which also relies on an understanding of self-similar sets — another very useful collection of objects for constructing examples! Specifically, we will prove Theorem 3.4.12 (stated earlier on page 49) which establishes that the Assouad spectrum can increase under Lipschitz maps.

*Proof of Theorem 3.4.12.* Given $s \in (0,1)$, we construct a compact set $F \subseteq \mathbb{R}^2$ such that for all $\theta \in (0,1)$

$$\dim_{\mathrm{A}} F = \overline{\dim}_{\mathrm{B}} F = \dim_{\mathrm{A}}^{\theta} F = s$$

but

$$\dim_{\mathrm{A}}^{\theta} \pi F = \min\left\{\frac{s}{1-\theta}, 1\right\}$$

where $\pi$ is projection onto the first coordinate. Let $c \in (0, 1/2)$ be such that $-\log 2/\log c = s$ and, for each integer $k \geqslant 1$, let

$$F'_k = \{S_{\boldsymbol{i}}(0) : \boldsymbol{i} \in \{1,2\}^k\}$$

where $S_1 : x \mapsto cx$ and $S_2 : x \mapsto cx + (1-c)$. The IFS $\{S_1, S_2\}$ generates a self-similar set $E \subseteq \mathbb{R}$ satisfying the SSC which has Assouad (and box and Hausdorff) dimension equal to $s$. Let $F_k \subseteq \mathbb{R}^2$ be a set consisting of $2^k$ points such that $\pi'(F_k) = F'_k$ and

$$\pi(F_k) = \{jc^k : j = 0, \ldots, 2^k - 1\}$$

where $\pi'$ is projection onto the second coordinate and $\pi$ is projection onto the first coordinate. The set $\pi(F_k)$ lies inside an interval of length $2^k c^k = c^{k(1-s)}$. There are many ways to choose $F_k$ and so, for concreteness, we assume it is the graph of an increasing function. Finally, let $t \in (2c, 1)$ and

$$F = \bigcup_{k \geqslant k_0} (F_k + t^k)$$

where $k_0$ is chosen such that $k \geqslant k_0$ guarantees $t^k - t^{k+1} > (2c)^{k+1}$, see Figure 10.1. Starting at $k_0$ ensures the pieces $(F_k + t^k)$ do not overlap. We abuse notation slightly here writing $F_k + t^k$ to denote $F_k$ translated in the horizontal direction by $t^k$.

The set $F$ satisfies the required properties, apart from being compact. However, we may obtain a compact set by taking the closure without affecting any of the dimensions we consider. In order to prove the result it suffices to show

(a) $\dim_{\mathrm{A}} F \leqslant s$
(b) $\dim_{\mathrm{A}}^{1-s} \pi F \geqslant 1$

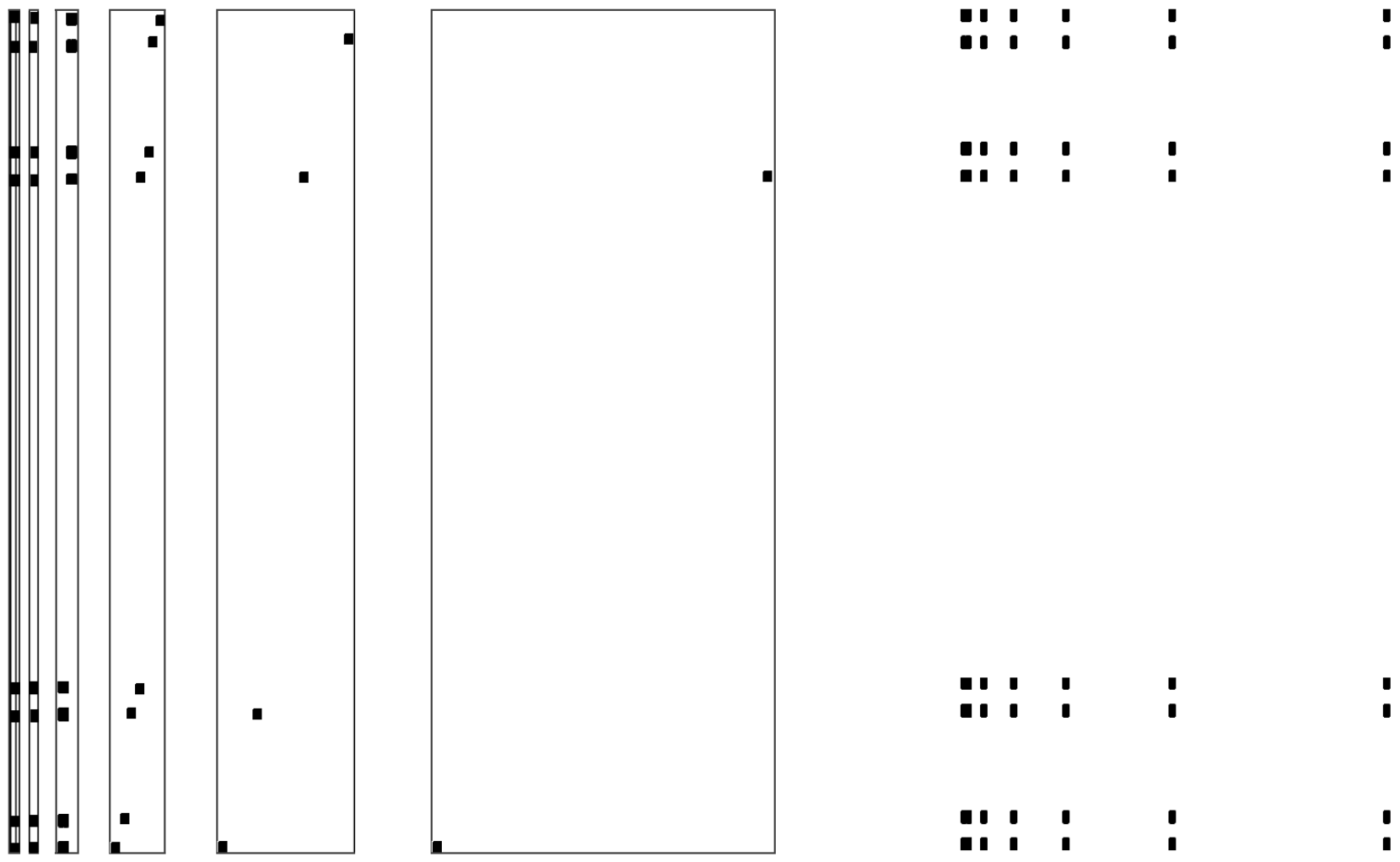

Figure 10.1 Left: the construction of $F$ with $c = 1/5$ and $t = 1/2$. Rectangles containing the finite sets $F_k + t^k$ are shown along with the points comprising $F_k + t^k$. Right: the associated set $Z \times E$, which contains a subset which is bi-Lipschitz equivalent to $F$.

and apply Lemma 3.4.4 and Corollary 3.4.6 to both $F$ and $\pi F$, recalling that the upper box dimension cannot increase under Lipschitz maps and so $\overline{\dim}_{\mathrm{B}} \pi F \leqslant \overline{\dim}_{\mathrm{B}} F \leqslant \dim_{\mathrm{A}} F$.

To prove (a), let $F_0$ be the set formed from $F$ by projecting each piece $F_k + t^k$ onto the vertical line passing through $(t^k, 0)$. The induced bijection between $F$ and $F_0$ is bi-Lipschitz, since the horizontal spacing of points in $F_k + t^k$ is $c^k$ and the vertical spacing is (at least) $c^k$. Moreover, $F_0 \subseteq Z \times E$ where $Z = \{t^k : k \geqslant 1\}$ and $E$ is the self-similar set mentioned above. Therefore, using Theorem 10.1.2 and Lemma 2.4.2,

$$\dim_{\mathrm{A}} F = \dim_{\mathrm{A}} F_0 \leqslant \dim_{\mathrm{A}} Z \times E \leqslant \dim_{\mathrm{A}} Z + \dim_{\mathrm{A}} E = s.$$

It is a short exercise to show that $\dim_{\mathrm{A}} Z = 0$, which we leave to the reader.

To prove (b), let $x = t^k \in \pi F$ and $r = c^k$. Since $r^{1-s} = c^{k(1-s)} = 2^k c^k$, $\pi(F_k + t^k) \subseteq B(x, r^{1-s})$ and therefore

$$N_r(B(x, r^{1-s}) \cap \pi F) \geqslant 2^{k-1} = \tfrac{1}{2} \frac{r^{1-s}}{r}.$$

Since we can take $k$ arbitrarily large, we conclude $\dim_{\mathrm{A}}^{1-s} \pi F \geqslant 1$. □

## 10.2 Orthogonal projections

### 10.2.1 Dimension theory of orthogonal projections

How dimension behaves under orthogonal projection is a classical problem in geometric measure theory. It was first considered by Marstrand in a seminal article from 1954 [201], although Besicovitch considered the Lebesgue measure of projections of planar sets with positive and finite 1-dimensional Hausdorff measure somewhat earlier, see for example [31]. Consider, for example, a line segment $L$ embedded in $\mathbb{R}^3$. Then for all but one 2-dimensional hyperplane $V \subseteq \mathbb{R}^3$, the orthogonal projection of $L$ onto $V$ is also a line segment and therefore the dimension of $L$ (which is 1) is typically preserved under projection onto 2-dimensional subspaces. Conversely, consider the unit sphere $\mathbb{S}^2 \subseteq \mathbb{R}^3$ (which is a 2-dimensional set). Then for all 1-dimensional subspaces $V \subseteq \mathbb{R}^3$, the projection of $\mathbb{S}^2$ onto $V$ is a line segment and so the dimension of $\mathbb{S}^2$ is not typically preserved under projection onto 1-dimensional subspaces, but does typically attain the largest value possible. Despite the simplicity of these examples, these phenomena occur much more generally.

For integers $k, d$ with $1 \leqslant k < d$, consider projections of $\mathbb{R}^d$ onto $k$-dimensional subspaces $V$. The $k$-dimensional subspaces of $\mathbb{R}^d$ come with a natural $k(d-k)$ dimensional manifold structure, and so come equipped with a $k(d-k)$ dimensional analogue of Lebesgue measure (this may be formally expressed in terms of the Haar measure on the $d$-dimensional orthogonal group $O(d)$). We identify each $k$-dimensional subspace $V$ with the orthogonal projection $\pi : \mathbb{R}^d \to V$ and denote the set of all such projections as $G(k, d)$. We can thus make statements about 'almost all' orthogonal projections $\pi \in G(k, d)$. The manifold $G(k, d)$ is usually called the *Grassmannian manifold*, see [205, Chapter 3]. We are abusing notation slightly, writing $G(k, d)$ for both the collection of $k$-dimensional subspaces of $\mathbb{R}^d$ and for the space of associated projections onto these subspaces.

The classical result concerning Hausdorff dimension (often referred to as *Marstrand's projection theorem* or the *Marstrand-Mattila projection theorem*), is that, for a fixed Borel set $F \subseteq \mathbb{R}^d$, almost all $\pi \in G(k, d)$ satisfy

$$\dim_{\mathrm{H}} \pi F \;=\; \min\{k, \dim_{\mathrm{H}} F\}. \tag{10.3}$$

This was first proved by Marstrand in the case $d = 2$ [201] and later by Mattila [204] in the general case, based on a potential-theoretic approach of Kaufman [160]. Similar results exist for the packing and upper and lower box dimensions in that for almost all $\pi \in G(k,d)$ the dimension of $\pi F$ is equal to a constant. These results were proved by Falconer and Howroyd [78]. Unlike the Hausdorff dimension result, the constants are not as simple as $\min\{k, \dim F\}$ and are given by *dimension profiles*. Dimension profiles were introduced by Howroyd [143] and Falconer-Howroyd [78] in the 1990s. We refer the reader to the survey articles [74, 206] for an overview of the dimension theory of projections.

### 10.2.2 The Assouad dimension of orthogonal projections

Compared to the Hausdorff, box and packing dimensions, a striking difference in the case of Assouad dimension is that the dimension of $\pi F$ need not be almost surely constant. This was first observed by Fraser and Orponen [106, Theorem 2.5] who constructed a compact set in the plane which projects to sets with two different Assouad dimensions in positively many directions $\pi \in G(1,2)$. We recall the construction below. Despite the Assouad dimension not necessarily taking a constant value almost surely, it turns out that it can be almost surely bounded from below in a meaningful way.

**Theorem 10.2.1** Let $F \subseteq \mathbb{R}^d$ be any non-empty set and $k < d$ be a positive integer. Then

$$\dim_{\mathrm{A}} \pi F \geqslant \min\{k, \dim_{\mathrm{A}} F\}$$

for almost all $\pi \in G(k,d)$.

The case when $k = 1$ and $d = 2$ was proved by Fraser and Orponen, see [106, Theorem 2.1], and the general case was proved by Fraser [89]. Interestingly, Theorem 10.2.1 is false if one replaces Assouad dimension by upper box dimension. The upper box dimension of $\pi F$ is almost surely constant, but this constant can be strictly smaller than the minimum of $k$ and the upper box dimension of $F$, see [77] or the surveys [74, 206].

Theorem 10.2.1 establishes that for a nonempty set $F \subseteq \mathbb{R}^d$, the set

$$\{\pi \in G(k,d) : \dim_{\mathrm{A}} \pi F < \min\{k, \dim_{\mathrm{A}} F\}\} \tag{10.4}$$

is a null set with respect to the Grassmannian measure. One can ask for more precise information here by, for example, estimating the Hausdorff dimension of the *exceptional set* (10.4). This problem has attracted much interest in the context of Marstrand's projection theorem for Hausdorff dimension, see [74, 206]. For example, given a Borel set $F \subseteq \mathbb{R}^d$ with

$\dim_{\mathrm{H}} F \leqslant k$ and $0 < s \leqslant \dim_{\mathrm{H}} F$,

$$\dim_{\mathrm{H}}\{\pi \in G(k,d) : \dim_{\mathrm{H}} \pi F < s\} \leqslant k(d-k) + s - k.$$

Moreover, Kaufman and Mattila [162] constructed a compact set $F \subseteq \mathbb{R}^d$ with any given Hausdorff dimension in $[0,k]$ such that

$$\dim_{\mathrm{H}}\{\pi \in G(k,d) : \dim_{\mathrm{H}} \pi F < \dim_{\mathrm{H}} F\} = k(d-k) + \dim_{\mathrm{H}} F - k$$

exhibiting sharpness of the estimate above in the case of largest $s$. In particular, the set of exceptions can be large in the sense of dimension, despite being a null set. Orponen [226] established a striking result concerning the exceptional set of dimensions in the context of Assouad dimension.

**Theorem 10.2.2** Let $F \subseteq \mathbb{R}^2$ be any non-empty set. Then

$$\dim_{\mathrm{H}}\{\pi \in G(1,2) : \dim_{\mathrm{A}} \pi F < \min\{1, \dim_{\mathrm{A}} F\}\} = 0.$$

We are able to give a full and self-contained proof of Theorem 10.2.1 below. Theorem 10.2.2 is rather harder to prove and relies on deep modern advances, including Shmerkin's inverse theorem for $L^q$-norms [254].

*Proof of Theorem 10.2.1:* We follow the proof from [89, Section 3.3.1]. For all $\pi \in G(k,d)$, $\pi(\overline{F}) \subseteq \overline{\pi(F)}$ and therefore we may assume that $F$ is closed. It follows from Theorem 5.1.3 that there exists a compact set $E \subseteq \mathbb{R}^d$ which is a weak tangent to $F$ such that $\dim_{\mathrm{H}} E = \dim_{\mathrm{A}} F$. Examining the proof of Theorem 5.1.3, we may assume that the similarity maps generating $E$ are homotheties. Therefore we may choose a sequence of homothetic similarity maps $T_k$ on $\mathbb{R}^d$ and a compact set $X \subseteq \mathbb{R}^d$ such that

$$T_k(F) \cap X \to E \tag{10.5}$$

in $d_{\mathcal{H}}$ as $k \to \infty$. For $k \geqslant 1$, write $T_k(x) = c_k x + t_k$ for a real constant $c_k > 0$ and a translation $t_k \in \mathbb{R}^d$.

In the following, we identify $\pi(\mathbb{R}^d)$ with $\mathbb{R}^k$ for each $\pi \in G(k,d)$. Define a map $\pi \circ T_k \circ \pi^{-1}$ from $\pi(\mathbb{R}^d)$ to itself by

$$\{\pi \circ T_k \circ \pi^{-1}(x)\} = \{\pi(T_k(y)) : \pi(y) = x\}.$$

Since $T_k$ is homothetic the set $\{\pi(T_k(y)) : \pi(y) = x\}$ is a singleton and so this map is well-defined. For $x \in \pi(\mathbb{R}^d)$,

$$\pi \circ T_k \circ \pi^{-1}(x) = \pi(c_k \pi^{-1}(x) + t_k) = c_k x + \pi(t_k)$$

and so $\pi \circ T_k \circ \pi^{-1}$ is also a homothetic similarity. Since $\pi : \mathcal{K}(\mathbb{R}^d) \to \mathcal{K}(\mathbb{R}^k)$ is continuous, it follows from (10.5) that

$$\Big( \left(\pi \circ T_k \circ \pi^{-1}\right) (\pi(F)) \Big) \cap \pi(X) \;=\; \pi(T_k(F)) \cap \pi(X) \supseteq \pi\left(T_k F \cap X\right) \to \pi(E) \qquad (10.6)$$

in $d_{\mathcal{H}}$ as $k \to \infty$. Note that $\pi(X)$ is a compact subset of $\pi(\mathbb{R}^d)$ and $\mathcal{K}(\pi(X))$ is compact and so we may assume, by taking a subsequence if necessary, that $(\pi \circ T_k \circ \pi^{-1})(\pi(F)) \cap \pi(X)$ converges to a compact set $E' \subseteq \pi(X)$ in $d_{\mathcal{H}}$ as $k \to \infty$. In particular, $E'$ is a weak tangent to $\pi(F)$ and it follows from (10.6) that $E' \supseteq \pi(E)$.

We have established that for *all* $\pi \in G(k,d)$, the set $\pi(E)$ is a subset of a weak tangent to $\pi(F)$. It therefore follows from Theorem 5.1.2 and the Marstrand-Mattila projection theorem that, for *almost all* $\pi \in G(k,d)$,

$$\begin{aligned} \dim_{\mathrm{A}} \pi(F) \geqslant \dim_{\mathrm{A}} \pi(E) \geqslant \dim_{\mathrm{H}} \pi(E) &= \min\{k, \dim_{\mathrm{H}} E\} \\ &= \min\{k, \dim_{\mathrm{A}} F\} \end{aligned}$$

completing the proof of Theorem 10.2.1. □

In order to construct a set $F$ such that $\dim_{\mathrm{A}} \pi F$ is not almost surely constant, we first consider projections of self-similar sets, recall Chapter 7. The following is a special case of [106, Theorems 2.3-2.4].

**Theorem 10.2.3** Let $F \subseteq \mathbb{R}^2$ be a self-similar set containing at least two points and suppose the defining maps are homotheties. Then, for all $\pi \in G(1,2)$,

(i) if $\mathcal{H}^{\dim_{\mathrm{H}} \pi F}(\pi F) > 0$, then $\dim_{\mathrm{A}} \pi F \;=\; \dim_{\mathrm{H}} \pi F$,
(ii) if $\mathcal{H}^{\dim_{\mathrm{H}} \pi F}(\pi F) = 0$, then $\dim_{\mathrm{A}} \pi F \;=\; 1$.

Moreover, provided $F$ is not contained in a line, there is at least one $\pi \in G(1,2)$ such that $\dim_{\mathrm{A}} \pi F = 1$.

*Proof* The first part of this result follows from Theorems 7.2.4 and 7.2.5 since the projections $\pi F$ are self-similar. This relies on the fact that the maps are homotheties. The observation that there is at least one $\pi$ such that $\dim_{\mathrm{A}} \pi F = 1$ was proved in [106, Theorem 2.4] using a variant of the example of Bandt and Graf, see (7.3). □

Theorem 10.2.3 provides a family of examples where $\dim_{\mathrm{A}} \pi F$ is not almost surely constant. The following was proved in [106, Theorem 2.5], and we sketch their argument below.

**Theorem 10.2.4** For any $s$ satisfying $\log_5 3 < s < 1$, there exists a compact set $F \subseteq \mathbb{R}^2$ with Hausdorff and Assouad dimensions equal to $s$ for which there are two non-empty disjoint open sets $I, J \subseteq G(1,2)$ such that:

(i) $\dim_A \pi F = s$ for all $\pi \in I$
(ii) $\dim_A \pi F = 1$ for almost all $\pi \in J$.

*Proof* Fix $c \in (1/5, 1/3)$ and let $F$ be the self-similar attractor of the IFS on $[0,1]^2$ given by $\{x \mapsto cx,\ x \mapsto cx+(0,1-c),\ x \mapsto cx+(1-c,0)\}$, see Figure 10.2. Observe that $F$ satisfies the OSC and so $\dim_H F = \dim_A F = -\log 3/\log c =: s \in (\log_5 3, 1)$. It follows from Theorem 10.2.3 that:

(i) $\dim_A \pi F = \dim_H \pi F$ if and only if $\mathcal{H}^{\dim_H \pi F}(\pi F) > 0$
(ii) $\dim_A \pi F = 1$ if and only if $\mathcal{H}^{\dim_H \pi F}(\pi F) = 0$.

Marstrand's projection theorem guarantees that $\dim_H \pi F = s < 1$ almost surely. Therefore, in light of the above dichotomy, it is sufficient to show that there is a non-empty open interval[1] of projections $\pi$ for which $\mathcal{H}^s(\pi F) > 0$, and a non-empty open interval of projections $\pi$ within which $\mathcal{H}^s(\pi F) = 0$ almost surely. The first of these tasks is straightforward because one can easily find a non-empty open interval of projections $\pi$ for which $\pi F$ is a self-similar set satisfying the OSC and therefore $\dim_H \pi F = s$ and $\mathcal{H}^s(\pi F) > 0$, see Section 6.4. The existence of such an interval relies on the assumption $c < 1/3$ since (after rescaling) the problem reduces to positioning three pairwise disjoint closed intervals of length $c$ inside the unit interval.

The second task is more delicate, but was solved by Peres, Simon and Solomyak [232]. They defined the set of *intersection parameters* $\mathcal{IP} = \{\pi : \pi \text{ is not injective on } F\}$ and proved that $\mathcal{IP}$ contains a non-empty interval provided $c > 1/5$ and, moreover, for almost every $\pi \in \mathcal{IP}$ we have $\mathcal{H}^s(\pi F) = 0$. This can be found in [232, Theorem 1.2(i) and Example 2.8], see Figure 10.2. □

The examples in Theorem 10.2.4 were the first examples demonstrating that the Assouad dimension need not be almost surely constant under projection. However, it is conjectured (but not proved, see Theorem 7.3.1 and Question 17.5.3), that for these examples, the quasi-Assouad dimension and Assouad spectrum *are* almost surely constant, since the projections are self-similar.

[1] We identify $G(1,2)$ with $[0,1)$ so we may talk about 'intervals' of projections $\pi$.

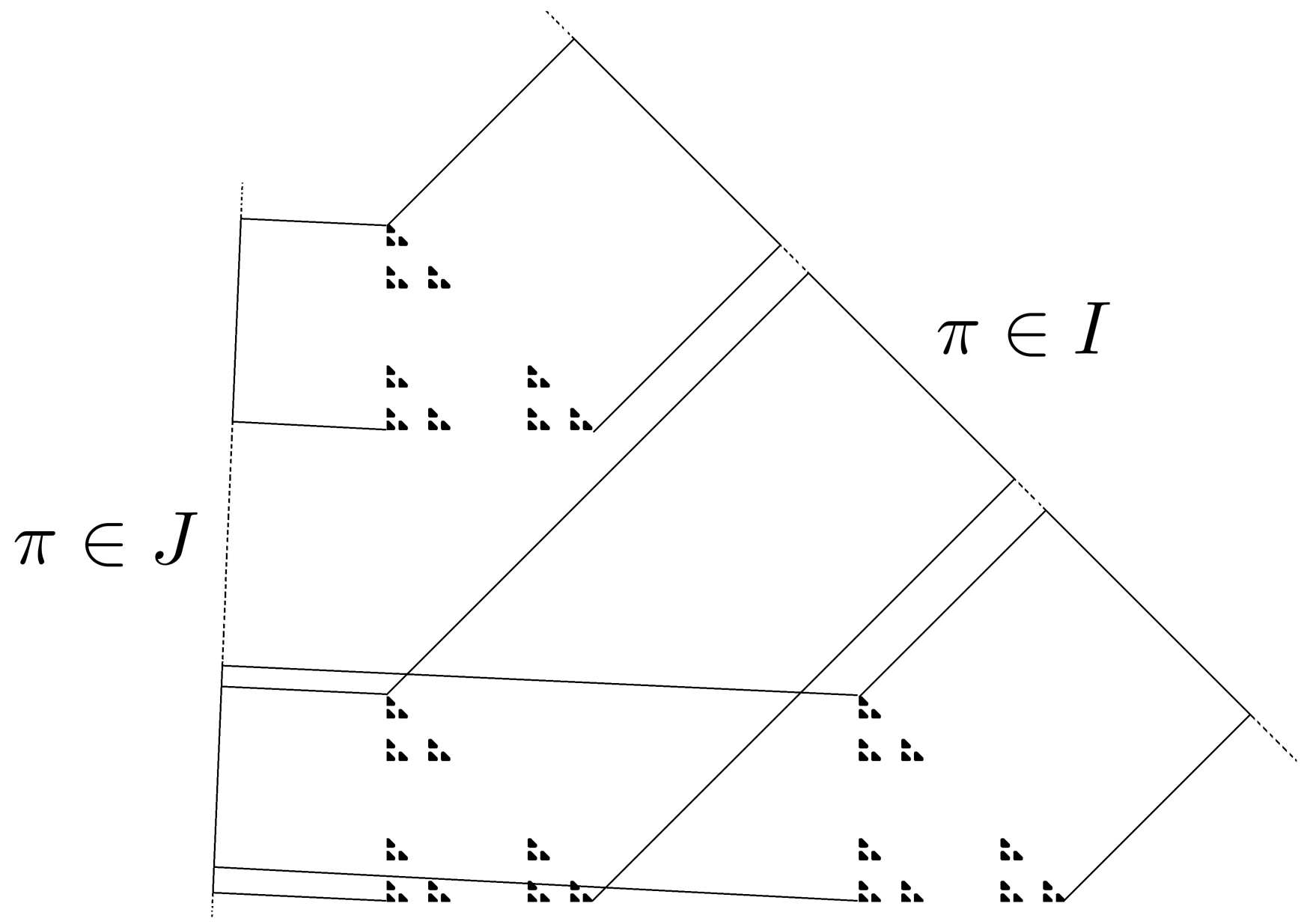

Figure 10.2 The set $F$ (with $c = 3/10$) and two projections. There is a small interval of projections $I$ for which the three scaled down copies of $F$ project into pairwise disjoint intervals, meaning that projections in these directions satisfy the OSC. There is also an interval $J$ within which the projection map is not injective.

In fact, the behaviour of the Assouad dimension under orthogonal projection can get *much* more complicated than simply taking on two distinct values with positive measure. The following construction appears in [103, Theorem 2.1].

**Theorem 10.2.5** Given an upper semi-continuous function $\phi : G(1, 2) \to [0, 1]$, there exists a compact set $F$ in the plane with $\dim_{\mathrm{A}} F = 0$ such that $\dim_{\mathrm{A}} \pi F = \phi(\pi)$ for all $\pi \in G(1, 2)$.

This result proves that $\dim_{\mathrm{A}} \pi F$ can take on any finite or countable number of distinct values with positive measure. Moreover, $\dim_{\mathrm{A}} \pi F$ can avoid all values almost surely. Again, the construction in [103] does not yield almost sure non-constancy for the quasi-Assouad dimension or Assouad spectrum.

### 10.2.3 An application to the box dimensions of projections

As mentioned above, there is a Marstrand-type theorem for the box dimensions. That is, for a bounded set $F \subseteq \mathbb{R}^d$, the dimensions $\overline{\dim}_{\mathrm{B}} \pi F$ take a common value for almost all $\pi \in G(k,d)$ with a similar statement for the lower box dimensions. The common value is given by a *dimension profile*, which is a little more complicated than the almost sure constant in Marstrand's theorem, which is simply $\min\{k, \dim_{\mathrm{H}} F\}$. Recall, that dimension profiles first appeared in [78, 143] and the almost sure constancy for the box dimensions of projections was established in [143].

A capacity approach to dimension profiles was introduced much more recently in [73]. Given the kernel

$$\phi_r^k(x) = \min\left\{1, \left(\frac{r}{|x|}\right)^k\right\} \tag{10.7}$$

the *capacity* of a bounded set $F$ with respect to this kernel is

$$\frac{1}{C_r^k(F)} = \inf_{\mu \in \mathcal{M}(F)} \int \int \phi_r^k(x-y) d\mu(x) d\mu(y).$$

Here $\mathcal{M}(F)$ denotes the collection of Borel probability measures supported by $F$. The double integral inside the infimum is called the *energy* of $\mu$ and the inner integral is called the *potential* of $\mu$ at $y$.

It was shown in [73] that the box dimension profiles are given by

$$\underline{\dim}_{\mathrm{B}}^k F = \liminf_{r \to 0} \frac{\log C_r^k(F)}{-\log r}, \quad \overline{\dim}_{\mathrm{B}}^k F = \limsup_{r \to 0} \frac{\log C_r^k(F)}{-\log r}$$

and, as mentioned above, these give the almost sure box dimensions of projections, see [73, Theorem 1.1].

**Theorem 10.2.6** Let $F \subseteq \mathbb{R}^d$ be a non-empty compact set. Then for almost all $\pi \in G(k,d)$

$$\underline{\dim}_{\mathrm{B}} \pi F = \underline{\dim}_{\mathrm{B}}^k F \quad \text{and} \quad \overline{\dim}_{\mathrm{B}} \pi F = \overline{\dim}_{\mathrm{B}}^k F.$$

Interestingly the dimension profiles are not necessarily given by a simple expression analogous to the Hausdorff dimension case. However, they always satisfy

$$\frac{\overline{\dim}_{\mathrm{B}} F}{1 + (1/k - 1/d)\overline{\dim}_{\mathrm{B}} F} \leqslant \overline{\dim}_{\mathrm{B}}^k F \leqslant \min\{k, \overline{\dim}_{\mathrm{B}} F\} \tag{10.8}$$

with the same bounds holding for $\underline{\dim}_{\mathrm{B}}^k F$ with lower box dimension replacing upper box dimension. Moreover, these bounds are the best possible, as shown in [77, 146].

Work of Falconer, Fraser and Shmerkin [76] demonstrated that the Assouad spectrum can be used to estimate the box dimension profiles and this leads to the following application in the context of the box dimensions of projections. More transparent corollaries will follow. The result is more efficiently expressed in terms of the upper Assouad spectrum, see Section 3.3.2.

**Theorem 10.2.7** Let $1 \leqslant k < d$, $F \subseteq \mathbb{R}^d$ bounded and $\theta \in (0,1)$. Then, for almost all $\pi \in G(k,d)$,

$$\underline{\dim}_{\mathrm{B}} \pi F \geqslant \underline{\dim}_{\mathrm{B}} F - \max\{0,\ \overline{\dim}_{\mathrm{A}}^{\theta} F - k,\ (\dim_{\mathrm{A}} F - k)(1-\theta)\}$$

and

$$\overline{\dim}_{\mathrm{B}} \pi F \geqslant \overline{\dim}_{\mathrm{B}} F - \max\{0,\ \overline{\dim}_{\mathrm{A}}^{\theta} F - k,\ (\dim_{\mathrm{A}} F - k)(1-\theta)\}.$$

*Proof* We follow the proof in [76]. By rescaling if necessary, we may assume that $|F| < 1/2$. Fix $\theta \in (0,1)$, $s > \overline{\dim}_{\mathrm{A}}^{\theta} F$, $t > \dim_{\mathrm{A}} F$ and let $C > 0$ be a constant such that, for all $0 < r < R < 1$ and $x \in F$,

$$N_r(B(x,R) \cap F) \leqslant C \left(\frac{R}{r}\right)^t$$

and, for all $0 < r \leqslant R^{1/\theta} < R < 1$ and $x \in F$,

$$N_r(B(x,R) \cap F) \leqslant C \left(\frac{R}{r}\right)^s.$$

Also let $u < \underline{\dim}_{\mathrm{B}} F$. Let $0 < r < 1$ be small and let $\{x_i\}_{i=1}^M$ be a maximal $r$-separated set of points in $F$. For small enough $r$, we can guarantee $M \geqslant r^{-u}$ by applying the definition of lower box dimension. Place a point mass of weight $1/M$ at each $x_i$ and let $\mu$ be the sum of these point masses, that is,

$$\mu = \sum_{i=1}^{M} \frac{1}{M} \delta_{x_i}.$$

Write $A = \lceil \log_2(2|F|r^{-1}) \rceil$ and $a = \lceil (1-\theta)\log_2(r^{-1}) \rceil$ and observe that,

for sufficiently small $r$, $a < A$. The potential of $\mu$ at $x_i$ is

$$\begin{aligned}\int \phi_r^k(x_i - y)d\mu(y) &\leqslant \sum_{n=0}^{A} 2^{-(n-1)k}\mu(B(x_i, 2^n r)) \qquad \text{(applying (10.7))}\\ &\leqslant \sum_{n=0}^{A} 2^{-(n-1)k}\frac{1}{M}N_r\big(B(x_i, 2^n r)\cap F\big)\\ &\leqslant 2^k r^u \Bigg(\sum_{n=a}^{A} 2^{-nk}C\Big(\frac{2^n r}{r}\Big)^s\\ &\qquad\qquad + \sum_{n=0}^{a-1} 2^{-nk}C\Big(\frac{2^n r}{r}\Big)^t\Bigg)\\ &\leqslant c\, r^u \max\Big\{1,\ r^{-(s-k)},\ r^{-(t-k)(1-\theta)}\Big\}\end{aligned}$$

where $c > 0$ is a constant independent of $r$. Summing over $i$, the energy of $\mu$ can be bounded as

$$\int\int \phi_r^k(x-y)d\mu(x)d\mu(y) \leqslant c\, r^u \max\Big\{1,\ r^{-(s-k)},\ r^{-(t-k)(1-\theta)}\Big\}$$

and therefore the capacity $C_r^k(F)$ satisfies

$$C_r^k(F) \geqslant c^{-1}r^{-u}\min\Big\{1,\ r^{(s-k)},\ r^{(t-k)(1-\theta)}\Big\}.$$

Thus, by Theorem 10.2.6, for almost all $\pi \in G(k, d)$,

$$\begin{aligned}\underline{\dim}_{\mathrm{B}}\pi F &= \underline{\dim}^k_{\mathrm{B}} F\\ &= \liminf_{r\to 0}\frac{\log C_r^k(F)}{-\log r}\\ &\geqslant \liminf_{r\to 0}\frac{\log\big(c^{-1}r^{-u}\min\{1,\ r^{(s-k)},\ r^{(t-k)(1-\theta)}\}\big)}{-\log r}\\ &= u - \max\{0,\ s-k,\ (t-k)(1-\theta)\}\end{aligned}$$

which completes the proof for lower box dimension. The proof for upper box dimension is similar and omitted. □

Theorem 10.2.7 yields simple conditions guaranteeing that the box dimension profiles are given by an expression analogous to that in the statement of the Marstrand-Mattila projection theorem. The following is an adaptation of [76, Corollary 1.3].

**Corollary 10.2.8** Let $1 \leqslant k < d$ and $F \subseteq \mathbb{R}^d$ be a non-empty bounded set. Then for almost all $\pi \in G(k, d)$

$$\overline{\dim}_{\mathrm{B}} \pi F \geqslant \overline{\dim}_{\mathrm{B}} F - \max\{0, \ \dim_{\mathrm{qA}} F - k\}.$$

In particular, if $\dim_{\mathrm{qA}} F \leqslant \max\{k, \overline{\dim}_{\mathrm{B}} F\}$, then for almost all $\pi \in G(k, d)$

$$\overline{\dim}_{\mathrm{B}} \pi F = \min\{k, \overline{\dim}_{\mathrm{B}} F\}.$$

*Proof* This follows immediately from Theorem 10.2.7 by letting $\theta \to 1$. □

A natural question is when Theorem 10.2.7 improves on the general lower bounds from (10.8). Based on knowledge of the Assouad dimension, we get the following improvement, which first noted in [76, Corollary 1.5].

**Corollary 10.2.9** Let $1 \leqslant k < d$, $F \subseteq \mathbb{R}^d$ be a non-empty bounded set and suppose $\overline{\dim}_{\mathrm{B}} F \leqslant k$ and

$$\max\{k, \overline{\dim}_{\mathrm{B}} F\} \leqslant \dim_{\mathrm{A}} F < \frac{(kd + 2\overline{\dim}_{\mathrm{B}} F (d-k))k}{kd + \overline{\dim}_{\mathrm{B}} F (d-k)}.$$

Then, for almost all $\pi \in G(k, d)$,

$$\overline{\dim}_{\mathrm{B}} \pi F \geqslant \overline{\dim}_{\mathrm{B}} F - \frac{(\dim_{\mathrm{A}} F - k)\overline{\dim}_{\mathrm{B}} F}{k} > \frac{\overline{\dim}_{\mathrm{B}} F}{1 + (1/k - 1/d)\overline{\dim}_{\mathrm{B}} F}.$$

*Proof* This follows immediately from Theorem 10.2.7. If $\overline{\dim}_{\mathrm{B}} F = k$, then the bound follows by letting $\theta \to 0$. If $\overline{\dim}_{\mathrm{B}} F < k$, then choose $\theta$ as large as possible such that $\overline{\dim}_{\mathrm{A}}^{\theta} F \leqslant k$. By Lemma 3.4.4 this can always be done for

$$\theta \geqslant 1 - \frac{\overline{\dim}_{\mathrm{B}} F}{k}$$

and so, for almost all $\pi \in G(k, d)$,

$$\overline{\dim}_{\mathrm{B}} \pi F \geqslant \overline{\dim}_{\mathrm{B}} F - (\dim_{\mathrm{A}} F - k)(1-\theta) \geqslant \overline{\dim}_{\mathrm{B}} F - \frac{(\dim_{\mathrm{A}} F - k)\overline{\dim}_{\mathrm{B}} F}{k}.$$

□

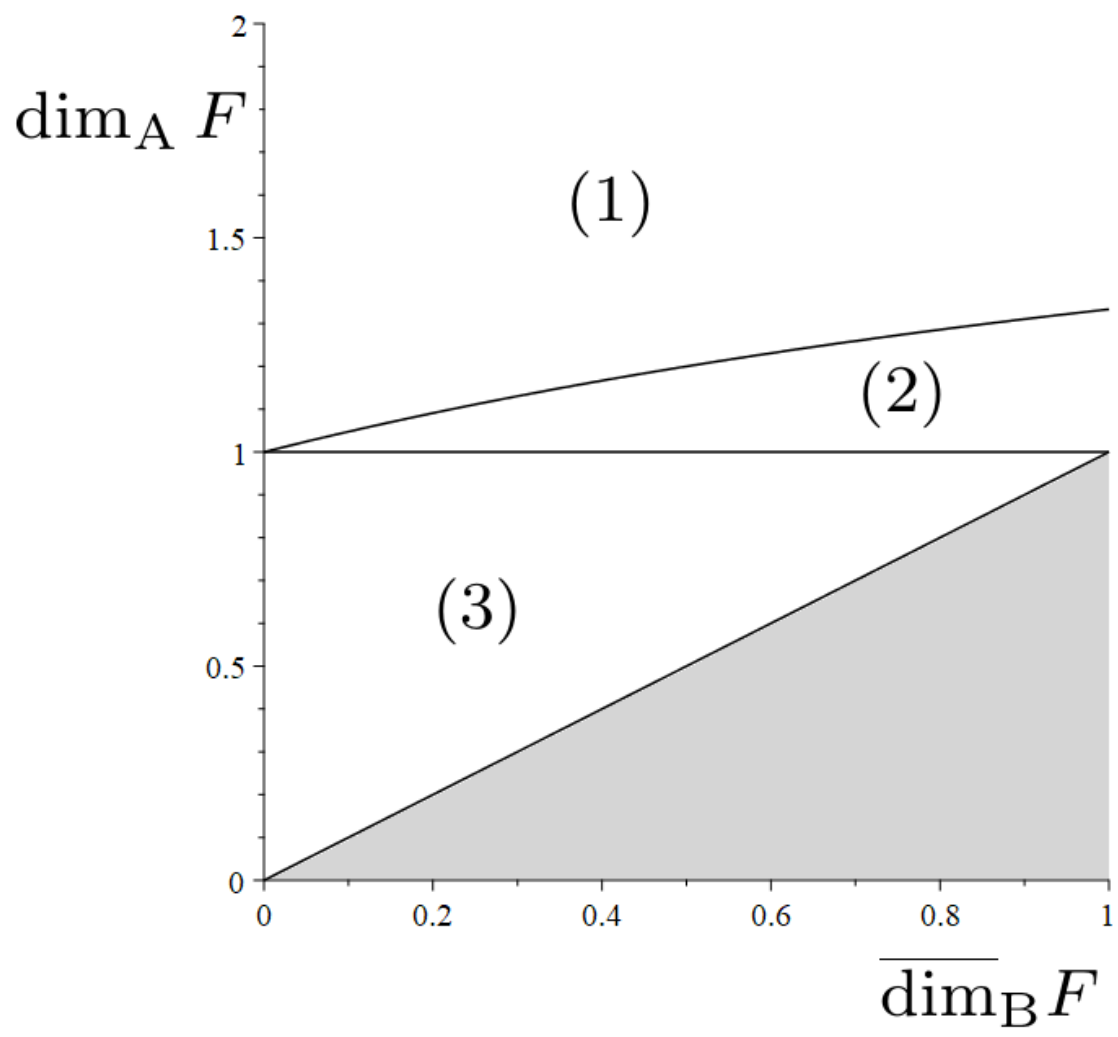

Figure 10.3 An example with $k = 1$ and $d = 2$. Depending on which region the pair $(\overline{\dim}_{\mathrm{B}} F, \dim_{\mathrm{A}} F)$ lies in, we may be able to apply one of Corollary 10.2.8 or Corollary 10.2.9: (1) no information; (2) apply Corollary 10.2.9 to gain improvements on general bounds from (10.8); (3) apply Corollary 10.2.8 to gain a sharp result, giving $\overline{\dim}_{\mathrm{B}} \pi F = \min\{1, \overline{\dim}_{\mathrm{B}} F\} = \overline{\dim}_{\mathrm{B}} F$ almost surely. The pair $(\overline{\dim}_{\mathrm{B}} F, \dim_{\mathrm{A}} F)$ never lies in the grey region, since $\dim_{\mathrm{A}} F \geqslant \overline{\dim}_{\mathrm{B}} F$. The curve bounding regions (1) and (2) is given by $y = (2x + 2)/(x + 2)$.

### 10.2.4 The lower dimension and projections

The lower dimension is very unstable with regard to orthogonal projections. Just as $\dim_{\mathrm{A}} \pi F$ is not almost surely constant, we will see that $\dim_{\mathrm{L}} \pi F$ and $\dim_{\mathrm{L}}^{\theta} \pi F$ are not almost surely constant. However, in this case the examples are much easier to construct. The reason for this is captured by the set $F \subseteq \mathbb{R}^2$ given by

$$F = ([0, 1] \times \{0\}) \cup \{(0, 1)\},$$

that is, a line segment together with a point not contained in the same affine hyperplane as the line segment. As such, $\dim_{\mathrm{L}} F = \dim_{\mathrm{L}}^{\theta} F = 0$ for all $\theta \in (0, 1)$. However, for a set of $\pi \in G(1, 2)$ of positive measure the projection of the isolated point is contained in the projection of the line segment. As such, the projection $\pi F$ will be a line segment, and will therefore satisfy $\dim_{\mathrm{L}} \pi F = \dim_{\mathrm{L}}^{\theta} \pi F = 1$ for all $\theta \in (0, 1)$. However, for another set of $\pi \in G(1, 2)$ of positive measure, the projection $\pi F$ will be a line segment together with an isolated point, and will therefore satisfy $\dim_L \pi F = \dim_L^\theta \pi F = 0$ for all $\theta \in (0,1)$. We can 'upgrade' this simple example to prove the following.

**Theorem 10.2.10** Let $X \subseteq [0,1]$ be a finite set or a countable set with the property that it can be enumerated as a strictly decreasing sequence $x_1 > x_2 > \dots$. Then there exists a compact set $F \subseteq \mathbb{R}^2$ such that for all $x \in X$

$$\dim_L \pi F = \dim_L^\theta \pi F = x$$

for a set of $\pi \in G(1,2)$ of positive measure.

*Proof* We consider only the countable case and observe that the finite case follows by a trivial modification. For each $i \in \mathbb{N}$, let $E_i \subseteq [0,1]$ be a self-similar set satisfying

$$\dim_L E_i = \dim_L^\theta E_i = \dim_H E_i = x_i.$$

Let

$$F = ([0,1/2] \times \{0\}) \cup \bigcup_{i \in \mathbb{N}} (2^{-i} E_i \times \{2^{-i}\})$$

where $2^{-i}E_i$ denotes $E_i$ scaled by a factor of $2^{-i}$, see Figure 10.4. By construction, for every $i$ there will be a set of projections $\pi$ of positive measure such that $\pi F$ consists of a line segment together with pairwise disjoint (similar) copies of the sets $\{E_j : 1 \leqslant j \leqslant i\}$. It follows from Lemma 3.4.9 that for such $\pi$

$$\dim_L \pi F = \dim_L^\theta \pi F = \min\{x_j : 1 \leqslant j \leqslant i\} = x_i$$

which proves the theorem. □

## 10.3 Slices and intersections

Relating the dimension of a set with the dimensions of its slices is a well-studied and difficult problem in geometric measure theory. A 'slice' of a set is the intersection of the set with an affine plane of lower dimension than the ambient space. Consider the unit square $F = [0,1]^2 \subseteq \mathbb{R}^2$ and given a 1-dimensional subspace $V \subseteq \mathbb{R}^2$ consider the family of lines in $\mathbb{R}^2$ orthogonal to $V$. That is $\{V^\perp + v : v \in V\}$, where $V^\perp$ is the orthogonal complement of $V$. Some of these lines intersect $F$ and some do not and it should be clear that for Lebesgue positively many $v \in V$ one has $\dim F \cap (V^\perp + v) = 1$. Perhaps we can expect this behaviour more generally. That is, given compact $F \subseteq \mathbb{R}^2$, then perhaps for 'most' $V \in G(1,2)$ and 'many' $v \in V$

$$\dim F \cap (V^{\perp} + v) \geqslant \dim F - 1.$$

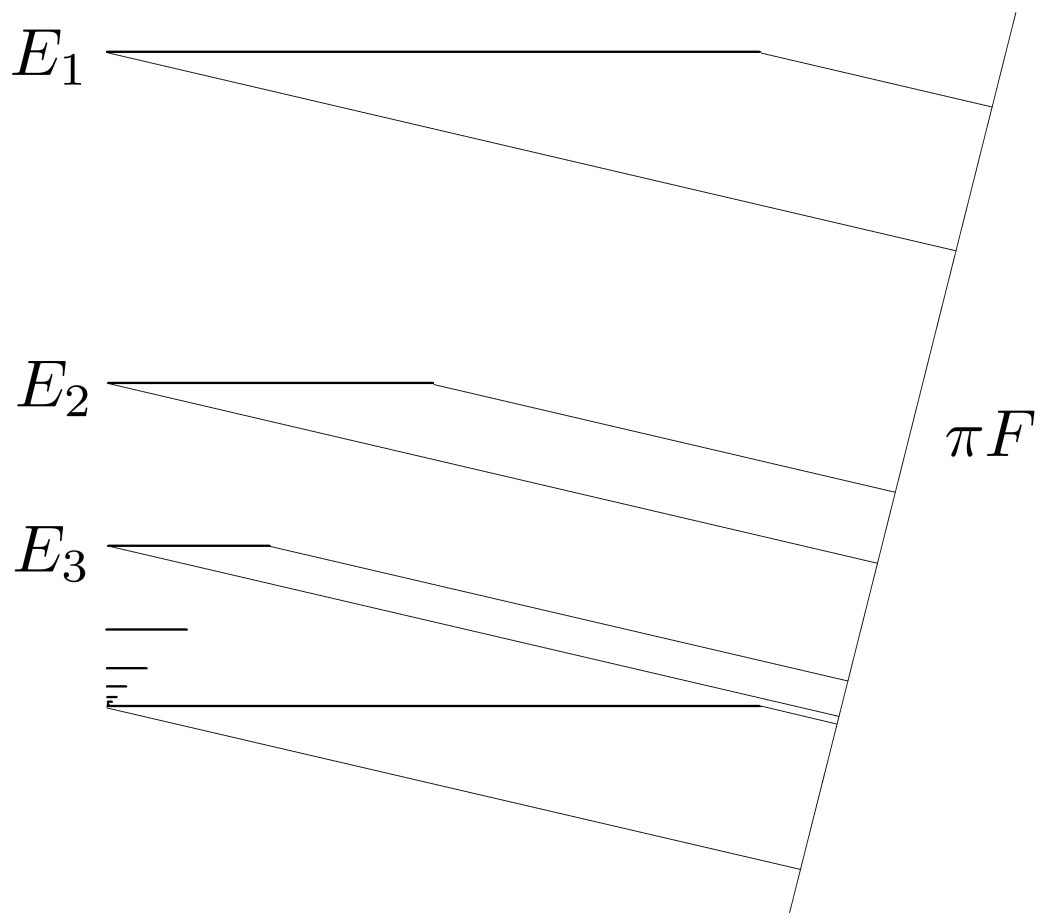

Figure 10.4 The set $F$ together with a projection $\pi F$ where only the sets $\pi E_1$, $\pi E_2$ and $\pi E_3$ are 'visible'. For the purpose of illustration, the self-similar sets have been replaced by line segments.

This is formalised by Marstrand's slicing theorem which was first established in the plane [201], and for the higher dimensional result see [205, 70]. One version of the theorem is as follows, see [205, Theorem 10.10].

**Theorem 10.3.1** Let $d > k \geqslant 1$ be integers and $F \subseteq \mathbb{R}^d$ be a Borel set with $\mathcal{H}^s(F) > 0$ for some $s > (d-k)$. Then, for almost all $V \in G(d-k, d)$,

$$\dim_{\mathrm{H}} \Big( F \cap (V^{\perp} + v) \Big) \geqslant s - (d-k)$$

for positively many $v \in V$.

Here "almost all" $V \in G(d-k, d)$ refers to the Grassmannian measure on $G(d-k, d)$ and "positively many" $v \in V$ refers to Lebesgue measure on $V$.

It turns out that this result is very much measure-theoretic and relies on the fact that the Hausdorff dimension cannot be large on a (measure theoretically) small set. This is not the case for the box dimensions, Assouad dimension and everything in between. We construct a set below with *full* box, quasi-Assouad and Assouad dimension such that *all* slices of the set have dimension 0 — we can even arrange for all the slices to be finite!

**Theorem 10.3.2** Given $d \geqslant 1$, there exists a compact set $F \subseteq \mathbb{R}^d$ such that $\underline{\dim}_{\mathrm{B}} F = d$ but for all integers $1 \leqslant k < d$, all $V \in G(d-k, d)$, and all $v \in V$, the slice $F \cap (V^\perp + v)$ consists of at most $k+1$ points.

Note that the upper bound of $k+1$ on the cardinality of the slices here is sharp since any set $F \subseteq \mathbb{R}^d$ with at least $k+1$ points intersects at least one $k$-dimensional affine plane in a set of at least $k+1$ points.

*Proof* Fix $d \geqslant 1$ and let $E_0 = \{1/\log(n) : n = 2, 3, 4, \dots\}$. It is straightforward to show that $\underline{\dim}_{\mathrm{B}} E_0 = 1$ and we leave the details as an exercise. For comparison, see Theorem 3.4.7 and note that $1/\log(n)$ is 'slower' than $1/n^p$ for any $p$. Let

$$E = (E_0)^d = E_0 \times \cdots \times E_0 \subseteq \mathbb{R}^d$$

which, by the formulae from Section 10.1, has lower box dimension $d$. Note that every point in $E$ is isolated and therefore, for every $x \in E$,

$$\varepsilon(x) = \inf\{|x-y| : y \in E, y \neq x\} > 0.$$

Let

$$E' = \bigcup_{x \in E} B(x, \varepsilon(x)/10)$$

which by the definition of $\varepsilon(\cdot)$ is a pairwise disjoint union of balls. Note that any set consisting of precisely one point from each of the balls $B(x, \varepsilon(x)/10)$ is bi-Lipschitz equivalent to $E$ and, therefore, has lower box dimension $d$. Our set $F$ will be (the closure of) one of these sets. Since $E$ is countable we can enumerate the balls in the construction of $E'$ as $\{B_1, B_2, \dots\}$ and we construct $F$ by placing one point inside each ball $B_1, B_2, \dots$ in turn. Place a point $x_1$ in $B_1$ arbitrarily and then, once we have placed the first $n$ points, the $(n+1)$th point $x_{n+1}$ is placed in $B_{n+1}$ such that, for all $1 \leqslant k < d$, $x_{n+1}$ does not lie in any $k$-dimensional affine plane $(V^\perp + v)$ which already contains $k+1$ points from the set $\{0\} \cup \{x_1, x_2, \dots, x_n\}$. To see why such placement of $x_{n+1}$ is possible, observe that there are only finitely many such affine planes because every collection of $k+1$ points is contained in a unique $k$-dimensional affine plane. Therefore, the union of all affine planes we must avoid has zero $d$-dimensional Lebesgue measure, whereas the set of possible choices for $x_{n+1}$, namely the set $B_{n+1}$, has positive $d$-dimensional Lebesgue measure. Let

$$F = \{0\} \cup \{x_1, x_2, \dots\}$$

which is clearly compact, has lower box dimension $d$, and by construction cannot intersect any $k$-dimensional affine plane $(V^{\perp} + v)$ in strictly more than $(k+1)$ points. □

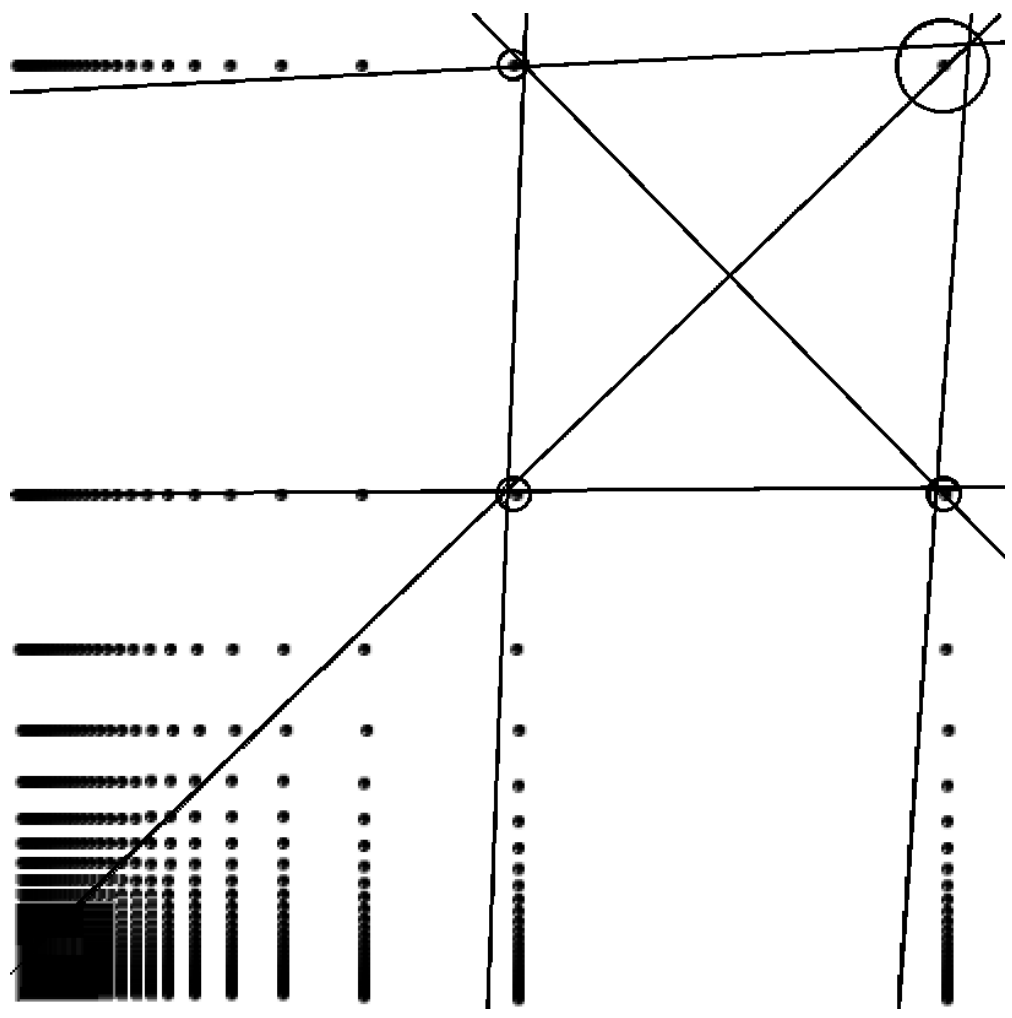

Figure 10.5 The construction of $F$ for $d = 2$: the product set $E = E_0 \times E_0$ is shown along with the construction balls $B(x, \varepsilon(x)/10)$. So far the points $x_1, \dots, x_4$ have been placed and when placing $x_5$ all lines containing at least two of $x_1, \dots, x_4$ must be avoided. There are $\binom{4}{2} = 6$ such lines and they are also shown. Note that they pass through the points $x_1, \dots, x_4$, which are not necessarily points in $E$.

# 11

# Two famous problems in geometric measure theory

In this chapter we discuss the Assouad dimension in the context of two famous conjectures in geometric measure theory: Falconer's distance set conjecture, and the Kakeya conjecture. In Section 11.1 we prove the Assouad dimension version of one variant of the distance set problem in the plane, see Theorem 11.1.3. In Section 11.2 we discuss the Kakeya conjecture and give an elementary proof that 'extended Kakeya sets' have maximal Assouad dimension, see Theorem 11.2.2.

## 11.1 Distance sets

Given a set $F \subseteq \mathbb{R}^d$, the *distance set* of $F$ is defined by

$$D(F) = \{|x-y| \ : \ x, y \in F\}.$$

The distance set problem, which stems from a seminal paper of Falconer [63], is to relate the dimensions of $F$ with the dimensions of $D(F)$. *Falconer's distance set conjecture* refers to several related conjectures, one version of which is as follows.

**Conjecture 11.1.1** Let $d \geqslant 2$ and $F \subseteq \mathbb{R}^d$ be an analytic set. If $\dim_{\mathrm{H}} F > d/2$, then $D(F)$ has positive Lebesgue measure (or even non-empty interior). In particular, $\dim_{\mathrm{H}} D(F) = 1$.

There are numerous partial results in support of this conjecture, but the problem is still open, even in the plane. See [207, Conjecture 4.5] and the subsequent discussion. That said, since 2010 the problem has been receiving a lot of attention in the literature, and significant progress is being made, see for example, [252, 253, 225, 167, 124].

If we only consider the dimension (not measure) of the distance set,

then dimension equal to $d/2$ should suffice. We may also replace the Hausdorff dimension with Assouad dimension, and obtain a different conjecture.

**Conjecture 11.1.2** Let dim denote either the Hausdorff or Assouad dimension, let $d \geqslant 2$ and $F \subseteq \mathbb{R}^d$ be an analytic set. If $\dim F \geqslant d/2$, then $\dim D(F) = 1$.

It may also be interesting to consider further variants of the problem with dim replaced by other notions of dimension or spectra. The planar Assouad dimension version of the conjecture was resolved in [89], apart from the case $\dim_A F = 1$. This remaining case was later resolved in [90]. To the best of our knowledge, all other versions of this conjecture are open for all values of $d \geqslant 2$.

**Theorem 11.1.3** If $F \subseteq \mathbb{R}^2$ with $\dim_A F > 1$, then $\dim_A D(F) = 1$.

A result of Shmerkin [253] resolves the planar Hausdorff dimension version of the problem with some additional assumptions.

**Theorem 11.1.4** If $F \subseteq \mathbb{R}^2$ is a Borel set with $\dim_H F = \dim_P F > 1$, then $\dim_H D(F) = 1$.

Here we show how to obtain Theorem 11.1.3 from Theorem 11.1.4, although this was not how Theorem 11.1.3 was originally proved because Shmerkin's result was not available at the time. The key technical observation is that one can pass questions on distance sets to weak tangents. The following lemma is from [89, Lemma 3.1].

**Lemma 11.1.5** Let $F \subseteq \mathbb{R}^d$ be a non-empty closed set and suppose $E$ is a weak tangent to $F$. Then

$$\dim_A D(F) \geqslant \dim_A D(E).$$

*Proof* We follow the proof from [89, Section 3.2.1]. Since $E$ is a weak tangent to $F$, we may find a non-empty compact set $X \subseteq \mathbb{R}^d$ with $|X| > 0$ and a sequence of similarity maps $T_k : \mathbb{R}^d \to \mathbb{R}^d$ such that

$$T_k(F) \cap X \to E \tag{11.1}$$

in $d_{\mathcal{H}}$ as $k \to \infty$. Write $c_k \in (0, \infty)$ for the similarity ratio of the map $T_k$. Consider the sequence of compact sets given by

$$c_k D(F) \cap [0, |X|]$$

where

$$c_k D(F) = \{c_k z : z \in D(F)\}.$$

Using compactness of $(\mathcal{K}([0,|X|]), d_{\mathcal{H}})$, we can find a strictly increasing infinite sequence of integers $(k_n)_{n \geqslant 1}$ such that $c_{k_n} D(F) \cap [0,|X|]$ converges to a compact set $B$ in the Hausdorff metric $d_{\mathcal{H}}$. In particular, the set $B$ is a weak tangent to $D(F)$ and therefore Theorem 5.1.2 gives

$$\dim_{\mathrm{A}} D(F) \geqslant \dim_{\mathrm{A}} B.$$

It therefore suffices to show that $D(E) \subseteq B$. Let $z = |x-y| \in D(E)$ for some $x, y \in E$. It follows from (11.1) that we can find sequences $x_k, y_k \in T_k(F) \cap X$ such that $x_k \to x$ and $y_k \to y$. For each $k$, $T_k^{-1}(x_k), T_k^{-1}(y_k) \in F$ and so

$$c_k^{-1}|x_k - y_k| = |T_k^{-1}(x_k) - T_k^{-1}(y_k)| \in D(F).$$

Moreover, $|x_k - y_k| \leqslant |X|$ and therefore

$$|x_k - y_k| \in c_k D(F) \cap [0,|X|].$$

It follows that

$$z = |x-y| = \lim_{k\to\infty} |x_k - y_k| = \lim_{n\to\infty} |x_{k_n} - y_{k_n}| \in B$$

as required. □

*Proof of Theorem 11.1.3:* Let $F \subseteq \mathbb{R}^2$ be a closed set with $\dim_{\mathrm{A}} F = s > 1$. We will deal with the non-closed case at the end of the proof. It follows from Theorem 5.1.3 that $F$ has a weak tangent $E$ such that

$$\dim_{\mathrm{H}} E = \dim_{\mathrm{P}} E = \dim_{\mathrm{A}} E = s > 1.$$

It follows that $\dim_{\mathrm{A}} D(F) \geqslant \dim_{\mathrm{A}} D(E) \geqslant \dim_{\mathrm{H}} D(E)$ by Lemma 11.1.5 and, moreover, since $\dim_{\mathrm{P}} E = \dim_{\mathrm{H}} E > 1$, Theorem 11.1.4 implies $\dim_{\mathrm{H}} D(E) = 1$.

All that remains is the case when $F$ is not closed. However, $\dim_{\mathrm{A}} \overline{F} = \dim_{\mathrm{A}} F > 1$ and $\overline{D(F)} \supseteq D(\overline{F})$. Therefore, using the result in the closed case,

$$\dim_{\mathrm{A}} D(F) = \dim_{\mathrm{A}} \overline{D(F)} \geqslant \dim_{\mathrm{A}} D(\overline{F}) = 1$$

as required. □

Clearly, $\dim_{\mathrm{A}} F > 1$ does not guarantee that $D(F)$ has non-empty interior since such sets may be countable, but it *does* guarantee the existence of a weak tangent to $D(F)$ with non-empty interior, by applying Theorem 5.1.5.

**Corollary 11.1.6** If $F \subseteq \mathbb{R}^2$ is any set with $\dim_A F > 1$, then $D(F)$ has a weak tangent with non-empty interior.

In the higher dimensional setting, Lemma 11.1.5 yields partial results and also shows that the Hausdorff dimension version of Conjecture 11.1.2 *implies* the Assouad dimension version, see [89, Theorem 2.5]. Moreover, partial results for Hausdorff dimension may be 'converted' into analogous partial results for Assouad dimension. See the survey [207] for some history and discussion of various partial results in higher dimensions, such as the well-known Falconer-Erdoğan-Wolff bounds. These bounds were recently improved in [55, 56].

**Theorem 11.1.7** Let $d \geqslant 2$ be an integer and $F \subseteq \mathbb{R}^d$ be any non-empty set. Then

$$\dim_A D(F) \geqslant \inf_{\substack{E \in \mathcal{K}(\mathbb{R}^d): \\ \dim_H E = \dim_A F}} \dim_A D(E).$$

For example, if $d = 3$ and $\dim_A F > 9/5$, then $\dim_A D(F) = 1$.

The non-quantitative part of Theorem 11.1.7 was proved in [89, Theorem 2.5] and provides a recipe for directly transferring *any* bounds concerning the Hausdorff dimension version of the distance set problem into analogous bounds concerning the Assouad dimension version. The quantitative example with $d = 3$ follows by applying the non-quantitative part to the Hausdorff dimension results from [55]. Higher dimensional analogues of these bounds can also be derived by considering [56]. It appears that the appropriate analogues of these bounds are *not* known to hold for other notions of dimension, such as box or packing dimension.

As for sets with small Assouad dimension, it was proved in [100, Theorem 2.9] that if $F \subseteq \mathbb{R}^d$ satisfies $0 < \dim_A F < d$, then

$$\dim_A D(F) > \dim_A F/d, \tag{11.2}$$

which improves on estimates from Theorem 11.1.7 when $\dim_A F$ is small.

Observe that the distance set $D(F)$ consists of the norms of points in the *difference set* $F - F$. By relating the difference set and an appropriate orthogonal projection of the product set $F \times F$, the results from Section 10.1 show

$$\overline{\dim}_B D(F) \leqslant \overline{\dim}_B (F \times F) \leqslant 2\overline{\dim}_B F$$

for any bounded $F \subseteq \mathbb{R}^d$. However, there exist compact sets in $\mathbb{R}^d$ with Hausdorff dimension 0, but whose distance set contains an interval. For example, [166, Proposition 3.1] provides a compact set $F \subseteq \mathbb{R}^d$ with Hausdorff dimension 0 but for which $F \cap (F+t) \neq \varnothing$ for all $t \in [0,1]^d$. This shows that $[0, \sqrt{d}] \subseteq D(F)$ and so $\dim_{\mathrm{A}} D(F) = 1$. That said, the product-projection argument above guarantees $\dim_{\mathrm{H}} D(F) \leqslant \dim_{\mathrm{H}} F + \dim_{\mathrm{P}} F$ for $F \subseteq \mathbb{R}^d$ and so the size of the distance set can be controlled by the size of $F$ in some sense. The situation for Assouad dimension is rather different since Assouad dimension can increase under orthogonal projection. In fact, arbitrarily small sets can have distance sets with full Assouad dimension. Such examples are provided by simple self-similar sets, requiring only an irrationality condition. This shows that for every $s \in (0,1)$, there exists an Ahlfors regular set with dimension $s$ but $\dim_{\mathrm{A}} D(F) = 1$. The following is from [89, Theorem 2.7].

**Theorem 11.1.8** Let $F \subseteq \mathbb{R}^d$ be a self-similar set which is not a single point and suppose that two of the defining similarity ratios are given by $\alpha, \beta \in (0,1)$ satisfying $\log \alpha / \log \beta \notin \mathbb{Q}$. Then $\dim_{\mathrm{A}} D(F) = 1$.

*Proof* We follow the argument from [89, Section 3.2.3]. Let $x, y \in F$ be distinct points and let $\Delta := |x-y| > 0$. Let $S_\alpha$ be a defining similarity map with similarity ratio $\alpha$ and $S_\beta$ be a defining similarity map with similarity ratio $\beta$. Using invariance under the defining IFS, for integers $m, n \geqslant 0$, the points $S_\alpha^m \circ S_\beta^n(x)$, $S_\alpha^m \circ S_\beta^n(y) \in F$. Therefore

$$|S_\alpha^m \circ S_\beta^n(x) - S_\alpha^m \circ S_\beta^n(y)| = \alpha^m \beta^n \Delta \in D(F).$$

The proof of Theorem 11.1.8 is now similar to Theorem 7.2.1. For integer $k \geqslant 1$, let $T_k : [0,1] \to [0,1]$ be the similarity defined by $T_k(x) = \Delta^{-1}\beta^{-k}x$ and, using compactness of $\mathcal{K}([0,1])$, extract a convergent subsequence of $T_k(D(F)) \cap [0,1]$ in the Hausdorff metric. The limit of this subsequence is a weak tangent to $D(F)$. Moreover,

$$\{0\} \cup \left\{\alpha^m \beta^n \ : \ m,n \in \mathbb{Z},\, m \geqslant 0, n \geqslant -k\right\} \cap [0,1] \ \subseteq \ T_k(D(F)) \cap [0,1]$$

for all $k \geqslant 1$ and therefore the weak tangent contains

$$\overline{\{\alpha^m \beta^n : m \in \mathbb{N}, n \in \mathbb{Z}\}} \cap [0,1] = [0,1]$$

with this equality following from Dirichlet's theorem and the irrationality of $\log \alpha / \log \beta$. □

By changing the contraction ratios at each stage in the construction, Theorem 11.1.8 can be used to find a compact set $F \subseteq [0,1]$ with $\dim_{\mathrm{A}} F = 0$ and $\dim_{\mathrm{A}} D(F) = 1$. See [89, Example 2.6] for the details.

Despite the previous examples, Theorem 11.1.3 is sharp in the sense that one cannot guarantee that the distance set has full Assouad dimension for sets $F$ with $\dim_A F > s$ for any $s < 1$. The next result follows [89, Example 2.8].

**Theorem 11.1.9** For every $s \in [0,1)$, there exists a compact set $F \subseteq [0,1]$ with $\dim_A F = \dim_H F \geqslant s$ but $\dim_A D(F) < 1$.

*Proof* Let $s \in [0,1)$, choose an integer $N \geqslant 2$ satisfying $(N+1)/2 - N^s > 1$, and another integer $K$ satisfying $N^s \leqslant K < (N+1)/2$. For $i \in \{0, \dots, K-1\}$, let $S_i : [0,1] \to [0,1]$ be the similarity defined by

$$S_i(x) = (x + 2i)/N$$

and $F \subseteq [0,1]$ be the self-similar set which is the attractor of the IFS $\{S_i\}_{i=0}^{K-1}$. By construction, this IFS satisfies the OSC, and therefore

$$\dim_A F = \dim_H F = \frac{\log K}{\log N} \geqslant s.$$

Let $\pi$ denote orthogonal projection onto the subspace of $\mathbb{R}^2$ spanned by the vector $(1,-1)$ and consider the set $\pi(F \times F)$. This set turns out to also be a self-similar set satisfying the OSC, this time defined by $2K-1$ similarities with common similarity ratio $1/N$. The map $g$ defined by $g(x) = |x|$ is bi-Lipschitz on $(-\infty, 0)$ and $[0, \infty)$ and so cannot increase Assouad dimension. The fact that $g$ is Lipschitz on the whole of $\mathbb{R}$ is not sufficient here since Assouad dimension is not stable under Lipschitz maps. Since $g(\pi(F \times F)) = D(F)$ we obtain

$$\dim_A D(F) \leqslant \dim_A \pi(F \times F) \leqslant \frac{\log(2K-1)}{\log N} < 1$$

as desired. □

There are many variants on the distance set problem. For example, [281] considers triples of distances between triples of points instead of distances between pairs of points. The resulting 'triangle set' is a subset of $\mathbb{R}^3$ and the problem is to understand the dimension of the triangle set, given the dimension of the original set. Partial results concerning the Assouad dimension version of the problem can be found in [281].

## 11.2 Kakeya sets

A *Kakeya set* is a compact set in $\mathbb{R}^d$ which contains a unit line segment in every possible direction. Besicovitch proved that such sets can have zero $d$-dimensional Lebesgue measure [30], but it is a major open problem in geometric measure theory if Kakeya sets can be even smaller than this, that is, have Hausdorff dimension less than $d$. Sometimes a Kakeya set with zero $d$-dimensional Lebesgue measure is referred to as a *Besicovitch set*.

**Conjecture 11.2.1** (Kakeya Conjecture) *If $K \subseteq \mathbb{R}^d$ is a Kakeya set, then $\dim_{\mathrm{H}} K = d$.*

Much progress has been made on this problem, but full resolution is still very much out of reach, see [158]. The case when $d = 1$ is trivial and the case when $d = 2$ was resolved by Davies in 1971 [51]. One can ask about other dimensions, and replacing the Hausdorff dimension by the upper box dimension or Assouad dimension should, in principle, make the conjecture easier to prove. However, progress seems to be guided by results on the Hausdorff dimension and so progress on the other dimensions is no faster.

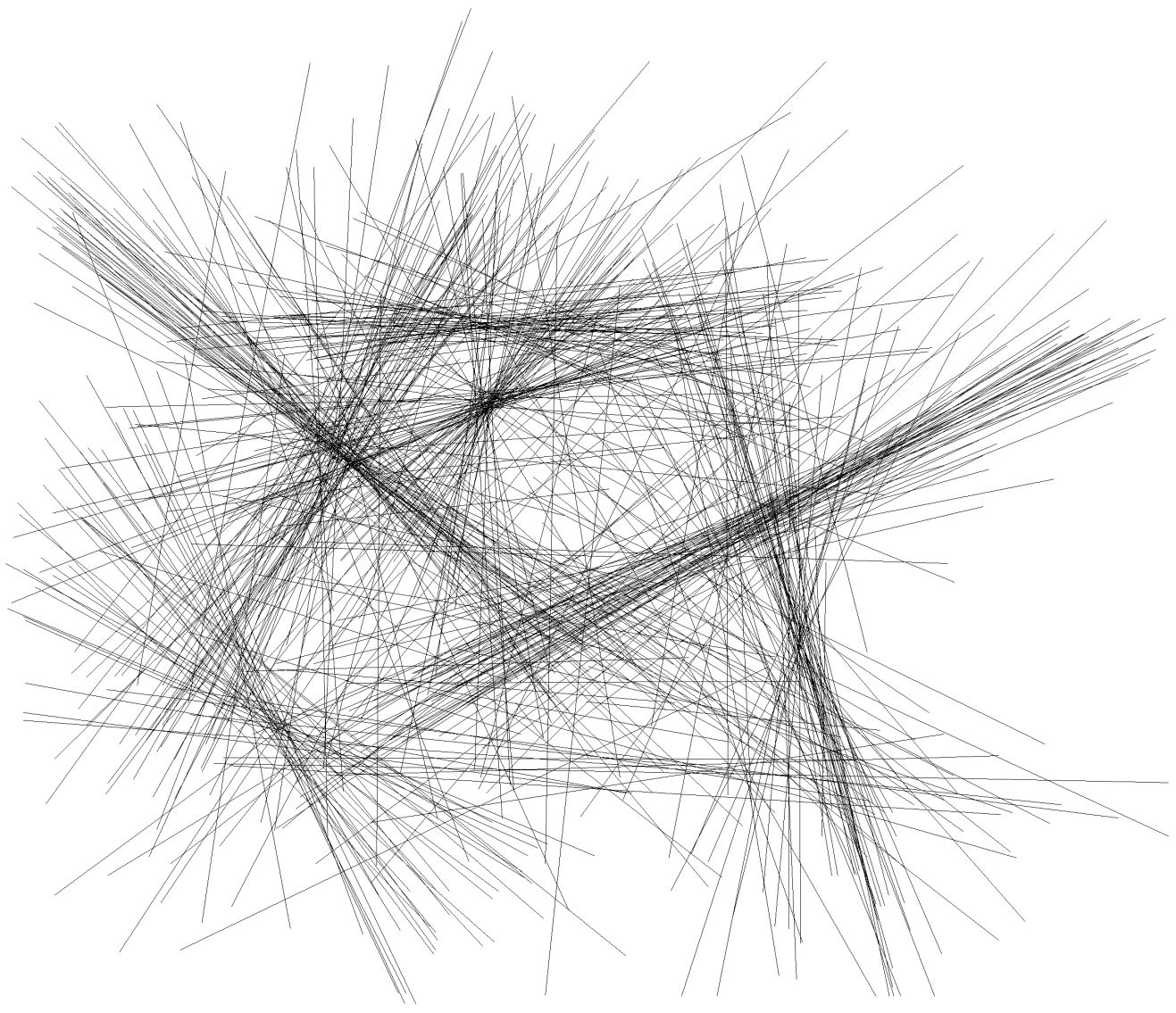

Figure 11.1 A Kakeya set in $\mathbb{R}^2$.

If one replaces 'unit line segment' with 'infinite line segment' in the definition of Kakeya set, then, again, the conjecture should be easier to prove. This approach was explored by Keleti [165]. *Keleti's line segment extension conjecture* [165] is that extending any collection of line segments in $\mathbb{R}^d$ to the corresponding collection of doubly infinite lines in the same directions does not alter the Hausdorff dimension. Keleti proved that a positive answer to this conjecture (for all $d$) would imply the Kakeya conjecture for upper box dimension, and a positive answer for a particular value of $d$ would prove that any Kakeya set in $\mathbb{R}^d$ has Hausdorff dimension at least $d-1$. Keleti proved that his conjecture is true when $d = 2$, see also [80, Theorem 3.2].

It turns out that the Assouad dimension of extended Kakeya sets is easy to compute. We say $K \subseteq \mathbb{R}^d$ is an *extended Kakeya set* if for all directions $\alpha \in \mathbb{S}^{d-1}$ there exists a translation $t_\alpha \in \mathbb{R}^d$ such that

$$\{\lambda\alpha + t_\alpha : \lambda \in [0, \infty)\} \subseteq K. \tag{11.3}$$

The following was proved in [105, Theorem 7.1] and we follow their proof below.

**Theorem 11.2.2** If $K \subseteq \mathbb{R}^d$ is an extended Kakeya set, then

$$\dim_{\mathrm{A}} K = d.$$

*Proof* We prove the theorem by showing that the unit ball $B(0,1) \subseteq \mathbb{R}^d$ is a weak tangent to $K$. For each integer $k \geqslant 1$, let $T_k : \mathbb{R}^d \to \mathbb{R}^d$ be defined by $T_k(x) = x/k$. Let $\Delta_n \subseteq \mathbb{S}^{d-1}$ be a sequence of finite sets which converge in the Hausdorff metric to $\mathbb{S}^{d-1}$, where $\mathbb{S}^{d-1} = \{x \in \mathbb{R}^d : \|x\| = 1\}$ is the unit sphere. Let

$$\Delta_n^* = \{t\alpha : \alpha \in \Delta_n,\, t \in [0,1]\} \subseteq B(0,1)$$

and observe that $\Delta_n^*$ converges to $B(0,1)$ in the Hausdorff metric. Let $\varepsilon > 0$ and choose $n$ sufficiently large to ensure that $\Delta_n^*$ is $\varepsilon$-close to $B(0,1)$ in the Hausdorff metric, and also choose $k$ sufficiently large, depending on $n$, to ensure that $|T_k(t_\alpha)| < \varepsilon/\sqrt{2}$ for all $\alpha \in \Delta_n$, where $t_\alpha$ is as in (11.3). We may do this since $\Delta_n$ is finite and therefore the set

$$\{t_\alpha : \alpha \in \Delta_n\}$$

is bounded, although the bound may depend on $n$. Using (11.3) the $\varepsilon$-neighbourhood of $T_k(K) \cap B(0,1)$ contains $\Delta_n^*$ and the $\varepsilon$-neighbourhood of $\Delta_n^*$ contains $B(0,1)$. This means that the $2\varepsilon$-neighbourhood of $T_k(K) \cap B(0,1)$ contains $B(0,1)$ and so the Hausdorff distance between $T_k(K) \cap$ $B(0,1)$ and $B(0,1)$ is bounded by $2\varepsilon$, which proves that $B(0,1)$ is a weak tangent to $K$. Theorem 5.1.2 therefore implies that

$$\dim_A K \geqslant d,$$

as required. □

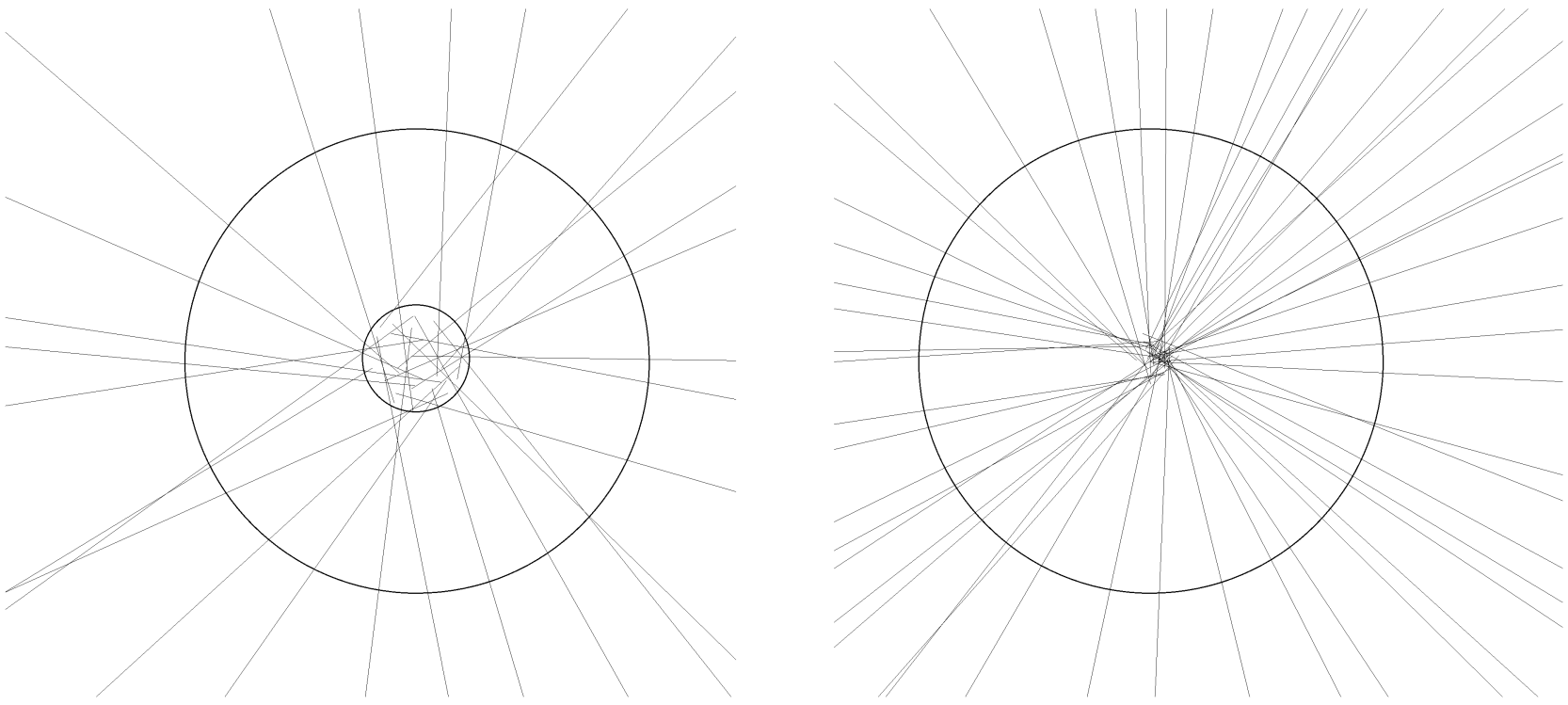

Figure 11.2 Two plots of $T_k(K)$ where $K$ is an extended Kakeya set with the unit ball for reference. On the left $\varepsilon$ is roughly 0.2 and on the right it is roughly 0.1.

# 12
# Conformal dimension

In this chapter we briefly discuss conformal dimension and highlight some particular areas of interest in connection with the topics covered in this book. The concept of conformal dimension considers how much a given notion of dimension can drop under quasi-symmetric deformation. This is often a very challenging problem and we refer the reader to Mackay and Tyson's book [195] for a thorough treatment. We give a proof that equicontractive self-similar sets have conformal Assouad dimension 0, see Theorem 12.1.1.

## 12.1 Lowering the Assouad dimension by quasi-symmetry

Consider a notion of dimension, dim, which in what follows will most often be Hausdorff or Assouad dimension. Then the conformal variant of this dimension, written $\mathcal{C}\dim$, measures how small the dimension can be made under quasi-symmetric deformation. This concept was introduced by Pansu [228] and has given rise to many useful invariants used in classifying, for example, boundaries of hyperbolic groups. See [195] for a detailed treatment.

More formally, we say a homeomorphism $\phi$ between two sets $E, F \subseteq \mathbb{R}^d$ is a *quasi-symmetry* if there exists a homeomorphism $\eta : [0, \infty) \to [0, \infty)$ such that

$$\frac{|\phi(x) - \phi(y)|}{|\phi(x) - \phi(z)|} \leqslant \eta\left(\frac{|x - y|}{|x - z|}\right)$$

for all $x, y, z \in E$ with $x \neq z$. Note that if $\phi$ is bi-Lipschitz, then it is easily seen to be a quasi-symmetry, but quasi-symmetries certainly need not be bi-Lipschitz. Given a notion of dimension dim, the conformal dimension of $F \subseteq \mathbb{R}^d$ is defined by

$$\mathcal{C}\dim F = \inf\{\dim \phi(F) \ : \ \phi \text{ is quasi-symmetric}\}.$$

Immediately we see $\mathcal{C}\dim F \leqslant \dim F$ but this inequality can be strict. For example, if $F \subseteq \mathbb{R}^d$ is a self-similar set satisfying the SSC, then the conformal Assouad dimension of $F$ is 0. This shows that there is no limit on how much the conformal dimension can drop from the original dimension. See Chapter 7 for background on self-similar sets.

**Theorem 12.1.1** If $F \subseteq \mathbb{R}^d$ is a self-similar set satisfying the SSC, then $\mathcal{C}\dim_{\mathrm{A}} F = 0$.

*Proof* We prove the result in the equicontractive case. That is, where $F$ is the attractor of an IFS of similarities $\{S_i\}_{i \in \mathcal{I}}$ which share a common similarity ratio $c \in (0,1)$, see (7.1). The general case follows from [195, Theorem 2.1.2], see also [270]. Let

$$\delta = \min_{i \neq j} \inf_{\substack{x \in S_i(F) \\ y \in S_j(F)}} |x-y| > 0$$

which is strictly positive since the SSC holds. Let $\lambda \in (0,1)$ and $F_\lambda$ be the self-similar set defined by the IFS $\{S_i^\lambda\}_{i \in \mathcal{I}}$ where $S_i^\lambda$ is a similarity with the same (unique) fixed point as $S_i$ but with similarity ratio $\lambda c$. Since the original IFS satisfies the SSC, the IFS defining $F_\lambda$ satisfies the SSC for sufficiently small $\lambda$. Moreover, for sufficiently small $\lambda$,

$$\min_{i \neq j} \inf_{\substack{x \in S_i^\lambda(F_\lambda) \\ y \in S_j^\lambda(F_\lambda)}} |x-y| \geqslant \delta/2.$$

We assume from now on that $\lambda$ is small enough for both of the previous statements to hold. By Corollary 6.4.4, the Assouad dimension of $F_\lambda$ is

$$\dim_{\mathrm{A}} F_\lambda = \frac{\log \# \mathcal{I}}{-\log(\lambda c)}$$

which can be made arbitrarily small by letting $\lambda \to 0$. We will show that $F$ and $F_\lambda$ are quasi-symmetric equivalent, that is, there is a quasi-symmetry $\phi : F \to F_\lambda$, which proves the result.

Let $\phi : F \to F_\lambda$ be defined by $\phi(\Pi(\boldsymbol{i})) = \Pi_\lambda(\boldsymbol{i})$ for $\boldsymbol{i} \in \mathcal{I}^{\mathbb{N}}$ where $\Pi$ is the symbolic coding map associated with the original IFS and $\Pi_\lambda$ is the symbolic coding map for $F_\lambda$, see Section 6.1. This map is well-defined since the SSC guarantees that the symbolic coding maps are bijective. For all $\boldsymbol{i}, \boldsymbol{j} \in \mathcal{I}^{\mathbb{N}}$ with $\boldsymbol{i} \neq \boldsymbol{j}$

$$c^{|\boldsymbol{i} \wedge \boldsymbol{j}|}\delta|F| \leqslant |\Pi(\boldsymbol{i}) - \Pi(\boldsymbol{j})| \leqslant c^{|\boldsymbol{i} \wedge \boldsymbol{j}|}|F| \tag{12.1}$$

and

$$(\lambda c)^{|\boldsymbol{i} \wedge \boldsymbol{j}|}(\delta/2)|F_\lambda| \leqslant |\Pi_\lambda(\boldsymbol{i}) - \Pi_\lambda(\boldsymbol{j})| \leqslant (\lambda c)^{|\boldsymbol{i} \wedge \boldsymbol{j}|}|F_\lambda| \tag{12.2}$$

where $\boldsymbol{i} \wedge \boldsymbol{j} \in \mathcal{I}^*$ is the longest common prefix of $\boldsymbol{i}$ and $\boldsymbol{j}$. Consider $x, y, z \in F$ with $x = \Pi(\boldsymbol{i})$, $y = \Pi(\boldsymbol{j})$, $z = \Pi(\boldsymbol{k})$ for $\boldsymbol{i}, \boldsymbol{j}, \boldsymbol{k} \in \mathcal{I}^{\mathbb{N}}$ with $\boldsymbol{i} \neq \boldsymbol{k}$. It follows from (12.1)–(12.2) that

$$\begin{aligned}\frac{|\phi(x) - \phi(y)|}{|\phi(x) - \phi(z)|} &\leqslant \frac{(\lambda c)^{|\boldsymbol{i} \wedge \boldsymbol{j}|}}{(\lambda c)^{|\boldsymbol{i} \wedge \boldsymbol{k}|}(\delta/2)} \\ &= 2\delta^{-1-\frac{\log(\lambda c)}{\log c}} \left(\frac{c^{|\boldsymbol{i} \wedge \boldsymbol{j}|}\delta}{c^{|\boldsymbol{i} \wedge \boldsymbol{k}|}}\right)^{\frac{\log(\lambda c)}{\log c}} \\ &\leqslant 2\delta^{-1-\frac{\log(\lambda c)}{\log c}} \left(\frac{|x - y|}{|x - z|}\right)^{\frac{\log(\lambda c)}{\log c}}.\end{aligned}$$

Therefore $\phi$ is a quasi-symmetry with $\eta : [0, \infty) \to [0, \infty)$ given by

$$\eta(t) = 2\delta^{-1-\frac{\log(\lambda c)}{\log c}} t^{\frac{\log(\lambda c)}{\log c}}$$

which completes the proof. □

The topological dimension, $\dim_{\mathrm{T}} F$, of a non-empty set $F \subseteq \mathbb{R}^d$ is the minimal integer $n \geqslant 0$ such that for every $r > 0$ there exists an $r$-cover of $F$ by open sets such that every $x \in F$ lies in at most $n + 1$ of the covering sets. As such, totally disconnected sets have topological dimension 0, and $\dim_{\mathrm{T}} \psi([0, 1]^k) = k$ for all $k$ and homeomorphisms $\psi : [0, 1]^k \to \mathbb{R}^d$. Moreover, $\dim_{\mathrm{T}} F \leqslant \dim_{\mathrm{H}} F$ for all sets $F$. Since the Assouad and Hausdorff dimensions are bounded below by the topological dimension, which is preserved under homeomorphism, one always has

$$\mathcal{C}\dim_{\mathrm{A}} F \geqslant \mathcal{C}\dim_{\mathrm{H}} F \geqslant \dim_{\mathrm{T}} F$$

and Theorem 12.1.1 shows that this lower bound can be realised in non-trivial situations. Conformal Assouad dimension behaves well in the context of weak tangents, due to the following result of Mackay and Tyson [195, Proposition 6.1.7], which should be compared with Theorem 5.1.2.

**Theorem 12.1.2** Let $F \subseteq \mathbb{R}^d$ be closed, $E \subseteq \mathbb{R}^d$ be compact, and suppose $E$ is a weak tangent to $F$. Then $\mathcal{C}\dim_{\mathrm{A}} F \geqslant \mathcal{C}\dim_{\mathrm{A}} E$.

This is particularly useful if a given set has easily understood tangents. For example, this result combined with Theorem 5.1.5 yields the following equivalence.

**Corollary 12.1.3** For $F \subseteq \mathbb{R}^d$, $\dim_{\mathrm{A}} F = d$ if and only if $\mathcal{C}\dim_{\mathrm{A}} F = d$.

*Proof* The non-trivial implication follows since $\dim_{\mathrm{A}} F = d$ guarantees that $[0,1]^d$ is a weak tangent to $F$ by Theorem 5.1.5. Therefore, $\mathcal{C}\dim_{\mathrm{A}} F \geqslant \dim_{\mathrm{T}}[0,1]^d = d$. □

Theorem 12.1.2 was used by Mackay [194] to compute the conformal Assouad dimension of Bedford-McMullen carpets; recall the notation and setup for these in Section 8.3. There turns out to be a simple dichotomy, showing that the conformal Assouad dimension of a Bedford-McMullen carpet is either the Assouad dimension or zero. See [154, Theorems A and B] for further results in this direction.

**Theorem 12.1.4** Let $F$ be a Bedford-McMullen carpet. If both $N_{\max} < n$ and $N_0 < m$, then

$$\mathcal{C}\dim_{\mathrm{A}} F = 0,$$

and, otherwise,

$$\mathcal{C}\dim_{\mathrm{A}} F = \dim_{\mathrm{A}} F = \frac{\log N_0}{\log m} + \frac{\log N_{\max}}{\log n}.$$

*Proof sketch.* In the first case, it is shown that $F$ is *uniformly disconnected*, and then the conformal Assouad dimension is 0 by [195, Theorem 2.1.2]. In the second case, a weak tangent is constructed with maximal Assouad dimension which is of the form $[0,1] \times E$ where $E$ is a self-similar set. It is known that the Assouad dimension of such product sets cannot be lowered by quasi-symmetric deformation, see [228, Lemma 6.3], and therefore Theorem 12.1.2 implies that the conformal Assouad and Assouad dimensions coincide. □

In general it is a notoriously difficult problem to compute the conformal dimension. Self-similar sets have attracted a large amount of attention and still very little is known outside of the SSC case. A celebrated result of Tyson and Wu [271] shows that the Sierpiński triangle has conformal Assouad and conformal Hausdorff dimension equal to 1 (its topological dimension), see Figure 12.1. Recall that the Hausdorff and Assouad dimensions of the Sierpiński triangle are both given by $\log 3/\log 2$. To prove the upper bound they construct an explicit family of IFS attractors which are all quasi-symmetric equivalent to each other and have dimensions approaching 1. Another famous example, this time for more nefarious reasons, is the Sierpiński carpet, see Figure 12.1. This is another self-similar set satisfying the OSC but this time the conformal (Assouad or Hausdorff) dimension is unknown. The Sierpiński carpet is the attractor of the IFS formed by dividing the unit square into a $3 \times 3$ grid, removing the middle square and taking orientation preserving similarities mapping the square onto the remaining 8 pieces. Therefore, for $F$ the Sierpiński carpet

$$\dim_{\mathrm{T}} F = 1 < \frac{\log 8}{\log 3} = \dim_{\mathrm{H}} F = \dim_{\mathrm{A}} F.$$

It is known that

$$1 + \frac{\log 2}{\log 3} \leqslant \mathcal{C}\dim_{\mathrm{H}} F \leqslant \mathcal{C}\dim_{\mathrm{A}} F < \frac{\log 8}{\log 3}.$$

The first bound follows since

$$[0,1] \times E \subseteq F$$

where $E$ is the middle third Cantor set and therefore by [228, Lemma 6.3] the conformal Hausdorff dimension of $F$ cannot be lowered below the dimension of this product set, see also [269]. The fact that $\mathcal{C}\dim_{\mathrm{A}} F < \dim_{\mathrm{A}} F$ is due to Keith and Laakso [164, Corollary 1.0.5], see also [195, Example 6.2.3]. For other more recent developments on this problem, see [169, 176].

Rossi and Suomala [243] investigated the conformal Hausdorff dimension of the Mandelbrot percolation sets, recall Section 9.4. They proved that there exists a $t$ such that almost surely, conditioned on non-extinction,

$$\mathcal{C}\dim_{\mathrm{H}} M = t < \dim_{\mathrm{H}} M.$$

In contrast, Theorem 9.4.1 and Corollary 12.1.3 show that almost surely, conditioned on non-extinction,

$$\mathcal{C}\dim_{\mathrm{A}} M = \dim_{\mathrm{A}} M = d.$$

Therefore Mandelbrot percolation sets are typically minimal for Assouad dimension but not for Hausdorff dimension.

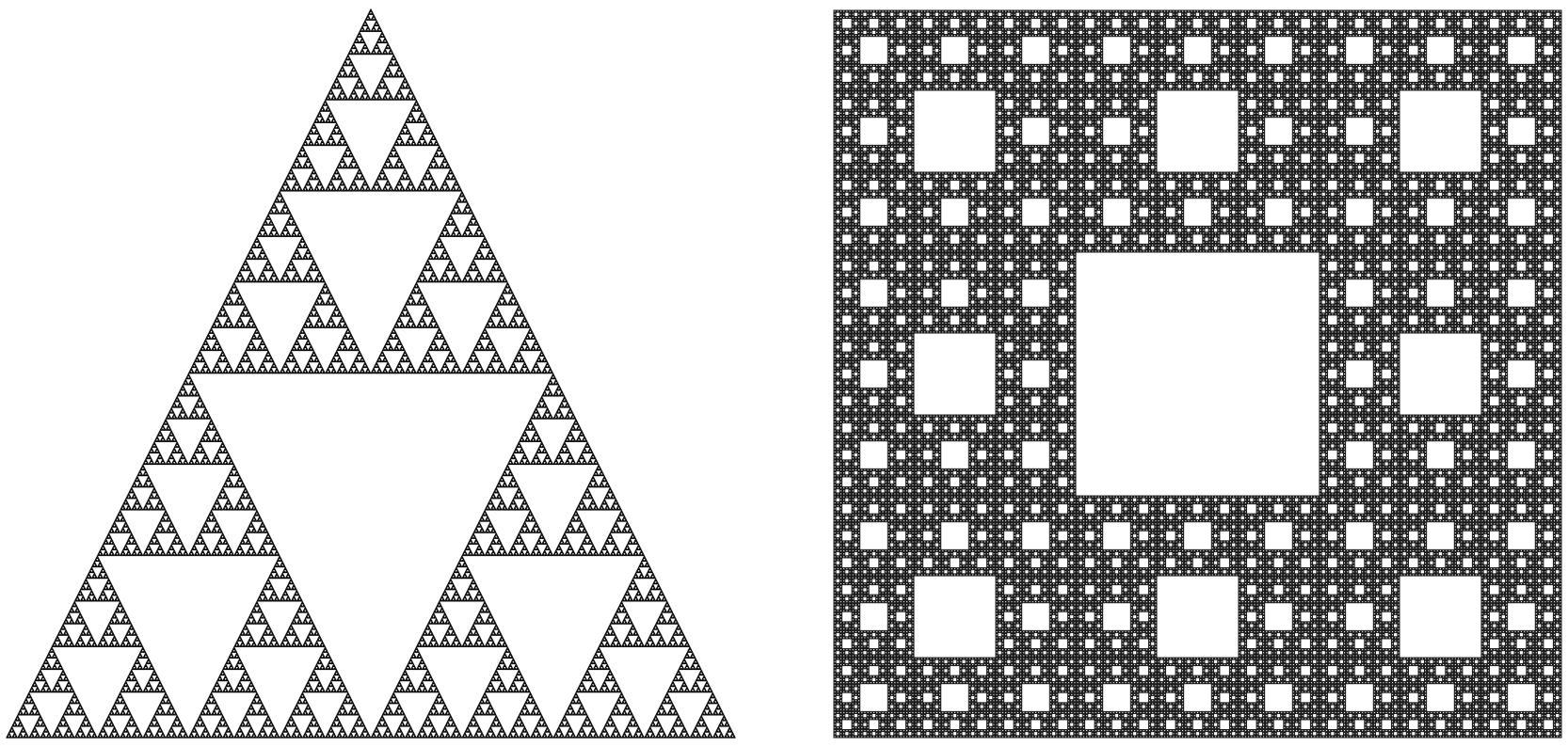

Figure 12.1 Left: the Sierpiński triangle which has conformal dimension 1. Right: the Sierpiński carpet, which has conformal dimension strictly in between its topological and Hausdorff dimension. The precise value is unknown.

# PART THREE

# APPLICATIONS

# 13
# Applications in embedding theory

Even though the Assouad dimension appeared in other contexts, see Section 2.3, it was the pioneering work of Assouad in the 1970s [6, 7, 8], which served to popularise the notion and establish its deep connections with several areas of mathematics, most notably embedding theory. In this chapter we discuss various problems in embedding theory, all of which can be understood as asking for conditions on a metric space $X$ guaranteeing the existence of an injection (the embedding) mapping $X$ into Euclidean space. We will often insist that the embeddings possess some sort of 'metric regularity', such as being bi-Lipschitz. In Section 13.1 we discuss Assouad's embedding theorem, see Theorem 13.1.3, which begins with a discussion of doubling and uniformly perfect metric spaces and the connection between these notions and the Assouad and lower dimensions in Section 13.1.1. In Section 13.2 we consider the 'spiral winding problem', following [85, 93, 159], which asks: given a planar spiral, how regular can a homeomorphism be which maps a line segment onto the spiral? Thus, the spiral winding problem is a specific example of an embedding problem, where the 'embedding' is the inverse of the solution; that is, we ask for regular embeddings of spirals into the line. The Assouad spectrum plays a key role here, see (13.7) and Figure 13.5. Finally, in Section 13.3 we discuss further applications of the Assouad dimension in embedding theory, as studied by Robinson [240]. Here we are interested in embeddings of subsets of Banach space into Euclidean space which are more regular than those provided by Assouad's embedding theorem.

## 13.1 Assouad's embedding theorem

### 13.1.1 Doubling and uniformly perfect metric spaces

Even though we work largely in Euclidean space throughout this book, much of the theory applies in arbitrary metric spaces. One of the fundamental facts concerning Assouad dimension is that it precisely characterises doubling metric spaces. Recall that a metric space $X$ is *doubling* if there exists a constant $K \geqslant 1$ such that any ball in $X$ may be covered by fewer than $K$ balls with half the radius. Moreover, the *Assouad dimension* of a metric space $X$ is defined by

$$\dim_{\mathrm{A}} X = \inf \Bigg\{ \alpha : \text{ there exists a constant } C > 0 \text{ such that, for all } 0 < r < R \text{ and } x \in X, \; N_r\big(B(x,R)\big) \leqslant C\left(\frac{R}{r}\right)^{\alpha} \Bigg\}$$

where $N_r(E)$ is the minimum number of balls of diameter $r$ required to cover $E \subseteq X$.

**Theorem 13.1.1** A metric space $X$ is doubling if and only if $\dim_{\mathrm{A}} X < \infty$.

*Proof* If $\dim_{\mathrm{A}} X < \infty$ then it follows immediately from the definition that $X$ is doubling. For the reverse implication, assume that $X$ is doubling. Therefore, for all $x \in X$ and $\rho > 0$ we can cover $B(x, \rho)$ by $K$ balls of radius $\rho/2$. Fix $0 < r < R$ and let $n$ be the unique integer such that $2^{1-n}R \leqslant r < 2^{2-n}R$. Seeking a cover of $B(x, R)$ by balls of diameter $r$, first cover $B(x, R)$ by $K$ balls of radius $R/2$. Then cover these balls by $K$ balls of radius $R/4$ and continue inductively until we obtain a cover by balls of radius $2^{-n}R \leqslant r/2$. This yields

$$N_r(B(x,R)) \leqslant K^n \leqslant K^2 \left(\frac{R}{r}\right)^{\frac{\log K}{\log 2}},$$

which proves that $\dim_{\mathrm{A}} X \leqslant \frac{\log K}{\log 2} < \infty$ as required. □

Recall Lemma 4.1.1 which states that a measure is doubling if and only if it has finite Assouad dimension. This also holds for measures supported on general metric spaces and combining Theorem 4.1.3 with the characterisation of doubling in terms of Assouad dimension, one sees that every doubling metric space supports a doubling measure. The special case of this fact which implies that every closed subset of Euclidean space supports a doubling measure was known as *Dynkin's conjecture* [59], until it was resolved by Konyagin and Vol'berg [173].

Just as the Assouad dimension characterises the doubling property, the lower dimension characterises the property of being *uniformly perfect.* Recall that a metric space is *perfect* if every point in the space is an accumulation point. Expressed differently, this means that for all points $x$ in $X$ and $0 < R$ we have

$$B(x,R) \setminus \{x\} \neq \varnothing.$$

Uniform perfectness is a quantified version of this where the accumulation can always be taken to be no faster than exponential. A metric space $X$ is *uniformly perfect* if there exists a constant $\kappa \in (0,1)$ such that, for all points $x$ in $X$ and $0 < R < |X|$,

$$B(x,R) \setminus B(x,\kappa R) \neq \varnothing$$

where $|X|$ denotes the diameter of $X$. The *lower dimension* of a metric space $X$ is defined by

$$\dim_{\mathrm{L}} X = \sup \Bigg\{ \alpha : \text{ there exists a constant } C > 0 \text{ such that, for all } 0 < r < R \leqslant |X| \text{ and } x \in X, \ N_r\big(B(x,R)\big) \geqslant C\left(\frac{R}{r}\right)^{\alpha} \Bigg\}$$

where $d$ is the metric on $X$ and $|X| = \sup_{x,y \in X} d(x,y)$ is the diameter of $X$. The following equivalence was proved in [153, Lemma 2.1].

**Theorem 13.1.2** A metric space $X$ is uniformly perfect if and only if $\dim_{\mathrm{L}} F > 0$.

*Proof* If $\dim_{\mathrm{L}} X > s > 0$, then there is a constant $C > 0$ such that for all $x \in X$ and $0 < R < |X|$

$$N_r(B(x,R)) \geqslant C\left(\frac{R}{r}\right)^s .$$

Choose $\kappa \in (0,1)$ such that $C\kappa^{-s} \geqslant 2$ and it follows that

$$B(x,R) \setminus B(x,(\kappa/2)R) \neq \varnothing$$

since otherwise

$$1 = N_{\kappa R}(B(x,R)) \geqslant C \left( \frac{R}{\kappa R} \right)^s \geqslant 2.$$

Therefore $X$ is uniformly perfect.

In the other direction, suppose $X$ is uniformly perfect with constant $\kappa$. Let $x_1 \in X$, $0 < r < R < |X|$ and let $n$ be the unique integer such that $(\kappa/3)^{n+1}R < r \leqslant (\kappa/3)^n R$. Since $X$ is uniformly perfect, there exists a point

$$x_2 \in B(x_1, R) \setminus B(x_1, \kappa R)$$

where, in particular, the balls $B(x_1, (\kappa/3)R)$ and $B(x_2, (\kappa/3)R)$ are disjoint. Repeat the above process within each of these balls to obtain points $x_1, x_2, x_3, x_4$ such that the balls $B(x_i, (\kappa/3)^2 R)$ are pairwise disjoint for $i = 1, 2, 3, 4$. Repeat this process inductively to obtain $2^n$ points $x_i$ such that the balls $B(x_i, (\kappa/3)^n R)$ are pairwise disjoint for $i = 1, \dots, 2^n$ . It follows that any set of diameter $r$ cannot contain more than one of the points $x_i$ and therefore

$$N_r(B(x_1, R)) \geqslant 2^n \geqslant (\kappa/3) \left( \frac{R}{r} \right)^{\frac{\log 2}{\log(3/\kappa)}}$$

which proves that $\dim_{\mathrm{L}} X \geqslant \frac{\log 2}{\log(3/\kappa)} > 0$ as required. □

### 13.1.2 Assouad's embedding theorem

The theorem which first threw the Assouad dimension into the limelight is *Assouad's embedding theorem*, proved during Assouad's doctoral studies [6, 7, 8]. The question Assouad was considering is: when can a given metric space be embedded into Euclidean space? The ideal notion of embedding would be via a bi-Lipschitz map — thus preserving all of the metric properties up to universal constants. It turns out that the Assouad dimension *almost* precisely characterises this problem.

Recall from Lemma 2.4.2 that Assouad dimension is preserved under bi-Lipschitz maps. This already makes the Assouad dimension very useful in studying embedding problems. For example, if we want a bi-Lipschitz embedding into $\mathbb{R}^d$ (or any space with finite Assouad dimension), then having finite Assouad dimension or, equivalently, satisfying the doubling property, see Theorem 13.1.1, is a necessary condition. This is true of other dimensions but, since Assouad dimension is the largest among those we consider, it is the most useful. For example, Troscheit used this approach to prove that the Brownian map cannot be embedded in $\mathbb{R}^d$ using a bi-Lipschitz (or even quasi-symmetric) embedding, see [267]. This was achieved by proving that the Assouad dimension of the Brownian map is almost surely infinite despite the Hausdorff dimension being almost surely 4. The Assouad dimension is also useful in studying purely Euclidean embedding problems. Recall Theorem 9.4.1 which implies that Mandelbrot percolation in $\mathbb{R}^d$ cannot be embedded in $\mathbb{R}^{d-1}$ even if the parameters are such that the Hausdorff and box dimension are almost surely very close to 0.

Assouad's embedding theorem says that being doubling is *almost* sufficient to guarantee a bi-Lipschitz embedding into Euclidean space.

**Theorem 13.1.3** If $(X, d)$ is a metric space with $\dim_A X < \infty$, then for all $\varepsilon \in (0, 1)$ the space $(X, d^\varepsilon)$ embeds into $\mathbb{R}^n$ via a bi-Lipschitz map for some $n \geqslant 1$ depending on $\dim_A X$ and $\varepsilon$.

Note that the metric space $(X, d^\varepsilon)$ is known as an $\varepsilon$-*snowflaking* of $X$ and is simply the set $X$ equipped with the snowflaked metric $d^\varepsilon$ defined by $d^\varepsilon(x, y) = d(x, y)^\varepsilon$. The requirement $\varepsilon < 1$ guarantees the metric $d^\varepsilon$ satisfies the triangle inequality. The snowflaking in Theorem 13.1.3 is necessary, unfortunately, and this is the reason that the Assouad dimension does not precisely characterise the bi-Lipschitz embedding problem. The Heisenberg group $\mathbb{H}$ equipped with the Carnot metric is an example of a doubling metric space which cannot be embedded in Euclidean space of any ambient dimension, see [227, 250]. Another well-known example of a doubling space which cannot embed into Euclidean space via a bi-Lipschitz mapping is a construction of Lang and Plaut [178] based on Laakso graphs, see [177].

Naor and Neiman [215] proved that for $\varepsilon \in (1/2, 1)$, one may choose the $n$ in Theorem 13.1.3 depending only on $\dim_A X$ and not on the snowflaking constant $\varepsilon$. The cost of this improvement was in the dependency on $\varepsilon$ of the distortion constant of the bi-Lipschitz embedding. Tao [264] considered the embedding properties of the Heisenberg group in detail, proving that the $\varepsilon$-snowflakings of $\mathbb{H}$ can be embedded into $\mathbb{R}^n$ with $n$ independent of $\varepsilon$ and with the optimal distortion constant.

## 13.2 The spiral winding problem

The *spiral winding problem* asks whether a line segment can be mapped bijectively onto a given spiral using a certain class of map, most often bi-Lipschitz maps, see [93, 85, 159]. Dimension theory can be of some service here. For example, since Assouad dimension is preserved under bi-Lipschitz maps,[1] if the spiral has Assouad dimension strictly larger than 1, then a line segment certainly cannot be mapped onto the spiral. We shall see that this applies to many polynomial spirals. A more general formulation of the winding problem concerns Hölder maps; that is: when can a line segment be mapped onto a given spiral using a bi-Hölder map and, if this can be done, then what are the optimal Hölder exponents? Dimension theory can still be used here provided one has information about a dimension which behaves reasonably well under Hölder maps. Notably, this rules out the Assouad dimension, recall the discussion in Section 3.4.3. The box dimension works well, but since the Assouad spectrum should be larger than the box dimension and we do have some, albeit complicated, formulae for how the Assouad dimension distorts under Hölder maps, this makes the Assouad spectrum the best candidate for informing this version of the winding problem, see Lemma 3.4.13.

Given a *winding function* $\phi : (1, \infty) \to (0, \infty)$, which we assume is continuous, strictly decreasing, and satisfies $\phi(x) \to 0$ as $x \to \infty$, the associated *spiral* is the set

$$\mathcal{S}(\phi) = \{\phi(x) \exp(\mathrm{i}x) : 1 < x < \infty\} \subseteq \mathbb{C}.$$

If $\phi(x) = e^{-cx}$ for some $c > 0$, then it *is* possible to map $(0, 1)$ onto the spiral $\mathcal{S}(\phi)$ via a bi-Lipschitz map. This example is known as the *logarithmic spiral* (associated to $c > 0$) and the bi-Lipschitz solution to the winding problem in this case goes back to Katznelson, Nag and Sullivan [159]. Moreover, if $\phi$ is sub-exponential, that is, if

$$\frac{\log \phi(x)}{x} \to 0 \qquad (x \to \infty),$$

then this cannot be done, thus illustrating sharpness of the logarithmic family for the bi-Lipschitz formulation of the winding problem, see [85]. This also suggests that the natural family of spirals to consider in the Hölder formulation is the polynomial family $\phi_p$ defined by $\phi_p(x) = x^{-p}$

[1] This is true of any of the other notions of dimension we discuss in this book, but the Assouad dimension is certainly the best choice since it is the greatest.

for $p > 0$. We write $\mathcal{S}_p$ for the spiral associated to $\phi_p$, so

$$\mathcal{S}_p = \mathcal{S}(\phi_p) = \{x^{-p} \exp(\mathrm{i}x) : 1 < x < \infty\} \tag{13.1}$$

for brevity. The spirals $\mathcal{S}_p$ are sometimes referred to as *generalised hyperbolic spirals* and certain special cases show up in various contexts. For example, the *hyperbolic spiral* corresponds to $p = 1$, and the *lituus* family includes the $p = 1/2$ case. A *lituus* is the locus of points $z \in \mathbb{C}$ preserving the area of the circular sector $\{w : \arg(w) \in (0, \arg(z)), |w| \leqslant |z|\}$, where we include multiplicity when $\arg(z) > 2\pi$.

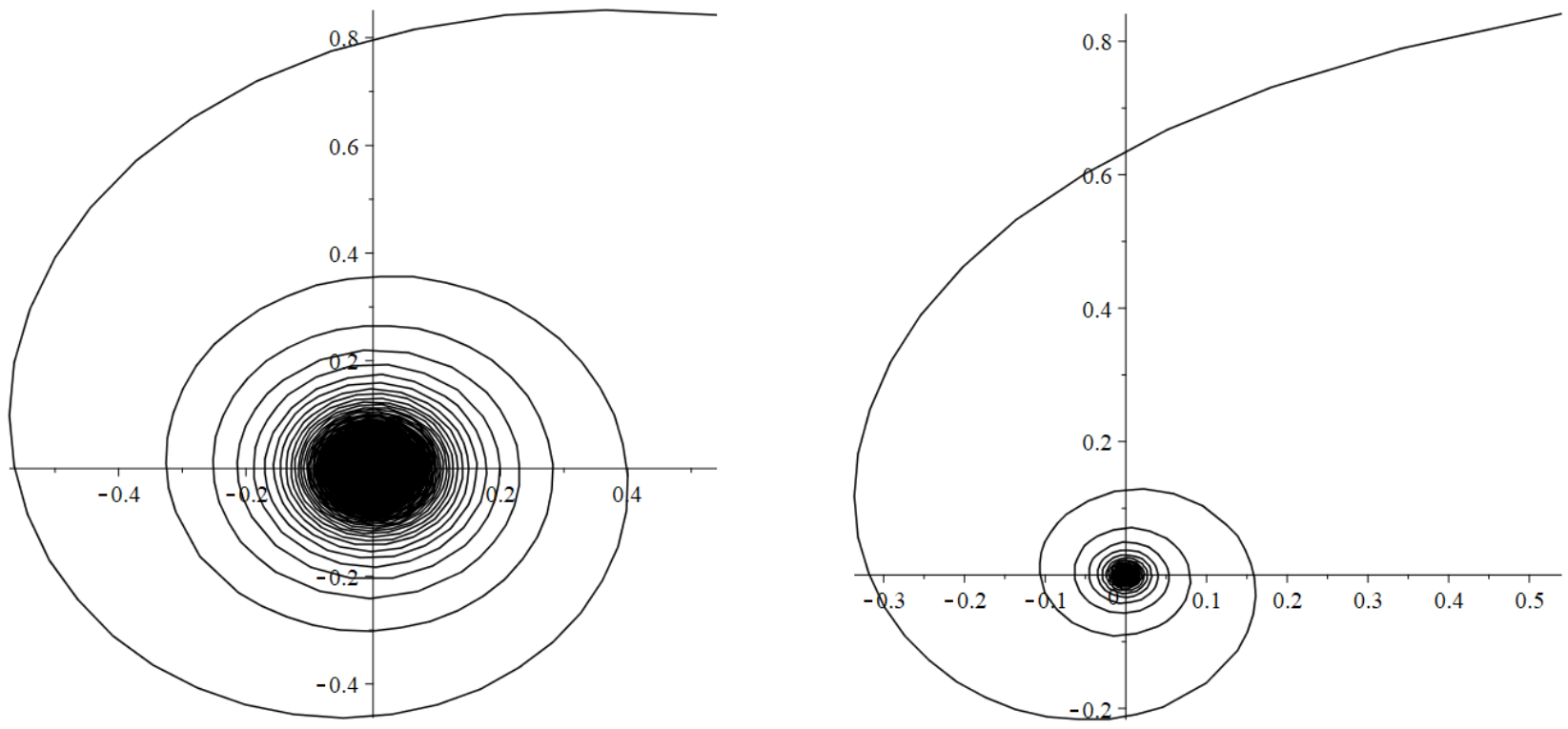

Figure 13.1 Two spirals: on the left $\phi(x) = x^{-1/2}$ (a *lituus*), and on the right $\phi(x) = x^{-1}$ (a hyperbolic spiral).

The Hölder version of the winding problem was resolved for the polynomial spirals $\mathcal{S}_p$ in [93].

**Theorem 13.2.1** If $f : (0, 1) \to \mathcal{S}_p$ is an $(\alpha, \beta)$-Hölder homeomorphism, then

$$\alpha < p$$

and

$$\beta \geqslant \frac{p\alpha}{p - \alpha}.$$

Moreover, for all $\alpha \in [\frac{p}{p+1}, p) \cap (0, 1]$, there exists an $(\alpha, \frac{p\alpha}{p-\alpha})$-Hölder homeomorphism between $(0, 1)$ and $\mathcal{S}_p$.

*Proof* The general bounds are fairly straightforward to obtain and we give the proof below. However, constructing sharp examples is more involved but can be achieved using the map $g_t : (0,1) \to \mathcal{S}_p$ given by

$$g_t(x) = x^{tp} \exp(\mathrm{i}/x^t)$$

where

$$t = \frac{\alpha}{p-\alpha}.$$

To obtain the bounds, decompose $\mathcal{S}_p$ into the disjoint union of 'full turns'

$$\mathcal{S}_p = \bigcup_{k \geqslant 1} \mathcal{S}_p^k$$

where

$$\mathcal{S}_p^k = \{x^{-p} \exp(\mathrm{i}x) : 1 + 2\pi(k-1) < x \leqslant 1 + 2\pi k\} \tag{13.2}$$

for each integer $k \geqslant 1$. Also, given a homeomorphism $f : (0,1) \to \mathcal{S}_p$, we decompose $(0,1)$ into the corresponding half-open intervals

$$\mathcal{I}^k = f^{-1}(\mathcal{S}_p^k). \tag{13.3}$$

It is convenient to extend $f$ continuously such that the domain includes 0 and 1 and we assume without loss of generality that $f(0) = 0$. We have

$$\frac{1}{8^p} k^{-p} \leqslant |\mathcal{S}_p^k| = |f(\mathcal{I}^k)| \leqslant C|\mathcal{I}^k|^\alpha \tag{13.4}$$

where $C \geqslant 1$ is the constant coming from the Hölder condition on $f$, see (3.10). Therefore

$$1 = \sum_{k=1}^{\infty} |\mathcal{I}^k| \geqslant \frac{1}{8^{p/\alpha} C^{1/\alpha}} \sum_{k=1}^{\infty} k^{-p/\alpha} \tag{13.5}$$

which forces $\alpha < p$.

For each integer $l \geqslant 1$, let

$$x_l = \sum_{k=l}^{\infty} |\mathcal{I}^k|$$

as in Figure 13.2. Combining this with (13.4) and (13.5) yields

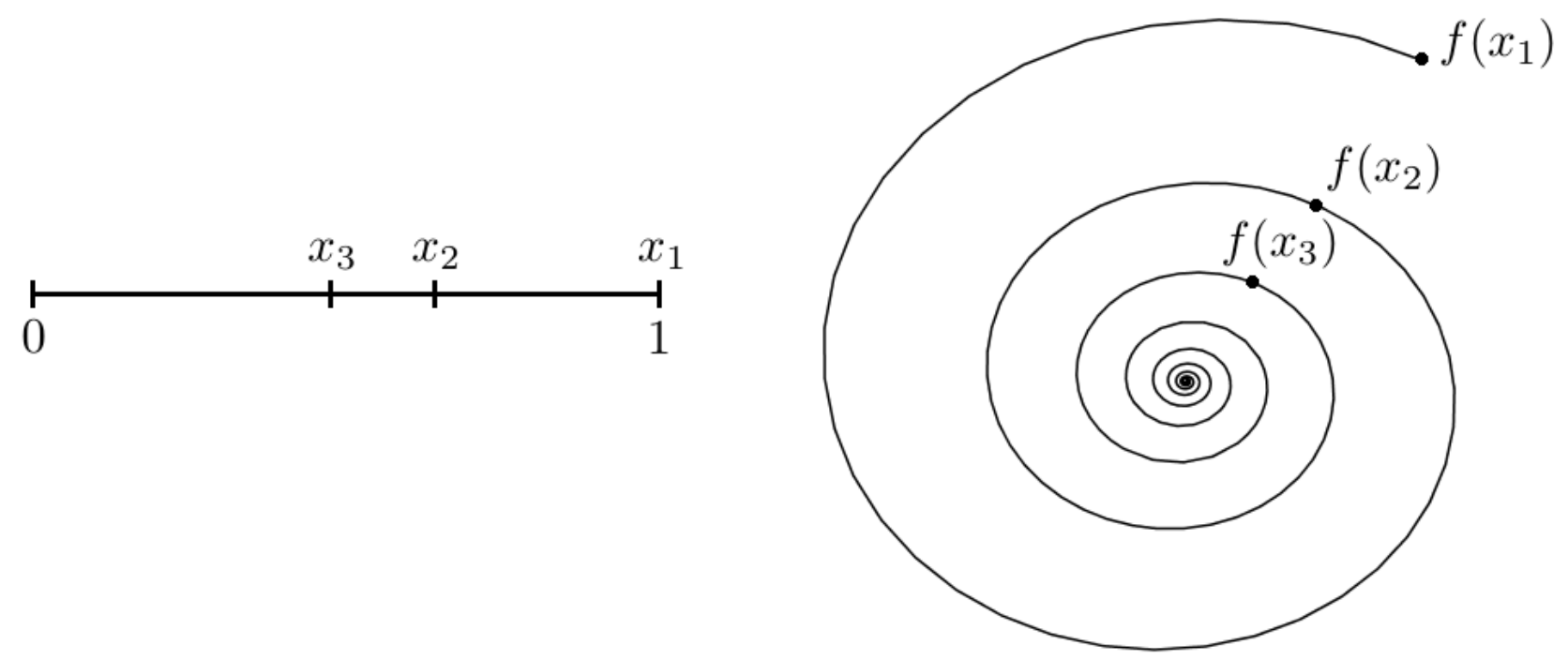

Figure 13.2 The location of the points $x_1, x_2, x_3 \in (0,1]$ and their images in $\mathcal{S}_p$.

$$\begin{aligned}
\frac{1}{C} \leqslant \frac{|f(x_l)-f(0)|}{|x_l|^\beta} = \frac{|f(x_l)|}{\left(\sum_{k=l}^\infty |\mathcal{I}^k|\right)^\beta} &\leqslant (8^{p/\alpha}C^{1/\alpha})^\beta \frac{l^{-p}}{\left(\sum_{k=l}^\infty k^{-p/\alpha}\right)^\beta} \\
&\leqslant 8^{p\beta/\alpha}C^{\beta/\alpha}\frac{l^{-p}}{\left(\int_{k=l}^\infty z^{-p/\alpha}dz\right)^\beta} \\
&= 8^{p\beta/\alpha}C^{\beta/\alpha}(p/\alpha-1)^\beta \frac{l^{-p}}{l^{(1-p/\alpha)\beta}} \\
&\to 0
\end{aligned}$$

as $l \to \infty$, if $-p-(1-p/\alpha)\beta < 0$. This forces

$$\beta \geqslant \frac{p\alpha}{p-\alpha}$$

as required.[2] □

We now turn our attention to the dimensions of spirals and the resulting applications to the winding problem. We begin with the box dimensions of $\mathcal{S}_p$. These are strictly greater than 1 for $p \in (0,1)$, which

[2] The sharp relationship between $\alpha$ and $\beta$ is somewhat reminiscent of *Sobolev conjugates*. For $1 \leqslant p < d$, the *Sobolev embedding theorem* states that $W^{1,p}(\mathbb{R}^d) \subseteq L^q(\mathbb{R}^d)$ where

$$q = \frac{dp}{d-p}$$

is the Sobolev conjugate of $p$. Here $W^{1,p}(\mathbb{R}^d)$ is the *Sobolev space* consisting of real-valued functions $f$ on $\mathbb{R}^d$ such that both $f$ and all weak derivatives of $f$ are in $L^p(\mathbb{R}^d)$.

therefore yields non-trivial information in this range. The following result can be found in [265, 274, 93]. See [284] for a treatment of the box dimensions of spirals in $\mathbb{R}^3$.

**Theorem 13.2.2** For $p \in (0,1)$

$$\dim_{\mathrm{B}} \mathcal{S}_p = \frac{2}{1+p}$$

and for $p \geqslant 1$, $\dim_{\mathrm{B}} \mathcal{S}_p = 1$.

*Proof* We only prove the result for $p \in (0,1)$ and refer the reader to [93] for the other cases. Let $r \in (0,1)$ and $k(r)$ be the unique positive integer satisfying

$$k(r)^{-(p+1)} \leqslant r < (k(r)-1)^{-(p+1)},$$

noting that $r^{-1/(p+1)} \leqslant k(r) < r^{-1/(p+1)} + 1 \leqslant 2r^{-1/(p+1)}$. The importance of $k(r)$ is that the sets $\mathcal{S}_p^k$ for $k \leqslant k(r)$ are not 'wound tightly' when viewed at scale $r$ and therefore the naïve covering strategy where one places balls of diameter $r$ along the spiral is likely to be optimal, see Figure 13.3. Here $\mathcal{S}_p^k$ is as in (13.2). Moreover, covering the curves $\mathcal{S}_p^k$ for $k > k(r)$ altogether will not be any easier than covering the ball $B(0, k(r)^{-p})$, which can be done with fewer than $c\left(k(r)^{-p}/r\right)^2$ many balls of diameter $r$ for an absolute constant $c$. This means

$$\begin{aligned} N_r(\mathcal{S}_p) &\leqslant N_r\left(B(0,k(r)^{-p})\right) + \sum_{k=1}^{k(r)} N_r\left(\mathcal{S}_p^k\right) \\ &\leqslant c\left(\frac{k(r)^{-p}}{r}\right)^2 + \pi \sum_{k=1}^{k(r)} \frac{k^{-p}}{r} \\ &\leqslant (c+\pi) r^{-\frac{2}{1+p}} \end{aligned}$$

which proves the upper bound. Moreover, the lower bound also follows since

$$N_r(\mathcal{S}_p) \geqslant N_r\left(\mathcal{S}_p \cap B(0,k(r)^{-p})\right) \geqslant c_0 \left(\frac{k(r)^{-p}}{r}\right)^2 \geqslant c_0 2^{-p} r^{-\frac{2}{1+p}}$$

for another absolute constant $c_0$. Here we have used the definition of $k(r)$ to ensure that we need a comparable number of $r$-balls to cover $\mathcal{S}_p \cap B(0, k(r)^{-p})$ as we need to cover $B(0, k(r)^{-p})$. □

Theorem 13.2.2 directly gives non-trivial information concerning the

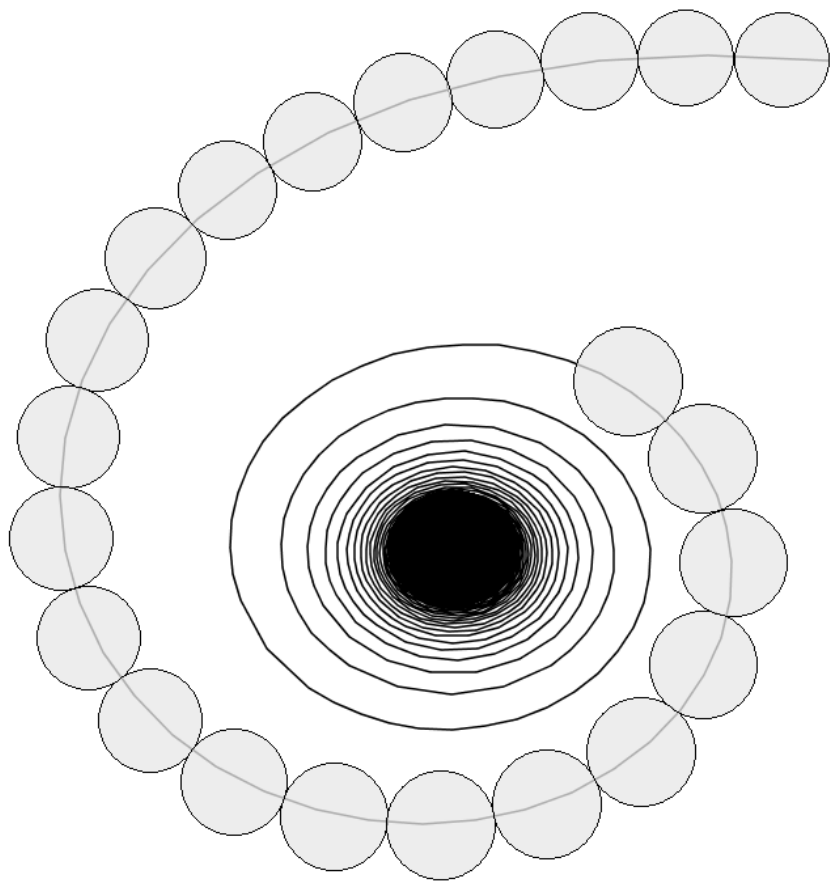

Figure 13.3 An efficient cover for part of a spiral which is not wound as tightly as the covering scale.

winding problem: if $f : (0,1) \to \mathcal{S}_p$ is an onto $\alpha$-Hölder map and $p \in (0,1)$, then

$$\alpha \leqslant \frac{p+1}{2} < 1 \tag{13.6}$$

and therefore $f$ cannot be Lipschitz.

Next we consider the Assouad spectrum. The following result was proved in [93], although the case when $p \in (0,1)$ is derivable from [111, Theorem 7.2], which also included more general sub-exponential spirals. Note that for $p \leqslant 1$, the Assouad spectrum is given by the general upper bound from Lemma 3.4.4 but, for $p > 1$, it is not.

**Theorem 13.2.3** For $p \in (0,1)$ and $\theta \in (0,1)$,

$$\dim_{\mathrm{A}}^{\theta} \mathcal{S}_p = \min\left\{\frac{2}{(1+p)(1-\theta)}, 2\right\}$$

and, for $p \geqslant 1$ and $\theta \in (0,1)$,

$$\dim_{\mathrm{A}}^{\theta} \mathcal{S}_p = \min\left\{1 + \frac{\theta}{p(1-\theta)}, 2\right\}.$$

*Proof* We prove the result for $p > 1$, and refer the reader to [93] for the $p \in (0,1]$ case. It suffices to consider $0 < \theta < \frac{p}{1+p}$ since, for $\theta' = \frac{p}{1+p}$,

$$1 + \frac{\theta'}{p(1-\theta')} = 2,$$

and so it follows by continuity of the Assouad spectrum that $\dim_{\mathrm{A}}^{\theta'} \mathcal{S}_p = 2$ and therefore $\dim_{\mathrm{A}}^{\theta} \mathcal{S}_p = 2$ for all $\theta > \frac{p}{1+p}$ by Corollary 3.3.3. We first provide the upper bound.

Let $r \in (0,1)$ and $L(r)$, $l(r)$ be the unique positive integers satisfying

$$L(r)^{-(p+1)} \leqslant r < (L(r)-1)^{-(p+1)}.$$

and

$$l(r)^{-p} \leqslant r^{\theta} < (l(r)-1)^{-p},$$

respectively. Certainly $L(r) > l(r)$ for all sufficiently small $r$ since $\theta < \frac{p}{p+1}$. Note that for some uniform constant $c \geqslant 1$

$$N_r\left(B(z, r^{\theta}) \cap \mathcal{S}_p\right) \;\leqslant\; cN_r\left(B(0, r^{\theta}) \cap \mathcal{S}_p\right)$$

for all $z \in \mathcal{S}_p$ and so it suffices to only consider $z = 0$. Arguing as in the proof of Theorem 13.2.2, where $\mathcal{S}_p^k$ is as in (13.2),

$$\begin{aligned} N_r\left(B(0, r^{\theta}) \cap \mathcal{S}_p\right) &\leqslant N_r\left(B(0, L(r)^{-p}) \cap \mathcal{S}_p\right) \;+\; \sum_{k=l(r)}^{L(r)} N_r\left(\mathcal{S}_p^k\right) \\ &\leqslant c\left(\frac{L(r)^{-p}}{r}\right)^2 \;+\; \pi \sum_{k=l(r)}^{L(r)} \frac{k^{-p}}{r} \\ &\leqslant c\left(\frac{r^{\theta}}{r}\right)^{\frac{2}{(1-\theta)(p+1)}} \;+\; 2\pi\left(\frac{r^{\theta}}{r}\right)^{1+\frac{\theta}{p(1-\theta)}} \end{aligned}$$

which proves that $\dim_{\mathrm{A}}^{\theta} \mathcal{S}_p$ is bounded above by the maximum of

$$\frac{2}{(1-\theta)(p+1)} \qquad \text{and} \qquad 1 + \frac{\theta}{p(1-\theta)}$$

which is always the latter since $p > 1$.

For the lower bound, if one follows the proof of Theorem 13.2.2, then one might use the estimate

$$\begin{aligned} N_r\left(B(0, r^{\theta}) \cap \mathcal{S}_p\right) \geqslant N_r\left(\mathcal{S}_p \cap B(0, L(r)^{-p})\right) &\geqslant c_0\left(\frac{L(r)^{-p}}{r}\right)^2 \\ &\geqslant c_0 2^{-p}\left(\frac{r^{\theta}}{r}\right)^{\frac{2}{(1-\theta)(p+1)}} \end{aligned}$$

which shows

$$\dim_{\mathrm{A}}^{\theta} \mathcal{S}_p \;\geqslant\; \frac{2}{(1-\theta)(p+1)}$$

which is *not* sharp. The sharp bound is achieved by considering the 'other part' of $\mathcal{S}_p$. Since any ball of diameter $r$ can intersect at most two of the sets $\mathcal{S}_p^k$ provided $l(r) \leqslant k \leqslant L(r)$, see (13.2),

$$N_r\left(B(0, r^\theta) \cap \mathcal{S}_p\right) \geqslant \frac{1}{2} \sum_{k=l(r)}^{L(r)} N_r\left(\mathcal{S}_p^k\right) \geqslant c_1 \sum_{k=l(r)}^{L(r)} \frac{k^{-p}}{r}$$

for a uniform constant $c_1$. This can be estimated from below for small $r$ by

$$\begin{aligned} c_1 r^{-1} \int_{l(r)}^{L(r)} z^{-p} dz &= \frac{c_1}{p-1} r^{-1} \left(l(r)^{1-p} - L(r)^{1-p}\right) \\ &\geqslant \frac{c_1}{2(p-1)} r^{-1-\frac{\theta(1-p)}{p}} \\ &= \frac{c_1}{2(p-1)} \left(\frac{r^\theta}{r}\right)^{1+\frac{\theta}{p(1-\theta)}} \end{aligned}$$

which gives the desired lower bound. □

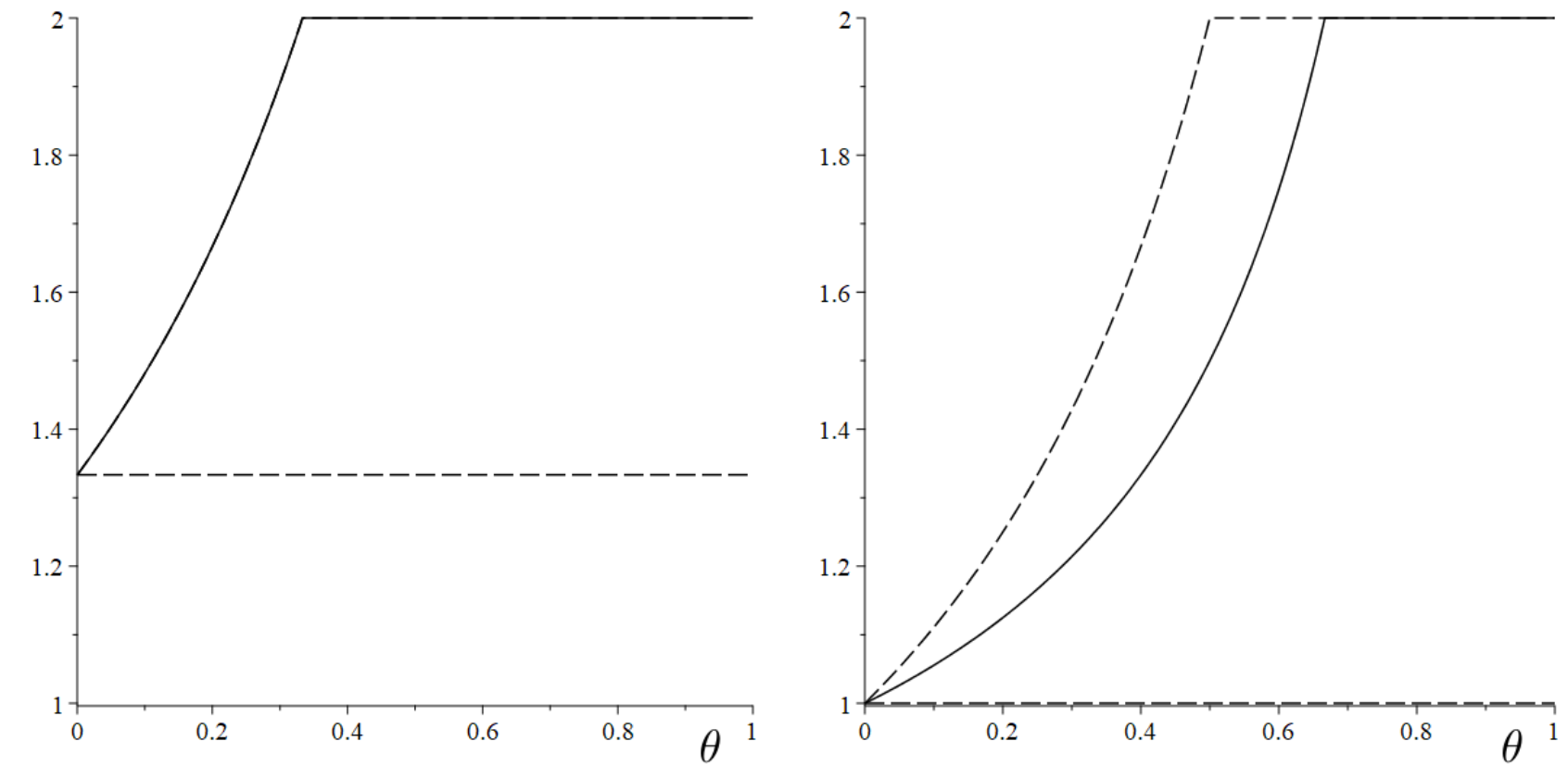

Figure 13.4 Plots of $\dim_{\mathrm{A}}^\theta \mathcal{S}_p$ (solid) along with the general upper and lower bounds from Lemma 3.4.4 (dashed) for comparison. On the left, $p = 1/2$, and on the right, $p = 2$.

The following Corollary follows immediately from Theorem 13.2.3, and was proved in the range $p \in (0, 1)$ in [111].

**Corollary 13.2.4** For all $p > 0$, $\dim_{\mathrm{A}} \mathcal{S}_p = \dim_{\mathrm{qA}} \mathcal{S}_p = 2$.

The Assouad dimension does not behave well under Hölder maps, and so we cannot derive any information regarding the spiral winding problem from knowledge of the Assouad dimension, despite it being as large as possible. However, the Assouad spectrum is more regular and *can* be controlled, recall Lemma 3.4.13 and Corollary 3.4.14.

Suppose $f : (0,1) \to \mathcal{S}_p$ is an $(\alpha, \beta)$-Hölder homeomorphism, in which case $f^{-1} : \mathcal{S}_p \to (0,1)$ is a $(\beta^{-1}, \alpha^{-1})$-Hölder homeomorphism, and

$$\rho = \inf\{\theta \in (0,1) : \dim_{\mathrm{A}}^{\theta} \mathcal{S}_p = \dim_{\mathrm{A}} \mathcal{S}_p\} = \frac{p}{p+1}.$$

Applying Corollary 3.4.14 we get

$$1 = \dim_{\mathrm{qA}} f^{-1}(\mathcal{S}_p) \geqslant \frac{(1-\rho)\dim_{\mathrm{qA}} \mathcal{S}_p}{\alpha^{-1} - \beta^{-1}\rho} = \frac{2(1-\frac{p}{p+1})}{\alpha^{-1} - \beta^{-1}\frac{p}{p+1}}$$

and rearranging gives

$$\beta \geqslant \frac{p\alpha}{1+p-2\alpha}$$

or, formulated in terms of $\alpha$,

$$\alpha \leqslant \frac{p\beta + \beta}{p + 2\beta}. \tag{13.7}$$

Comparing (13.6), (13.7), and the sharp result from Theorem 13.2.1 which, formulated in terms of $\alpha$, is

$$\alpha \leqslant \frac{p\beta}{p+\beta}$$

we see that neither (13.6) nor (13.7) are sharp but the Assouad spectrum 'outperforms' the box dimension. However, as we relax the restrictions on the inverse map by letting $\beta \to \infty$, the estimates obtained from the Assouad spectrum approach those obtained by the box dimension.

## 13.3 Almost bi-Lipschitz embeddings

Recall Assouad's embedding theorem, Theorem 13.1.3. It turns out that if the metric space we want to embed is a subset of Banach space and the Assouad dimension of the *difference set* is finite, then we can upgrade the conclusion of Theorem 13.1.3 to yield more regular embeddings. This sort of problem is considered in detail in [240]. Given a subset $X$ of a Banach space $\mathcal{B}$, the *difference set* is defined by

$$X - X = \{x - y : x, y \in X\} \subseteq \mathcal{B}.$$

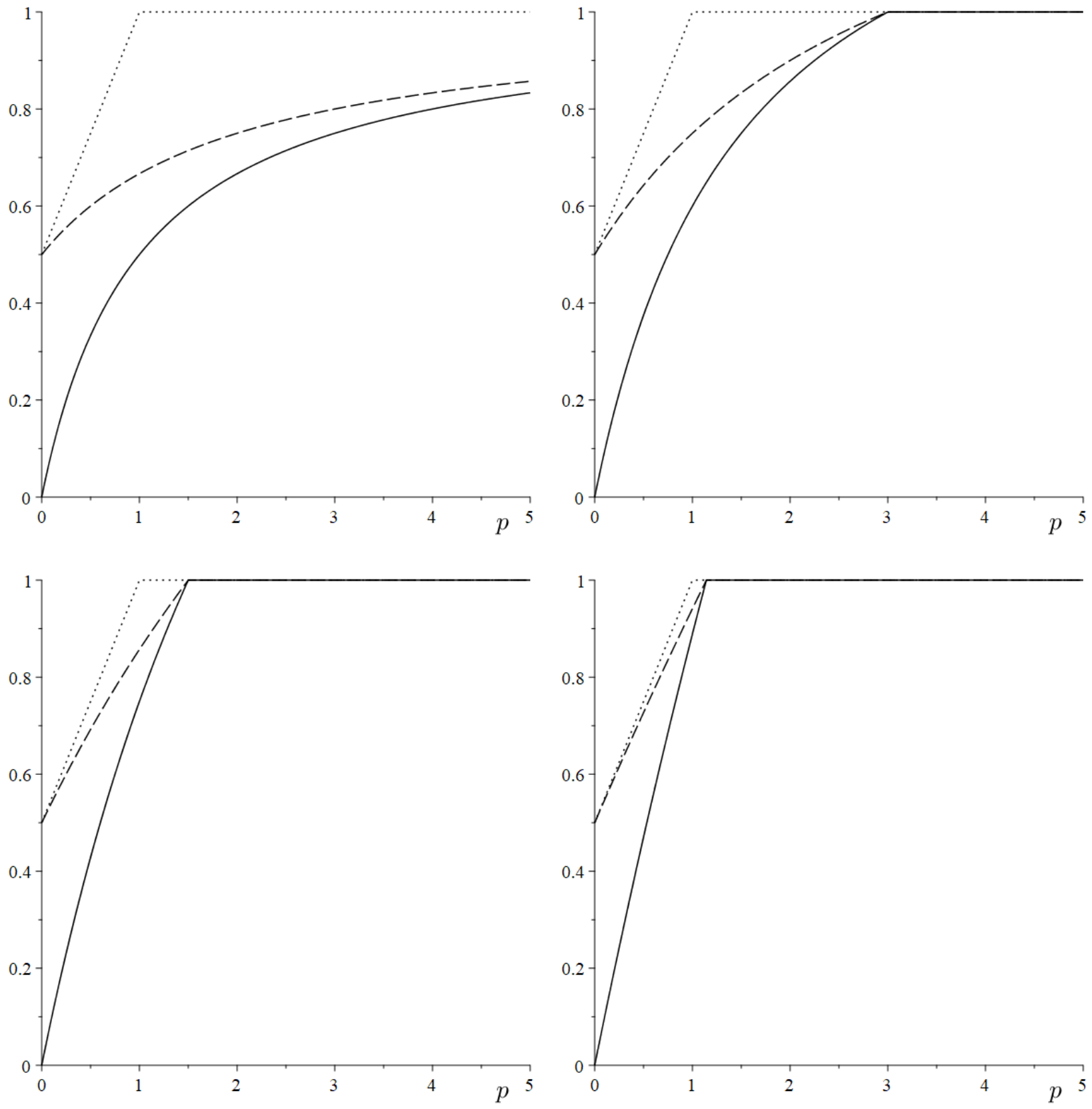

Figure 13.5 Suppose there exists an $(\alpha, \beta)$-Hölder homeomorphism between $(0,1)$ and $\mathcal{S}_p$. We plot upper bounds for $\alpha$ as a function of $p$ for fixed $\beta$. The sharp bounds are plotted as a solid line, the bounds obtained from the Assouad spectrum are plotted as a dashed line, and the bounds obtained from the box dimension are plotted as a dotted line. On the top left, we fix $\beta = 1$, on the top right, $\beta = 3/2$, on the bottom left, $\beta = 3$, and on the bottom right, $\beta = 8$.

The following embedding theorem for subsets of Banach space was proved by Robinson [241], see also [235].

**Theorem 13.3.1** Let $X$ be a compact subset of a real Banach space $\mathcal{B}$ and suppose that

$$\dim_{\mathrm{A}}(X - X) < s < d$$

with $d \in \mathbb{N}$ and

$$\gamma > \frac{d+1}{d-s}.$$

Then there exists a linear map $f : \mathcal{B} \to \mathbb{R}^d$ which is injective and $\gamma$-almost bi-Lipschitz on $X$; that is, there exists $c, \rho > 0$ such that

$$c^{-1} \frac{\|x-y\|}{|\log \|x-y\||^\gamma} \leqslant |f(x) - f(y)| \leqslant c\|x-y\|$$

for all $x, y \in X$ with $\|x-y\| \leqslant \rho$.

In fact, this theorem can be strengthened in several ways, for example to *almost homogeneous* $X - X$ and the existence of the embedding $f$ can be upgraded to a generic statement, see [240, Theorem 9.18]. Note that the difference set $X - X$ is an orthogonal projection, see Section 10.2, of the product $X \times X$ and if $\dim_A X < \infty$, then $\dim_A X \times X = 2 \dim_A X < \infty$ by results in Section 10.1. However, despite orthogonal projections being Lipschitz, the Assouad dimension of $X - X$ may be larger than that of $X \times X$. In fact, [240, Lemma 9.12] constructs an example of a set $X$ in a Banach space such that $\dim_A X = 0$ but $\dim_A(X - X) = \infty$. This highlights the fact that the assumption that the difference set is doubling in Theorem 13.3.1 is much stronger than the assumption that $X$ is doubling in Theorem 13.1.3. Then again, the conclusion that there is a linear almost bi-Lipschitz embedding into Euclidean space is stronger than the conclusion of Theorem 13.1.3 that every snowflaking embeds via a bi-Lipschitz map.

Margaris [199, 200] later observed that the assumption that $X - X$ is doubling in Theorem 13.3.1 can be weakened in a meaningful way. In fact, all that is needed is that $X - X$ be 'doubling near the origin'. Given $s > 0$, we say $X - X$ is *$s$-doubling near the origin* if there exists $C \geqslant 1$ such that for all $0 < r < R$

$$N_r\Big(B(0,R) \cap (X - X)\Big) \leqslant C \left(\frac{R}{r}\right)^s .$$

**Theorem 13.3.2** Let $X$ be a compact subset of a real Banach space $\mathcal{B}$ and suppose $X - X$ is $s$-doubling near the origin. Then, given any $\gamma > 1$, there exists $d \in \mathbb{N}$ and a linear map $f : \mathcal{B} \to \mathbb{R}^d$ which is injective and $\gamma$-almost bi-Lipschitz on $X$, that is, there exists $c > 0, \rho > 0$ such that

$$c^{-1} \frac{\|x-y\|}{|\log \|x-y\||^\gamma} \leqslant |f(x) - f(y)| \leqslant c\|x-y\|$$

for all $x, y \in X$ with $\|x-y\| \leqslant \rho$.

# 14
# Applications in number theory

Given the origins and motivation for studying the Assouad dimension, applications in embedding theory are perhaps expected. In this chapter we discuss a collection of, perhaps more surprising, applications to several distinct problems in number theory. In Section 14.1 we consider the problem of finding arithmetic progressions inside sets of integers. Recall, for example, the Erdős conjecture on arithmetic progressions, of which the celebrated Green-Tao Theorem [123] is a special case. It turns out that the Assouad dimension can be used to prove a weak 'asymptotic' version of this conjecture, first proved in [110]. In Section 14.2 we discuss some problems in Diophantine approximation, the general area concerned with how well real numbers can be approximated by rationals. Here there is a well-developed connection between estimating the Hausdorff dimension of sets of badly approximable numbers and the lower dimension, see work of Das, Fishman, Simmons and Urbański [45]. Finally, in Section 14.3 we discuss an elegant use of the Assouad dimension in problems of definability of the integers due to Hieronymi and Miller [131].

## 14.1 Arithmetic progressions

Arithmetic progressions are among the most natural and well-studied mathematical objects. Determining whether a given set contains an arithmetic progression is often a subtle problem and has been studied in many contexts. Perhaps the most famous example of this is *Szemerédi's theorem*, which states that a set of positive integers with positive upper natural density contains arbitrarily long arithmetic progressions. The *upper natural density* of a set $X \subseteq \mathbb{Z}^+$ is defined as $\limsup_{n\to\infty} \frac{\#X\cap\{1,\ldots,n\}}{n}$. The *Erdős conjecture on arithmetic progressions* suggests that more should be true.

**Conjecture 14.1.1** If $X \subseteq \mathbb{Z}^+$ is such that

$$\sum_{x\in X} 1/x = \infty,$$

then $X$ should contain arbitrarily long arithmetic progressions.

This conjecture is sometimes attributed to Erdős-Turán, following [62]. A celebrated breakthrough in this direction is the Green-Tao theorem [123], which asserts that the *primes* contain arbitrarily long arithmetic progressions. Recall that the reciprocals of the primes form a divergent series, and so this is a special case of the Erdős conjecture. It turns out that the Assouad dimension is intimately related to 'approximate' arithmetic structure, and in this section we show how it can be used to prove an approximate version of the Erdős conjecture, following [110].

An arithmetic progression of *length* $k \geqslant 1$ and *gap size* $\Delta > 0$ is a set of the form

$$\{x,\, x+\Delta,\, x+2\Delta,\, \ldots,\, x+(k-1)\Delta\},$$

for some $x \in \mathbb{R}$, that is, a collection of $k$ points each separated from the next by a common distance $\Delta$. We consider the following weakening of the property of containing arithmetic progressions, which was introduced and studied in [110, 107].

**Definition 14.1.2** A set of positive integers $X \subseteq \mathbb{Z}^+$ *gets arbitrarily close to arbitrarily long arithmetic progressions* if, for all $k \in \mathbb{N}$ and $\varepsilon > 0$, there exists an arithmetic progression $P$ of length $k$ and gap size $\Delta > 0$ such that

$$\sup_{p\in P} \inf_{x\in X} |p - x| \leqslant \varepsilon\Delta.$$

This definition should be understood as saying that, for arbitrarily large $k$ and arbitrarily small $\varepsilon > 0$, $X$ gets within $\varepsilon$ of an arithmetic progression of length $k$. The fact that $\varepsilon\Delta$ appears instead of $\varepsilon$ is the necessary normalization. This is based on the observation that all arithmetic progressions of length $k$ are essentially the same: they are all equal to $\{1, 2, \ldots, k\}$ upon rescaling and translation.

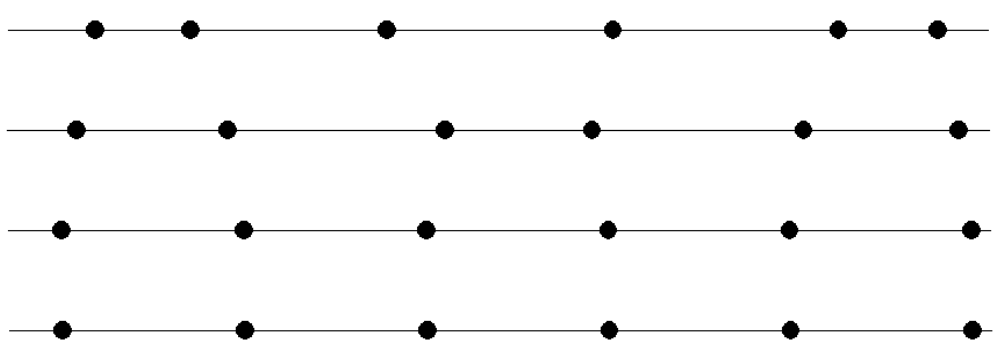

Figure 14.1 From top row to bottom row: three approximations to an arithmetic progression of length 6, where $\varepsilon$ is 1/3, 1/10, 1/150, respectively, followed by a genuine arithmetic progression of length 6. At this resolution the $\varepsilon = 1/150$ approximation is indistinguishable from the genuine arithmetic progression.

**Theorem 14.1.3** If $X \subseteq \mathbb{Z}^+$ is such that

$$\sum_{x \in X} 1/x = \infty,$$

then $X$ gets arbitrarily close to arbitrarily long arithmetic progressions.

To the lay audience this theorem may have the same aesthetic appeal as resolution of the full conjecture, however, this result is actually rather straightforward to prove. The theorem was proved in [110, Theorem 2.11], see also [91]. We begin by demonstrating that reciprocals forming a divergent series is enough to guarantee full Assouad dimension — even maximal Assouad spectrum. The fact that the Assouad dimension is 1 was shown in [110, Lemma 2.10], and this is all that is required to prove Theorem 14.1.3. However, we go slightly further here by computing the Assouad spectrum. We write $1/X = \{1/x : x \in X\}$.

**Lemma 14.1.4** If $X \subseteq \mathbb{Z}^+$ is such that

$$\sum_{x \in X} 1/x = \infty,$$

then, for all $\theta \in (0,1)$,

$$\dim_{\mathrm{A}}^{\theta}(1/X) = \min\left\{\frac{1/2}{1-\theta},\ 1\right\}.$$

Therefore, $\overline{\dim}_{\mathrm{B}}(1/X) = 1/2$, $\dim_{\mathrm{A}}(1/X) = \dim_{\mathrm{qA}}(1/X) = 1$, and $\dim_{\mathrm{A}} X = 1$.

*Proof* Let $X = \{x_n\}_{n \in \mathbb{N}} \subseteq \mathbb{N}$ with $x_n$ strictly increasing. To prove the result, it suffices to show that $\dim_{\mathrm{A}}^{1/2}(1/X) \geqslant 1$, noting that the result follows from this together with Lemma 3.4.6 and the fact that $\overline{\dim}_{\mathrm{B}}(1/X) \leqslant \overline{\dim}_{\mathrm{B}}(1/\mathbb{N}) = 1/2$, recall Theorem 2.1.1. The fact that $\dim_A X = 1$ follows since the Assouad dimension is preserved under the inversion $x \mapsto 1/x$, see [192, Theorem A.10 (1)].

In order to reach a contradiction, suppose that this is false, that is, $\dim_A^{1/2}(1/X) < s < 1$ for some $s$. It follows that there exists $C > 0$ such that, for all integer $k \geqslant 0$,

$$N_{4^{-k}}\left(B(0,2^{-k}) \cap 1/X\right) \leqslant C\left(\frac{2^{-k}}{4^{-k}}\right)^s = C2^{ks}.$$

If $x_n, x_{n+1} \leqslant 2^{k+1}$ for some $k$, then

$$1/x_n - 1/x_{n+1} > 4^{-(k+1)}$$

and so no open set of diameter $4^{-k}$ can cover the reciprocals of more than four distinct points in $X_k := X \cap [2^k, 2^{k+1}]$. It follows that

$$|X_k| \leqslant 4N_{4^{-k}}\left(B(0,2^{-k}) \cap 1/X\right) \leqslant 4C\,2^{ks}.$$

Therefore

$$\sum_{x \in X} 1/x \leqslant \sum_{k=0}^{\infty} \sum_{n: x_n \in X_k} 1/x_n \leqslant \sum_{k=0}^{\infty} |X_k| 2^{-k} \leqslant 4C \sum_{k=0}^{\infty} 2^{(s-1)k} < \infty,$$

a contradiction. □

Having full Assouad dimension is a strong property, which, in particular, guarantees that $[0,1]$ is a weak tangent to $X$ by Theorem 5.1.5. This is enough information to approximate arithmetic structure.

**Lemma 14.1.5** If $[0,1]$ is a weak tangent to $X \subseteq \mathbb{Z}^+$, then $X$ gets arbitrarily close to arbitrarily long arithmetic progressions.

*Proof* Let $k \geqslant 1$ be an integer and $\varepsilon \in (0,1)$. Since $[0,1]$ is a weak tangent to $X$ we can find an interval $[n, n+k\Delta]$ for integers $n, \Delta \geqslant 1$ such that

$$d_{\mathcal{H}}\left((k\Delta)^{-1}[n, n+k\Delta] \cap X - (k\Delta)^{-1}n,\ [0,1]\right) \leqslant \varepsilon/k$$

and therefore, upon rescaling and translating,

$$d_{\mathcal{H}}\left([n, n+k\Delta] \cap X,\ [n, n+k\Delta]\right) \leqslant \varepsilon\Delta.$$

It follows that for each $p \in P = \{n, n+\Delta, n+2\Delta, \ldots, n+k\Delta\}$ there must exist a point $x \in X$ such that $|p - x| \leqslant \varepsilon\Delta$, which proves the result. □

Theorem 14.1.3 then follows from Lemma 14.1.4, Theorem 5.1.5, and Lemma 14.1.5.

One may bound the Assouad dimension of a set above in terms of how far away it is from approximating arithmetic progressions of a given length. This is explored in detail in [107]. A quantitative improvement of Theorem 14.1.3 was proved in [113] where the error $o(\Delta)$ is replaced with $O(\Delta^\alpha)$ for any $\alpha \in (0,1)$.

### 14.1.1 Discrete Kakeya sets

Lemma 14.1.5 shows that full Assouad dimension guarantees the existence of almost arithmetic progressions. This also holds for arbitrary subsets of $\mathbb{R}$ and even in higher dimensions or finite dimensional Banach spaces with 'arithmetic progression' replaced by 'arithmetic patch', see [110]. Rather than looking for arithmetic patches, one may consider a discrete analogue of the Kakeya problem, recall the (continuous) Kakeya problem which was discussed in Section 11.2. For example, say $K \subseteq \mathbb{R}^d$ is a *discrete approximate Kakeya set* if for all $\varepsilon > 0$, $k \geqslant 1$ and directions $e \in \mathbb{S}^{d-1}$, there exists an isometric embedding $S : \mathbb{R} \to \mathbb{R}^d$ and an arithmetic progressions $P \subseteq \mathbb{R}$ of length $k$ and gap $\Delta$ such that

$$\sup_{p \in P} \inf_{x \in K} |S(p) - x| \leqslant \varepsilon \Delta$$

and $S(\mathbb{R})$ is in direction $e$. We say an isometric embedding $V \subseteq \mathbb{R}^d$ of $\mathbb{R}$ is in *direction* $e \in \mathbb{S}^{d-1}$ if $V$ is a translate of the span of $e$. One may ask, see [107, 245], whether or not being a discrete approximate Kakeya set guarantees full Assouad dimension, that is, $\dim_A K = d$. However, Saito [245] proved that this is false. In fact, Saito constructed a compact discrete approximate Kakeya set $K \subseteq \mathbb{R}^d$ with $\dim_A K = 1$. This is clearly sharp, since such sets must have a line segment as a weak tangent.

## 14.2 Diophantine approximation

Diophantine approximation concerns how well real numbers can be approximated by rationals. One of the fundamental results is Dirichlet's theorem, which states that for all $x \in \mathbb{R}^d$ there exist infinitely many $\mathbf{p} = (p_1, \ldots, p_d) \in \mathbb{Z}^d$ and $q \in \mathbb{N}$ such that

$$\|x - \mathbf{p}/q\| < \frac{1}{q^{1+1/d}} \tag{14.1}$$

where $\|\cdot\|$ is the supremum norm on $\mathbb{R}^d$. Importantly, the entries in the rational approximants $\mathbf{p}/q$ share a common denominator $q \in \mathbb{N}$ which quantifies the level of approximation. Certain numbers can be approximated much better than (14.1), but others cannot. The *badly approximable numbers* are the collection of points $x \in \mathbb{R}^d$ such that there is a constant $c > 0$ such that for all $\mathbf{p} = (p_1, \ldots, p_d) \in \mathbb{Z}^d$ and $q \in \mathbb{N}$

$$\|x - \mathbf{p}/q\| \geqslant \frac{c}{q^{1+1/d}}.$$

We denote the set of badly approximable numbers in $\mathbb{R}^d$ by $\mathrm{BA}(d)$. Some authors include the rational numbers in $\mathrm{BA}(d)$ and others do not.[1] Our definition does not include them. The set $\mathrm{BA}(d)$ has been widely studied in a variety of contexts. For example, it is known to be a set of zero $d$-dimensional Lebesgue measure, but yet has full Hausdorff dimension, see [248, Chapter III].

A common problem is to investigate the Hausdorff dimension of the intersection of $\mathrm{BA}(d)$ with a given set. The hope is that, since $\mathrm{BA}(d)$ is large from a dimension point of view, for a given set $F \subseteq \mathbb{R}^d$

$$\dim_{\mathrm{H}} F \cap \mathrm{BA}(d) = \dim_{\mathrm{H}} F. \tag{14.2}$$

This clearly cannot be true for arbitrary $F$. For example, take $F = \mathbb{R}^d \setminus \mathrm{BA}(d)$, which is a set of full Hausdorff dimension which does not intersect $\mathrm{BA}(d)$ at all. Alternatively, take $F = [0,1] \times \{0\} \subseteq \mathbb{R}^2$. Again, this set does not intersect $\mathrm{BA}(d)$ at all, which can be seen by applying Dirichlet's theorem in dimension 1. However, if $F$ is somehow chosen in a more 'representative' way, for example without deliberately trying to avoid $\mathrm{BA}(d)$, then perhaps (14.2) will hold. This kind of problem has been studied in two main contexts: when $F$ is chosen to be a smooth submanifold of $\mathbb{R}^d$, or when $F$ is a fractal set. The setting where $F$ is a manifold is motivated by work of Davenport [50][2] and has developed into a broad research area, see [27] for some background. The second context, where $F$ is a fractal set, has gained a lot of attention in the literature, see [45, 46, 171]. It is in this second context that a surprising connection to the lower dimension has emerged.

[1] Rational numbers are surprisingly hard to approximate by rational numbers (other than themselves).

[2] In his 1964 paper, Davenport writes 'I do not know whether there is a set of $\alpha$ with the cardinal of the continuum such that the pair $(\alpha, \alpha^2)$ is badly approximable for simultaneous approximation'. Badziahin and Velani [10] answered this specific question in the affirmative.

A set $F \subseteq \mathbb{R}^d$ is *hyperplane diffuse* if there exists $\beta > 0$ such that, for all $R \in (0,1)$, $x \in F$ and affine hyperplanes $V \subseteq \mathbb{R}^d$,

$$B(x, R) \cap F \setminus V_{\beta R} \neq \varnothing$$

where $V_\varepsilon$ denotes the $\varepsilon$ neighbourhood of $V$. The idea behind this concept is that hyperplane diffuse sets uniformly avoid all affine hyperplanes at all locations and scales. For example, any differentiable curve or smooth manifold fails to be hyperplane diffuse. The following result was obtained by Das, Fishman, Simmons and Urbański [45, Corollary 2.6].

**Theorem 14.2.1** If $F \subseteq \mathbb{R}^d$ is closed and hyperplane diffuse, then

$$\dim_{\mathrm{H}} F \cap \mathrm{BA}(d) \geqslant \dim_{\mathrm{L}} F.$$

We give the proof of Theorem 14.2.1 below, at least the part of the proof where the lower dimension comes into play. There are two subtly different ways of applying this result to concrete examples. Firstly, we can look for hyperplane diffuse sets which already have equal Hausdorff and lower dimension.

**Corollary 14.2.2** If $F \subseteq \mathbb{R}^d$ is a self-similar set, a self-conformal set or the limit set of a Kleinian group with no parabolic elements which, in addition, is not contained in a hyperplane, then

$$\dim_{\mathrm{H}} F \cap \mathrm{BA}(d) = \dim_{\mathrm{H}} F.$$

One might also hope to apply this theorem directly to Kleinian limit sets where $k_{\min} \geqslant \delta(\Gamma)$, since here $\dim_{\mathrm{L}} L(\Gamma) = \dim_{\mathrm{H}} L(\Gamma)$, see Theorem 9.3.3. However, digging into the proof of this fact in [92] demonstrates that $[0,1]^{k_{\min}}$ (embedded in $\mathbb{R}^d$) is a weak tangent to such sets, which prevents them from being hyperplane diffuse if $k_{\min} < d$. The theorem *does* apply to Kleinian groups with parabolic elements if $k_{\min} = d$, and to any Fuchsian limit set.

The second approach is to upgrade the lower bound by replacing the lower dimension with modified lower dimension. Although here one must be careful that the subsets of $F$ which are used to witness the modified lower dimension are themselves hyperplane diffuse. One then looks for sets with equal modified lower and Hausdorff dimension, which are more common and include Bedford-McMullen carpets, for example, see Section 8.3.

**Corollary 14.2.3** If $F \subseteq \mathbb{R}^2$ is a Bedford-McMullen carpet constructed such that at least two columns and at least two rows are used, then

$$\dim_{\mathrm{H}} F \cap \mathrm{BA}(d) = \dim_{\mathrm{H}} F.$$

*Proof* This follows immediately from Theorem 14.2.1 and Lemma 8.3.7 together with the observation that the subsystems constructed in the proof of Lemma 8.3.7 are hyperplane diffuse under the additional assumptions regarding used columns and rows. □

Before proving Theorem 14.2.1, we need to introduce *Schmidt games*, which often play a role in Diophantine approximation. These games were introduced by Schmidt in [247] and are played in a complete metric space $(X, d)$. For us the metric space will always be a closed subset of Euclidean space. The game is played by two players, Alice and Bob. Given $0 < \alpha, \beta < 1$, Alice and Bob play the $(\alpha, \beta)$-*game* as follows:

(i) Bob begins by choosing $r_0 > 0$ and $x_0 \in X$ and we write $B_0 = B(x_0, r_0)$ and $r_n = (\alpha\beta)^n r_0$ for all $n \in \mathbb{N}$.
(ii) On Alice's $n$th turn, she chooses $y_n$ such that $d(x_n, y_n) + \alpha r_n \leqslant r_n$.
(iii) On Bob's $(n+1)$st turn, he chooses $x_{n+1}$ such that $d(y_n, x_{n+1}) + \alpha\beta r_n \leqslant \alpha r_n$.
(iv) The inequalities above ensure that the closed balls

$$B_n = B(x_n, r_n)$$

and

$$A_n = B\left(y_n, \alpha r_{(n-1)}\right)$$

form an interlaced decreasing sequence

$$B_0 \supseteq A_1 \supseteq B_1 \supseteq A_2 \supseteq \cdots \supseteq B_n \supseteq A_{n+1} \supseteq B_{n+1} \supseteq \cdots, \tag{14.3}$$

and hence intersect at a unique point which is called the *outcome* of the game.

Given a set $S \subseteq X$, if Alice has a strategy guaranteeing that the outcome lies in $S$, then $S$ is called $(\alpha, \beta)$-*winning.* If for some fixed $\alpha$, the set $S$ is $(\alpha, \beta)$-winning for all $0 < \beta < 1$, then $S$ is called $\alpha$-*winning.* If $S$ is $\alpha$-winning for some $0 < \alpha < 1$, then $S$ is called *winning.*

We quote the following result from [45, Proposition 2.3] which follows by piecing together results from [38].

**Lemma 14.2.4** Let $F \subseteq \mathbb{R}^d$ be closed and hyperplane diffuse. Then $\mathrm{BA}(d) \cap F$ is winning on $F$.

First, [38, Theorem 2.5] establishes that $\mathrm{BA}(d)$ is hyperplane absolute winning, which just means winning for a different variant of Schmidt's game. Second, [38, Proposition 4.7] proves that if $F$ is hyperplane diffuse and $E$ is hyperplane absolute winning on $F$, then $E$ is winning on $F$.

*Proof of Theorem 14.2.1.* We follow the proof of [45, Proposition 2.5]. Let $s = \dim_{\mathrm{L}} F$, and $\varepsilon > 0$. Therefore there is a constant $C > 0$ such that, for all $x \in F$, $\beta \in (0, 1/4)$ and $R \in (0, 1]$, we can find a $3\beta R$-separated subset of $B(x, (1-\beta)R) \cap F$ of size at least

$$N := C\beta^{-(s-\varepsilon)}.$$

Let $\alpha \in (0, 1)$ and $\beta \in (0, 1/4)$ be chosen such that $\mathrm{BA}(d) \cap F$ is $(\alpha, \beta)$-winning on $F$. We know such an $\alpha$ exists since $\mathrm{BA}(d) \cap F$ is winning on $F$ by Lemma 14.2.4 and $\beta$ can be arbitrarily small.

For each ball $A = B(x, R)$, let $y_i(A)$ (for $i = 1, \dots, N$) be a $3\beta R$-separated subset and write $f_i(A) = B(y_i(A), \beta R)$. Suppose Alice is playing a winning strategy for the $(\alpha, \beta)$-game such that the outcome of the game always ends up in $\mathrm{BA}(d) \cap F$. In response to Alice playing the ball $A_n$, Bob can play any of the balls $f_i(A_n)$, which are legal moves by construction. Moreover, suppose Bob begins with an arbitrary ball $B_0 = B(x_0, 1)$. Since the balls $f_i(A_n)$ are pairwise disjoint, this array of plausible strategies for Bob gives rise to a Cantor set

$$E = \bigcup_{(i_0, i_1, \dots) \in \{1, \dots, N\}^{\mathbb{N}}} \bigcap_n f_{i_n}(A_n)$$

of plausible outcomes for the game, all of which must lie in $\mathrm{BA}(d) \cap F$. Let $\mu$ be the uniform measure on $E$ chosen such that, for $i \in \{1, \dots, N\}$, $\mu(f_i(A_n) \cap E) = N^{-n}$. Since the balls $f_i(A_n)$ have radius $(\alpha\beta)^n$ and are separated by $(\alpha\beta)^n$, there is a constant $c > 0$ such that for any $r \in (0, 1)$ and $x \in E$

$$\mu(B(x, r)) \leqslant c r^{\frac{\log N}{-\log(\alpha\beta)}}$$

which, by the mass distribution principle Lemma 3.4.2, shows that

$$\dim_{\mathrm{H}} E \geqslant \frac{\log N}{-\log(\alpha\beta)} = \frac{(s-\varepsilon)\log\beta - \log C}{\log\beta + \log\alpha} \to (s - \varepsilon)$$

as $\beta \to 0$. Letting $\varepsilon \to 0$, proves $\dim_{\mathrm{H}} \mathrm{BA}(d) \cap F \geqslant \dim_{\mathrm{H}} E \geqslant s = \dim_{\mathrm{L}} F$. □

Finally, we briefly highlight another connection between Diophantine approximation and the lower dimension, this time of measures. Following [28], a measure $\mu$ is said to be *weakly absolutely $\alpha$-decaying* for $\alpha > 0$ if there exist constants $C, r_0 > 0$ such that for all $\varepsilon > 0$

$$\mu(B(x, \varepsilon r)) \leqslant C\varepsilon^{\alpha} \mu(B(x, r))$$

for all $x \in \operatorname{supp}(\mu)$ and $r \leqslant r_0$. It was proved in [140, Proposition 2.6] that

$$\dim_{\mathrm{L}} \mu = \sup\{\alpha \geqslant 0 : \mu \text{ is weakly absolutely } \alpha\text{-decaying}\}.$$

In particular, there exists weakly absolutely $\alpha$-decaying measures with support equal to $F$ if and only if $\dim_{\mathrm{L}} F > 0$ by virtue of Theorem 4.1.3. In [28] the notion of weakly absolutely $\alpha$-decaying was used to study Diophantine properties of limit sets of Kleinian groups. For example, [28, Corollary 1] proves that if $\Gamma$ is a non-elementary, geometrically finite Kleinian group and $\mu$ is a weakly absolutely $\alpha$-decaying measure supported on $K \subseteq L(\Gamma)$, then $K$ is $\mu$-extremal. Here *$\mu$-extremal* refers to a Diophantine property which, roughly speaking, means that the $\mu$ measure of the set of very well-approximable points is zero. In the Euclidean setting $x \in \mathbb{R}^d$ is called *very well-approximable* if there exists $\varepsilon > 0$ such that there exist infinitely many $\mathbf{p} = (p_1, \ldots, p_d) \in \mathbb{Z}^d$ and $q \in \mathbb{N}$ such that

$$\|x - \mathbf{p}/q\| < \frac{1}{q^{1+1/d+\varepsilon}} \tag{14.4}$$

where $\|\cdot\|$ is the supremum norm on $\mathbb{R}^d$. This should be compared with (14.1). In the setting of limit sets of Kleinian groups, ‘approximation by rationals’ is replaced by ‘approximation by $\Gamma$-orbits’; see [28] for the details and [47, 48] for more variants on ‘decaying’ properties of measures and application to Diophantine approximation.

## 14.3 Definability of the integers

In this section we very briefly remark upon an application of the Assouad dimension in an area rather far from the topic of this book: definability in first-order expansions of the real field. We made no attempt to make this section self-contained and refer the reader to [54, 131, 196] for more details and [137] for the necessary background in model theory. In what follows $\overline{\mathbb{R}} = (\mathbb{R}, \cdot, +)$ denotes the field of real numbers.

In first order logic, the subsets of $\mathbb{R}^d$ which are *definable* in $\overline{\mathbb{R}}$ are the *semialgebraic sets*, that is, the sets which can be written as a finite union of sets expressed in terms of a finite collection of polynomial equations or (strict) inequalities. In particular, the semialgebraic sets in $\mathbb{R}$ are precisely the sets which are a finite union of singletons and open intervals. Therefore, $\mathbb{N}$ is *not* definable in $\overline{\mathbb{R}}$. This is representative of the simplicity or 'tameness' of $\overline{\mathbb{R}}$. More generally, and following [131], we think of *expansions* $\mathcal{R}$ of $\overline{\mathbb{R}}$ as being 'tame' if they do *not* define $\mathbb{N}$ and 'wild' if they *do* define $\mathbb{N}$ as a first-order structure. This should be thought of in the context of *o-minimality*, a notion which characterises the simplest expansions. An expansion $\mathcal{R}$ of $\overline{\mathbb{R}}$ is *o-minimal* if the only definable subsets of $\mathbb{R}$ are the semialgebraic sets and, as such, expansions which define $\mathbb{N}$ are certainly not o-minimal. The following theorem was proved in [131, Theorem A] and uses the Assouad dimension to characterise which expansions of $\overline{\mathbb{R}}$ fail to define $\mathbb{N}$ and are therefore 'tame' in the above sense.

**Theorem 14.3.1** An expansion $\mathcal{R}$ of $\overline{\mathbb{R}}$ does not define $\mathbb{N}$ if and only if $\dim_{\mathrm{T}} E = \dim_{\mathrm{A}} E$ for all images $E$ of closed definable sets under definable continuous functions.

Note here that a function is definable if its graph is a definable set. One direction of the equivalence in Theorem 14.3.1 is straightforward since $\dim_{\mathrm{A}} \mathbb{N} = 1 > 0 = \dim_{\mathrm{T}} \mathbb{N}$. The Assouad dimension is favoured in this analysis (instead of upper box dimension, say), since it is the largest available notion of metric dimension and therefore provides the most striking result. However, beyond the aesthetic appeal, the fact that the Assouad dimension appears specifically has useful consequences in other definability questions, see [212] and [275].

Following [131], given a set $E \subseteq \mathbb{R}$, $\mathcal{R} = (\overline{\mathbb{R}}, E)$ denotes a minimal expansion which defines $E$ and all of the definable sets in $\overline{\mathbb{R}}$. For example, $(\overline{\mathbb{R}}, \mathbb{N})$ necessarily defines all Borel sets, see [163, (37.6)]. This is yet further evidence of the 'wildness' of expansions which define $\mathbb{N}$: once you have $\mathbb{N}$ you also have all Borel sets. Interesting further examples along these lines include $(\overline{\mathbb{R}}, \{2^{-k} : k \in \mathbb{N}\})$ which does *not* define $\mathbb{N}$ (this was proved by van den Dries [53]) and $(\overline{\mathbb{R}}, \{2^{-\sqrt{k}} : k \in \mathbb{N}\})$ which *does* define $\mathbb{N}$ (by Theorem 14.3.1). It is a straightforward but worthwhile exercise to verify that

$$\dim_{\mathrm{A}}\{2^{-k} : k \in \mathbb{N}\} = \dim_{\mathrm{T}}\{2^{-k} : k \in \mathbb{N}\} = 0,$$

but

$$\dim_{\mathrm{A}}\{2^{-\sqrt{k}} : k \in \mathbb{N}\} = 1 > 0 = \dim_{\mathrm{T}}\{2^{-\sqrt{k}} : k \in \mathbb{N}\}.$$

# 15
# Applications in probability theory

Understanding how the dimensions of sets change under distortion by a random function is an old problem in dimension theory. There is a general principle that 'Brownian motion doubles dimension', that is, given a compact set $F \subseteq \mathbb{R}$ with Hausdorff dimension at most $1/2$, the Hausdorff dimension of $F$ will almost surely double under distortion by Brownian motion, see (15.1). Motivated by a series of 'uniform dimension results' by Kaufman [161] and others, Balka and Peres proved that this dimension doubling phenomenon can be made uniform over all Borel subsets of a given set, provided the quasi-Assouad dimension of the given set is small enough, see Theorem 15.1.1. In Section 15.2, we consider the Assouad dimensions of sets and measures generated by random functions, following [141, 140].

## 15.1 Uniform dimension results for fractional Brownian motion

Brownian motion[1] is a fundamental example of a random process and gives rise to many interesting (random) fractals via graphs, images and level sets. It can be characterised as the unique random function $B : [0, \infty) \to \mathbb{R}$ satisfying

(i) $B(0) = 0$
(ii) $B$ is almost surely continuous
(iii) $B$ has independent increments, that is, $B(t_1) - B(s_1)$ and $B(t_2) - B(s_2)$ are independent for $0 \leqslant s_1 < t_1 < s_2 < t_2$

[1] Some authors refer to the precise mathematical formulation of this process as the *Wiener process*, reserving the term *Brownian motion* to describe the trajectory of a piece of pollen suspended in water.

(iv) the increments are normally distributed with

$$B(t) - B(s) \sim \mathcal{N}(0, |s-t|).$$

Fractional Brownian motion $B_\alpha$ (with index $\alpha \in (0,1)$) is a natural variant of Brownian motion.[2] This process was introduced by Mandelbrot and van Ness [198] and can be characterised as the unique random function $B_\alpha : [0,\infty) \to \mathbb{R}$ satisfying

(i) $B_\alpha(0) = 0$
(ii) $B_\alpha$ is a Gaussian process
(iii) $B_\alpha$ is almost surely continuous
(iv) $B_\alpha$ has stationary increments with the covariance condition

$$\mathbb{E}(B_\alpha(s)B_\alpha(t)) = (|s|^{2\alpha} + |t|^{2\alpha} - |s-t|^{2\alpha})/2$$

(v) the increments are normally distributed with

$$B_\alpha(t) - B_\alpha(s) \sim \mathcal{N}(0, |s-t|^{2\alpha}).$$

When $\alpha = 1/2$, fractional Brownian motion recovers classical Brownian motion and as $\alpha$ decreases the process becomes 'more fractal', see [1, 70, 214, 156] for more details. Many interesting dimension questions arise.

Given a fixed Borel set $F \subseteq [0,1]$, Kahane [156] proved that almost surely

$$\dim_H B_\alpha(F) = \min\left\{1, \frac{\dim_H F}{\alpha}\right\}. \tag{15.1}$$

In fact, $B_\alpha$ is almost surely $\gamma$-Hölder for any $\gamma \in (0,\alpha)$ and so the upper bound follows from (3.11). A similar result holds for packing dimension but with $\dim_H F$ replaced with the $\alpha$-dimensional *packing profile* of $F$, see Section 10.2. This result is due to Xiao [277]. Higher dimensional analogues of these results are also given in [156, 277].

An interesting and challenging problem is to try to obtain these results *uniformly*; that is, almost surely the conclusion holds for *all* Borel sets $F$ simultaneously. However, it is easy to see that this cannot be true since the zero set of $B_\alpha$, that is, the set

$$\{x \in [0,1] : B_\alpha(x) = 0\}$$

is known to have Hausdorff dimension $1-\alpha$ almost surely. Therefore, almost surely there is always a Borel set with positive dimension which gets mapped to a set with dimension 0 under $B_\alpha$. Surprisingly, one can obtain a uniform dimension result for planar Brownian motion, that is $\mathbf{B} : \mathbb{R} \to \mathbb{R}^2$ defined by $\mathbf{B}(x) = (B^1(x), B^2(x))$ where $B^1$ and $B^2$ are independent Brownian motions. It was proved by Kaufman [161] that almost surely

$$\dim_{\mathrm{H}} \mathbf{B}(F) = 2 \dim_{\mathrm{H}} F \qquad \text{and} \qquad \dim_{\mathrm{P}} \mathbf{B}(F) = 2 \dim_{\mathrm{P}} F$$

simultaneously for all Borel sets $F \subseteq [0,1]$, see also [214]. Balka and Peres [12] investigated the question of when a uniform result can be recovered in the 1-dimensional setting, where a surprising and elegant connection to the Assouad dimension was unearthed. The idea is to replace $[0,1]$ with a 'small enough' set $D \subseteq [0,1]$ such that a uniform dimension result can be obtained for all subsets of $D$.

[2] $\alpha$ is often called the 'Hurst parameter' or 'Hurst index'.

**Theorem 15.1.1** If $D \subseteq [0,1]$ satisfies $\dim_{\mathrm{qA}} D \leqslant \alpha$, then almost surely

$$\dim_{\mathrm{H}} B_\alpha(F) = \frac{\dim_{\mathrm{H}} F}{\alpha}$$

simultaneously for all $F \subseteq D$.

Balka and Peres also obtain this result with Hausdorff dimension replaced by packing dimension and their result is stronger than stated since they use a modification of the quasi-Assouad dimension, which is bounded above by the quasi-Assouad dimension. In fact their notion of dimension turns out to be sharp for the packing dimension problem. First, they define the *modified Assouad dimension* by

$$\dim_{\mathrm{MA}} F = \lim_{\delta \to 1} \inf \left\{ \sup_i h_{F_i}(\delta) \ : \ F = \bigcup_i F_i \right\}$$

where $h_F(\delta)$ is the function used to define the quasi-Assouad dimension in Section 3.2. To get a better feel for this definition, notice the similarity between this and how the upper box dimension is 'modified' in order to obtain the packing dimension, see (1.1). In particular, $\dim_{\mathrm{P}} F \leqslant \dim_{\mathrm{MA}} F \leqslant \dim_{\mathrm{qA}} F$. The following is [12, Theorem 1.8].

**Theorem 15.1.2** For an analytic set $D \subseteq [0,1]$, the following are equivalent:

(i) $\dim_{\mathrm{MA}} D \leqslant \alpha$
(ii) almost surely

$$\dim_{\mathrm{P}} B_\alpha(F) = \frac{\dim_{\mathrm{P}} F}{\alpha}$$

simultaneously for all $F \subseteq D$.

Note the philosophical connection between this result and Corollary 10.2.8. Both results say that if the quasi-Assouad dimension is small enough, then the dimension profile under consideration is maximal.

Many other related results are proved and discussed in [12]. Balka and Peres specifically ask if $\dim_{\mathrm{MA}} F = \dim_{\mathrm{H}} F$ for any self-similar set $F$, see [12, Problem 1.16]. This would follow from a positive answer to the question of whether $\dim_{\mathrm{qA}} F = \dim_{\mathrm{H}} F$ for every self-similar set $F$, see Question 17.5.3.

## 15.2 Dimensions of random graphs

Another interesting random fractal arising from (fractional) Brownian motion is the *graph*, defined by

$$G(B_\alpha) = \{(x, B_\alpha(x)) : x \in [0,1]\}.$$

This set is known to have Hausdorff, packing and box dimension almost surely given by $2 - \alpha$, see [1] or [70, Section 16.3]. It was shown in [141] that the Assouad dimension of $G(B_\alpha)$ is almost surely 2 (independent of $\alpha$). Notice the parallel between this result and the dimension results for Mandelbrot percolation, see Theorem 9.4.1. The fact that the graph of classical Brownian motion has Assouad dimension 2 can be derived from results in [43].

Howroyd [140, Theorem 2.5] proved that for all doubling measures $\mu$ supported on $[0,1]$, the pushforward of $\mu$ onto the graph $G(B_\alpha)$ is almost surely not doubling. By virtue of Theorem 4.1.3 the graph $G(B_\alpha)$ necessarily supports doubling measures, but the result of Howroyd shows that almost surely these measures cannot be the pushforward of a deterministic measure on $[0,1]$. This is rather different from the Hausdorff dimension case, where the pushforward of Lebesgue measure onto $G(B_\alpha)$ is almost surely a measure of Hausdorff dimension $2 - \alpha$, that is, a measure which realises the Hausdorff dimension of the graph.

# 16
# Applications in functional analysis

In this final applications chapter, we consider a range of problems in functional and harmonic analysis. In Section 16.1 we begin with a well-established connection between the Assouad and lower dimensions and Hardy inequalities, following work of Lehrbäck [186] and others. The basic connection is that conditions on a domain, in terms of the Assouad and lower dimensions of the codomain, determine whether the domain admits a Hardy inequality. In Section 16.2 we explore a different problem involving maximal operators averaged over spheres, considered by Anderson, Hughes, Roos and Seeger [4], where the Assouad dimension plays a role in determining whether certain $L^p$ improving estimates are satisfied. Moreover, the Assouad spectrum makes an appearance in examining the sharpness of the Assouad dimension results. Finally, in Section 16.3 we consider the relationship between $L^p$ smoothness of measures and the Assouad and lower spectra, following [109].

## 16.1 Hardy inequalities

*Hardy's inequality* states that for a sequence of non-negative real numbers $a_n$ and $p > 1$

$$\sum_{n=1}^{\infty} \left( \frac{a_1 + a_2 + \cdots + a_n}{n} \right)^p \leqslant \left( \frac{p}{p-1} \right)^p \sum_{n=1}^{\infty} a_n^p.$$

Various continuous analogues of this inequality are known to hold. For example, we say an open domain $\Omega \subseteq \mathbb{R}^d$ admits a *p-Hardy inequality* if there exists a constant $C \geqslant 1$ such that

$$\int_\Omega \left( \frac{|f(x)|}{d(x, \Omega^c)} \right)^p dx \leqslant C \int_\Omega \mathrm{Lip}_x(f)^p dx$$

for all Lipschitz $f : \mathbb{R}^d \to \mathbb{R}$ with compact support contained in $\Omega$, where $d(x,\Omega^c) = \inf_{y \notin \Omega} |x-y|$ is the distance from $x$ to the codomain $\Omega^c = \mathbb{R}^d \setminus \Omega$ and

$$\mathrm{Lip}_x(f) = \sup_{y \neq x} \frac{|f(x)-f(y)|}{|x-y|}$$

is the (upper) Lipschitz constant of $f$ at $x \in \mathbb{R}^d$. Often the family of Lipschitz functions is replaced by smooth functions with compact support, in which case $\mathrm{Lip}_x(f)$ is replaced by $|\nabla f(x)|$.

Not all domains admit a $p$-Hardy inequality and it turns out that the Assouad and lower dimensions of the codomain play a key role. The following result can be found in [186, Corollary 1.3].

**Theorem 16.1.1** If $\Omega \subseteq \mathbb{R}^d$ is an open set and either $\dim_{\mathrm{A}} \Omega^c < d-p$ or $\Omega^c$ is unbounded and $\dim_{\mathrm{L}} \Omega^c > d-p$, then $\Omega$ admits a $p$-Hardy inequality.

Theorem 16.1.1 follows as a corollary to a more general result concerning weighted Hardy inequalities, which we do not consider here but refer the reader to [186] for more information, including an extensive bibliography. There are many other connections between the Assouad dimension and inequalities in functional analysis. We refer the reader to the survey [187] for more details but briefly mention a few highlights here. Connections between the Assouad dimension and *Sobolev inequalities* are considered in [189, 57]. In [188] it is established that (for subsets of Euclidean space) the Assouad and *Aikawa* dimensions coincide, see [2]. This latter result connects the Assouad dimension to yet another family of integral inequalities, see [188, Definition 3.2]. Finally, connections between the Assouad dimension and *Muckenhoupt weights* are given in [58], see also [187].

Returning to Theorem 16.1.1, it was proved in [174] that, for open domains $\Omega \subseteq \mathbb{R}^d$, $\dim_{\mathrm{L}} \Omega^c > 0$ is equivalent to $\Omega$ admitting a $d$-Hardy inequality. We prove one direction of this equivalence to give some intuition as to how the lower dimension plays a role.

**Theorem 16.1.2** If $\Omega \subseteq \mathbb{R}^d$ is an open domain with $\dim_{\mathrm{L}} \Omega^c = 0$, then $\Omega$ does not admit a $d$-Hardy inequality.

*Proof* Suppose $\dim_{\mathrm{L}} \Omega^c = 0$, in which case $\Omega^c$ is not uniformly perfect by Theorem 13.1.2. Therefore, for an arbitrary $\kappa \in (0,1/2)$ which we fix, we can find $y \in \Omega^c$ and $0 < R < 1/d$ such that $B(y,R) \setminus B(y,\kappa R) \subseteq \Omega$. Consider $f$ defined by $f(x) = \max\{1-|y-x|, 0\}$ for $x \in \Omega$ and $f(x) = 0$

otherwise. Then

$$\begin{aligned}\int_{\Omega}\left(\frac{|f(x)|}{d(x,\Omega^c)}\right)^d dx &\geqslant \int_{B(y,R/2)\setminus B(y,\kappa R)}\left(\frac{1-|y-x|}{|x-y|}\right)^d dx\\ &= \int_{B(0,R/2)\setminus B(0,\kappa R)}\left(\frac{1-|x|}{|x|}\right)^d dx\\ &= \omega_d\int_{\kappa R}^{R/2}\left(\frac{1-x}{x}\right)^d x^{d-1}dx \qquad (16.1)\end{aligned}$$

where $\omega_d > 0$ is a constant depending only on $d$. This final equality is via a standard reduction for the integral of a radial function in spherical coordinates. Continuing, we may bound (16.1) below by

$$c_d\int_{\kappa R}^{R/2} x^{-1}dx = c_d\log(1/(2\kappa))$$

for another constant $c_d > 0$. In particular, this lower bound grows without bound as $\kappa \to 0$. However,

$$\int_{\Omega}\mathrm{Lip}_x(f)^d dx \leqslant \int_{B(y,1)} dx < \infty$$

is bounded independently of $\kappa$ and so $\Omega$ does not admit a $d$-Hardy inequality. □

Sometimes in the literature one can find analogues of Theorem 16.1.1 where the codomain $\Omega^c$ is replaced by the boundary $\partial\Omega$. Essentially, this subtle difference only plays a role in the case of *weighted* Hardy inequalities, which we do not discuss here. Note that if $\dim_{\mathrm{A}}\partial\Omega < d-1$, then $\Omega^c = \partial\Omega$. To see this, observe that if $\Omega^c$ has interior, then we can find $x, y \in \mathbb{R}^d$ and $r > 0$ such that $B(x,r) \subseteq \Omega$ and $B(y,r) \subseteq \Omega^c$. Therefore, the projection of $\partial\Omega$ onto the $(d-1)$-dimensional hyperplane perpendicular to the line joining $x$ and $y$ has non-empty interior, which forces $\dim_{\mathrm{H}}\partial\Omega \geqslant d-1$. In the context of Theorem 16.1.1, this observation means that the alternative assumption that $\dim_{\mathrm{A}}\partial\Omega < d-p$ is actually equivalent to the (*a priori* stronger) assumption that $\dim_{\mathrm{A}}\Omega^c < d-p$. The situation for lower dimension is more subtle, see [185].

## 16.2 $L^p \to L^q$ bounds for spherical maximal operators

A connection between the Assouad dimension and certain $L^p \to L^q$ inequalities was established in [4]. Moreover, the Assouad spectrum is used in order to examine sharpness of the results. Let $d \geqslant 2$ and $f : \mathbb{R}^d \to \mathbb{R}$ be a locally integrable function and consider the average of $f$ over the sphere centred at $x \in \mathbb{R}^d$ with radius $t > 0$, that is,

$$A_t f(x) = \int f(x - ty) d\sigma$$

where $\sigma$ is the normalised surface measure on the unit sphere in $\mathbb{R}^d$. Given a set $E \subset [1, 2]$, consider the *maximal function*

$$M_E f(x) = \sup_{t \in E} |A_t f(x)|.$$

A common problem in harmonic analysis is to examine when a given maximal function is an embedding from $L^p$ into $L^q$ for some $p, q \geqslant 1$. The maximal function is called $L^p$ *improving* if $q > p$. The following $L^p \to L^q$ embedding result was proved in [4, Theorem 1].

**Theorem 16.2.1** If $(1/p, 1/q) \in \mathcal{R}_d(E)$, then there is a constant $C \geqslant 1$ such that

$$\|M_E f\|_q \leqslant C \|f\|_p$$

for all locally integrable $f : \mathbb{R}^d \to \mathbb{R}$, where $\mathcal{R}_d(E)$ is a quadrilateral region determined by $\dim_A E$, $\overline{\dim}_B E$ and $d$. Specifically, $\mathcal{R}_d(E)$ is the interior of the convex hull of the four points:

$$(0,0), \qquad \left( \frac{d-1}{d-1+\overline{\dim}_B E}, \frac{d-1}{d-1+\overline{\dim}_B E} \right),$$

$$\left( \frac{d-\overline{\dim}_B E}{d-\overline{\dim}_B E+1}, \frac{1}{d-\overline{\dim}_B E+1} \right),$$

and

$$\left( \frac{d(d-1)}{d^2+2\dim_A E-1}, \frac{d-1}{d^2+2\dim_A E-1} \right)$$

together with the diagonal line segment joining the first two points with the first point included and the second point not included.

If $d = 2$ and $\dim_A E > 1/2$, then the proof of the above theorem is not given in [4] but deferred to a later paper. The portion of $\mathcal{R}_d(E)$ intersecting the diagonal $p = q$ only depends on $\overline{\dim}_B E$ and $d$ and this connection with the upper box dimension was observed earlier [249].

Theorem 16.2.1 is pretty close to being sharp and the sharpness was investigated in [4, Theorem 2] using the Assouad spectrum. An example of what is proved in [4] is that if $E$ is such that

$$\dim_{\mathrm{A}}^{\theta} E = \min\left\{\frac{\overline{\dim}_{\mathrm{B}} E}{1-\theta}, \dim_{\mathrm{A}} E\right\}, \tag{16.2}$$

for all $\theta \in (0,1)$ — that is, the Assouad spectrum of $E$ is as large as possible given the basic bounds in Lemma 3.4.4 — then the conclusion of Theorem 16.2.1 is false for $(1/p, 1/q) \notin \overline{\mathcal{R}(E)}$. This can be interpreted as showing that, for $E$ satisfying (16.2), Theorem 16.2.1 is sharp up to boundary points.

## 16.3 Connection with $L^p$-norms

There is a simple connection between the Assouad and lower spectra and $L^p$ norms, which was observed in [109]. Recall a measure $\mu$ on $\mathbb{R}$ is *absolutely continuous* if all Lebesgue null sets are given zero $\mu$ measure. If $\mu$ is absolutely continuous, then there is a Lebesgue integrable function $f$, the *density* or *Radon-Nikodym derivative*, such that

$$\mu(E) = \int_E f(x)\, dx$$

for all Borel sets $E$. Given $p \geqslant 1$, we write $L^p[0,1]$ for the set of all integrable functions $g : [0,1] \to \mathbb{R}$ such that

$$\|g\|_p := \left(\int_0^1 g(x)^p\, dx\right)^{1/p} < \infty$$

and we say $\mu \in L^p[0,1]$ if $\mu$ is absolutely continuous with density $f \in L^p[0,1]$. As such, we are restricting our attention to measures with support in $[0,1]$. We write $f \in L^{-p}[0,1]$ if the set $E = \{x \in [0,1] : f(x) = 0\}$ is of Lebesgue measure zero and

$$\left(\int_{[0,1]\setminus E} 1/f(x)^p dx\right)^{1/p} < \infty,$$

and $\mu \in L^{-p}[0,1]$ if $\mu$ is absolutely continuous with density $f \in L^{-p}[0,1]$.

The smoothness of $\mu$ can then be characterised by which $L^p[0,1]$ and $L^{-p}[0,1]$ spaces the measure belongs to — the larger the exponents $p$, the smoother the measure is. Absolutely continuous measures on $\mathbb{R}$ immediately satisfy $\dim_{\mathrm{H}} \mu = 1$ and it is natural to consider other notions of dimension. The main result in [109, Theorem 2.1] establishes a sharp one-way implication between smoothness and regularity, as characterised by $L^p$ properties and dimensions respectively. More general results than those we present below are available in [109].

**Theorem 16.3.1** Suppose $\mu$ is an absolutely continuous measure and $p_1, p_2 \geqslant 1$ are such that $\mu \in L^{p_1}[0,1] \cap L^{-p_2}[0,1]$. Then, for all $\theta \in (0,1)$,

$$1 - \frac{\theta p_1 + p_2}{p_1 p_2(1-\theta)} \leqslant \dim_{\mathrm{L}}^{\theta} \mu \leqslant \dim_{\mathrm{A}}^{\theta} \mu \leqslant 1 + \frac{p_1 + \theta p_2}{p_1 p_2(1-\theta)}.$$

*Proof* Recall *Hölder's inequality* for measurable functions $f, g$:

$$\|fg\|_1 \leqslant \|f\|_p \|g\|_q,$$

where $q = p/(p-1)$ is the Hölder conjugate of $p$, and the *reverse Hölder inequality*, assuming in addition that $f(x) \neq 0$ for almost every $x$,

$$\|fg\|_1 \geqslant \|f\|_{-p} \|g\|_{p/(p+1)}.$$

Here $\|f\|_{-p}$ and $\|g\|_{p/(p+1)}$ are not norms but convenient notation for

$$\|f\|_{-p} = \left( \int_0^1 f(x)^{-p} dx \right)^{-1/p}$$

and

$$\|g\|_{p/(p+1)} = \left( \int_0^1 f(x)^{p/(p+1)} dx \right)^{(p+1)/p}.$$

Write $f$ for the density of $\mu$ and $\chi_A$ for the indicator function on a set $A$. Fix $\theta \in (0,1)$, $x \in \mathrm{supp}(\mu)$ and $r \in (0,1)$. Write $q_1 \in (1,\infty)$ for the Hölder conjugate of $p_1$. Applying Hölder's inequality and the reverse Hölder inequality,

$$\begin{aligned} \frac{\mu(B(x,r^\theta))}{\mu(B(x,r))} = \frac{\|f\chi_{B(x,r^\theta)}\|_1}{\|f\chi_{B(x,r)}\|_1} &\leqslant \frac{\|f\|_{p_1} \|\chi_{B(x,r^\theta)}\|_{q_1}}{\|f\|_{-p_2} \|\chi_{B(x,r)}\|_{p_2/(1+p_2)}} \\ &= \left( \frac{\|f\|_{p_1}}{\|f\|_{-p_2}} \right) \frac{(2r^\theta)^{1/q_1}}{(2r)^{1+1/p_2}} \\ &= \left( \frac{\|f\|_{p_1} 2^{1/q_1}}{\|f\|_{-p_2} 2^{1+1/p_2}} \right) \left( \frac{r^\theta}{r} \right)^{\frac{\theta/q_1 - 1 - 1/p_2}{\theta - 1}}, \end{aligned}$$

which proves

$$\dim_{\mathrm{A}}^{\theta} \mu \leqslant \frac{\theta/q_1 - 1 - 1/p_2}{\theta - 1} = 1 + \frac{p_1 + \theta p_2}{p_1 p_2(1-\theta)},$$

as required. We omit the proof of the estimate for $\dim_{\mathrm{L}}^{\theta} \mu$. □

The inequalities in Theorem 16.3.1 are sharp, see [109], and estimates are also available for $p_1, p_2 = \infty$.

# 17
# Future directions

In this chapter we briefly outline what we believe are some of the key questions for the future and discuss some of the challenges which this book has unearthed. We pose several open problems and, in order for this list to remain up to date, an active list will be maintained on the author's personal homepage. In 2020 this page was hosted here: `http://www.mcs.st-andrews.ac.uk/~jmf32/`. The questions are presented roughly in the same order as the relevant topics are encountered in the book.

## 17.1 Finite stability of modified lower dimension

Modified lower dimension is known to be finitely stable in two situations: when the sets in question are closed, see Lemma 3.4.10, or when the sets in question are 'properly separated', see [111, Proposition 4.3]. The general case remains unknown.

**Question 17.1.1** Is it true that, for arbitrary $E, F \subseteq \mathbb{R}^d$,

$$\dim_{\mathrm{ML}} E \cup F \;=\; \max\{\dim_{\mathrm{ML}} E,\ \dim_{\mathrm{ML}} F\}?$$

It follows from monotonicity and Lemma 3.4.10 that

$$\max\{\dim_{\mathrm{ML}} E,\ \dim_{\mathrm{ML}} F\} \leqslant \dim_{\mathrm{ML}} E \cup F \leqslant \max\left\{\dim_{\mathrm{ML}} \overline{E},\ \dim_{\mathrm{ML}} \overline{F}\right\}. \quad (17.1)$$

## 17.2 Dimensions of measures

Theorem 4.1.3 shows that one can approximate the Assouad and lower dimensions of a compact set by the Assouad and lower dimensions of doubling measures supported on the set (from above and below, respectively). However, it is not always possible to find a measure with Assouad dimension precisely equal to that of the set, see [152, 173]. The analogous problem for lower dimension remains open, and was posed explicitly in [152, Question 5.4].

**Question 17.2.1** Does there exist a compact set $F \subseteq \mathbb{R}^d$ such that $\dim_{\mathrm{L}} \mu < \dim_{\mathrm{L}} F$ for all doubling measures fully supported on $F$?

In the context of box dimension, it was proved in [75, Proposition 3.4] that there exists a compact set $F \subseteq \mathbb{R}^d$ such that

$$\overline{\dim}_{\mathrm{B}} F < \overline{\dim}_{\mathrm{B}} \mu$$

for all doubling measures supported on $F$. This should be compared to Theorem 4.2.3. It was proved in [75, Theorem 2.1] that the word 'doubling' cannot be dropped here.

## 17.3 Weak tangents

The fact that the Assouad dimension is always realised by the Hausdorff dimension of a weak tangent is a powerful theoretical tool, recall Theorem 5.1.3. It is natural to ask if this can be strengthened by replacing the Hausdorff dimension with something smaller. It is easy to construct examples showing that the Hausdorff dimension cannot be replaced by the *lower* dimension (or even the lower spectrum) in general. Consider the middle third Cantor set and add a point in the middle of every complementary interval. The resulting set $F$ is compact (in fact it is an inhomogeneous self-similar set, see [87, 152, 255]) and it is a short exercise to show that $\dim_{\mathrm{A}} F = \frac{\log 2}{\log 3}$ but that every weak tangent has an isolated point. Therefore, all weak tangents $E$ of $F$ satisfy $\dim_{\mathrm{L}} E = \dim_{\mathrm{L}}^{\theta} E = 0$ for all $\theta \in (0, 1)$. However, it is not known if the *modified lower* dimension can be used.

**Question 17.3.1** Is it true that for all closed sets $F \subseteq \mathbb{R}^d$ there exists a weak tangent $E$ of $F$ such that $\dim_{\mathrm{ML}} E = \dim_{\mathrm{A}} F$?

A slightly weaker version of this question would be to ask if there exist weak tangents with modified lower dimension arbitrarily close to $\dim_A F$. On the other hand, a slightly stronger version of the question would be to ask if there exists a weak tangent which contains an Ahlfors-regular subset of dimension $\dim_A F$. An affirmative answer to the strong version of the question would have powerful applications. For example, it would allow Theorem 11.1.3 to be extended to the $\dim_A F = 1$ case, see Question 17.9.1 below. This would follow by applying Lemma 11.1.5 and the main result of [225].

Given a closed set $F \subseteq \mathbb{R}^d$, recall the set of Hausdorff dimensions attained at the level of weak tangents is

$$\Delta(F) = \{\dim_H E : E \in \mathcal{W}(F)\},$$

where $\mathcal{W}(F)$ denotes the set of 'reasonable' weak tangents of $F$, see Section 5.2. It is known that $\Delta(F)$ always contains its maximum and minimum and, subject to this constraint, any $\mathcal{F}_\sigma$ set can be realised as $\Delta(F)$ for some $F$. However, not much more is known about the set-theoretic complexity of $\Delta(F)$.

**Question 17.3.2** Given a closed set $F \subseteq \mathbb{R}^d$, is the set $\Delta(F)$ always analytic or Borel and, if so, does it always belong to a particular finite Borel class?

It is certainly not true that for *any* set $\Delta \subseteq [0, d]$ containing its maximum and minimum there is a compact set $F$ such that $\Delta(F) = \Delta$. This can be seen by a simple cardinality argument, since the collection of possible $\Delta$ has cardinality $\#\mathcal{P}(\mathbb{R})$ (that is, the cardinality of the power set of $\mathbb{R}$), but the collection of possible choices for $F$ has cardinality $\#\mathbb{R}$ (that is, the cardinality of $\mathbb{R}$), which is strictly smaller.

## 17.4 Further questions of measurability

Motivated by the discussion in Section 17.3, we note that a simpler measurability question was considered in [88, 111]. It was shown that the function $\dim : \mathcal{K}(\mathbb{R}^d) \to \mathbb{R}$ is Borel measurable (even of *Baire class* precisely 2) if dim is the Assouad or lower dimension and $\mathcal{K}(\mathbb{R}^d)$ is the set of non-empty compact subsets of $\mathbb{R}^d$ equipped with the Hausdorff metric, see (5.1). Similar measurability results hold for the Assouad and lower spectra, see [111, Theorems 5.1 and 5.2]. For comparison, the Hausdorff and upper and lower box dimensions were shown to be measurable (and of Baire class 2) in [208] but there it was also proved that the packing dimension is *not* measurable. It is not known if the modified lower dimension is measurable, see [111, Question 9.4].

**Question 17.4.1** Is the function $\dim_{\mathrm{ML}} : \mathcal{K}(\mathbb{R}^d) \to \mathbb{R}$ Borel measurable? If so, does it belong to a particular Baire class?

## 17.5 IFS attractors

The Assouad and lower dimension are well-understood in the case of planar self-affine sets with a grid structure. Bedford-McMullen and Lalley-Gatzouras type carpets were considered in [194], Barański type carpets were considered in [88], higher dimensional sponges in [97, 139, 45], and restricted families of planar self-affine sets without a grid structure were considered in [21, 101]. It seems the most natural examples to consider next are either more general planar self-affine carpets without a grid structure or higher dimensional versions of these sets. Alternatively, the Assouad and lower spectra are not as well understood, even in the plane. The following question is motivated by Theorem 8.4.1 proved in [21] as well as the typical formula for the Assouad dimension of carpet type self-affine sets.

**Question 17.5.1** Let $F \subseteq \mathbb{R}^2$ be a self-affine set satisfying the SSC and which is not self-similar. Does there always exist $\pi \in G(1,2)$ such that

$$\dim_{\mathrm{A}} F = \dim_{\mathrm{A}} \pi F + \max_{x \in \pi F} \dim_{\mathrm{A}} (\pi^{\perp} + x) \cap F?$$

It is interesting to examine the relationship between the various different notions of dimension in a given setting. For a given class of sets one can ask *what relationships are possible between the dimensions and spectra?* For example, all the dimensions and spectra are necessarily equal for Ahlfors regular sets. Interestingly, the configuration $\dim_{\mathrm{L}} F < \dim_{\mathrm{B}} F = \dim_{\mathrm{A}} F$ is possible for self-affine carpets, but not for self-similar sets (even with overlaps) and the configuration $\dim_{\mathrm{L}} F = \dim_{\mathrm{B}} F < \dim_{\mathrm{A}} F$ is not possible for self-affine carpets, but is possible for self-similar sets with overlaps. Roughly speaking, the reason for this is as follows. The non-uniform scaling present in self-affine carpets allows one to 'spread' the set out making certain locations easier to cover and thus making the lower dimension drop. Whereas, one can use overlaps to 'pile' the set up at certain locations making them harder to cover and thus raising the Assouad dimension. It would be interesting to add Hausdorff dimension to this discussion, but there are some configurations which remain mysterious to us.

| Configuration | carpet | self-affine | self-similar | Fuchsian | Kleinian |
|---|---|---|---|---|---|
| L = H = B = A | ✓ | ✓ | ✓ | ✓ | ✓ |
| L = H = B < A | × | ✓ | ✓ | ✓ | ✓ |
| L = H < B = A | × | ? | × | × | × |
| L < H = B = A | × | ✓ | × | × | ✓ |
| L = H < B < A | × | ? | × | × | × |
| L < H = B < A | × | ? | × | × | ✓ |
| L < H < B = A | × | ? | × | × | × |
| L < H < B < A | ✓ | ✓ | × | × | × |

Table 17.1 *This table summarises which configurations of dimension are possible in certain families of fractal sets.* Carpet *refers to Bedford-McMullen carpets, see Section 8.3,* self-affine *refers to general self-affine sets, see Chapter 8,* self-similar *refers to general self-similar sets, see Chapter 7,* Fuchsian *refers to limit sets of geometrically finite Fuchsian groups, and* Kleinian *refers to limit sets of geometrically finite Kleinian groups, see Section 9.3. The configurations refer to the relationships between lower, Hausdorff, box and Assouad dimensions, with the obvious labelling. The symbol ✓ means the configuration is known to be possible within the given class, × means the configuration is known to be impossible within the given class, and ? means we do not know if the configuration is possible within the given class.*

**Question 17.5.2** Are any of the configurations marked with a question mark in Table 17.1 possible in the relevant class of sets? The entries marked with a tick or a cross have all been explained in this book.

We still do not have a deep understanding of the Assouad spectrum of overlapping self-similar sets. In the situations where the quasi-Assouad dimension is known, it coincides with the upper box dimension, rendering the Assouad spectrum constant.

**Question 17.5.3** Is it true that, for all self-similar sets $F \subseteq \mathbb{R}^d$,

$$\dim_{\mathrm{qA}} F = \dim_{\mathrm{B}} F?$$

An affirmative answer to this question would also answer a question posed by Balka and Peres [12, Problem 1.16], which asks the same question but with quasi-Assouad dimension replaced by their *modified Assouad dimension* (which is bounded above by quasi-Assouad dimension).

In all cases where the quasi-Assouad and Assouad dimensions are both known for an attractor of an IFS satisfying the OSC or SSC, they coincide. On the other hand, they are often distinct in the overlapping case, recall Theorem 7.3.1, or in the random case, recall Theorems 9.4.1-9.4.2.

**Question 17.5.4** Does there exist a self-affine set $F \subseteq \mathbb{R}^d$ which satisfies the SSC and is such that $\dim_{\mathrm{qA}} F < \dim_{\mathrm{A}} F$?

Note that the answer to this question is *yes* if one does not assume separation conditions, for example, the self-similar sets from Theorem 7.2.1. Moreover, it is *no* if one replaces self-affine with self-similar or self-conformal (whilst maintaining the SSC).

Theorem 6.3.1 proved that the lower spectrum of the attractor of a non-trivial IFS of bi-Lipschitz contractions is always strictly positive, at least for some range of $\theta$. We do not know if the lower dimension is necessarily strictly positive. Recall, it was proved in [278] that the lower dimension of any self-*affine* set is strictly positive and the example (6.4) on page 90 shows that the lower dimension can be zero if the contractions are not bi-Lipschitz.

**Question 17.5.5** Let $F \subseteq \mathbb{R}^d$ be the attractor of an IFS consisting of bi-Lipschitz contractions and assume $F$ is not a singleton. Is it true that $\dim_{\mathrm{L}} F > 0$?

## 17.6 Random sets

Recall Theorem 9.4.4 which established a dichotomy for determining the almost sure value of $\dim_{\mathrm{A}}^{\phi} M$ for different functions $\phi$, where $M$ is the limit set of Mandelbrot percolation, see Sections 3.3.3 and 9.4. We condition on non-extinction of $M$. The situation is understood when

$$\frac{\log(R/\phi(R))}{\log|\log R|} \tag{17.2}$$

either converges to 0 or diverges as $R \to 0$, but otherwise remains unclear.

**Question 17.6.1** Suppose (17.2) neither converges to 0 nor diverges. Then what is the almost sure value of $\dim_{\mathrm{A}}^{\phi} M$?

Note that $\dim_{\mathrm{A}}^{\phi} M$ may be expressed in terms of tail events (we only care about small scale behaviour), and so there is an almost sure value of $\dim_{\mathrm{A}}^{\phi} M$ by Kolmogorov's 0-1 law. Moreover, for any $\phi$, the almost sure value of $\dim_{\mathrm{A}}^{\phi} M$ must be in the interval $[s, d]$ where $s$ is the almost sure Hausdorff dimension of $M$ and $d$ is the ambient spatial dimension.

**Question 17.6.2** Is it true that for all $t \in [s, d]$, there exists $\phi$ such that almost surely $\dim_{\mathrm{A}}^{\phi} M = t$?

## 17.7 General behaviour of the Assouad spectrum

An example was constructed in [111, Section 8] demonstrating that the Assouad spectrum can be strictly decreasing on infinitely many disjoint intervals accumulating at $\theta = 0$. However, we believe this behaviour is not possible at $\theta = 1$.

**Question 17.7.1** Is it true that for any set $F \subseteq \mathbb{R}^d$, there exists $\theta_F \in (0, 1)$ such that $\dim_{\mathrm{A}}^{\theta} F$ is non-decreasing on the interval $(\theta_F, 1)$?

Recall that $\dim_{\mathrm{A}}^{\theta} F$ is often constant in a neighbourhood of $\theta = 1$, but may be strictly increasing in certain special cases, see Theorem 3.3.5.

We have seen that the Assouad spectrum can exhibit many different features. However, for sets which exhibit a certain 'rough homogeneity', such as dynamical invariance or statistical self-similarity, then the Assouad spectrum tends to take on one of two forms. It is either constant (recall Mandelbrot percolation or self-similar sets with overlaps) or consists of a convex increasing part followed by a constant part (recall self-affine carpets, spirals, polynomial sequences). It turns out that there is even more similarity between the examples in this latter situation than meets the eye. Suppose the phase transition between the increasing and constant parts of the spectrum occurs at $\theta = \rho \in (0, 1)$. Following some non-trivial algebraic manipulation[1], we find that for $F$ a (non-uniform fibre) Bedford-McMullen carpet (Theorem 8.3.3), a polynomial spiral

[1] "but isn't *all* algebraic manipulation trivial?" - Anon

(Theorem 13.2.3), or a polynomial sequence (Theorem 3.4.7),

$$\dim_{\mathrm{A}}^{\theta} F = \min \left\{ \dim_{\mathrm{B}} F + \frac{(1-\rho)\theta}{(1-\theta)\rho} \left( \dim_{\mathrm{qA}} F - \dim_{\mathrm{B}} F \right), \ \dim_{\mathrm{qA}} F \right\}. \tag{17.3}$$

We note that in each of these examples $\dim_{\mathrm{qA}} F = \dim_{\mathrm{A}} F$ but we chose to use $\dim_{\mathrm{qA}}$ in the formula since it is perhaps more representative of what is really going on. The formula (17.3) shows that in each of these examples the Assouad spectrum is determined by the box dimension, the (quasi-)Assouad dimension, and the phase transition $\rho$. It is perhaps disappointing for an object with so many degrees of freedom to be determined by such a small set of values, but we believe this is indicative of the simplicity of these examples rather than the simplicity of the Assouad spectrum as a concept. Lemma 3.4.4 guarantees that

$$\rho \geqslant 1 - \frac{\overline{\dim}_{\mathrm{B}} F}{\dim_{\mathrm{qA}} F}. \tag{17.4}$$

For the polynomial sequences (17.4) is always an equality, for the polynomial spirals (17.4) is an equality if and only if $p \leqslant 1$, and for the Bedford-McMullen carpets (17.4) is never an equality.

It seems from this observation that the phase transition $\rho$ merits deeper investigation. Indeed, in these three examples $\rho$ is a constant which holds particular geometric significance for the set in question, capturing some fundamental scaling property. For carpets, the $k$th level rectangles in the standard construction of $F$ are of size $m^{-k} \times n^{-k}$ and therefore $\rho = \frac{\log m}{\log n}$ is the 'logarithmic eccentricity'. For spirals, the *$k$th revolution*, given by

$$\{x^{-p} \exp(ix) : 1 + 2\pi(k-1) < x \leqslant 1 + 2\pi k\},$$

has diameter comparable to $k^{-p}$, while the distance between the end points (or, outer radius minus inner radius) is comparable to $k^{-(p+1)}$. These are fundamental measurements considered in the winding problem, see [93], and measure how big the $k$th revolution is and how tightly it is wound, respectively. Again the 'logarithmic eccentricity' is

$$\frac{\log\left(k^{-p}\right)}{\log\left(k^{-(p+1)}\right)} = \frac{p}{p+1} = \rho.$$

It is probably naïve to expect (17.3) to be the 'rule of thumb' for sets exhibiting 'rough homogeneity', but it would certainly be interesting to identify other natural classes of set for which this formula holds for a particular choice of 'fundamental ratio' $\rho$. Natural examples to consider include self-affine sponges, where there are many 'logarithmic eccentricities' to deal with, or Lalley-Gatzouras carpets where the logarithmic eccentricities vary across the set. As such, it seems extremely unlikely that (17.3) holds generally in these classes.

Finally, recall that Lemma 3.4.4 and Corollary 3.3.4 show

$$\left(\frac{1-\rho}{1-\theta}\right)\dim_{\mathrm{qA}} F \leqslant \dim_{\mathrm{A}}^{\theta} F \leqslant \frac{\overline{\dim}_{\mathrm{B}} F}{1-\theta} \tag{17.5}$$

and the range over which at least one of these bounds is non-trivial is $(0,\rho)$. The lower bound is necessarily attained at $\theta=\rho$ and the upper bound is attained at $\theta=0$ (extending by continuity). If we take the natural weighted average of these bounds across the relevant interval, we get

$$\left(\frac{\theta}{\rho}\right)\left(\frac{1-\rho}{1-\theta}\right)\dim_{\mathrm{qA}} F + \left(1-\frac{\theta}{\rho}\right)\frac{\overline{\dim}_{\mathrm{B}} F}{1-\theta}$$
$$= \overline{\dim}_{\mathrm{B}} F + \frac{(1-\rho)\theta}{(1-\theta)\rho}\left(\dim_{\mathrm{qA}} F - \overline{\dim}_{\mathrm{B}} F\right)$$

and so recover the increasing part of (17.3).

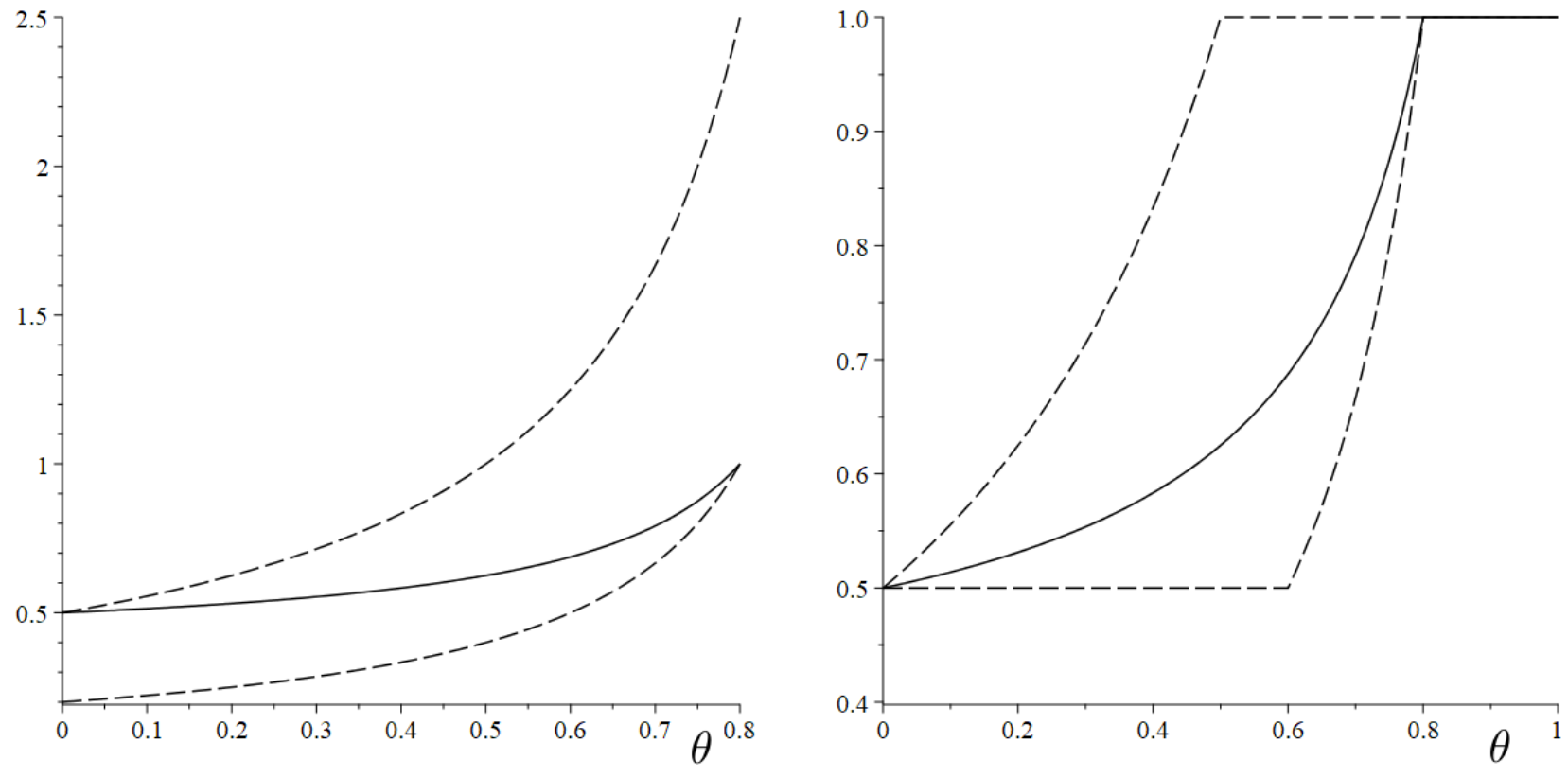

Figure 17.1 Left: the upper and lower bounds from (17.5) (dashed) and the weighted average (solid). Right: plots of the bounds on the Assouad spectrum (dashed) from Lemma 3.4.4 and Corollary 3.3.4 and the formula for the Assouad spectrum given by (17.3) (solid). In both cases $\overline{\dim}_{\mathrm{B}} F = 1/2$, $\dim_{\mathrm{qA}} F = 1$, and $\rho = 4/5$.

## 17.8 Projections

Theorem 10.2.1 establishes an almost sure lower bound for the Assouad dimension of projections of sets in $\mathbb{R}^d$ onto $k$-dimensional subspaces. Theorem 10.2.2 upgrades this fact in the case $d = 2$ and $k = 1$ by establishing that the lower bound in fact holds for all projections outside of a set of exceptions of Hausdorff dimension 0. It is natural to ask about the higher dimensional setting. One cannot hope for a 0-dimensional set of exceptions in general. Suppose $F$ is contained in a 1-dimensional subspace $V \subseteq \mathbb{R}^d$ and consider projections of $F$ onto $k$-dimensional subspaces. Then there is a $k(d-1-k)$-dimensional family of $k$-dimensional subspaces (those transversal to $V$) onto which the projection of $F$ is singleton.

**Question 17.8.1** Is it true that, for any non-empty set $F \subseteq \mathbb{R}^d$ and $1 \leqslant k < d$,

$$\dim_{\mathrm{H}}\{\pi \in G(k,d) : \dim_{\mathrm{A}} \pi F < \min\{k, \dim_{\mathrm{A}} F\}\} \leqslant k(d-1-k)?$$

There is also the question of possible refinements of the conclusion that the Hausdorff dimension of the exceptional set is 0. The following question was asked in [226], and there are obvious analogous questions in higher dimensions.

**Question 17.8.2** Is it true that, for any non-empty set $F \subseteq \mathbb{R}^2$,

$$\dim_{\mathrm{P}}\{\pi \in G(1,2) : \dim_{\mathrm{A}} \pi F < \min\{1, \dim_{\mathrm{A}} F\}\} = 0?$$

One could also consider different gauge functions to further refine the Hausdorff dimension result.

Theorem 10.2.5 revealed that the function $\pi \mapsto \dim_{\mathrm{A}} \pi F$ may realise *any* upper semi-continuous function $\phi : G(1,2) \to [0,1]$ upon careful selection of a compact set $F \subseteq \mathbb{R}^2$. Note, however, that $\pi \mapsto \dim_{\mathrm{A}} \pi F$ is not necessarily upper semi-continuous (for example, consider a line segment). On the other hand, $\pi \mapsto \dim_{\mathrm{A}} \pi F$ cannot realise *all* functions $\phi : G(1,2) \to [0,1]$. This can be seen by another simple cardinality argument. The set of possible functions $\phi$ has cardinality $\#\mathcal{P}(\mathbb{R})$ (that is, the cardinality of the power set of $\mathbb{R}$), but the collection of possible choices for compact $F$ has cardinality $\#\mathbb{R}$ (that is, the cardinality of $\mathbb{R}$), which is strictly smaller. Recall the connection between this issue and that following Question 17.3.2.

**Question 17.8.3** Given a compact set $F \subseteq \mathbb{R}^2$, is the function $\pi \mapsto \dim_{\mathrm{A}} \pi F$ Borel measurable? If so, does it lie in a finite Borel class?

## 17.9 Distance sets

The planar version of the distance set problem for Assouad dimension was resolved in [90], but progress in higher dimensions so far is tied to the Hausdorff dimension problem, which is still wide open. The first question we ask is very concrete, and the second is a little more open ended, motivated by the fact that the Assouad dimension problem is 'ahead' of the Hausdorff dimension problem in the plane.

**Question 17.9.1** Is it true that if $F \subseteq \mathbb{R}^d$ with $d \geqslant 3$, then

$$\dim_{\mathrm{A}} D(F) \geqslant \frac{2}{d} \dim_{\mathrm{A}} F?$$

**Question 17.9.2** For sets $F \subseteq \mathbb{R}^d$ with $d \geqslant 3$, can we move the Assouad dimension version of the distance set problem 'ahead' of the Hausdorff or box dimension version? For example, can we prove that $\dim_{\mathrm{A}} F \geqslant s$ guarantees $\dim_{\mathrm{A}} D(F) = 1$ for some $s$ for which we do not yet have the corresponding Hausdorff or box dimension result?

A popular variant of the distance set problem is to consider *pinned* distance sets. Given $F \subseteq \mathbb{R}^d$ and $x \in F$, the pinned distance set (of $F$ at $x$) is defined by

$$D_x(F) = \{|x - y| : y \in F\}.$$

Note that pinned distance sets are subsets of the distance set. Many natural questions arise, but generally speaking if the distance set is large, we expect there to be a large pinned distance set; or even *many* large pinned distance sets, see [253]. For example, we can ask for a pinned version of Theorem 11.1.3.

**Question 17.9.3** Is it true that if $F \subseteq \mathbb{R}^2$ is a set with $\dim_{\mathrm{A}} F > 1$, then there exists $x \in F$ such that $\dim_{\mathrm{A}} D_x(F) = 1$?

Finally, a more general question, which equally applies to any of the Hausdorff, box, packing or Assouad dimensions.

**Question 17.9.4** Let dim denote one of the Hausdorff, box, packing or Assouad dimensions. Given $s \in (0, d]$, what is

$$\inf\{\dim D(F) : \dim F = s\}?$$

Is the function

$$s \mapsto \inf\{\dim D(F) : \dim F = s\}$$

independent of the specific choice of dim?

## 17.10 The Hölder mapping problem and dimension

Motivated by the spiral winding problem discussed in Section 13.2, consider the following general problem. Fix bounded sets $X, Y \subseteq \mathbb{R}^d$ which are homeomorphic to each other, for example, the line segment $(0,1)$ and the spiral $\mathcal{S}_p$. First, ask if there exists a bi-Hölder homeomorphism between $X$ and $Y$. If the answer to this is yes, then ask for sharp estimates on the forward and backwards Hölder exponents $\alpha$ and $\beta$. Second, compute the dimensions and dimension spectra of $X$ and $Y$ and consider the estimates on $\alpha$ and $\beta$ which come directly from this information (by applying results from Section 3.4.3, for example). It seems interesting to develop a thorough understanding of how good this information is: when do the sharp estimates come from dimension theory and when do they not? Moreover, it would also be interesting to know which dimensions or dimension spectra generally provide better information. For example, the polynomial spirals we considered in Section 13.2 provide examples where sharp information is *not* provided by dimension theory, and where the Assouad spectrum performs best. This discussion is deliberately rather vague so that it can be applied in many situations. Rather than get too hung up on a particular version of the problem, a better overall strategy is to unearth general phenomena using elegant examples.

## 17.11 Dimensions of graphs

One of the earliest examples of what is now known as a fractal is the graph of the *Weierstrass function*. Discovered in 1872 by Carl Weierstrass, it was the first example of a continuous real-valued function which is nowhere differentiable.[2] Given integer $b \geqslant 2$ and $a \in (1/b, 1)$, the Weierstrass function $W_{a,b} : \mathbb{R} \to \mathbb{R}$ is defined by

$$W_{a,b}(x) = \sum_{n=0}^{\infty} a^n \cos(2\pi b^n x).$$

It was a long-standing open problem to compute the Hausdorff dimension of the graph of $W_{a,b}$, that is, the set[3]

$$G(W_{a,b}) = \{(x, W_{a,b}(x)) : x \in [0,1]\},$$

[2] Bolzano had an example of such a function as early as 1830 but it was only published posthumously in 1922.

[3] the domain is not important here, but it is convenient to assume it is bounded.

see, for example, [32]. The conjecture was that the Hausdorff dimension should equal

$$2 + \frac{\log a}{\log b} \in (1, 2). \tag{17.6}$$

This was resolved by Shen [251] in full generality and for certain ranges of parameters a little earlier by Barański, Bárány and Romanowska [16]. The box and packing dimensions of the graph are also given by (17.6), but this was known somewhat earlier, see for example [144]. We do not know the Assouad or quasi-Assouad dimension of the graph, but there is some heuristic evidence, provided by Yu, to suggest that it is not given by (17.6) in general. More precisely, a related question was considered by Yu [279] where the cosine function appearing in the definition of the Weierstrass function is replaced by the 'distance to the nearest integer' function $T : \mathbb{R} \to \mathbb{R}$ given by $T(x) = \min\{|x - n| : n \in \mathbb{Z}\}$. The only important difference here is that the smooth oscillations provided by cos are replaced with sawtooth type oscillations. The resulting functions are known as *Takagi functions* and are easier to handle than the Weierstrass function, see the survey [3], and the Hausdorff, box and packing dimensions are also known to be given by the formula (17.6). Yu showed that the Assouad dimension of the graph of the Takagi function can exceed the box dimension [279, Theorem 1.1]. He was able to bound the Assouad dimension from below by $1 + 1/k$ if $ab > 1$ is the root of a Littlewood polynomial of degree $k - 1$. Recall that a Littlewood polynomial is a polynomial all of whose coefficients are $\pm 1$. In particular, fixing $ab$ and thus $k$, the box dimension of the graph can be made arbitrarily close to 1 by increasing $b$ (and decreasing $a$ accordingly).

**Question 17.11.1** What are the Assouad dimensions, quasi-Assouad dimensions, and Assouad spectra of the graphs of the Weierstrass and Takagi functions?

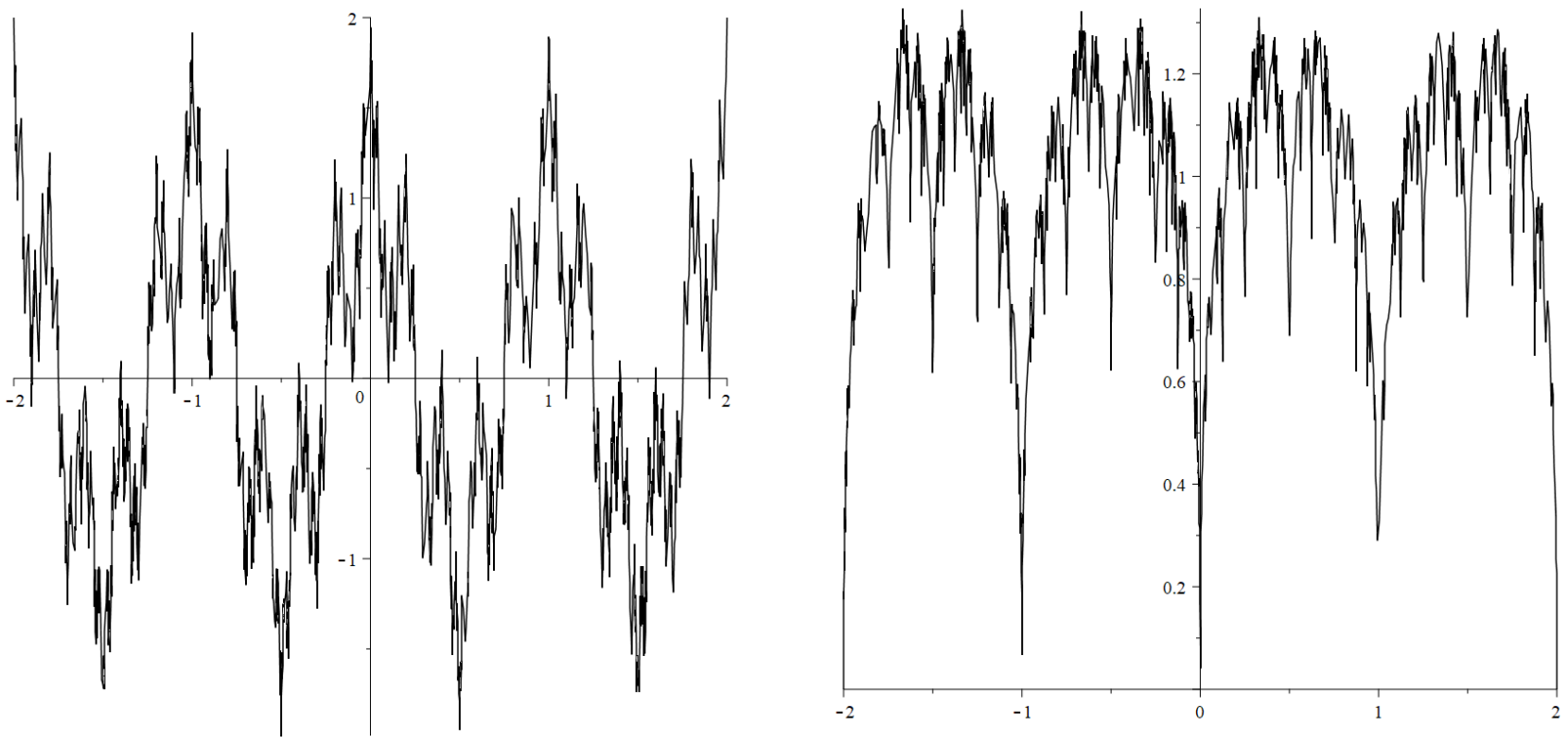

Figure 17.2 Left: graph of the Weierstrass function with $b = 5$ and $a = 1/2$. Right: graph of the Takagi function with $b = 2$ and $a = 3/4$.

# List of notation

**Dimensions and spectra**

| | |
|---|---|
| $\dim_{\mathrm{A}}^{\phi}$ | $\phi$-Assouad dimension |
| $\dim$ | an unspecified dimension |
| $\dim_{\mathrm{A}}$ | Assouad dimension |
| $\dim_{\mathrm{A}}^{\theta}$ | Assouad spectrum |
| $\overline{\dim}_{\mathrm{A}}^{\theta}$ | upper (Assouad) spectrum |
| $\underline{\dim}_{\mathrm{B}}$ | lower box dimension |
| $\overline{\dim}_{\mathrm{B}}$ | upper box dimension |
| $\dim_{\mathrm{B}}$ | box dimension |
| $\mathcal{C}\dim_{\mathrm{A}}$ | conformal Assouad dimension |
| $\mathcal{C}\dim_{\mathrm{H}}$ | conformal Hausdorff dimension |
| $\dim_{\mathrm{H}}$ | Hausdorff dimension |
| $\dim_{\mathrm{L}}$ | lower dimension |
| $\dim_{\mathrm{L}}^{\theta}$ | lower spectrum |
| $\dim_{\mathrm{ML}}$ | modified lower dimension |
| $\dim_{\mathrm{P}}$ | packing dimension |
| $\dim_{\mathrm{qA}}$ | quasi-Assouad dimension |
| $\dim_{\mathrm{qL}}$ | quasi-lower dimension |
| $\dim_{\mathrm{T}}$ | topological dimension |

**General notation**

| | |
|---|---|
| $\lvert x\rvert$ | the absolute value of a real or complex number $x$ |
| $B(x,r)$ | a (usually closed) ball centred at $x$ with radius $r$ |
| $\#F$ | cardinality of a set $F$, which may be finite or infinite |
| $E \times F$ | Cartesian product of sets $E$ and $F$ |
| $\overline{F}$ | closure of a set $F$ |
| $F^c$ | complement of a set $F$ (usually with $\mathbb{R}^d$ as the ambient space) |
| $\mathbb{C}$ | the complex numbers, often identified with the complex plane |
| $\mathrm{Con}(d)$ | isometry group of $\mathbb{D}^{d+1}$ |
| $\lvert F\rvert$ | the diameter $\lvert F\rvert = \sup_{x,y\in F} \lvert x-y\rvert$ of a set $F$, which may be finite or infinite |
| $\mathbb{D}^{d+1}$ | Poincaré disk model of $d+1$ dimensional hyperbolic space |
| $D(F)$ | the distance set $D(F) = \{\lvert x-y\rvert : x, y \in F\}$ of a set $F$ |
| $G(k,d)$ | Grassmannian manifold consisting of $k$-dimensional subspaces of $\mathbb{R}^d$ |
| $\mathcal{H}^s$ | $s$-dimensional Hausdorff measure |
| $\mathcal{H}^s_r$ | $r$-approximate $s$-dimensional Hausdorff measure |
| $d_{\mathcal{H}}$ | Hausdorff metric |
| $\mathcal{K}(X)$ | set of all non-empty compact subsets of a set $X$ |
| $L(\Gamma)$ | the limit set of a Kleinian group $\Gamma$ |
| $\mathcal{L}^d$ | $d$-dimensional Lebesgue measure |
| $\mu$ | a locally finite Borel measure |
| $M_r(F)$ | maximum number of points in an $r$-separated subset of $F$ |
| $\mathbb{N}$ | the natural numbers (which do not include 0) |
| $N_r(F)$ | minimum number of open sets of diameter $r$ required to cover $F$ |
| $O(d)$ | $d$-dimensional orthogonal group |
| $\pi$ | an element of $G(k,d)$, which may be identified with an orthogonal projection |
| $\pi^{\perp}$ | orthogonal complement of $\pi \in G(k,d)$ |
| $\mathbb{Q}$ | the rational numbers |
| $\mathbb{R}$ | the real numbers, often identified with the real line |
| $\mathbb{R}^d$ | $d$-dimensional Euclidean space |
| $SO(d)$ | $d$-dimensional special orthogonal group |
| $\mathrm{supp}(\mu)$ | support $\mathrm{supp}(\mu) = \{x : \mu(B(x,r)) > 0 \text{ for all } r > 0\}$ of a measure $\mu$ |
| $\mathbb{Z}$ | the integers |

**Notation relating to iterated function systems**

| | |
|---|---|
| $I$ | the identity map |
| $\mathcal{I}$ | a finite index set for an iterated function system (IFS) |
| $\mathcal{I}^k$ | the set of all sequences of length $k$ with entries in $\mathcal{I}$ |
| $\mathcal{I}^*$ | the set of all finite sequences with entries in $\mathcal{I}$ |
| $\mathcal{I}^{\mathbb{N}}$ | the set of all infinite sequences with entries in $\mathcal{I}$ |
| $\boldsymbol{i}^-$ | the sequence $(i_1, i_2, \ldots, i_{k-1})$ for $\boldsymbol{i} = \big(i_1, i_2, \ldots, i_k\big) \in \mathcal{I}^k$ |
| $[\boldsymbol{i}]$ | the cylinder $\{\boldsymbol{j} \in \mathcal{I}^{\mathbb{N}} : \boldsymbol{i} \prec \boldsymbol{j}\}$ for $\boldsymbol{i} \in \mathcal{I}^k$ |
| IFS | iterated function system |
| OSC | open set condition |
| $\Pi$ | symbolic coding map associated with an IFS |
| $\prec$ | 'is a prefix of', e.g. $\boldsymbol{i} \prec \boldsymbol{j}$ for $\boldsymbol{i}, \boldsymbol{j} \in \mathcal{I}^*$ |
| $Q(\mathbf{d}, r)$ | approximate square centred at $\mathbf{d}$ with radius $r$ |
| $S_{\boldsymbol{i}}$ | the composition $S_{i_1} \circ S_{i_2} \circ \cdots \circ S_{i_k}$ for $\boldsymbol{i} = \big(i_1, i_2, \ldots, i_k\big) \in \mathcal{I}^k$ |
| SSC | strong separation condition |
| VSSC | very strong separation condition |
| WSP | weak separation property |

# Index